\documentclass[twoside,12pt]{amsart}

\usepackage{amssymb}
\usepackage{amsthm}
\usepackage{amsfonts}
\usepackage{amsmath}
\usepackage[T1]{fontenc}
\usepackage{lscape}
\usepackage{url}
\usepackage{longtable}
\usepackage{graphicx}
\usepackage{latexsym}
\usepackage{colortbl}
\usepackage[all]{xy}

\newtheorem{thm}{Theorem}[section]

\newtheorem*{mtheorem*}{Main Theorem}
\newtheorem{lem}[thm]{Lemma}
\newtheorem*{lem*}{Lemma}
\newtheorem*{cor*}{Corollary}
\newtheorem{dfn}[thm]{Definition}
\newtheorem{prop}[thm]{Proposition}
\newtheorem{cor}[thm]{Corollary}

\newtheorem{rem}[thm]{Remark}
\newtheorem{hypo}[thm]{Hypothesis}
\newenvironment{prf}{{\sc Proof.}}{\hfill $\Box$}

\def\Aut{ \text{\rm Aut} }
\def\Inn{ \text{\rm Inn} }
\def\Out{ \text{\rm Out} }
\def\Hom{ \text{\rm Hom} }
\def\Irr{ \text{\rm Irr} }

\def\IBr{ \text{\rm IBr} }

\def\Ind{ \text{\rm Ind} }

\def\Res{ \text{\rm Res} }

\def\Spin{ \text{\rm Spin} }
\def\Stab{ \text{\rm Stab} }

\def\diag{ \text{\rm diag} }
\def\GL{ \text{\rm GL} }

\def\Sp{ \text{\rm Sp} }

\def\PSp{ \text{\rm PSp} }
\def\SL{ \text{\rm SL} }
\def\PSL{ \text{\rm PSL} }
\def\PGL{ \text{\rm PGL} }
\def\SU{ \text{\rm SU} }
\def\PSU{ \text{\rm PSU} }

\def\rk{ \text{\rm rk} }
\newcommand{\CJxxx}{\text{\rm ad}}

\def\myhline{}

\newcommand{\R}{\mathcal{R}}
\newcommand{\A}{{\rm A}}
\newcommand{\B}{{\rm B}}
\newcommand{\C}{{\rm C}}

\newcommand{\mycol}{\mathbin{:}}
\newcommand{\Le}{{L_\varepsilon}}

\newcommand{\Syl}{\mathop{\rm Syl}\nolimits}

\newcommand{\cR}{{\mathcal R} }
\newcommand{\cE}{{\mathcal E} }

\newcommand{\cW}{{\mathcal W} }
\newcommand{\ibr}{{\rm {}IBr} }
\newcommand{\midd}{\mid}

\newcommand{\id}{{\rm {}id}}

\renewcommand\thetable{\Roman{table}}

\begin{document}

%\frontmatter

\title{The inductive blockwise Alperin weight condition for the Chevalley groups $F_4(q)$}

\author{Jianbei An}
\address{Department of Mathematics, The University of Auckland,
Private Bag 92019, Auckland 1142, New Zealand}
\email{j.an@auckland.ac.nz}

\author{Gerhard Hiss}
\address{Lehrstuhl f{\"u}r Algebra und Zahlentheorie, RWTH Aachen University,
52056 Aachen, Germany}
\email{gerhard.hiss@math.rwth-aachen.de}

\author{Frank L\"ubeck}
\address{Lehrstuhl f{\"u}r Algebra und Zahlentheorie, RWTH Aachen University,
52056 Aachen, Germany}
\email{frank.luebeck@math.rwth-aachen.de}



\date{May 12, 2022}

\begin{abstract}
We verify the inductive blockwise Alperin weight condition in odd
characteristic~$\ell$ for the finite exceptional Chevalley groups 
$F_4(q)$ for~$q$ not divisible by~$\ell$.
\end{abstract}

\maketitle

%\mainmatter

\tableofcontents

\markright{INDUCTIVE BLOCKWISE ALPERIN WEIGHT CONDITION FOR $F_4(q)$}
\section{Introduction}

Alperin's weight conjecture, published in~\cite{Ap} in 1987, is one of the 
famous intriguing conjectures in the representation theory of finite groups.
It postulates the coincidence of two invariants of a finite group, which are
defined from global, respectively local data, and which seem unrelated at 
first sight.
Given a prime $\ell$ and a finite group~$G$, the global invariant in Alperin's
weight conjecture is the number of conjugacy classes of~$G$ of elements of order
prime to~$\ell$. This is the same as the number of absolutely irreducible 
$\ell$-modular characters of~$G$. The local invariant is the number of conjugacy 
classes of $\ell$-weights of~$G$. An $\ell$-weight is a pair $(Q,\chi)$, 
where~$Q$ is a finite $\ell$-subgroup of~$G$, and $\chi$ is an irreducible 
character of $N_G( Q )/Q$ with the property that $\ell$ does not divide 
$|N_G(Q)/Q|/\chi( 1 )$. If $(Q,\chi)$ is an $\ell$-weight, we also call~$Q$
a weight subgroup. In his original paper~\cite{Ap}, Alperin 
proved the truth of his conjecture in many instances, for example in case~$G$ is 
a finite group of Lie type of characteristic~$\ell$. Soon after its appearance, 
the Alperin weight conjecture has been verified for various series of finite 
groups, in particular the alternating and symmetric groups and their Schur 
covering groups, as well as the general linear groups; see \cite{AlpFo} 
and~\cite{MiOl}.

Beginning with a paper by Kn{\"o}rr and Robinson \cite{KR}, the Alperin weight 
conjecture has generated a huge new field of research, including various 
reformulations, generalizations and reductions. The number of contributions
in this direction is too large to be listed here. %Staczewski, Dade ...

Since the classification of the finite simple groups, it has become a common
approach to reduce questions on general finite groups to related, often more
complex questions on finite simple groups or their universal covering groups.
A reduction theorem to Alperin's weight conjecture in this spirit has been 
established by Navarro and Tiep~\cite{nt}. The authors formulated a set of 
conditions on a non-abelian finite simple group and all of its perfect covering 
groups. If this set of conditions is satisfied, the corresponding finite simple 
group is called AWC-good. The main theorem of \cite{nt} states that if all 
non-abelian finite simple groups are AWC-good, then the Alperin weight 
conjecture holds for all finite groups. 

There is a blockwise version of Alperin's weight conjecture, already present 
in~\cite{Ap}, and a reduction theorem for this blockwise version was proved by 
Sp{\"a}th in~\cite{Spaeth13}. Refining the conditions for AWC-goodness, 
Sp{\"a}th introduced a collection of properties for a non-abelian finite simple 
group, called the \textit{inductive blockwise Alperin weight condition} 
(inductive BAW condition for short). Naturally, this set of conditions, which is 
formulated for a fixed prime~$\ell$ and a fixed $\ell$-block, is even more 
intricate than the one of~\cite{nt}. If the condition is satisfied for a 
prime~$\ell$ and all $\ell$-blocks of a non-abelian finite simple group~$G$ and 
its covering groups, we say the inductive BAW condition holds for~$G$ at the 
prime~$\ell$. According to the main theorem of~\cite{Spaeth13}, the blockwise 
Alperin weight conjecture for the prime~$\ell$ holds for all finite groups if 
the inductive BAW condition is satisfied for all non-abelian finite simple 
groups at the prime~$\ell$. A simplified version of the inductive BAW 
conditions, adapted to our purposes, is presented in Hypothesis~\ref{defBWC} 
below.

The inductive BAW condition has been verified for various series of simple 
groups, for example some groups of Lie type of small rank already by Sp{\"a}th 
in~\cite{Spaeth13}. Schulte shows in~\cite{ESch} that the inductive 
BAW condition holds for the exceptional Chevalley groups $G_2( q )$ and the 
Steinberg triality groups. Malle verified the inductive BAW condition for the
Suzuki and Ree groups in \cite[Theorem~$5.1$]{MalleAb}. For blocks with cyclic
defect groups the inductive Alperin Weight condition holds by a theorem of
Koshitani and Sp{\"a}th; see \cite[Theorem~$1.1$]{KS16}. Also, if $G$ is a finite
simple group of Lie type of characteristic $\ell$, then the inductive BAW condition
for~$G$ holds at the prime~$\ell$ by a result of Sp{\"a}th; 
see~\cite[Theorem~C]{Spaeth13}. Our article is a contribution to the
programme of verifying the inductive BAW condition for the non-abelian finite 
simple groups.

\begin{mtheorem*}\label{thm}
Let $\ell$ be an odd prime and let $G = F_4(q)$ for a prime power~$q$ coprime 
to~$\ell$. Then the inductive blockwise Alperin weight condition holds for every 
$\ell$-block of $G$.
\end{mtheorem*}

Let~$G$ and $\ell$ be as in the main theorem. The proof of this relies on a 
careful analysis of the $\ell$-blocks of~$G$ and their
invariants, as well as the candidates for the weight subgroups.
The major part of our paper is devoted to the case
$\ell = 3$. If $\ell > 3$, the Sylow $\ell$-subgroups of~$G$ are abelian, and 
substantial results towards the main theorem are already contained 
in~\cite{MalleAb} by exhibiting natural bijections between weights and
absolutely irreducible $\ell$-modular characters, and it remains to establish 
equivariance of
these bijections with respect to outer automorphisms of~$G$.
The classification of the $\ell$-blocks follows the route
laid out by the paper \cite{BrouMi} by Brou{\'e} and Michel, i.e.\ it is based on
the classification of the semisimple conjugacy classes of elements of order 
prime to~$\ell$. The semisimple conjugacy classes are grouped into
finitely many class types, and our results can be proved uniformly for all
elements inside each class type. The class types and properties of the 
corresponding elements, in particular their centralizers are determined and
enumerated in~\cite{LL}.

Our task is simplified to some extent 
%as~$G$ is isomorphic to its dual group,
as the Schur multiplier of~$G$, with one
exception, is trivial, as the outer automorphism group of~$G$ is cyclic and as
all proper Levi subgroups of~$G$ are of classical type. The case $\ell = 3$ 
presents an interesting example for the verification of the inductive BAW 
condition, as the Sylow $3$-subgroups of~$G$ are non-abelian. The results for
the principal $3$-block reveal a distinctive different behavior in the cases when
$9$ divides or does not divide $q^2 - 1$. Our investigations
for $\ell = 3$ are largely supported by the fact that~$3$ is good for all proper 
Levi subgroups of~$G$, and that the
radical $3$-subgroups of~$G$, which are the candidates for the weight subgroups
have been classified in~\cite{AH2} and~\cite{AnDF4}. The missing classification
of the radical $2$-subgroups of~$G$ prevents us from extending our results
to the case $\ell = 2$. Many of our investigations
are highly assisted by deep results of Bonnaf{\'e}, Dat and Rouquier
\cite{BoDaRo}, as well as by recent work of Boltje and 
Perepelitsky~\cite{BoPe}.

Let us now comment on the contents of the individual sections of our article.
Section~$2$ is devoted to the introduction of notation and background material 
on groups and representations, occasionally refined and extended for our 
purpose. Subsection~\ref{InductiveBlockwiseAlperin} contains the version of 
the inductive BAW 
condition relevant to our investigations. In Section~$3$ we recall the principal 
concepts and results on finite groups of Lie type needed later on. In 
particular, we summarize the above mentioned theorems of Bonnaf{\'e}, Dat and 
Rouquier and various others in Theorem~\ref{BoRoEtAl} to have a convenient 
reference. Some consequences of this major theorem are derived. We also collect 
some auxiliary results useful in our later study. Section~$4$ introduces the 
group $F_4(q)$ as the group of fixed points under a Steinberg morphism of a 
simple algebraic group~$\mathbf{G}$ of type~$F_4$ over the algebraic closure of 
the field with~$q$ elements. We establish our notation for the corresponding 
root system and the Weyl group. We also introduce class types, and a duality
of Levi subgroups of~$G$ arising from the fact that~$G$ is isomorphic to its dual 
group. We then investigate in great detail the structure of some Levi subgroups
of~$G$. This yields a first new result in our paper, 
Corollary~\ref{MaximalExtendability}, which states that the $e$-split Levi 
subgroups of~$G$ satisfy the maximal extendibility condition. This is a crucial
ingredient in the proof of our main theorem and might be of independent 
interest. Section~$5$ is devoted to the description of the $\ell$-blocks of~$G$
and some of their invariants for primes~$\ell$ not dividing~$q$. 
The description for the primes $\ell > 3$ was known before and is due to 
Brou{\'e}, Malle and Michel~\cite{BrouMaMi}, as well as 
Brou{\'e} and Michel~\cite{BrouMi}; see also~\cite{MalleAb}. The main effort
here is spent to handle the prime $\ell = 3$ for the non-unipotent blocks. 
The unipotent blocks have been treated by Enguehard in~\cite{En00}.
Although not pursued furthermore, we also include a subsection for the prime 
$\ell = 2$. The results are presented in form of tables in the Appendix; see
Tables~\ref{1}--\ref{19} and~\ref{21}. The action of the outer automorphisms
of~$G$ on the set of absolutely irreducible $\ell$-modular characters of an 
invariant block is determined in Proposition~\ref{ActionAutomorphisms}.

Section~$6$ is dedicated solely to the prime~$3$. We recall the construction 
of the radical $3$-subgroups of~$G$ established in~\cite{AnDF4} and~\cite{AH2}
to some detail, as this will be important later on. The defect groups of the
$3$-blocks of~$G$ are, in particular, radical $3$-subgroups. We describe these
defect groups by identifying the corresponding conjugacy class of radical 
$3$-subgroups. A preliminary result states that the Sylow $3$-subgroups of
the centralizers of semisimple elements are all radical $3$-subgroups. It would
be interesting to find an a priori reason for this observation, which might 
extend to other groups of Lie type. The weight subgroups are
radical $3$-subgroups, and we determine the candidates for the weight subgroups
among the radical $3$-subgroups.

Section~$7$ contains a proof of one half of the inductive BAW condition by 
giving a bijection between the absolutely irreducible 
$\ell$-modular characters of a block and the conjugacy classes of weights 
associated to the block; see Theorem~\ref{thmAWC}. It is worth remarking that 
the blocks in~$G$ are split in the sense that the associated canonical 
characters extend to their inertia subgroups.

Finally, Section~$8$ proves the equivariance of the bijections established in 
Theorem~\ref{thmAWC} with respect to outer automorphisms. This equivariance is 
established in Theorem~\ref{thmmain} after a long series of rather technical 
and involved preparations. It is unfortunate that we were not able to 
find a uniform approach to these results. Instead, we develop numerous ad 
hoc methods for specific situations. On the other hand, the variety of the 
methods introduced here might facilitate analogous investigation for other
exceptional groups of Lie type.

The referee of the first version of this article has pointed us to the
preprint~\cite{FLZ}. This reduces the proof of the equivariance condition for 
the non-quasi-isolated blocks of~$G$ to the verification of the inductive 
blockwise Alperin weight condition for quasi-isolated blocks of simple groups 
involved in~$G$. At this time, to the best of our knowledge, this verification
has not been established completely for all the cases relevant to our work, 
which includes groups such as $\Sp_6( q )$ or $\Spin_7( q )$. We have therefore 
decided to retain the original presentation of Section~$8$, as this is 
comprehensive and self-contained.

%\markboth{Notation and preliminaries}
\section{Notation and preliminaries}
\label{Preliminaries}

Throughout this section~$G$ denotes a group, which is assumed to be finite, 
except in Subsection~\ref{Groups} below.

\subsection{Groups}\label{Groups}
Most of our notation for groups is standard.
If $g, x \in G$ we write $g^x := x^{-1} g x$ for the right conjugation of~$g$ 
by~$x$. This notation is extended to subsets of~$G$.
If $K, H \leq G$, then we write  $K \leq_G H$ whenever there exists
$x \in G$ with $K^x \leq H$. Analogously, we define $H =_G K$. 
For $H \leq G$, we write
\[\Out_G (H) := N_G (H)/HC_G (H).\]
Notice that if $H \leq M \leq G$, then $\Out_M( H )$ naturally embeds into
$\Out_G(H)$. The commutator subgroup of~$G$ is denoted by $[G,G]$.

If~$p$ is a prime and~$G$ is a finite abelian $p$-group, we write 
$\Omega_1( G )$ for the subgroup of~$G$ generated by the elements of order~$p$. 
Then $\Omega_1(G)$ is an elementary abelian $p$-group.

If $n, m$ are positive integers, then $[n]$ denotes a cyclic group of order $n$,
and $[n]^m$ the direct product of~$m$ copies of~$[n]$. If $n \in \{ 2, 3, 6 \}$,
we also omit the outer brackets for simplicity of notation.

Our notation for simple groups and related groups follows one of the standard
conventions from the literature. 

Recall the notation $\SL^\varepsilon_n(q)$ and
$\GL^\varepsilon_n(q)$ for $\varepsilon \in \{ 1, - 1\}$: if 
$\varepsilon = 1$, then these are the special
linear and general linear groups of degree~$n$ over the field
$\mathbb{F}_q$; if $\varepsilon = -1$, then these are the corresponding 
special unitary and unitary groups, respectively, defined over
$\mathbb{F}_{q^2}$. Recall that $\SL_2(q) = \SL_2^{-1}(q) = \Sp_2(q)$.

The notation for group extensions follows the Atlas \cite{Atlas} convention,
i.e.\ $A.B$ is a group with a normal subgroup isomorphic to~$A$ and corresponding
factor group isomorphic to~$B$; see \cite[Page~xx]{Atlas}.
If not indicated by brackets, we read group extensions $A.B.C$ from left to 
right, that is, $A.B.C=(A.B).C$.

\subsection{Characters and modular systems}
Let~$\ell$ be a prime. Fix an $\ell$-modular system 
$(\mathcal{K}, \mathcal{O}, \overline{\mathcal{O}})$ for~$G$, 
where~$\mathcal{O}$ is a complete discrete valuation ring of characteristic~$0$
with residue class field $\overline{\mathcal{O}}$ of characteristic~$\ell$ and
field of fractions~$\mathcal{K}$, which is large enough for~$G$,
i.e.\ $\mathcal{K}$ contains a $|G|$th root of unity.
In the following, the term $\ell$-block refers to a block of 
$\overline{\mathcal{O}}G$ or its lift to $\mathcal{O}G$. 
If~$\theta$ is a $\mathcal{K}$-valued class function of~$G$, we write $\check{\theta}$ for 
the restriction of~$\theta$ to the set of $\ell$-regular elements of~$G$.
The set of ordinary irreducible characters (i.e.\ irreducible $\mathcal{K}$-characters) 
of~$G$ is denoted by $\Irr(G)$, and the subset of $\Irr(G)$ of $\ell$-defect 
zero characters by $\Irr^0(G)$. We also write $\IBr_\ell(G)$ for the set of 
irreducible Brauer characters of~$G$ (with respect to 
$(\mathcal{K},\mathcal{O},\overline{\mathcal{O}})$). If~$B$ is a union of $\ell$-blocks,
we use the notation $\Irr(B)$ and $\IBr(B)$ for the sets
of irreducible ordinary, respectively Brauer characters of~$B$, and we write
$\mathbb{Z}[\Irr(B)]$, respectively $\mathbb{Z}[\IBr(B)]$ for the corresponding
sets of generalized characters.

Let $H \leq G$, and let $\chi$ and $\psi$ denote $\mathcal{K}$-valued class 
functions on~$G$, respectively~$H$. Then $\Res^G_H( \chi )$ and 
$\Ind_H^G( \psi )$ denote the restriction of~$\chi$ to~$H$, respectively the
class function of~$G$ obtained by inducing~$\psi$ to~$G$.
If $H$ is a finite group and $\chi$ and $\psi$ class function of~$G$,
respectively~$H$, with values in~$\mathcal{K}$, we write $\chi \boxtimes \psi$ 
for the outer product of~$\chi$ and~$\psi$. This is a $\mathcal{K}$-valued class 
function of $G \times H$.

\subsection{Actions of automorphisms}
We collect some miscellaneous results on the action of automorphisms on characters.

\addtocounter{thm}{3}
\begin{lem}
\label{StableInLinearCase}
Let $N \unlhd G$ and let $\sigma \in \Aut(G)$ stabilize~$N$. Assume that every 
coset of~$N$ in~$G$ contains a $\sigma$-stable element.

Let $\lambda \in \Irr( G )$ with $\lambda( 1 ) = 1$ such that
$\Res^G_N( \lambda )$ is $\sigma$-stable. Then $\lambda$ is $\sigma$-stable.
\end{lem}
\begin{prf}
This is a straightforward computation using the fact that~$\lambda$ is a 
homomorphism.
%Let $x \in G$. Then $x$ can be written as $x = ng$ with $n \in N$ and $g \in G$
%such that $\sigma(g) = g$. Thus
%        \begin{eqnarray*}
%		{^\sigma\!\lambda}( x ) & = & {^\sigma\!\lambda}( n g ) \\
%                & = & \lambda( \sigma^{-1}(n) ) \lambda( \sigma^{-1}( g ) ) \\
%                & = & {^\sigma\!\lambda}( n ) \lambda( g ) \\
%                & = & \lambda( n ) \lambda( g ) \\
%                & = & \lambda( x ).
%        \end{eqnarray*}
\end{prf}

\begin{lem}
\label{FixingReduction}
Let $N \unlhd G$ with $G/N$ a group of prime order.
Further, let $\sigma \in \Aut(G)$ stabilize~$N$ and let $\chi \in \Irr( G )$.
Suppose that $\Res^G_N( \chi )$ is $\sigma$-stable. If $\Res^G_N( \chi )$
is irreducible, assume that $0 \neq \chi( g ) = {^\sigma\!\chi}( g )$ for some 
$g \in G \setminus N$. Then $\sigma$ fixes~$\chi$.
\end{lem}
\begin{prf}
Suppose first that $\Res^G_N( \chi )$ is reducible, and let $\psi
\in \Irr( N )$ be an irreducible constituent of $\Res^G_N( \chi )$.
By our assumption, $\chi = \Ind_N^G( \psi )$, and ${^\sigma\!\chi} = 
\Ind_N^G( {^\sigma\!\psi} ) = \chi$, as ${^\sigma\!\psi} $ is an irreducible 
constituent of $\Res^G_N( \chi )$.

Suppose next that $\Res^G_N( \chi )$ is irreducible. Then ${^\sigma\!\chi}$
is an extension of $\Res^G_N( \chi )$ to $G$. All such extensions are of the
form $\chi \lambda$, where $\lambda$ is an irreducible character of~$G/N$.
Hence $\chi( g ) = {^\sigma\!\chi}( g ) = \chi( g ) \lambda( g )$ for some 
$\lambda \in \Irr( G/N )$. As $\chi( g ) \neq 0$, this implies 
$\lambda( g ) = 1$ and hence $\lambda$ is the trivial character, as~$G/N$ has 
prime order.
\end{prf}

\medskip

\noindent
Let us record a corollary, which is relevant in our applications.
\begin{cor}
\label{FixingReductionCor}
Let $M, N \unlhd G$ such that $G/N$ is a group of prime order and such that $G = MN$.
Let $\sigma \in \Aut(G)$ stabilize~$M, N$ and let $\chi \in \Irr( G )$. Assume that
$\Res^G_N( \chi )$ is reducible or that $\Res^G_{M\cap N}( \chi )$ is irreducible. 
Assume finally that $\Res^G_M( \chi )$ and $\Res^G_N( \chi )$ are $\sigma$-stable. 
Then~$\chi$ is $\sigma$-stable.
\end{cor}
\begin{prf}
By Lemma~\ref{FixingReduction} we may assume that $\Res^G_{M \cap N}( \chi )$ is 
irreducible. As $G/N \cong M/(M\cap N)$, there is $g \in M \setminus N$ such 
that $\chi( g ) \neq 0$. Since $\Res_M^G( \chi )$ is $\sigma$-invariant, the 
claim follows from Lemma~\ref{FixingReduction}.
\end{prf}

\medskip

\noindent
We will need the following variant of Brauer's permutation lemma.

\begin{lem}
\label{BPL}
Let $m \leq n$ be positive integers, $U \in \mathcal{K}^{m \times n}$ of
full rank, $M \in \GL_n( \mathcal{K} )$. Suppose that there are permutation
matrices $P, P', Q$ of the appropriate sizes such that $UQ = PU$ and
$MUQ = P'MU$. Then $P' = MPM^{-1}$, and thus $P$ and $P'$ have the same trace.
\end{lem}
\begin{prf}
We have $P'MU = MUQ = MPU$, and hence $P'M = MP$, as~$U$ has full rank.
\end{prf}

\medskip
\noindent 
We indicate an application of this lemma.

\begin{lem}
\label{OrbitCount}
Let $\sigma$ be an automorphism of~$G$ and let~$B$ be a $\sigma$-stable union 
of $\ell$-blocks of~$G$. Suppose that $\mathcal{U}$ is a $\sigma$-stable set
of $\mathcal{K}$-valued class functions of~$G$ with $|\mathcal{U}| = |\IBr(B)|$, 
such that $\check{\mathcal{U}} := \{ \check{\theta} \mid \theta \in \mathcal{U} \}$ 
is a $\mathcal{K}$-basis of $\mathcal{K} \otimes_{\mathbb{Z}} 
\mathbb{Z}[\IBr( B )]$. Then the number of~$\sigma$-stable elements of $\IBr(B)$ 
equals the number of~$\sigma$-stable elements of~$\mathcal{U}$.
\end{lem}
\begin{prf}
Write~$m$ and~$n$ for the number of irreducible Brauer characters in~$B$ and
the number of $\ell$-regular classes of~$G$, respectively. Let 
$U \in \mathcal{K}^{m \times n}$ denote the character table of 
$\check{\mathcal{U}}$ and let~$Q$ denote the $(n \times n)$-permutation matrix 
arising from the permutation of~$\sigma$ on the set of $\ell$-regular classes. 
As~$\sigma$ fixes $\mathcal{U}$, hence~$\mathcal{\check{U}}$, there is a 
permutation matrix~$P$ such that $PU = UQ$. By assumption, there is 
$M \in \GL_m( \mathcal{K} )$ such that $MU$ is the Brauer character table 
of~$B$. In particular,~$U$ has full rank. As~$\sigma$ fixes~$B$, it stabilizes 
$\IBr(B)$, and thus there is a permutation matrix~$P'$ such that $P'MU = MUQ$. 
By Lemma~\ref{BPL}, the permutation matrices~$P$ and~$P'$ have the same trace, 
and hence~$\sigma$ has the same number of fixed points on~$\check{U}$
and on $\IBr(B)$. As~$U$ has full rank, it does not contain duplicate rows. 
Hence $\theta \in \mathcal{U}$ is fixed by~$\sigma$, if and only 
if~$\check{\theta}$ is fixed by $\sigma$. This proves our claim.
\end{prf}

\addtocounter{subsection}{5}
\subsection{Central products}\label{CentralProducts}
Let $H_1, H_2 \leq G$ with $[H_1,H_2] = 1$, put $Z := H_1 \cap H_2$ and
$H := H_1H_2$. Then $Z \leq Z( H )$ and $H$ is a central product of~$H_1$
and~$H_2$ over~$Z$, written as $H = H_1 \circ_Z H_2$. Let $U_i \leq H_i$
with $Z \leq U_i$ for $i = 1, 2$. Then $U := U_1U_2 = U_1 \circ_Z U_2$
and $N_H( U ) = N_{H_1}( U_1 )N_{H_2}( U_2 ) = N_{H_1}( U_1 ) \circ_Z 
N_{H_2}( U_2 )$. If $H_1$ is abelian, we also have $C_H (U ) = 
H_1C_{H_2}( U_2 ) = H_1 \circ_Z C_{H_2}( U_2 )$ and thus $\Out_H( U ) 
\cong \Out_{H_2}( U_2 )$.
%%%%%%%%%%%%%%%%%%%%%%%%%%%%%%%%%%%%%%%%%%%%%%%%%%%%%%%%%%%%%%%%%%%%%%%%%%%%%%%%
%%%%%%%%%%%%%%%%%%%%%%%%%%%%%%%%%%%%%%%%%%%%%%%%%%%%%%%%%%%%%%%%%%%%%%%%%%%%%%%%
%%
%% For details see Sheets A, B, C of 14.08.2019.
%%
%%%%%%%%%%%%%%%%%%%%%%%%%%%%%%%%%%%%%%%%%%%%%%%%%%%%%%%%%%%%%%%%%%%%%%%%%%%%%%%%
%%%%%%%%%%%%%%%%%%%%%%%%%%%%%%%%%%%%%%%%%%%%%%%%%%%%%%%%%%%%%%%%%%%%%%%%%%%%%%%%

Every $\chi \in \Irr(H)$ can uniquely be 
written as $\chi = \chi_1 \chi_2$ for $\chi_i \in \Irr( H_i )$, $i = 1, 2$
with $\chi_2(1)\Res^{H_1}_Z( \chi_1 ) = \chi_1(1)\Res^{H_2}_Z( \chi_2 )$. 
(Under the latter 
assumption, the product $\chi_1\chi_2 : H \rightarrow \mathcal{K}, h_1h_2 
\mapsto \chi_1( h_1 )\chi_2( h_2 )$ is well defined.) If $U_i \unlhd H_i$ with
$Z \leq U_i$ for $i = 1, 2$, and $\vartheta_i \in \Irr( U_i )$ for $i = 1, 2$
such that $\vartheta := \vartheta_1 \vartheta_2 \in \Irr( U_1 \circ_Z U_2 )$ is 
invariant in $H = H_1 \circ_Z H_2$, then each $\vartheta_i$ is invariant 
in~$H_i$ for $i = 1, 2$; moreover, $\vartheta$ extends to~$H$, if and only if
$\vartheta_i$ extends to $H_i$ for $i = 1, 2$.
%%%%%%%%%%%%%%%%%%%%%%%%%%%%%%%%%%%%%%%%%%%%%%%%%%%%%%%%%%%%%%%%%%%%%%%%%%%%%%%%
%%%%%%%%%%%%%%%%%%%%%%%%%%%%%%%%%%%%%%%%%%%%%%%%%%%%%%%%%%%%%%%%%%%%%%%%%%%%%%%%
%%
%% For details see Sheets 1, 2 of 13.11.2019.
%%
%%%%%%%%%%%%%%%%%%%%%%%%%%%%%%%%%%%%%%%%%%%%%%%%%%%%%%%%%%%%%%%%%%%%%%%%%%%%%%%%
%%%%%%%%%%%%%%%%%%%%%%%%%%%%%%%%%%%%%%%%%%%%%%%%%%%%%%%%%%%%%%%%%%%%%%%%%%%%%%%%

\subsection{Blocks and weights}
\label{BlocksAndWeights}
For easier reference we summarize a few well known results from modular 
representation theory of finite groups, thereby introducing our notation.

An $\ell$-subgroup~$R$ of~$G$ is {\it radical} if $O_\ell(N_G(R)) = R$.
We denote by~$\cR_\ell(G)$ the set of all radical $\ell$-subgroups of~$G$, and
by~$\cR_\ell(G)/G$ the set of $G$-conjugacy classes of $\cR_\ell(G)$.

Let~$b$ be an $\ell$-block of~$G$. We write $D(b)$ for the set of defect groups 
of~$b$. We will use the concept of Brauer pairs in the following. An excellent
reference for this notion, originally introduced by Alperin and Brou{\'e} under 
the name of subpairs in \cite{AlpBro}, is \cite[Sections~$2$,~$3$]{Kessar}. We
will always implicitly assume that the Brauer pairs are defined with respect 
to~$\ell$, or, more precisely, with respect to the field~$\overline{\mathcal{O}}$. Let
$(R,b_R)$ be a {\it Brauer pair}, i.e.~$R$ is an $\ell$-subgroup of~$G$ 
and~$b_R$ is an $\ell$-block of~$C_G(R)$. For every subgroup $Q \leq R$, there 
is a unique Brauer pair $(Q,b_Q)$ with $(Q,b_Q) \leq (R,b_R)$. Recall 
that~$(R,b_R)$ is called a $b$-Brauer pair, if $(\{1\},b) \leq (R,b_R)$.  
Let $(D, b_D)$ denote a maximal $b$-Brauer pair. Then $Z(D) \in D(b_D)$ and 
$D \in D(b)$. We call $(R,b_R)$ {\it centric}, if $Z(R) \in D(b_R)$. 

Let $(Q,\varphi)$ be a {\it weight} of~$G$, i.e.~$Q$ is an $\ell$-subgroup 
of~$G$ and~$\varphi$ is an $\ell$-defect~$0$ character of $N_G(Q)/Q$. Let~$b_Q$ 
be a block of~$C_G(Q)$ covered by the block of $N_G(Q)$ containing (the 
inflation of)~$\varphi$. Then $(Q,b_Q)$ is a centric Brauer pair and 
$Q \in \cR_\ell(G)$. If~$(Q,b_Q)$ 
is a $b$-Brauer pair, we say that $(Q,\varphi)$ is a $b$-weight. In this case,
we also have
\begin{equation}
\label{CenterContainment}
Z(D) \leq Z(R) \leq R \leq D 
\end{equation}
for some conjugate~$R$ of~$Q$; see, 
e.g.\ \cite[Chapter~$5$, Theorem~$5.21$]{NaTsu}. This fact will be used 
frequently in the following. If~$Q$ is an $\ell$-subgroup such that there 
exists a ($b$-)weight~$(Q,\varphi)$, we call~$Q$ a ($b$-){\it weight subgroup}.

Let $R \in \R_\ell(G)$ and let $b_R$ be an $\ell$-block of~$C_G(R)$. We write
$N_G( R, b_R )$ for the stabilizer of~$b_R$ in~$N_G( R )$ and put
$$
\Out_G(R, b_R) := N_G(R, b_R)/C_G(R)R.
$$
Recall that $\Out_G( D, b_D )$ is an $\ell'$-group if, as above,~$D$ is a defect
group of the $\ell$-block~$b$ of~$G$ and $(D,b_D)$ is a maximal $b$-Brauer pair.

Now assume that $(R,b_R)$ is centric. In this case, we denote by $\theta_R \in 
\Irr( C_G( R ) )$ the canonical character of~$b_R$, i.e.\ the unique ordinary 
irreducible character in~$b_R$ with $Z(R)$ in its kernel. We may and will also
view~$\theta_R$ as a character of $C_G( R )R$ via inflation over~$R$. Write 
$N_G( R, \theta_R )$ for the stabilizer of~$\theta_R$ in~$N_G( R )$. This notion 
is independent of whether we view~$\theta_R$ as a character of $C_G(R)$ or as 
one of $C_G(R)R$. Clearly, $N_G( R, \theta_R) = N_G( R, b_R )$ and thus
$\Out_G( R, b_R ) = N_G( R, \theta_R)/C_G(R)R$.  Put
\[
\begin{array}{l}
\Irr^0(N_G( R, \theta_R) \mid \theta_R) := \\
\{ \zeta \in \Irr(N_G( R, \theta_R)) \mid  
\zeta(1)_\ell=|N_G( R, \theta_R){:} R|_\ell, \hbox {  and  }
\zeta \hbox { covers } \theta_R \},
\end{array}
\]
and
\begin{equation}
\label{NumberOfWeights0}
\cW(R, b_R) := |\Irr^0(N_G( R, \theta_R) \mid \theta_R)|.
\end{equation}
By Schur's theory of projective characters, the fact that $\theta_R$ is 
invariant under~$N_G( R, \theta_R)$ yields an element $\alpha \in 
H^2( \Out_G( R, \theta_R ), \mathcal{K}^\times )$, called the 
{\em K{\"u}ls\-ham\-mer-Puig} class associated to the centric Brauer pair 
$(R,b_R)$. If $\mathcal{K}_\alpha( \Out_G( R, \theta_R ) )$
denotes the corresponding twisted group algebra, we have
$$\cW( R, b_R ) = |\Irr^0( \mathcal{K}_\alpha \Out_G( R, \theta_R ) )|.$$
In particular, if $\theta_R$, viewed as a character of~$C_G(R)R$, extends to 
$N_G( R, \theta_R)$, then 
\begin{equation}
\label{NumberOfWeights}
\cW(R, b_R) = |\Irr^0( \Out_G( R, b_R ) )|
\end{equation} 
by Clifford theory.

If $(R,b_R)$ is a $b$-Brauer pair and $\zeta \in 
\Irr^0(N_G( R, \theta_R) \mid \theta_R)$, then 
$\Ind_{N_G( R,  \theta_R )}^{N_G( R )}( \zeta )$ is a $b$-weight. By choosing a
maximal $b$-Brauer pair $(D,b_D)$ and applying this construction to all
centric $b$-Brauer pairs $(R,b_R)$ with $(R,b_R) \leq (D,b_D)$, we obtain a set 
of representatives for the $G$-conjugacy classes of $b$-weights, the number of
which is denoted by $\cW(b)$  in the following. 

\addtocounter{thm}{2}
\begin{prop}
\label{WeightCount}
Fix an $\ell$-block~$b$ of~$G$ and a maximal $b$-Brauer pair $(D,b_D)$.
Then
\begin{equation}
\label{SumsOfWeights}
\cW(b) = \sum_{R} \cW(R, b_R), 
\end{equation}
where $R \in \cR_\ell(G)$ runs through a set of $G$-conjugacy class 
representatives satisfying $R \leq D$ and the unique Brauer pair 
$(R, b_R)$ with $(R,b_R) \leq (D,b_D)$ is centric.

Suppose that~$D$ is abelian. Then the only summand 
in {\rm Equation~(\ref{SumsOfWeights})} is the one for $R = D$, and thus 
$\cW(b) = \cW(D, b_D)$.
\end{prop}
\begin{prf}
The first statement is due to Alperin and Fong; see the discussion in 
\cite[Page~$3$]{AlpFo}. 
%It is also contained in 
%\cite[Proposition~$5.4$]{Kessar}, although in a slightly different language.
If~$D$ is abelian, then $R = D$ by Equation~(\ref{CenterContainment}), and 
the second statement follows. 
\end{prf}

\smallskip

\noindent
Assume the situation and notation of Proposition~\ref{WeightCount} and that~$D$ 
is abelian. Then~$\theta_D$ extends to $N_G( D, \theta_D )$ if 
$\Out_G( D, b_D )$ is cyclic, or if~$\theta_D$ is linear and~$C_G( D )$ has a 
complement in $N_G( D, \theta_D )$.

Alperin's weight conjecture for~$b$ postulates the equality 
$\cW( b ) = |\IBr( b )|$.

\addtocounter{subsection}{1}
\subsection{The inductive blockwise Alperin weight condition}
\label{InductiveBlockwiseAlperin}

The following hypotheses constitute a simplified version of Conditions~(i) 
and~(ii) of \cite[Definition~$3.2$]{KS16}, adapted to our purpose. The latter,
in turn, is a specialization of \cite[Definitions~$4.1$, $5.17$]{Spaeth13}.

\addtocounter{thm}{1}
\begin{hypo}\label{defBWC} 
{\rm
Let~$\ell$ be a prime,~$G$ a finite non-abelian simple group with trivial Schur 
multiplier, and let~$b$ be an $\ell$-block of~$G$. Write $N_{\Aut(G)}( b )$ for 
the stabilizer in $\Aut( G )$ of~$b$. Assume that the following conditions are 
satisfied.

\smallskip

\noindent
(i) For every $Q \in \cR_\ell(G)$, there exists a subset $\ibr(b\midd Q)$ 
of $\ibr(b)$ satisfying the following conditions.

\begin{itemize}
\item[(1)] If $\alpha\in N_{\Aut(G)}(b)$, then $\ibr(b\midd Q)^\alpha=
\ibr(b\midd Q^\alpha)$.

\item[(2)] $\ibr(b)$ is the disjoint union of all $\ibr(b\midd Q)$, where~$Q$ 
runs through some set of representatives of~$\cR_\ell(G)/G$.
\end{itemize} 

\noindent
(ii) For every $Q \in \cR_\ell(G)$, there is a bijection
\[
\Omega_Q^G\colon \ibr(b\midd Q)\to \Irr^0(N_G(Q)/Q, b),
\] 
such that for all $\alpha\in N_{\Aut(G)}(b)$ and $\varphi\in\ibr(b\midd Q)$ we 
have that $\Omega_Q^G(\varphi)^\alpha=\Omega_{Q^\alpha}^G(\varphi^\alpha)$. Here,
$\Irr^0(N_G(Q)/Q, b)$ consists of all characters $\varphi\in \Irr^0(N_G(Q)/Q)$
such that $(Q, \varphi)$ is a $b$-weight.
}
\end{hypo}

\begin{rem}\label{remBWC}
{\rm
Let~$G$ and~$\ell$ be as in Hypothesis~\ref{defBWC}. 
Suppose in addition that $\Out(G)$ cyclic. If Hypothesis~\ref{defBWC} is 
satisfied for every $\ell$-block~$b$ of~$G$, then~$G$ satisfies the 
\textit{inductive blockwise Alperin weight condition at the prime~$\ell$}. 
This is proved in 
\cite[Lemma~$6.1$]{Spaeth13} under the stronger condition that~$G$ is AWC-good, 
but the proof of this lemma only uses Hypothesis~\ref{defBWC}. We are grateful 
to Britta Sp{\"a}th for pointing out the relevance of 
\cite[Lemma~$6.1$]{Spaeth13} for this reduction, and for clarifying remarks 
regarding this issue.\hfill{$\Box$}
}
\end{rem}

\noindent 
We formulate a set of conditions which simplifies the verification of 
Hypothesis~\ref{defBWC}.

\begin{hypo}
\label{defBWC_Simplified}
{\rm
Let~$G$,~$\ell$ and~$b$ be as in Hypothesis~\ref{defBWC}. Put $A :=
N_{\Aut(G)}( b )$. Let $Q_1 , \ldots , Q_n$ denote a set of
representatives for the $G$-conjugacy classes in the set
$ \{ Q \in \cR_\ell( G ) \mid \Irr^0( N_G( Q )/Q, b ) \neq \emptyset \}$.

(1) There are pairwise disjoint subsets $\IBr( b \midd Q_i )$, 
$1 \leq i \leq n$ of $\IBr( b )$ such that
$$\bigcup_{i = 1}^n \IBr(b \midd Q_i ) = \IBr( b ),$$
and there are bijections 
$$\Omega_i : \IBr( b \midd Q_i ) \rightarrow \Irr^0( N_G( Q_i )/Q_i, b )$$
for $1 \leq i \leq n$.

(2) Suppose that (1) holds and, in addition, that
$$\IBr( b \midd Q_i )^\alpha = \IBr( b \midd Q_j )\quad\text{and\ }
\Omega_j( \varphi^\alpha ) = \Omega_i( \varphi )^\alpha$$
for all $1 \leq i, j \leq n$, all $\alpha \in A$ with $Q_i^\alpha = Q_j$ and
all $\varphi \in \IBr( b \midd Q_i )$.\hfill{$\Box$}
}
\end{hypo}

\begin{rem}
\label{remBWC_Simplified}
{\rm 
Let~$G$,~$\ell$ and~$b$ be as in Hypothesis~\ref{defBWC}, and put 
$A := N_{\Aut(G)}( b )$.

(a) Assume that both conditions of Hypothesis~\ref{defBWC_Simplified} hold. Let
$Q \in \cR_\ell( G )$ be $G$-conjugate to $Q_i$ for some $1 \leq i \leq n$, say 
$Q = Q_i^g$ for some $g \in G$. Define 
$\IBr( b \midd Q ) := \IBr( b \midd Q_i )^g$, and 
$\Omega_Q^G : \IBr( b \midd Q ) \rightarrow \Irr^0( N_G( Q )/Q, b ),
\varphi^g \mapsto \Omega_i( \varphi )^g$ for $\varphi \in \IBr(b \midd Q_i )$.
Then $\IBr( b \midd Q )$ and $\Omega_Q^G$ are well-defined, and the collection 
of theses sets and maps satisfies Hypothesis~\ref{defBWC}.
%%%%%%%%%%%%%%%%%%%%%%%%%%%%%%%%%%%%%%%%%%%%%%%%%%%%%%%%%%%%%%%%%%%%%%%%%%%%%%%%
%%%%%%%%%%%%%%%%%%%%%%%%%%%%%%%%%%%%%%%%%%%%%%%%%%%%%%%%%%%%%%%%%%%%%%%%%%%%%%%%
%%
%% See sheets 1 - 3 of 15.10.2020.
%%
%%%%%%%%%%%%%%%%%%%%%%%%%%%%%%%%%%%%%%%%%%%%%%%%%%%%%%%%%%%%%%%%%%%%%%%%%%%%%%%%
%%%%%%%%%%%%%%%%%%%%%%%%%%%%%%%%%%%%%%%%%%%%%%%%%%%%%%%%%%%%%%%%%%%%%%%%%%%%%%%%

(b) Condition~(1) of~\ref{defBWC_Simplified} is satisfied if and only if the 
Alperin weight conjecture holds for~$b$. Suppose this is the case. If the 
$G$-orbit of some~$Q_i$ equals its $A$-orbit, e.g.\ if $n = 1$ which occurs 
if~$b$ has abelian defect, then Condition~(2) amounts to an equivariance 
condition with respect to the action of $N_A( Q_i )$.}\hfill{$\Box$}
\end{rem}

\medskip
\noindent
In the case of our interest, $G = F_4(q)$, every $G$-orbit of every~$Q_i$ 
is $A$-invariant, except for one instance where the $A$-orbit of some~$Q_i$ 
splits into two $G$-orbits.

%\markboth{Finite reductive groups}
\section{Finite reductive groups}
\label{FiniteReductiveGroups}
We continue by recalling some basic concepts and results from the theory of 
finite reductive groups and their representations, to the extent needed later 
on. 

\subsection{Notation}
Let $p$ be a prime number and let $\mathbb{F}$ denote an algebraic closure of
the finite field with~$p$ elements. Let~$\mathbf{G}$ be a connected reductive 
algebraic group over~$\mathbb{F}$. We also
let~$F$ denote a Steinberg endomorphism of~$\mathbf{G}$. Let~$\mathbf{H}$ be a
closed subgroup of~$\mathbf{G}$. We then write $\mathbf{H}^\circ$ for
the connected component of~$\mathbf{H}$ containing~$1$, and if~$\mathbf{H}$ is 
$F$-stable we write $H := \mathbf{H}^F$ for the 
finite group of $F$-fixed points of~$\mathbf{H}$. Thus $G = \mathbf{G}^F$ is a
finite reductive group. We also let~$\mathbf{G}^*$ denote a group dual 
to~$\mathbf{G}$ (with respect to a fixed $F$-stable maximally split torus 
of~$\mathbf{G}$), endowed with a dual Steinberg endomorphism~$F^*$.

We say that~$\mathbf{G}$ is of classical type, if every minimal $F$-stable
semisimple component~$\mathbf{H}$ of $[\mathbf{G},\mathbf{G}]$ is of Dynkin
type~$A$,~$B$,~$C$ or~$D$, and if $\mathbf{H}^F$ is not isomorphic to
${^3\!D}_4(q)$ for some power~$q$ of~$p$.

Finally, let~$\ell$ be a prime different from~$p$.

\subsection{Recollections and preliminary results}
If~$\mathbf{G}$ has connected center, the centralizers of semisimple elements in 
$\mathbf{G}^*$ are connected. In this case, two semisimple elements of~$G^*$ are 
conjugate in~$G^*$, if and only if they are conjugate in~$\mathbf{G}^*$.

An $F$-stable Levi subgroup of~$\mathbf{G}$ is called a {\it regular subgroup} 
of~$\mathbf{G}$. The regular subgroups of~$\mathbf{G}$ are exactly the 
centralizers of $F$-stable tori. If~$\mathbf{G}$ has connected center, so has 
any regular subgroup of~$\mathbf{G}$;
see \cite[Proposition~$8.1.4$]{C2}. If $\mathbf{L} \leq \mathbf{G}$ is  
regular in~$\mathbf{G}$, and $\mathbf{M} \leq \mathbf{L}$ is regular 
in~$\mathbf{L}$, then~$\mathbf{M}$ is regular in~$\mathbf{G}$.
If~$\mathbf{M} \leq \mathbf{G}$ is a regular subgroup, we write
$W_{\mathbf{G}}( \mathbf{M} ) := N_{\mathbf{G}}( \mathbf{M} )/\mathbf{M}$
for the relative Weyl group of~$\mathbf{M}$ in~$\mathbf{G}$. Then
$W_{\mathbf{G}}( \mathbf{M} )^F = N_G( \mathbf{M} )/M$.

The following lemma will be used to identify centralizers
of $3$-elements in $F_4(q)$.

\addtocounter{thm}{2}
\begin{lem}
\label{NonRegularCentralizers}
Assume that $Z( \mathbf{G} )$ is connected.
Let $\mathbf{H}$ be an $F$-stable closed subgroup of~$\mathbf{G}$ such that
$\mathbf{H}$ is reductive and $Z( \mathbf{H} )$ is not connected.
Let $s \in Z( \mathbf{H} )^F$ be of order coprime to
$|( Z( \mathbf{H} )/ Z( \mathbf{H} )^\circ )^F|$.
Then $\mathbf{H} \lneq C_{\mathbf{G}}( s )$.
\end{lem}
\begin{prf}
Our assumption implies that $s \in (Z( \mathbf{H} )^\circ)^F \leq
Z( \mathbf{H} )^\circ$.
As $\mathbf{H}^\circ$ is reductive, $Z( \mathbf{H} )^\circ$ is a torus, and thus
$C_{\mathbf{G}}( Z( \mathbf{H} )^\circ )$ is a regular subgroup of~$\mathbf{G}$,
hence has connected center. It follows that $\mathbf{H} \lneq
C_{\mathbf{G}}( Z( \mathbf{H} )^\circ ) \leq C_{\mathbf{G}}( s )$, proving
our claim.
\end{prf}

\medskip
\noindent
We will also need the following slight generalization of
\cite[Proposition~$4.2$]{GeHi1}.

\begin{prop}
\label{Proposition42}
Let $\mathcal{S}_\ell( G )$ and $\mathcal{S}_\ell( G^* )$ denote sets of
representatives for the conjugacy classes of $\ell$-elements of~$G$,
respectively~$G^*$.
Suppose that $\ell$ does not divide the determinant of the Cartan matrix
of the root system of~$\mathbf{G}$ and that centralizers of $\ell$-elements
in~$\mathbf{G}$ and $\mathbf{G}^*$ are connected.

Then there is a bijection $\mathcal{S}_\ell( G ) \rightarrow
\mathcal{S}_\ell( G^* )$, $t \mapsto t'$ such that the following holds.
If $C_{\mathbf{G}}( t )$ is a regular subgroup of~$\mathbf{G}$, then
$C_{\mathbf{G}^*}( t' )$ is a regular subgroup of~$\mathbf{G}^*$, and
there is an $F$-equivariant isomorphism $C_{\mathbf{G}}( t ) \rightarrow
C_{\mathbf{G}^*}( t' )^*$.
\end{prop}
\begin{prf}
This is contained in the proof of \cite[Proposition~$4.2$]{GeHi1}.
\end{prf}

\medskip
\noindent 
Assume the hypotheses and the notation of the above proposition. 
If~$C_{\mathbf{G}}( t )$ is not a regular subgroup of~$\mathbf{G}$, its dual 
group is not necessarily the centralizer of a semisimple element. An example
is provided by $\mathbf{G} = G_2(\mathbb{F})$, the simple group of type~$G_2$
over~$\mathbb{F}$ when $p \neq 3$. Then $\mathbf{G} \cong \mathbf{G}^*$, and
the hypotheses of Proposition~\ref{Proposition42} are satisfied. There is an
element~$s$ of order~$3$ whose centralizer is isomorphic to $\SL_3( \mathbb{F} )$;
see \cite[Table~$4.7.1$]{Gor}. If~$F$ is such that $G = \mathbf{G}^F = G_2(q)$
with $3 \mid q - 1$, then~$G$ contains a representative~$t$ of the 
$\mathbf{G}$-conjugacy class of~$s$, so that $C_G( t ) \cong \SL_3(q)$.
However, $C_{\mathbf{G}}( t )^* \cong \PGL_3( \mathbb{F} )$ has trivial center.

\addtocounter{subsection}{2}
\subsection{Character groups and cocharacter groups of maximal tori}
\label{CharacterAndCocharacterGroups}
Let~$\mathbf{T}$ be an $F$-stable maximal torus of~$\mathbf{G}$. Then 
$X( \mathbf{T} ) := \Hom( \mathbf{T}, \mathbb{F}^* )$ and $Y( \mathbf{T} ) := 
\Hom( \mathbb{F}^*, \mathbf{T} )$ denote the \textit{character group} and the
\textit{cocharacter group} of~$\mathbf{T}$, respectively. By 
$\langle \ , \ \rangle$ we denote the natural pairing between $X( \mathbf{T} )$ 
and $Y( \mathbf{T} )$. We choose an isomorphism
\begin{equation}
\label{iota}
\iota : \mathbb{F}^* \rightarrow \mathbb{Q}_{p'}/\mathbb{Z}
\end{equation}
(see \cite[Proposition~$3.1.3$]{C2}).
This gives rise to an isomorphism of abelian groups $Y( \mathbf{T} ) \otimes \mathbb{F}^*
\cong Y( \mathbf{T} ) \otimes \mathbb{Q}_{p'}/\mathbb{Z}$, and, in turn, to an
isomorphism
\begin{equation}
	\label{YotimesF}
	Y( \mathbf{T} ) \otimes \mathbb{Q}_{p'}/\mathbb{Z} \rightarrow \mathbf{T}
\end{equation}
(see \cite[Proposition~$3.1.2$]{C2}). Under this isomorphism,~$T$ corresponds to
the kernel of $F - 1$ on $Y( \mathbf{T} ) \otimes \mathbb{Q}_{p'}/\mathbb{Z}$.
This yields a further isomorphism
\begin{equation}
	\label{TF}
	Y( \mathbf{T} )/(F - 1)Y( \mathbf{T} ) \rightarrow T
\end{equation}
(see \cite[Proposition~$3.2.2$]{C2} or \cite[Proposition~$11.1.7$(ii)]{DiMi2}). 
Similarly, there is an isomorphism
\begin{equation}
	\label{XotimesF}
	\left( \text{kernel of\ } F - 1 \text{\ on\ } X( \mathbf{T} ) \right)
	\otimes \mathbb{Q}_{p'}/\mathbb{Z} \rightarrow \Irr( T )
\end{equation}
(see \cite[Proposition~$3.2.4$]{C2}). Now let $\chi \in X( \mathbf{T} )$ and
$a \in \mathbb{Q}_{p'}/\mathbb{Z}$ such that 
$\chi \otimes a \in X( \mathbf{T} ) \otimes \mathbb{Q}_{p'}/\mathbb{Z}$
lies in the kernel of $F - 1$, and let $\lambda \in \Irr(T)$ 
correspond to $\chi \otimes a$ under the isomorphism~(\ref{XotimesF}).
Next, let $\gamma \in Y( \mathbf{T} )$, and let $t \in T$ correspond to
$\gamma + (F-1)Y( \mathbf{T} ) \in Y( \mathbf{T} )/(F - 1)Y( \mathbf{T} )$
under the isomorphism~(\ref{TF}). Then
\begin{equation}
	\label{EvaluationOfCharacter}
	\lambda( t ) = \exp( 2\pi\sqrt{-1} \langle \chi, \gamma \rangle a ),
\end{equation}
where $\exp$ is the exponential function of $\mathbb{C}$, and where the $|G|$th 
roots of unity of~$\mathbb{C}$ are identified with the $|G|$th roots of unity 
of $\mathcal{K}$ (recall that $\mathcal{K}$ is our large enough field of
characteristic~$0$ containing the character values).

\subsection{Lusztig induction}
If~$\mathbf{M}$ is a regular subgroup of~$\mathbf{G}$, we write 
$R_{\mathbf{M}}^{\mathbf{G}}$ for the Lusztig induction map from the class 
functions of~$M$ to the class functions of~$G$. Strictly speaking, this map also
depends upon a parabolic subgroup~$\mathbf{P}$ containing~$\mathbf{M}$ as a 
Levi complement, so that we should write 
$R_{\mathbf{M} \leq \mathbf{P}}^{\mathbf{G}}$. We will always implicitly assume 
that the Mackey formula holds for~$\mathbf{G}$, in which case 
$R_{\mathbf{M} \leq \mathbf{P}}^{\mathbf{G}}$ is independent of such~$\mathbf{P}$;
see \cite[Theorem~$3.3.8$]{GeMa}. We therefore omit the~$\mathbf{P}$ from the 
notation. By~\cite{BoMi}, the Mackey formula holds for all connected reductive
groups relevant to this work. Let~$\mathbf{T}$ be an 
$F$-stable maximal torus of~$\mathbf{G}$ and $\theta \in \Irr(T)$. Suppose that
the pair $(\mathbf{T}^*,s)$, where~$\mathbf{T}^*$ is an $F^*$-stable maximal 
torus of~$\mathbf{G}^*$ and $s \in T^*$ corresponds to $(\mathbf{T},\theta)$ 
via duality; see \cite[Proposition~$11.1.16$]{DiMi2}. We then also write 
$R_{\mathbf{T}^*}^{\mathbf{G}}(s)$ for $R_{\mathbf{T}}^{\mathbf{G}}(\theta)$, as 
in \cite[p.~$167$]{DiMi2} or \cite[Remark~$8.22$(i)]{CaEn}.

\addtocounter{thm}{2}
\begin{lem}
\label{GenerationLemma}
Let $s \in G^*$ be a semisimple $\ell'$-element,~$t$ an $\ell$-element in 
$C_{G^*}( s )$ and put $\mathbf{C}^* := C_{\mathbf{G}^*}( st )$. Suppose that 
one of the following conditions is satisfied:

{\rm (a)} The subgroup $\mathbf{C}^* \leq \mathbf{G}^*$ is regular.

{\rm (b)} There is a regular subgroup $\mathbf{M}^* \leq \mathbf{G}^*$ with
$\mathbf{C}^* \leq \mathbf{M}^*$ such that~$\ell$ is good for~$\mathbf{M}$
and $\ell \nmid |Z( \mathbf{M} )/Z( \mathbf{M} )^\circ|$, where~$\mathbf{M}$ 
is a regular subgroup of~$\mathbf{G}$ dual to~$\mathbf{M}^*$.

{\rm (c)} The center of~$\mathbf{G}$ is connected and every unipotent 
character of~$C^*$ is uniform. 

\noindent
Then the class functions $\check{\theta}$ for 
$\theta \in \mathcal{E}( G, st )$ lie in the $\mathbb{Z}$-span of 
$\{ \check{\chi} \mid \chi \in \mathcal{E}( G, s ) \}$ if~{\rm (a)} 
or~{\rm (b)} holds, and in the $\mathbb{Q}$-span of this set, otherwise.
\end{lem}
\begin{prf}
(a) If~$\mathbf{C}^*$ is regular, the claim follows with exactly the same proof 
as that of \cite[Theorem~$3.1$]{GeHi1} (with~$L'$ replaced by~$\mathbf{C}^*$).

(b) Lusztig induction $R_{\mathbf{M}}^{\mathbf{G}}$ establishes a bijection, 
up to a sign, between $\mathcal{E}( M, st )$ and $\mathcal{E}( G, st )$; see 
\cite[($11.4.3$(ii))]{DiMi2}. Let $\theta \in \mathcal{E}( G, st )$. Then 
there is $\mu \in \pm \mathcal{E}( M, st )$ such that $\theta = 
R_{\mathbf{M}}^{\mathbf{G}}( \mu )$. Denote by $\gamma$ the characteristic 
function on the set of $\ell$-regular elements of~$G$ (or any subgroup of~$G$). 
Then $\gamma \cdot \theta = \gamma \cdot R_{\mathbf{M}}^{\mathbf{G}}( \mu ) = 
R_{\mathbf{M}}^{\mathbf{G}}( \gamma \cdot \mu )$; see 
\cite[Proposition~$10.1.6$]{DiMi2}. As~$\ell$ is good for~$\mathbf{M}$ and
does not divide $|Z( \mathbf{M} )/Z( \mathbf{M} )^\circ|$, we find 
elements $\nu_1, \ldots , \nu_d$ in $\mathcal{E}( M, s )$ and integers 
$z_1, \ldots , z_d$ such that $\gamma \cdot \mu = 
\sum_{i=1}^d z_i( \gamma \cdot \nu_i )$; see \cite[Theorem~A]{GeBSII}. It 
follows that $\gamma \cdot \theta = 
\sum_{i=1}^d z_i( \gamma \cdot R_\mathbf{M}^\mathbf{G}( \nu_i ) )$.
As the irreducible constituents of $R_{\mathbf{M}}^{\mathbf{G}}( \nu_i )$ lie in
$\mathcal{E}( G, s )$ for all $1 \leq i \leq d$ (see 
\cite[Proposition~$15.7$]{CaEn}), our claim follows.

(c) If every unipotent character of~$C^*$ is uniform, the same is true for the
elements of $\mathcal{E}( G,  st )$. Indeed, the assumption implies that the
matrix of scalar products between the elements of $\mathcal{E}( C^*, 1 )$ and 
the Deligne--Lusztig characters $R_{\mathbf{T}^*}^{\mathbf{C}^*}( 1_{T*} )$, 
where~$\mathbf{T}^*$ runs through a set of representatives of the 
$C^*$-conjugacy classes of $F^*$-stable maximal tori of~$\mathbf{C}^*$, is 
square and non-singular. By the compatibility of Lusztig's Jordan decomposition 
of characters with Deligne--Lusztig induction (see, e.g.\ 
\cite[Theorem~$15.8$]{CaEn}), this implies that every element 
of~$\mathcal{E}( G,  st )$ is a $\mathbb{Q}$-linear combination of 
Deligne--Lusztig characters $R_{\mathbf{T}^*}^{\mathbf{G}}(st)$ for $F$-stable 
maximal tori~$\mathbf{T}$ of~$\mathbf{G}$ such that~$st$ is contained in a dual 
$\mathbf{T}^*$ of~$\mathbf{T}$. It follows from \cite[Proposition~$10.1.6$]{DiMi2}
that ${R^\mathbf{G}_{\mathbf{T}^*}(st)}$ and
${R^\mathbf{G}_{\mathbf{T}^*}(s)}$ agree on $\ell$-regular elements of~$G$, 
yielding our assertion.
\end{prf}

\medskip

\noindent
Notice that the second condition on~$\ell$ in~(b) above is satisfied if
$Z( \mathbf{G} )$ is connected. 

\addtocounter{subsection}{1}
\subsection{Some major results and their consequences}
We will also use the following combination of major results of 
Bonnaf{\'e}-Dat-Rouquier, Boltje-Perepelitsky, Cabanes-Enguehard and Puig. The 
first of these papers generalizes results of Bonnaf{\'e}-Rouquier \cite{BoRo} 
and Kessar-Malle \cite{KeMa}. 

\addtocounter{thm}{1}
\begin{thm}
\label{BoRoEtAl}
Let $s \in G^*$ be a semisimple $\ell'$-element and assume that 
$C_{\mathbf{G}^*}( s )$ is connected. Put $\mathbf{M}^* := 
C_{\mathbf{G}^*}( Z( C_{\mathbf{G}^*}( s ) )^\circ )$. Then~$\mathbf{M}^*$ is 
a regular subgroup of~$\mathbf{G}^*$ minimal with the property that
$C_{\mathbf{G}^*}( s ) \leq \mathbf{M}^*$. Let~$\mathbf{M}$ be a
regular subgroup of~$\mathbf{G}$ dual to~$\mathbf{M}^*$.

{\rm (a) (Bonnaf{\'e}, Dat, Rouquier \cite[Theorems~$1.1$, $7.7$]{BoDaRo})}
There is a bijection $b \mapsto b'$ between the $\ell$-blocks contained in 
$\mathcal{E}_\ell( G, s )$ and those in $\mathcal{E}_\ell( M, s )$
such that $b$ and $b'$ are Morita equivalent and have a common defect group
$D \leq M$.

{\rm (b) (Bonnaf{\'e}, Dat, Rouquier \cite[Theorems~$1.1$, $7.7$]{BoDaRo}),
Puig, \cite[Theorem~$19.7$]{Puig}, 
Boltje, Perepelitsky~\cite[Theorems~$11.2$,~$13.4$]{BoPe}).} 
Let $b$ and $b'$ be blocks corresponding as in~{\rm (a)} and let $D \leq M$
be a common defect group of~$b$ and~$b'$. Then there is a maximal $b$-Brauer
pair $(D,b_D)$ and a maximal $b'$-Brauer pair $(D,b_D')$, such that the 
following statements hold.

There is a bijection $( R, b_R ) \mapsto ( R, b_R' )$ between the $b$-Brauer 
pairs $( R, b_R ) \leq (D,b_D)$ and the $b'$-Brauer pairs $( R, b_R' ) \leq 
(D,b_D')$ such that $\Out_G( R, b_R ) \cong \Out_M( R, b_R' )$. Moreover, $(R,b_R)$ 
is centric if and only if $(R,b_R')$ is. In the latter case, the 
K{\"u}lshammer-Puig classes $\alpha$ and $\alpha'$ associated to $(R,b_R)$ and 
$(R,b_R')$, respectively, correspond to each other under the isomorphism of
cohomology groups induced by the isomorphism $\Out_G( R, b_R ) \cong 
\Out_M( R, b_R' )$. In particular, $\cW_G( R, b_R ) = \cW_M( R, b_R' )$.

{\rm (c) (Cabanes, Enguehard \cite{CaEn97})} If $C_{\mathbf{G}^*}( s ) \lneq 
\mathbf{M}^*$, assume that~$\mathbf{M}$ has connected center, is of classical 
type and that~$\ell$ is odd. Then there is a bijection between the $\ell$-blocks 
contained in $\mathcal{E}_\ell( M, s )$ and those in 
$\mathcal{E}_\ell( C_{G^*}( s ), 1 )$. Moreover, corresponding blocks have the 
same number of irreducible $\ell$-modular characters and isomorphic defect 
groups. 
\end{thm}
\begin{prf}
Part~(a) is directly taken from the cited reference, so it suffices to prove~(b)
and~(c).

(b) The Morita equivalence in~(a) is induced by a splendid tilting complex~$C$;
see the proof of \cite[Theorem~$7.7$]{BoDaRo}. Let $(D,b_D)$ be a maximal 
$b$-Brauer pair. By~\cite[Theorem~$19.7$]{Puig}, there is a maximal $b'$-Brauer 
pair and an equivalence between the fusion systems associated to $(D,b_D)$ 
and $(D,b_D')$.  This yields all but the last statement of~(b). 

To obtain the last assertion, we follow~\cite{BoPe} (which also yields the other
statements of~(b)). The indecomposable direct summands of the complex~$C$ have
vertices contained in the diagonal $\Delta(D) \leq G \times M$ by 
\cite[Corollary~$3.8$]{BoDaRo}. Hence $C$ yields an $\ell$-permutation 
equivalence between~$b$ and~$b'$ in the sense of \cite[Definition~$9.8$]{BoPe}. 
The existence of $(D,b_D)$ and $(D,b_D')$ follows from 
\cite[Theorem~$10.11$]{BoPe}. The equivalence between the fusion systems 
associated to $(D,b_D)$ and $(D,b_D')$ follows from \cite[Theorem~$11.2$]{BoPe}. 
The correspondence of the K{\"u}lshammer-Puig classes $\alpha$ and $\alpha'$ 
associated to $(R,b_R)$ and $(R,b_R')$, respectively, follows from 
\cite[Theorem~$13.4$]{BoPe}. This implies that the twisted group algebras 
$\mathcal{K}_{\alpha} \Out_G( R, b_R )$ and 
$\mathcal{K}_{\alpha'} \Out_M( R, b_R' )$ are isomorphic. As $\cW_G( R, b_R ) = 
|\Irr^0( \mathcal{K}_{\alpha} \Out_G( R, b_R ) )|$, the last claim follows.

(c) If $C_{\mathbf{G}^*}( s ) = \mathbf{M}^*$, the result follows from~(a), as
the blocks of $\mathcal{E}_\ell( C_{G^*}( s ), s )$ 
and $\mathcal{E}_\ell( C_{G^*}( s ), 1 )$ are Morita equivalent. The Morita
equivalence is given by tensoring with a linear character (see, e.g.\ 
\cite[Proposition~$11.4.8$(ii)]{DiMi2}) and preserves defect groups.

Now assume that $C_{\mathbf{G}^*}( s ) \lneq \mathbf{M}^*$. Then~$\mathbf{M}^*$
is of classical type and~$\ell$ is odd by assumption. As the Mackey formula 
holds for regular subgroups of~$\mathbf{M}$ by~\cite{BoMi}, the Jordan 
decomposition of characters for~$M$ is compatible with Lusztig induction by 
\cite[Proposition~$5.3$]{En13}; see also \cite[Theorem~$4.7.2$]{GeMa}. 
The classification of $\ell$-blocks of~$M$ 
by Cabanes and Enguehard in \cite[Theorem~$3.3$]{CaEn97} then implies that there 
is a bijection $B \mapsto b$ between the $\ell$-blocks~$B$ in 
$\mathcal{E}_\ell( M, s )$ and the $\ell$-blocks~$b$ in 
$\mathcal{E}_\ell( C_{M^*}( s ), 1 ) = \mathcal{E}_\ell( C_{G^*}( s ), 1 )$
such that $|\Irr(B) \cap \mathcal{E}_\ell( M, s )|
= |\Irr(b) \cap \mathcal{E}_\ell( C_{M^*}( s ), 1 )|$.
By \cite[Theorem~$5.1$]{GeHi1}, this implies that~$B$ and~$b$ have the same 
number of irreducible $\ell$-modular characters, as $\mathbf{M}$ has connected 
center. Finally, by \cite[Remark~$3.6$]{CaEn97}, the bijection $B \mapsto b$
preserves isomorphism types of defect groups.
\end{prf}

\medskip

\noindent 
The following generalizations of a result of Brou{\'e} and Michel
\cite[Th{\'e}o\-r{\`e}\-me~$3.2$]{BrouMi} are extremely useful in identifying 
canonical characters of blocks. It holds under much weaker hypotheses, but its
derivation from \cite[Theorem~$4.14$]{BoDaRo} is simplified under the stronger
conditions given.

\begin{lem}[Bonnaf{\'e}, Dat, Rouquier \cite{BoDaRo}]
\label{ConjugateSemisimpleLabels}
Let $s \in G^*$ be a semisimple $\ell'$-element such that 
$C_{\mathbf{G}^*}( s )$ is connected and let $b \subseteq 
\mathcal{E}_\ell( G, s )$ be an $\ell$-block of~$G$. Let $(R,b_R)$ be a 
$b$-Brauer pair with $R$ abelian 
such that $C_{\mathbf{G}}( R )$ and $C_{\mathbf{G}}( y )$ are regular subgroups 
of $\mathbf{G}$ for all $y \in R$. Let $C_{\mathbf{G}}( R )^*$ denote a regular 
subgroup of~$\mathbf{G}^*$ dual to $C_{\mathbf{G}}( R )$.

Let $\theta \in \Irr(b_R)$ and suppose that $\theta \in 
\mathcal{E}( C_G(R), t )$ for some semisimple $\ell'$-element 
$t \in (C_{\mathbf{G}}( R )^*)^{F^*}$. Then $t$ is conjugate to $s$ 
in~$G^*$. If $C_{\mathbf{G}}( R )^* = C_{\mathbf{G}^*}( R^\diamond )$
for some subgroup $R^\diamond \leq G^*$, there is a $G^*$ conjugate
$R^\dagger \leq G^*$ such that $R^\dagger \leq C_{G^*}( s )$.
\end{lem}
\begin{prf}
This follows from \cite[Theorem~$4.14$]{BoDaRo}.
\end{prf}
%%%%%%%%%%%%%%%%%%%%%%%%%%%%%%%%%%%%%%%%%%%%%%%%%%%%%%%%%%%%%%%%%%%%%%%%%%%%%%%%
%%%%%%%%%%%%%%%%%%%%%%%%%%%%%%%%%%%%%%%%%%%%%%%%%%%%%%%%%%%%%%%%%%%%%%%%%%%%%%%%
%%
%% See also Lemma 2.9 of alperinF4_29032019 for a direct proof.
%%
%%%%%%%%%%%%%%%%%%%%%%%%%%%%%%%%%%%%%%%%%%%%%%%%%%%%%%%%%%%%%%%%%%%%%%%%%%%%%%%%
%%%%%%%%%%%%%%%%%%%%%%%%%%%%%%%%%%%%%%%%%%%%%%%%%%%%%%%%%%%%%%%%%%%%%%%%%%%%%%%%

\medskip
\noindent 
Lemma~\ref{ConjugateSemisimpleLabels} shows the relevance of regular 
centralizers of $\ell$-sub\-groups. We give a criterion for this property 
to hold.

\begin{lem}
\label{CentralizersAreRegular}
Assume that $\ell \nmid |Z( \mathbf{G}^* )/Z( \mathbf{G}^* )^\circ|$. Let
$R \leq G$ be an abelian
$\ell$-subgroup. Suppose that there is $z \in R$ such that $\mathbf{L} :=
C_{\mathbf{G}}(z)$ is a regular subgroup of~$\mathbf{G}$ and such that~$\ell$
is good for~$\mathbf{L}$. Then $C_{\mathbf{G}}(R)$ is regular in~$\mathbf{G}$.
\end{lem}
\begin{prf}
As~$R$ is abelian, we have $R \leq C_{\mathbf{G}}(z)$. Moreover,
$C_{\mathbf{G}}( R ) \leq \mathbf{L}$, and thus $C_{\mathbf{G}}( R ) =
C_{\mathbf{L}}( R )$. If~$\mathbf{L}^*$ denotes a regular subgroup
of~$\mathbf{G}^*$ dual to~$\mathbf{L}$, we have $\ell \nmid
|Z( \mathbf{L}^* )/Z( \mathbf{L}^* )^\circ|$; see, e.g.\
\cite[Proposition~$4.2$]{CeBo2}. By replacing~$\mathbf{G}$ with~$\mathbf{L}$,
we may thus assume that~$\ell$ is good for~$\mathbf{G}$. In this case we prove
the assertion by induction on~$|G|$.

If $R \leq Z( \mathbf{G} )$, there is nothing to prove. Thus assume that there
is $y \in R \setminus Z( \mathbf{G} )$, and put $\mathbf{M} :=
C_{\mathbf{G}}( y )$. Then~$\mathbf{M}$ is a regular subgroup of~$\mathbf{G}$,
as~$y$ is an $\ell$-element, $\ell \nmid
|Z( \mathbf{G}^* )/Z( \mathbf{G}^* )^\circ|$, and~$\ell$ is
good for~$\mathbf{G}$; see \cite[Corollary~$2.6$]{GeHi1}. It follows that
$\ell \nmid |Z( \mathbf{M}^* )/Z( \mathbf{M}^* )^\circ|$, and that~$\ell$ is
good for~$\mathbf{M}$. Since
$y \not\in Z(G)$, we have $|M| < |G|$. Hence~$C_{\mathbf{M}}( R )$ is a regular
subgroup of~$\mathbf{M}$ by induction.  As $C_{\mathbf{G}}( R ) =
C_{\mathbf{M}}( R )$, we are done.
\end{prf}

\medskip

\noindent
By letting $z = 1$, the conclusion of the above lemma holds in particular if
$Z( \mathbf{G}^* )$ is connected and~$\ell$ is good for~$\mathbf{G}$.

We record a further useful consequence of Theorem~\ref{BoRoEtAl}.
%%%%%%%%%%%%%%%%%%%%%%%%%%%%%%%%%%%%%%%%%%%%%%%%%%%%%%%%%%%%%%%%%%%%%%%%%%%%%%%%
%%%%%%%%%%%%%%%%%%%%%%%%%%%%%%%%%%%%%%%%%%%%%%%%%%%%%%%%%%%%%%%%%%%%%%%%%%%%%%%%
%%
%% \marginpar{The following Lemma~\ref{NewOutLemma} was included on February~$12$,
%% 2019. It is intended to replace some parts of Lemma~\ref{lemmain}.}
%%
%%%%%%%%%%%%%%%%%%%%%%%%%%%%%%%%%%%%%%%%%%%%%%%%%%%%%%%%%%%%%%%%%%%%%%%%%%%%%%%%
%%%%%%%%%%%%%%%%%%%%%%%%%%%%%%%%%%%%%%%%%%%%%%%%%%%%%%%%%%%%%%%%%%%%%%%%%%%%%%%%
\begin{lem}
\label{NewOutLemma}
Let $s \in G^*$ be semisimple such that $\mathbf{L}^* := C_{\mathbf{G}^*}( s )$
is a regular subgroup of $\mathbf{G}^*$. Choose a regular subgroup $\mathbf{L}
\leq \mathbf{G}$ dual to~$\mathbf{L}^*$.

Let $b \subseteq \cE_\ell( G , s )$ be an $\ell$-block of~$G$ and $b' \subseteq
\cE_\ell( L, s )$ the corresponding block according to
{\rm Theorem~\ref{BoRoEtAl}(a)}. Assume that $b' = \hat{s} \otimes b_0$,
where~$b_0$ is the principal block of~$L$ and $\hat{s}$ is a linear character 
of~$L$ corresponding to~$s$ via duality; see \cite[(8.19)]{CaEn}. Let 
$D \leq L$ be a common defect group of~$b$ and $b'$. Let $(R,b_R) \leq (D,b_D)$ 
be a centric $b$-Brauer pair. Then the following statements hold.
	
{\rm (a)} We have $\Out_G(R,b_R) = \Out_L(R)$.
If $\cW( R, b_R ) \neq 0$, then~$R$ is a radical subgroup of~$L$.
(For the definition of $\cW( R, b_R )$ see~{\rm (\ref{NumberOfWeights0})}).

{\rm (b)} The canonical character $\theta_R$ of $b_R$ extends to
$N_G( R, b_R )$.

{\rm (c)} If $N_L( R )$ fixes~$b_R$, then the homomorphism
$$N_L( R ) \rightarrow N_G( R, b_R )/RC_G(R) = \Out_G( R, b_R )$$
is surjective with kernel $RC_L( R )$.

{\rm (d)} If~$R$ is abelian and $\mathbf{K} :=
C_{\mathbf{L}}(R)$ is a regular subgroup of~$\mathbf{L}$, then
$\Out_G( R, b_R ) = W_{\mathbf{L}}( \mathbf{K} )^F$.
\end{lem}
\begin{prf}
By our hypothesis, the group~$\mathbf{M}^*$ of Theorem~\ref{BoRoEtAl} is equal
to $\mathbf{L}^* = C_{\mathbf{G}^*}( s )$, and~$\mathbf{M}$ may be chosen to 
be equal to~$\mathbf{L}$.

(a) By Theorem~\ref{BoRoEtAl}(b), there is a centric $b'$-Brauer pair 
$(R,b'_R)$ of $L$ such that $\Out_G(R,b_R) = \Out_L(R,b'_R)$. By
%
%In 
%particular,~$R$ is a defect group of some $\ell$-block of~$L$ and thus 
%a radical subgroup of~$L$. By
our assumption on~$b'$, the canonical character~$\theta_R'$ of~$b'_R$ equals
$\Res_{C_L(R)}^L( \hat{s} )$. In particular, $N_L( R, \theta_R' ) = N_L( R )$ 
and~$\theta_R'$ extends to $N_L( R, \theta_R' )$.
We obtain
\begin{eqnarray*}
\Out_G(R,b_R) & = & \Out_L( R, b_R' ) \\
& = & N_L(R,b_R')/RC_L(R) \\
& = & N_L( R )/RC_L( R ) \\
& = & \Out_L(R),
\end{eqnarray*}
giving our first claim. 

To prove the second, let $Q := O_\ell(N_L(R))$, so that $R \leq Q$ and 
$Q \cap RC_L(R) = R$ as $b_R'$ is centric. Now $O_\ell( {\Out}_L(R) )
\cong O_\ell( \Out_G( R, b_R ) ) \cong \{ 1 \}$, the first isomorphism arising 
from the first claim, and the second one from~(\ref{NumberOfWeights}) and our 
hypothesis $\cW(R,b_R) \neq 0$. This implies that the image of~$Q$ in 
$\Out_L(R)$ is trivial. Hence $Q=R$ and~$R$ is a radical $\ell$-subgroup of~$L$.

(b) As $\Res_{C_L(R)}^L( \hat{s} )$ extends to $N_L( R )$, 
the K{\"u}ls\-ham\-mer-Puig class associated to $(R,b_R')$ is trivial. By
Theorem~\ref{BoRoEtAl}(b), the K{\"u}lshammer-Puig class associated to
$(R,b_R)$ is trivial as well, hence $\theta_R$ extends to~$N_G( R, b_R)$.

(c) If $N_L( R )$ fixes~$b_R$, the map $N_L( R ) \rightarrow 
N_G( R, b_R )/RC_G(R)$ is well defined and has kernel $RC_L(R)$, hence is 
surjective by~(a).

(d) Now suppose that $R$ is abelian and that $\mathbf{K} = C_{\mathbf{L}}(R)$ 
is a regular subgroup of~$\mathbf{L}$. As $R \leq Z( K )$ is a defect group 
of~$b_R'$, we have $R = O_\ell( Z( K ) )$ and thus $N_L( K ) = N_L( R ) \leq
N_L( C_{\mathbf{L}}( R ) ) = N_L( \mathbf{K} ) \leq N_L( K )$. In particular, 
$N_L( K ) = N_L( \mathbf{K} )$ and thus 
$$\Out_L( R ) = N_L( R )/C_L( R ) =
N_{\mathbf{L}}( \mathbf{K} )^F/\mathbf{K}^F = W_{\mathbf{L}}( \mathbf{K} )^F.$$
This completes our proof.
\end{prf}

\addtocounter{subsection}{4}
\subsection{Non-regular centralizers}
We now work towards a variant of Lemma~\ref{NewOutLemma} in the case 
when $C_{\mathbf{G}^*}( s )$ is not regular, so that $C_{\mathbf{G}^*}( s )^*$ 
cannot be embedded into~$\mathbf{G}$.

\addtocounter{thm}{1}
\begin{lem}
\label{RelativeWeylGroups}
Let~$\mathbf{M}$ and~$\mathbf{M}^*$ be a pair of dual regular subgroups
of~$\mathbf{G}$ and~$\mathbf{G}^*$, respectively. Then there is an 
$F$-$F^*$-equivariant group isomorphism
$$W_{\mathbf{G}}( \mathbf{M} ) \stackrel{\alpha}{\longrightarrow} 
W_{\mathbf{G}^*}( \mathbf{M}^* ),\quad w \mapsto w^*,$$
satisfying the following condition. 
	
Let~$\mathbf{S}$ and~$\mathbf{S}^*$ be $F$-stable, respectively $F^*$-stable,
maximal tori of~$\mathbf{M}$, respectively~$\mathbf{M}^*$, 
let $\theta \in \Irr( S )$ and $s \in S^*$, such that the $M$-conjugacy
class of the pair $(\mathbf{S},\theta)$ corresponds to the $M^*$-conjugacy
class of the pair $(\mathbf{S}^*,s)$ under the duality of $\mathbf{M}$ and 
$\mathbf{M}^*$; see, e.g.\ \cite[Proposition~$11.1.16$]{DiMi2}. Further, let
$w \in W_{\mathbf{G}}( \mathbf{M} )^F$ and let $x \in N_G( \mathbf{M} )$
and $y \in N_{G^*}( \mathbf{M}^* )$ denote inverse images of~$w$
and~$w^*$, respectively. Then
$$R_{\mathbf{S}^*}^{\mathbf{M}}( s )^x =
R_{(\mathbf{S}^*)^{y}}^{\mathbf{M}}( s^y ).$$
\end{lem}
\begin{prf}
For the first assertion see the concluding remarks of \cite[Section~$8.2$]{CaEn}.
The displayed equation is \cite[$(4.3)$ Lemma]{DipFlI}.
\end{prf}

\begin{lem}\label{lemmain}
Let~$\mathbf{M}$ and~$\mathbf{M}^*$ be a pair of dual regular subgroups 
of~$\mathbf{G}$ and~$\mathbf{G}^*$, respectively, and let $s \in M^*$ be 
semisimple.
Let $x \in N_{G}( \mathbf{M} )$ and $y \in N_{G^*}( \mathbf{M}^* )$ be such 
that the images of~$x$ and~$y$ in $W_{\mathbf{G}}( \mathbf{M} )^F$ respectively 
$W_{\mathbf{G}^*}( \mathbf{M}^* )^{F^*}$ correspond under the isomorphism given
in {\rm Lemma~\ref{RelativeWeylGroups}}. 

Then $\cE(M, s)^x = \cE(M, s^y)$. In particular, if $y \in 
C_{G^*}( s )$, then~$x$ stabilizes $\cE(M, s)$.
\end{lem}
\begin{prf}
Let $\mathbf{S}^* \leq \mathbf{M}^*$ be an $F^*$-stable maximal torus with
$s \in S^*$.
By Lemma~\ref{RelativeWeylGroups} we have $R_{\mathbf{S}^*}^\mathbf{M}(s)^x 
= R_{{(\mathbf{S}^*)}^y}^\mathbf{M}(s^y)$. As
all irreducible constituents of $R_{{(\mathbf{S}^*)}^y}^\mathbf{M}(s^y)$ lie in
$\cE( M, s^y )$, the claim follows.
\end{prf}

\medskip

\noindent
Part~(a) of the next proposition is a corollary to Lemma \ref{lemmain}.
%%%%%%%%%%%%%%%%%%%%%%%%%%%%%%%%%%%%%%%%%%%%%%%%%%%%%%%%%%%%%%%%%%%%%%%%%%%%%%%%
%%%%%%%%%%%%%%%%%%%%%%%%%%%%%%%%%%%%%%%%%%%%%%%%%%%%%%%%%%%%%%%%%%%%%%%%%%%%%%%%
%%
%% \marginpar{Proposition~\ref{cormain} was simplified on April 16, 2020.
%% The previous version is Proposition~\ref{cormainOld} below.}
%%
%%%%%%%%%%%%%%%%%%%%%%%%%%%%%%%%%%%%%%%%%%%%%%%%%%%%%%%%%%%%%%%%%%%%%%%%%%%%%%%%
%%%%%%%%%%%%%%%%%%%%%%%%%%%%%%%%%%%%%%%%%%%%%%%%%%%%%%%%%%%%%%%%%%%%%%%%%%%%%%%%

\begin{prop}\label{cormain}
Let $R \leq G$ and $R^\dagger \leq G^*$ be abelian radical $\ell$-subgroups
such that $\mathbf{M} := C_{\mathbf{G}}(R)$ is a regular subgroup
of~$\mathbf{G}$ and $\mathbf{M}^* := C_{\mathbf{G}^*}( R^\dagger )$ is dual
to $\mathbf{M}$. Then $\Out_G(R) = W_{\mathbf{G}}( \mathbf{M} )^F$ and
$\Out_G(R^\dagger) = W_{\mathbf{G}^*}( \mathbf{M}^* )^{F^*}$.
	
Let $s \in M^*$ be a semisimple $\ell'$-element and put
$\mathbf{L}^* := C_{\mathbf{G}^*}( s ) \leq \mathbf{G}^*$. Assume that
$\mathbf{L}^*$ is connected.

{\rm (a)} Let $b_R \subseteq \cE_\ell(C_G(R), s)$ be an $\ell$-block of $C_G(R)$ 
and assume that $N_G( R, b_R ) = 
\{ x \in N_G( R ) \mid x \text{\ stabilizes\ } \cE(C_G(R), s) \}$. Then
the isomorphism 
$$
W_{\mathbf{G}}( \mathbf{M} )^F \stackrel{\alpha}{\longrightarrow} 
W_{\mathbf{G}^*}( \mathbf{M}^* )^{F^*}
$$
arising from {\rm Lemma~\ref{RelativeWeylGroups}} maps $\Out_G( R, b_R )$ to 
$\Out_{L^*}(R^\dagger)$, i.e.\ it induces
an isomorphism
$$
\Out_G( R, b_R ) \cong \Out_{L^*}(R^\dagger).
$$

{\rm (b)} Suppose that $(R,b_R)$ is centric and that
$\mathbf{K}^* := C_{\mathbf{L}^*}( R^\dagger )$ is a regular subgroup 
of~$\mathbf{L}^*$ and of~$\mathbf{M}^*$. Suppose further that 
$$\ell \nmid 
|Z( \mathbf{M}^* )/Z( \mathbf{M}^* )^\circ|
|Z( \mathbf{M} )/Z( \mathbf{M} )^\circ|.$$
Then
\[
\Out_{L^*}(R^\dagger) = W_{\mathbf{L}^*}( \mathbf{K}^* )^{F^*}.
\]
%%%%%%%%%%%%%%%%%%%%%%%%%%%%%%%%%%%%%%%%%%%%%%%%%%%%%%%%%%%%%%%%%%%%%%%%%%%%%%%%
%%%%%%%%%%%%%%%%%%%%%%%%%%%%%%%%%%%%%%%%%%%%%%%%%%%%%%%%%%%%%%%%%%%%%%%%%%%%%%%%
%%
%% \marginpar{Sufficient conditions for the last statement to hold are:
%% $R^\dagger = O_\ell( Z( K^* ) )$, or $R^\dagger$ is a radical subgroup of~$L^*$.}
%%
%%%%%%%%%%%%%%%%%%%%%%%%%%%%%%%%%%%%%%%%%%%%%%%%%%%%%%%%%%%%%%%%%%%%%%%%%%%%%%%%
%%%%%%%%%%%%%%%%%%%%%%%%%%%%%%%%%%%%%%%%%%%%%%%%%%%%%%%%%%%%%%%%%%%%%%%%%%%%%%%%
%%%%%%%%%%%%%%%%%%%%%%%%%%%%%%%%%%%%%%%%%%%%%%%%%%%%%%%%%%%%%%%%%%%%%%%%%%%%%%%%
%%%%%%%%%%%%%%%%%%%%%%%%%%%%%%%%%%%%%%%%%%%%%%%%%%%%%%%%%%%%%%%%%%%%%%%%%%%%%%%%
%%
%% See notes on Sheet A of 30.04.2020.
%%
%%%%%%%%%%%%%%%%%%%%%%%%%%%%%%%%%%%%%%%%%%%%%%%%%%%%%%%%%%%%%%%%%%%%%%%%%%%%%%%%
%%%%%%%%%%%%%%%%%%%%%%%%%%%%%%%%%%%%%%%%%%%%%%%%%%%%%%%%%%%%%%%%%%%%%%%%%%%%%%%%
\end{prop}
\begin{prf}
As~$R$ is an abelian radical subgroup of~$G$, we have
$R = O_\ell( C_G ( R ) )$ and thus $N_G( R ) = N_G( C_G( R ) ) = N_G(M)$ and
$N_{G}( R ) = N_{G}( \mathbf{M} )$. Analogously, $N_{G^*}( R^\dagger ) = 
N_{G^*}( M^* ) = N_{G^*}( \mathbf{M}^* )$, giving our first assertion. 
	
(a) We show that $\alpha$ maps $N_G( R, b_R )/M \leq N_G(M)/M = 
W_{\mathbf{G}}( \mathbf{M} )^F$ to $\Out_{L^*}(R^\dagger)$, naturally
embedded into $\Out_{G^*}( R^\dagger ) = 
W_{\mathbf{G}^*}( \mathbf{M}^* )^{F^*}$.

Let $x \in N_{G}( \mathbf{M} )$ and $y \in N_{G^*}( \mathbf{M}^* )$ be such
$\alpha( M{x} ) = M^*{y}$. Suppose first that $x \in N_G( R, b_R )$. Then~$x$ 
stabilizes $\cE_\ell( C_G ( R ), s )$ and thus 
$\cE_\ell( C_G ( R ), s ) = \cE_\ell( C_G ( R ), s )^x = 
\cE_\ell( C_G ( R ), s^y )$ by Lemma~\ref{lemmain}. Hence~$s^y$ 
is conjugate to~$s$ in~$M^*$, i.e.\ $ty \in L^*$ for some $t \in M^*$.
Hence~$\alpha$ maps $N_G( R, b_R )/M$ into $\Out_{L^*}( R^\dagger )$.

Now assume that $M^*{y} \in \Out_{L^*}( R^\dagger )$. By multiplying~$y$
with a suitable element of~$M^*$, we may assume that~$y \in L^*$. With
Lemma~\ref{lemmain} we conclude that~$x$ stabilizes $\cE(C_G(R), s)$.
Our assumption now implies that $x \in N_G(R, b_R)$.

(b) This is very similar to the last part of the proof of 
Lemma~\ref{NewOutLemma}. We claim that 
$R^\dagger = O_\ell( Z( C_{L^*}( R^\dagger ) ) )$. Once this
claim is proved, we can conclude $N_{L^*}( C_{L^*}( R^\dagger ) ) = 
N_{L^*}( R^\dagger ) = N_{L^*}( C_{\mathbf{L}^*}( R^\dagger ) )$,
yielding our assertion.

To prove the claim, first notice that $\mathbf{K}^* = 
\mathbf{L}^* \cap \mathbf{M}^* = C_{\mathbf{M}^*}( s )$. As $(R,b_R)$ is
centric, the defect group of~$b_R$ equals~$R$, and as~$R$ is contained in
$Z( M )$, we must have $R = O_\ell( Z( M ) )$. By
\cite[Proposition~$4.4.5$]{C2} and our assumption on~$\ell$, we conclude that
$R^\dagger = O_\ell( Z( M^* ) )$. If $s \in Z( M^* )$, then $M^* = 
C_{L^*}( R^\dagger )$ and our claim follows. Otherwise, $C_{\mathbf{M}^*}( s )$
is a proper regular subgroup of~$\mathbf{M}^*$. Let $\mathbf{M}_s$ denote a
regular subgroup of~$\mathbf{M}$ dual to~$C_{\mathbf{M}^*}( s )$. By
Theorem~\ref{BoRoEtAl}(a), the block $b_R \subseteq \cE_\ell( M, s )$
corresponds to a block of $\cE_\ell( M_s, s )$ with a defect group~$R_s$
conjugate to~$R$ in~$M$. As $Z( M ) \leq Z( M_s )$, we get $R = 
O_\ell( Z( M_s ) )$. The assumption on~$\ell$ descends to the regular
subgroups $\mathbf{M}_s$ and $C_{\mathbf{M}^*}( s )$ of~$\mathbf{M}$ and
$\mathbf{M}^*$, respectively; see, e.g.\ \cite[Proposition~$4.2$]{CeBo2}.
We thus obtain $R^\dagger = O_\ell( Z( C_{M^*}( s ) ) )$.
As $C_{M^*}( s ) = C_{L^*}( R^\dagger )$, this yields our claim.
\end{prf}

\section{The group $F_4(q)$ and some of its Levi subgroups}
\label{LeviSubgroups}

In this section we introduce the group $F_4(q)$ and investigate some of its 
subgroups. The main new result is contained in 
Corollary~\ref{MaximalExtendability}: Let $e$ be a positive integer. Then every 
irreducible character of an $e$-split Levi subgroup of~$F_4(q)$ extends to its 
inertia subgroup. This generalizes the result of Sp{\"a}th in 
\cite[Theorem~$1.1$]{Spaeth10} in case of the group $F_4(q)$. For the definition of
$e$-split Levi subgroups see, e.g.\ \cite[$3.5.1$]{GeMa}.

\subsection{Setup, notation and preliminaries}
\label{SetupF4}
Let~$p$ be a prime,~$f$ a positive integer and $q = p^f$. Further, 
let~$\mathbf{G}$ denote a simple algebraic group of type~$F_4$ over an algebraic 
closure~$\mathbb{F}$ of~$\mathbb{F}_p$, such that~$\mathbf{G}$ has a standard 
Frobenius morphism~$F_1$ with $\mathbf{G}^{F_1} = F_4(p)$. We put $F := F_1^f$, 
so that $G = \mathbf{G}^F = F_4( q )$. We have
$$|G| = q^{24}\Phi_1(q)^4\Phi_2(q)^4\Phi_3(q)^2\Phi_4(q)^2\Phi_6(q)^2
       \Phi_8(q)\Phi_{12}(q),$$
where $\Phi_i$ denotes the $i$th cyclotomic polynomial.

We choose an $F$-stable maximal torus~$\mathbf{T}_0$ of~$\mathbf{G}$ contained 
in an $F$-stable Borel subgroup~$\mathbf{B}$ of~$\mathbf{G}$, so 
that~$\mathbf{T}_0$ is maximally split and $F(t) = t^q$ for all 
$t \in \mathbf{T}_0$. Write $W := W_{\mathbf{G}}( \mathbf{T}_0 ) := 
N_{\mathbf{G}}( \mathbf{T}_0 )/\mathbf{T}_0$ for the Weyl group 
of~$\mathbf{G}$ (with respect to~$\mathbf{T}_0$). Let~$\mathbf{U}$ denote the 
unipotent radical of $\mathbf{B}$. The root system of~$\mathbf{G}$ is 
denoted by~$\Sigma$, the root subgroup giving rise to~$\alpha \in \Sigma$
by $\mathbf{U}_\alpha$, and $u_\alpha: \mathbb{F} \rightarrow \mathbf{U}_\alpha$ 
the corresponding isomorphism of algebraic groups. The choice of~$\mathbf{B}$ 
determines the 
set~$\Sigma^+$ of positive roots and the corresponding base $\alpha_i$, 
$i = 1, \ldots , 4$, numbered as in the following Dynkin diagram:

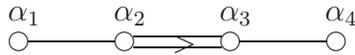
\begin{center}
\begin{picture}(123.5,60)(0,-22)
\put(   0,   0){\circle{7}}
\put(  40,   0){\circle{7}}
\put(  80,   0){\circle{7}}
\put( 120,   0){\circle{7}}
\put( 3.5,   0){\line(1,0){ 33}}
\put(43.5,   2){\line(1,0){ 33}}
\put(43.5,  -2){\line(1,0){ 33}}
\put(83.5,   0){\line(1,0){ 33}}
\put(58.0, -4){$>$}
\put(   -4,+8){$\alpha_1$}
\put(   36,+8){$\alpha_2$}
\put(   76,+8){$\alpha_3$}
\put(  116,+8){$\alpha_4$}
\end{picture}
\end{center}
Thus~$\alpha_1$, $\alpha_2$ are the long simple roots, and $\alpha_3$,
$\alpha_4$ the short ones. 

As $(\mathbf{G},F)$ is split, $F$ acts trivially on $W$, so that $W = W^F = 
N_G( \mathbf{T}_0 )/T_0$. Moreover, for each $\alpha \in \Sigma$, the root 
subgroup~$\mathbf{U}_\alpha$ is $F$-stable, and $F( u_\alpha( t ) ) = 
u_\alpha( t^q )$ for all $t \in \mathbb{F}$; in particular, $U_\alpha = 
\mathbf{U}_\alpha^F = \{ u_\alpha( t ) \mid t \in \mathbb{F}_q \}$.
For $\alpha \in \Sigma$, write $n_\alpha := u_{\alpha}(1)u_{-\alpha}(-1)
u_{\alpha}(1)$. Then $n_\alpha \in N_{\mathbf{G}}( \mathbf{T}_0 )^{F_1}$, 
and we write $s_{\alpha}$ for the image of~$n_{\alpha}$ in~$W$.
Let $w_0$ denote the longest element of~$W$.
Then $w_0 \in Z(W)$. If $\dot{w}_0 \in N_{\mathbf{G}}( \mathbf{T}_0 )$ is an 
inverse image of $w_0$, then~$\dot{w}_0$ inverts the elements of~$\mathbf{T}_0$ 
and $\mathbf{U}_\alpha^{\dot{w}_0} = \mathbf{U}_{-\alpha}$ for all $\alpha \in 
\Sigma$. In the following, we will refer to some computations in~$\Sigma$ using
CHEVIE~\cite{chevie}. For easier reference, Table~\ref{PositiveRoots} gives
the numbering of the roots of $\Sigma^+$ as in CHEVIE. In the column 
headed~$\alpha$ we list the CHEVIE number of~$\alpha \in \Sigma^+$,
and in the column headed ``Root'' the expansion of~$\alpha$ in the base
$\{ \alpha_1, \ldots , \alpha_4 \}$. The significance of the column headed
$\alpha^\dagger$ will be explained in Subsection~\ref{Duality} below.
The negative roots are numbered as $\alpha_{24 + i} := -\alpha_i$ for 
$1 \leq i \leq 24$. If $\alpha = \alpha_i$ for some $1 \leq i \leq 48$, we
write $n_i := n_{\alpha_i}$, $s_i := s_{\alpha_i}$ and $u_i := u_{\alpha_i}$. 

As already introduced in Subsection~\ref{CharacterAndCocharacterGroups}, we 
write $X := X( \mathbf{T}_0 ) = \Hom( \mathbf{T}_0, \mathbb{F}^* )$ and 
$Y := Y( \mathbf{T}_0 ) = \Hom(  \mathbb{F}^*, \mathbf{T}_0 )$ for the character 
group and the cocharacter group of~$\mathbf{T}_0$, respectively. Then 
$X \cong Y \cong 
\mathbb{Z}^4$ as abelian groups. The homomorphism $h_\alpha : \mathbb{F}^* 
\rightarrow \mathbf{T}_0$ in~$Y$ associated to $\alpha \in \Sigma$ as in 
\cite[Theorem~$12.1.1$]{C1} is called the coroot corresponding to~$\alpha$, 
and will be denoted by $\alpha^\vee$. For every subset $\Gamma \subseteq 
\Sigma$ write $\Gamma^\vee := \{ \alpha^\vee \mid \alpha \in \Gamma \}$. Then 
$\{ \alpha_1, \ldots, \alpha_4 \}  \subseteq X$ and 
$\{ \alpha_1^\vee, \ldots \alpha_4^\vee \} \subseteq Y$ are $\mathbb{Z}$-bases 
of~$X$ and~$Y$ respectively. The root datum of $\mathbf{G}$ is 
$(X, \Sigma, Y, \Sigma^\vee )$.

\begin{table}[t]
\caption{\label{PositiveRoots} The positive roots of $\Sigma$}
$
\begin{array}{rcr}  \hline\hline
\alpha & \text{Root} & \alpha^\dagger\rule[- 3pt]{0pt}{ 16pt} \\ \hline\hline
	1 &[ 1, 0, 0, 0 ] & 4 \rule[ 0pt]{0pt}{ 13pt} \\
	2 &[ 0, 1, 0, 0 ] & 3 \\
	3 &[ 0, 0, 1, 0 ] & 2 \\
	4 &[ 0, 0, 0, 1 ] & 1 \\
	5 &[ 1, 1, 0, 0 ] & 7 \\
	6 &[ 0, 1, 1, 0 ] & 9 \\
	7 &[ 0, 0, 1, 1 ] & 5 \\
	8 &[ 1, 1, 1, 0 ] & 16
\rule[- 2pt]{0pt}{ 5pt} \\ \hline\hline
\end{array}
$
$
\begin{array}{rcr}  \hline\hline
\alpha & \text{Root} & \alpha^\dagger \rule[- 3pt]{0pt}{ 16pt} \\ \hline\hline
	 9 &[ 0, 1, 2, 0 ] & 6 \rule[ 0pt]{0pt}{ 13pt} \\
	10 &[ 0, 1, 1, 1 ] & 11 \\
	11 &[ 1, 1, 2, 0 ] & 10 \\
	12 &[ 1, 1, 1, 1 ] & 18 \\
	13 &[ 0, 1, 2, 1 ] & 14 \\
	14 &[ 1, 2, 2, 0 ] & 13 \\
	15 &[ 1, 1, 2, 1 ] & 20 \\
	16 &[ 0, 1, 2, 2 ] & 8 
\rule[- 2pt]{0pt}{ 5pt} \\ \hline\hline
\end{array}
$
$
\begin{array}{rcc}  \hline\hline
\alpha & \text{Root} & \alpha^\dagger \rule[- 3pt]{0pt}{ 16pt} \\ \hline\hline
	17 &[ 1, 2, 2, 1 ] & 22  \rule[ 0pt]{0pt}{ 13pt} \\
	18 &[ 1, 1, 2, 2 ] & 12 \\
	19 &[ 1, 2, 3, 1 ] & 23 \\
	20 &[ 1, 2, 2, 2 ] & 15 \\
	21 &[ 1, 2, 3, 2 ] & 24 \\
	22 &[ 1, 2, 4, 2 ] & 17 \\
	23 &[ 1, 3, 4, 2 ] & 19 \\
	24 &[ 2, 3, 4, 2 ] & 21 
\rule[- 2pt]{0pt}{ 5pt} \\ \hline\hline
\end{array}
$
\end{table}

%[ 4, 3, 2, 1, 7, 9, 5, 16, 6, 11, 10, 18, 14, 13, 20, 8, 22, 12, 23, 15, 24,
%  17, 19, 21, 28, 27, 26, 25, 31, 33, 29, 40, 30, 35, 34, 42, 38, 37, 44, 32,
%  46, 36, 47, 39, 48, 41, 43, 45 ]

Let $\Gamma \subseteq \Sigma$. By $\bar{\Gamma}$ we denote the smallest closed 
subsystem of~$\Sigma$ containing~$\Gamma$. We put $\mathbf{L}_\Gamma := 
\langle \mathbf{T}_0, \mathbf{U}_\alpha \mid \alpha \in \bar{\Gamma} \rangle$ 
and $\mathbf{K}_\Gamma := 
\langle \mathbf{U}_\alpha \mid \alpha \in \bar{\Gamma} \rangle$. Then 
$\mathbf{L}_\Gamma$ and $\mathbf{K}_\Gamma$ are connected reductive algebraic 
groups, $\mathbf{K}_\Gamma$ is semisimple and $\mathbf{K}_\Gamma = 
[\mathbf{L}_\Gamma, \mathbf{L}_\Gamma]$; see, e.g.\ \cite[Theorem~$13.6$]{MaTe}. 
It follows from \cite[Proposition~$12.14$]{MaTe} that $\mathbf{K}_{\Gamma}$ is 
simply connected  if $\Gamma \subseteq \{ \alpha_1, \ldots , \alpha_4 \}$, i.e.\ 
if $\bar{\Gamma}$ is a parabolic subsystem. Notice that $n_\alpha \in K_\Gamma = 
\mathbf{K}_\Gamma^F$ if $\alpha \in \bar{\Gamma}$. Moreover, $\mathbf{L}_\Gamma$ 
and $\mathbf{K}_\Gamma$ are $F$-stable and $\mathbf{T}_0$, respectively 
$\mathbf{T}_0 \cap \mathbf{K}_\Gamma$, are maximally split maximal tori of 
$\mathbf{L}_\Gamma$ and $\mathbf{K}_\Gamma$, respectively; the corresponding 
root system is $\bar{\Gamma}$ and we have $\bar{\Gamma}^+ = 
\bar{\Gamma} \cap \Sigma^+$. The root datum of $\mathbf{L}_\Gamma$ with respect 
to~$\mathbf{T}_0$ equals $(X, \bar{\Gamma}, Y , \bar{\Gamma}^\vee )$. Let us 
also put $W_\Gamma := \langle s_\alpha \mid \alpha \in \Gamma \rangle$. Then the 
Weyl group of~$\mathbf{L}_\Gamma$ with respect to~$\mathbf{T}_0$ 
equals~$W_{\bar{\Gamma}}$. If~$\Gamma$ is a base of $\bar{\Gamma}$, then 
$W_\Gamma = W_{\bar{\Gamma}}$, but, in general, 
$W_\Gamma \lneq W_{\bar{\Gamma}}$.

If $\Gamma, \Delta \subseteq \Sigma$ are disjoint, closed subsystems of~$\Sigma$
such that $\Delta \cup \Gamma$ is closed, then
$[\mathbf{K}_\Delta, \mathbf{K}_\Gamma] = 1$ by the commutator relations for
root subgroups.

\addtocounter{thm}{1}
\begin{lem}
\label{AdjustmentOfInverseImages}
Suppose that $\beta_1, \ldots , \beta_m$ are $m$ distinct roots of~$\Sigma$ 
contained in some base (in particular, $m \leq 4$). Then for any $a_1, \ldots , 
a_m \in \mathbb{F}^*$, there is $t \in \mathbf{T}_0$ with $\beta_i( t ) 
= a_i$ for $1 \leq i \leq m$. Moreover,~$t$ may be chosen in $T_0$ if
$a_1, \ldots , a_m \in \mathbb{F}_q^*$.
\end{lem}
\begin{prf}
As $\mathbf{G}$ is adjoint, a base of~$\Sigma$ is a $\mathbb{Z}$-basis of~$X$. 
The result follows from $\mathbf{T}_0 \cong \Hom( X, \mathbb{F}^* )$; see 
\cite[Proposition~$3.1.2$(i)]{C2}. As the latter isomorphism is $F$-equivariant,
the final statement also follows.
\end{prf}

\medskip
\noindent
Of particular importance is the following construction, a special case of
twisting; see Subsection~\ref{Twisting} below. If $g \in \mathbf{G}$ such that 
$F(g)g^{-1}$ normalizes~$\mathbf{T}_0$ and maps to~$w_0$ under the natural 
epimorphism, then $\mathbf{T}_0^g$, $\mathbf{L}^g_\Gamma$ and 
$\mathbf{K}^g_\Gamma$ are $F$-stable for every $\Gamma \subseteq \Sigma$,
and $F$ acts trivially on
$N_{\mathbf{G}}( \mathbf{T}_0^g )/\mathbf{T}_0^g$. As~$w_0$ acts as $-\id$ 
on~$X$, we have $F(t) = t^{-q}$ for all $t \in \mathbf{T}_0^g$. A torus which 
is $G$-conjugate to~$\mathbf{T}_0$, respectively to~$\mathbf{T}_0^g$ is
called $1$-$F$-split, respectively $2$-$F$-split, where we omit the~$F$ if
clear from the context.
%%%%%%%%%%%%%%%%%%%%%%%%%%%%%%%%%%%%%%%%%%%%%%%%%%%%%%%%%%%%%%%%%%%%%%%%%%%%%%%%
%%%%%%%%%%%%%%%%%%%%%%%%%%%%%%%%%%%%%%%%%%%%%%%%%%%%%%%%%%%%%%%%%%%%%%%%%%%%%%%%
%%
%% See sheet 1 of 13.02.2020.
%%
%%%%%%%%%%%%%%%%%%%%%%%%%%%%%%%%%%%%%%%%%%%%%%%%%%%%%%%%%%%%%%%%%%%%%%%%%%%%%%%%
%%%%%%%%%%%%%%%%%%%%%%%%%%%%%%%%%%%%%%%%%%%%%%%%%%%%%%%%%%%%%%%%%%%%%%%%%%%%%%%%

\addtocounter{subsection}{1}
\subsection{A lift of the longest element}
\label{LiftOfLongestElement}

The existence of a lift of~$w_0$ of order~$2$ is indicated in 
\cite[Definition~$(2.23)$]{Griess}. We will make use of a
particular such lift with further properties.

\addtocounter{thm}{1}
\begin{lem}
\label{LiftOfLongestElementLem}
Let $w_0$ denote the longest element of~$W$. Then there is a lift
$\gamma \in N_{G}( \mathbf{T}_0 )$ of $w_0$ such that
$$u_\alpha( t )^\gamma = u_{-\alpha}( -t )$$
for all $\alpha \in \{ \pm \alpha_1, \ldots , \pm \alpha_4 \}$
and all $t \in \mathbb{F}$. In particular, $\gamma^2 = 1$.

Moreover, $n_1 , \ldots , n_4$ commute with $\gamma$.
\end{lem}
\begin{prf}
We start with a particular reduced word for $w_0$, and let $\gamma$ denote the 
corresponding product of the $n_{\alpha}$. 
%(According to~\cite{Steinberg},
%the element $\gamma$ is independent of the chosen reduced word for~$w_0$.)
We then use CHEVIE and 
\cite[Lemma~$7.2.1$(i)]{C1} to verify that $\gamma u_\alpha( t ) \gamma^{-1} = 
u_{-\alpha}( -t )$ for all $\alpha \in \{ \alpha_1, \ldots , \alpha_4 \}$ and
all $t \in \mathbb{F}$. The corresponding relation for the negative roots then 
follows from \cite[Proposition~$6.4.3$]{C1}. Now $\gamma^2$ fixes 
$\pm \alpha_1( t ), \ldots , \pm \alpha_4( t )$ for all $t \in \mathbb{F}$, as 
well as the elements of~$\mathbf{T}_0$, and, as 
$\mathbf{G} = \langle \mathbf{T}_0, \pm \alpha_1( t ), \ldots , \pm \alpha_4( t )
\mid t \in \mathbb{F} \rangle$, we get $\gamma^2 \in Z( \mathbf{G} ) 
= \{ 1 \}$.

Finally,~$\gamma$ acts as inverse-transpose automorphism on 
$\mathbf{L}_{\{ \alpha_j \}} \cong \SL_2( \mathbb{F} )$ in the natural 
representation of $\SL_2( \mathbb{F} )$, for $j = 1, \ldots , 4$. Thus~$\gamma$ 
fixes $n_j$ for $j = 1, \ldots , 4$.
\end{prf}

\medskip
\noindent
The above proof provides an example of a computation in the \textit{extended 
Weyl group} of~$\mathbf{G}$. This is a group associated to a Coxeter system, 
introduced and investigated by Tits in~\cite{Tits}. Let 
$\hat{W} \leq N_{\mathbf{G}}( \mathbf{T}_0 )$ denote the subgroup generated by 
$n_j$, $1 \leq j \leq 24$; for the definition of these elements 
see Subsection~\ref{SetupF4}. We constructed the~$n_j$ as matrices in the adjoint
representation of~$\mathbf{G}$, as described in \cite[Lemma~$4.3.1$]{C1}. 
One checks that the elements $n_1, \ldots ,n_4$
satisfy the relations exhibited in~\cite[$4.6$, Equations (1)--(4)]{Tits}, and 
that $\hat{W}$ has order $2^4|W|$ if~$p$ is odd. Thus in this case,~$\hat{W}$ is 
indeed the extended Weyl group as defined in~\cite[D{\'e}finition~$2.2$]{Tits}. 
If $p = 2$, we have $\hat{W} = W$. If~$p$ is odd, the map $\hat{W} \rightarrow W$ 
defined by sending~$n_j$ to $s_j$ for $j = 1, \ldots , 4$, is surjective with 
kernel of order $2^4$, generated by $n_j^2$, $j = 1, \ldots , 4$.

\addtocounter{subsection}{1}
\subsection{Automorphisms}
\label{Automorphisms}
Recall from Subsection~\ref{SetupF4} that~$F_1$ is the Steinberg endomorphism of 
$\mathbf{G}$ such that $\mathbf{G}^{F_1} = F_4(p)$. 
If~$p$ is odd, let $\sigma_1 := F_1$. If $p = 2$, let~$\sigma_1$ denote the
endomorphism of~$\mathbf{G}$ constructed in \cite[Theorem~$28$]{Steinberg},
such that $\sigma_1^2 = F_1$. Then~$\sigma_1$ is an automorphism of abstract
groups. Following \cite[Definition~$1.15.1$]{Gor}, we write $\Aut_1(\mathbf{G})$
for the set of automorphisms~$\psi$ of the abstract group~$\mathbf{G}$, such
that $\psi$ or $\psi^{-1}$ is an endomorphism of the algebraic
group~$\mathbf{G}$. By the results summarized in \cite[Subsection~$1.15$]{Gor},
we find
$\Aut_1( \mathbf{G} ) = \Inn( \mathbf{G} ) \rtimes \langle \sigma_1 \rangle$.
As $Z( \mathbf{G} )$ is trivial we may identify $\Inn( \mathbf{G} )$
with~$\mathbf{G}$ and $\Aut_1( \mathbf{G} )$ with
$\mathbf{G} \rtimes \langle \sigma_1 \rangle$.

\subsection{Twisting}\label{Twisting}
For $g \in \mathbf{G}$ let $\CJxxx_g : \mathbf{G} \rightarrow \mathbf{G}$,
$x \mapsto g^{-1} x g$ denote conjugation by~$g$.
Let $n \in \mathbf{G}$.
As usual, we write $Fn$ for the Steinberg morphism of~$\mathbf{G}$
defined by $F n := \CJxxx_{n} \circ F$, i.e.\
$$Fn: \mathbf{G} \rightarrow \mathbf{G}, x \mapsto n^{-1} F(x) n.$$
(Notice that, although we compose endomorphisms of~$\mathbf{G}$ from right to
left, as indicated by the symbol~$\circ$, conjugation by~$g$ is ``conjugation 
from the right''. The reason is that in GAP \cite{GAP04}, which we use for numerous 
computations in the Weyl group of~$\mathbf{G}$, this is the default way of 
conjugating in groups. As a consequence of this convention, $\CJxxx_g \circ \CJxxx_h 
= \CJxxx_{hg}$ for $g, h \in \mathbf{G}$) .

Let $\mathbf{N}_0 \unlhd \mathbf{M}_0$ denote closed $F$-stable subgroups 
of~$\mathbf{G}$ normalized by~$n$, so that~$\mathbf{N}_0$ and~$\mathbf{M}_0$ 
are $Fn$-stable. By the Lang-Steinberg theorem, there is
$g \in \mathbf{G}$ with $F(g)g^{-1} = n$; choose one such~$g$ and put $\mathbf{N}
:= \mathbf{N}_0^g$ and $\mathbf{M} := \mathbf{M}_0^g$. We say that~$\mathbf{M}$
is obtained from~$\mathbf{M}_0$ by \textit{($F$-)twisting with}~$n$. Notice 
that~$\mathbf{M}$ is $F$-stable and
$$\CJxxx_g: \mathbf{M}_0^{Fn} \rightarrow \mathbf{M}^F$$
is an isomorphism mapping $\mathbf{N}_0^{Fn}$ to $\mathbf{N}^F$. Let 
$\omega \in \Aut_1( \mathbf{G} )$ such that~$\omega$ stabilizes~$\mathbf{N}$ 
and~$\mathbf{M}$. Then 
$\CJxxx_g^{-1} \circ \omega \circ \CJxxx_g$
stabilizes~$\mathbf{N}_0$ and~$\mathbf{M}_0$ and we obtain the following 
commutative diagram of groups and automorphisms
$$
\begin{xy}
\xymatrix@C+30pt{
\mathbf{M}/\mathbf{N} \ar[r]^{\omega} & \mathbf{M}/\mathbf{N} \\
\mathbf{M}_0/\mathbf{N}_0 \ar[u]^{\CJxxx_g} \ar[r]_{\CJxxx_g^{-1} \circ \omega \circ \CJxxx_g} & 
\mathbf{M}_0/\mathbf{N}_0 \ar[u]_{\CJxxx_g}
}
\end{xy}
$$
where we use~$\omega$ and $\CJxxx_g$ to also denote the induced maps on the
factor groups.

Now assume in addition that~$\omega$ commutes with~$F$. Observe that this is the 
case, if and only if $\CJxxx_g^{-1} \circ \omega \circ \CJxxx_g$ commutes 
with~$Fn$.
%%%%%%%%%%%%%%%%%%%%%%%%%%%%%%%%%%%%%%%%%%%%%%%%%%%%%%%%%%%%%%%%%%%%%%%%%%%%%%%%
%%%%%%%%%%%%%%%%%%%%%%%%%%%%%%%%%%%%%%%%%%%%%%%%%%%%%%%%%%%%%%%%%%%%%%%%%%%%%%%%
%%
%% See sheet 1 of 13.02.2020.
%%
%%%%%%%%%%%%%%%%%%%%%%%%%%%%%%%%%%%%%%%%%%%%%%%%%%%%%%%%%%%%%%%%%%%%%%%%%%%%%%%%
%%%%%%%%%%%%%%%%%%%%%%%%%%%%%%%%%%%%%%%%%%%%%%%%%%%%%%%%%%%%%%%%%%%%%%%%%%%%%%%%
We obtain an analogous diagram for the groups of~$F$-fixed points.
$$
\begin{xy}
\xymatrix@C+30pt{
(\mathbf{M}/\mathbf{N} )^F\ar[r]^{\omega} & (\mathbf{M}/\mathbf{N})^F \\
(\mathbf{M}_0/\mathbf{N}_0)^F \ar[u]^{\CJxxx_g} \ar[r]_{\CJxxx_g^{-1} \circ \omega \circ \CJxxx_g} &
(\mathbf{M}_0/\mathbf{N}_0)^F \ar[u]_{\CJxxx_g}
}
\end{xy}
$$

\subsection{Subgroups of maximal rank}
\label{MaximalRank}
Centralizers of semisimple elements in~$\mathbf{G}$ are connected, reductive
and contain a maximal torus; see, e.g.\ 
\cite[Theorems~$3.5.6$, $3.5.4$, $3.5.3$(i)]{C2}. Subgroups with these 
properties are called \textit{connected reductive subgroups of maximal rank 
of}~$\mathbf{G}$. In particular, the $\mathbf{G}$-conjugates of the 
groups $\mathbf{L}_\Gamma$, where $\Gamma \subseteq \Sigma$ is a closed 
subsystem, are subgroups of maximal rank. Once more by 
\cite[Theorem~$3.5.3$(i)]{C2}, centralizers of semisimple elements 
in~$\mathbf{G}$ are of this latter form. 
We will, therefore, restrict the following considerations to subgroups 
conjugate to $\mathbf{L}_\Gamma$ for closed subsystems~$\Gamma$.

Let $\Gamma \subseteq \Sigma$ be a closed subsystem and put $\mathbf{L} = 
\mathbf{L}_{\Gamma}$. As $N_{\mathbf{G}}( \mathbf{L} ) = 
(N_{\mathbf{G}}( \mathbf{T}_0 ) \cap N_{\mathbf{G}}( \mathbf{L} ) )\mathbf{L}$, 
we have
\begin{eqnarray*}
W_{\mathbf{G}}( \mathbf{L} ) & = & 
(N_{\mathbf{G}}( \mathbf{T}_0 ) \cap 
	N_{\mathbf{G}}( \mathbf{L} ) )\mathbf{L}/\mathbf{L} \\
& \cong & (N_{\mathbf{G}}( \mathbf{T}_0 ) \cap 
	N_{\mathbf{G}}( \mathbf{L} ) )/N_{\mathbf{L}}( \mathbf{T}_0 ) \\
& \cong & \left[(N_{\mathbf{G}}( \mathbf{T}_0 ) \cap 
	N_{\mathbf{G}}( \mathbf{L} ))/\mathbf{T}_0\right]/\left[N_{\mathbf{L}}( \mathbf{T}_0 )/\mathbf{T}_0\right] 
\end{eqnarray*}
Let $\Gamma_0$ denote a base of~$\Gamma$. Then $W_{\mathbf{G}}( \mathbf{L} ) = 
\Stab_W( \Gamma_0 )$. Moreover, $N_{\mathbf{L}}( \mathbf{T}_0 )/\mathbf{T}_0 = 
W_{\Gamma}$, and thus the image of
$N_{\mathbf{G}}( \mathbf{T}_0 ) \cap N_{\mathbf{G}}( \mathbf{L} )$ in~$W$ is the
split extension $W_{\Gamma}.\Stab_W( \Gamma_0 )$. As all bases of $\Gamma$ are
conjugate in $W_{\Gamma}$, we also have
$W_{\Gamma}.\Stab_W( \Gamma_0 ) = \Stab_W( \Gamma )$.

The $G$-conjugacy classes of $F$-stable subgroups of~$\mathbf{G}$ which are 
$\mathbf{G}$-conjugate to~$\mathbf{L}$ are in bijection with the conjugacy classes 
of $W_{\mathbf{G}}( \mathbf{L} ) = \Stab_W( \Gamma_0 )$. (Recall that~$F$ acts
trivially on~$W$.)
This bijection is determined as follows. Let~$C$ be a conjugacy class of 
$\Stab_W( \Gamma_0 )$. Choose an element 
$v \in W_{\Gamma}.\Stab_W( \Gamma_0 )$, whose image in $\Stab_W( \Gamma_0 )$
lies in~$C$. Next, choose an inverse image
$\dot{v} \in N_{\mathbf{G}}( \mathbf{T}_0 ) \cap N_{\mathbf{G}}( \mathbf{L} )$ 
of~$v$, and an
element $g \in \mathbf{G}$ such that $F(g)g^{-1} = \dot{v}$.
Then~$C$ is mapped to the $G$-conjugacy class of the group $\mathbf{L}^g$. 
If $\mathbf{M}$ is an $F$-stable
$\mathbf{G}$-conjugate of $\mathbf{L}$ that corresponds to 
$v \in W_{\Gamma}.\Stab_W( \Gamma_0 )$ in the above sense, we call 
$(\Gamma,[v])$ the $F$-\textit{type} of~$\mathbf{M}$, where~$[v]$ denotes the
conjugacy class in $\Stab_W( \Gamma_0 )$ of the coset $W_{\Gamma}v$. We omit 
the~$F$ from the 
notation if it is clear from the context. Notice that~$\mathbf{M}$ is a regular 
subgroup of~$\mathbf{G}$, if and only if $\Gamma$ is a parabolic subsystem 
of~$\Sigma$, and that~$\mathbf{M}$ is a maximal torus, if and only if~$\Gamma$ 
is the empty set.

\addtocounter{thm}{3}
\begin{lem}
\label{TypesOfMaximalRankSubgroupsUnderPowersOfSigma}
Let~$\mathbf{M}$ denote an $F$-stable connected reductive subgroup
of~$\mathbf{G}$ of maximal rank which is $\mathbf{G}$-conjugate 
to~$\mathbf{L}_\Gamma$ for some closed subsystem $\Gamma \subseteq \Sigma$. 
Suppose that the $F$-type of~$\mathbf{M}$ equals $(\Gamma,[v])$.

Then, for any positive integer~$m$, the $F^m$-type of~$\mathbf{M}$ equals
$(\Gamma, [v^m] )$.
\end{lem}
\begin{prf}
We have $\mathbf{M} = \mathbf{L}_{\Gamma}^g$ for some 
$g \in \mathbf{G}$ such that $\dot{v} := F(g)g^{-1} \in 
N_{\mathbf{G}}( \mathbf{T}_0 ) \cap N_{\mathbf{G}}( \mathbf{L}_{\Gamma} )$
maps to~$v$ under the natural homomorphism. Clearly, $\mathbf{T}_0^g$ 
and~$\mathbf{M}$
are $F^m$-stable and thus $F^m(g)g^{-1} \in N_{\mathbf{G}}( \mathbf{T}_0 )
\cap N_{\mathbf{G}}( \mathbf{L}_{\Gamma } )$.
Using the fact that $F( \dot{v} ) \in \mathbf{T}_0\dot{v}$, one
shows by induction on~$m$ that $F^m(g)g^{-1} = \dot{v}^mt$ for some
$t \in \mathbf{T}_0$, which proves our assertion.
\end{prf}

\addtocounter{subsection}{1}
\subsection{Class types}
\label{ClassTypes}
Our description of the blocks is based upon a classification of the semisimple 
conjugacy classes of~$G$. These were first parameterized in~\cite{Shoji}, but we 
need more precise information on centralizers of semisimple elements, as given
in~\cite{LL}. Let~$s$ and~$s'$ be two semisimple elements of~$G$. We say 
that~$s$ and~$s'$ belong to the same $G$-class type, respectively 
$\mathbf{G}$-class type, if and only if~$C_{\mathbf{G}}(s)$ and
$C_{\mathbf{G}}(s')$ are conjugate in~$G$, respectively $\mathbf{G}$. In the
tables of~\cite{LL} corresponding to~$F_4(q)$, the $G$-class types are
labeled by triples $(i,j,k)$, where the first index,~$i$, distinguishes the
$\mathbf{G}$-class types, and the second index,~$j$, always takes the value~$1$,
owing to the fact that centralizers of semisimple elements
in~$\mathbf{G}$ are connected. In the following, we will omit the index~$j$, and
talk of the $G$-class type $(i,k)$ instead of $(i,1,k)$. 

The first index~$i$ runs from~$1$ to~$20$, and the second index~$k$ depends 
on~$i$. Not every pair $(i,k)$ of indices occurring in the tables in~\cite{LL} 
corresponds to a semisimple element of~$G$ of class type~$(i,k)$. For example, 
if~$q$ is even, there are no elements of class type~$(i,k)$ for 
$i \in \{ 2, 3, 5, 8, 11, 12, 16 \}$. On the other hand, for each~$(i,k)$, 
there is a power~$q'$ of a prime~$p'$, possibly $p' \neq p$, such that 
$F_4( q' )$ has a semisimple element of class type $(i,k)$.

This indicates that there is a generic description of the class types. Indeed, 
first notice the root datum $(X, \Sigma, Y, \Sigma^\perp)$ is generic, 
i.e.\ independent of~$p$, up to isomorphism of root data. (In fact, 
$(X, \Sigma, Y, \Sigma^\perp)$ is part of the \textit{generic finite reductive 
group} $(X, \Sigma, Y, \Sigma^\perp, W\text{\id})$ as defined in 
\cite[Definition in 1.A]{BrouMaMi}; see also \cite[Definition~$1.6.10$]{GeMa}). 
Up to conjugation in~$W$, there are~$19$ subsets $\Gamma_i \leq \Sigma$, 
$2 \leq i \leq 20$ such that $\bar{\Gamma}_i^\perp \neq 0$, where 
$\bar{\Gamma}_i^\perp = \{ \gamma \in Y \mid 
\langle \alpha, \gamma \rangle = 0 \text{\ for all\ } 
\alpha \in \bar{\Gamma}_i \}$, and such that $\Gamma_i$ is a base of 
$\bar{\Gamma}_i$. We choose notation such that $\Gamma_{20} = 
\emptyset$ and put $\Gamma_1 = \{ \alpha_1, \ldots , \alpha_4 \}$. After the 
choice of a 
prime~$p$, i.e.\ the group~$\mathbf{G} = F_4( \bar{\mathbb{F}}_p )$, and a 
maximal torus~$\mathbf{T}_0 \leq \mathbf{G}$ giving rise to the root datum
$(X, \Sigma, Y, \Sigma^\perp)$, one can construct the subgroups
$\mathbf{L}_{\Gamma_i} \leq \mathbf{G}$ as in Subsection~\ref{SetupF4}. Then 
$C_{\mathbf{G}}( s )$ is $\mathbf{G}$-conjugate to one of the 
$\mathbf{L}_{\Gamma_i}$ for every semisimple element $s \in \mathbf{G}$.
However, depending on~$p$, the group $\mathbf{L}_{\Gamma_i}$ can have trivial 
center.

Given $i \in \{ 1, \ldots , 20 \}$, the index~$k$ numbers the conjugacy 
classes of $\Stab_W( {\Gamma}_i )$. After the choice of a power~$q$ of~$p$, 
there is a Frobenius morphism~$F$ of~$\mathbf{G}$ such that $\mathbf{T}_0$ is 
contained in an $F$-stable Borel subgroup, and $G = \mathbf{G}^F = F_4(q)$. By
the results summarized in Subsection~\ref{MaximalRank}, the $G$-conjugacy 
classes of the $F$-stable $\mathbf{G}$-conjugates of the groups 
$\mathbf{L}_{\Gamma_i}$ are labeled by the pairs $(i,k)$. Let us write 
$\mathbf{M}_{i,k}$ for a representative of the corresponding $G$-conjugacy 
class of connected reductive subgroups of~$\mathbf{G}$ of maximal rank, 
adopting the convention that $\mathbf{M}_{i,1} := \mathbf{L}_{\Gamma_i}$. Then 
$C_{\mathbf{G}}( s )$ is $G$-conjugate to one of the $\mathbf{M}_{i,k}$ for
every semisimple element $s \in G$, in which case we say that~$s$ has class type 
$(i,k)$. However, depending on~$q$, the 
group $Z( \mathbf{M}_{i,k} )$ does not necessarily contain elements of~$G$ with 
centralizer~$\mathbf{M}_{i,k}$. Notice that, unless $k = 1$, the group 
$\mathbf{M}_{i,k}$ depends on~$F$, although the index $(i,k)$ does not. Notice 
also that a
semisimple element $s \in G$ of class type $(i,k)$ can have a different class 
type when viewed as element of $\mathbf{G}^{F^m}$ for a positive integer~$m$; 
see Lemma~\ref{TypesOfMaximalRankSubgroupsUnderPowersOfSigma}. We therefore
sometimes speak of the $G$-class type of~$s$, respectively the
$\mathbf{G}^{F^m}$-class type of~$s$ to be precise.

The quasi-isolated semisimple elements of~$G$ are exactly those corresponding
to $i = 1, \ldots , 5$. The centralizers of the other $F$-stable semisimple
elements are regular, unless $i \in \{ 8, 11, 12, 16 \}$. The trivial element
is of class type~$1$, and the regular semisimple elements are of 
$\mathbf{G}$-class type~$20$.

In the first column of Table~\ref{CT} we list the sets $\Gamma_i$ for 
$2 \leq i \leq 19$; these sets are the same as in \cite{LL}, up to two 
modifications for $i = 14, 15$. For each pair $(i,k)$, we also give 
representatives $v \in W_{\bar{\Gamma}_i}.\Stab_W( \Gamma_i )$ for the conjugacy
class of $\Stab_W( \Gamma_i ) =  
W_{\bar{\Gamma}_i}.\Stab_W( \Gamma_i )/W_{\bar{\Gamma}_i}$
with label~$(i,k)$. This labeling is the same as in \cite{LL}. We usually give 
several values of $v \in W_{\bar{\Gamma}_i}.\Stab_W( \Gamma_i )$ for a given 
$(i,k)$, to have more flexibility in the proofs of 
Section~\ref{ActionOnWeights}. We take one of the groups thus constructed (after 
a choice of a lift $\dot{v}$ and an element $g \in \mathbf{G}$ with $F(g)g^{-1} 
= \dot{v}$) as our representative~$\mathbf{M}_{i,k}$. The $W$-conjugacy classes 
of the elements in the coset $W_{\bar{\Gamma}_i}v$ determine the $G$-conjugacy 
classes of the $F$-stable maximal tori which have some representative 
in~$\mathbf{M}_{i,k}$.

\subsection{Duality}\label{Duality}
We will identify the dual group $\mathbf{G}^*$ with~$\mathbf{G}$, and thus~$G^*$
with~$G$. For later purposes, we will choose a specific identification. The dual 
root datum $( Y, \Sigma^\vee , X, \Sigma )$ is isomorphic to 
$( X, \Sigma , Y, \Sigma^\vee )$ via an isomorphism $\delta : X \rightarrow Y$ 
satisfying $\delta( \alpha_i) = \alpha_{5-i}^\vee$ for $i = 1, \ldots , 4$.
As~$\delta$ maps~$\Sigma$ to~$\Sigma^\vee$, this yields a permutation 
$\alpha \mapsto \alpha^\dagger$ of $\Sigma$ such that $\delta( \alpha^\dagger )
= \alpha^\vee$ for all $\alpha \in \Sigma$. This permutation, in fact an 
involution, is easily determined with CHEVIE. The
CHEVIE number of the image of a positive root $\alpha$ under this permutation
is given in Table~\ref{PositiveRoots} under the heading $\alpha^\dagger$.
The map $s_{\alpha} \mapsto s_{\alpha^\dagger}$ extends to an automorphism
$w \mapsto w^\dagger$ of~$W$.
By the isomorphism theorem \cite[Theorem~$9.6.2$]{springer}, there is an 
$F$-equivariant isomorphism $\mathbf{G} \rightarrow 
\mathbf{G}^*$ inducing the isomorphism~$\delta$ of root data as in 
\cite[$9.6.1$]{springer}, which we use to identify~$\mathbf{G}$ 
with~$\mathbf{G}^*$.

In view of~(\ref{iota}),~(\ref{TF}) and (\ref{XotimesF}), the isomorphism 
$\delta : X \rightarrow Y$ gives rise to a $W$-equivariant isomorphism 
$T_0 \rightarrow \Irr( T_0 ), s \mapsto \hat{s}$.

\addtocounter{thm}{2}
\begin{lem}
\label{KernelsOfDualCharacters}
Let $s \in T_0$, and let $\hat{s} \in \Irr(T_0)$ denote the irreducible
character arising from duality. Further, let $\alpha \in \Sigma$ such that 
$\alpha^\dagger(s) = 1$. Then $\alpha^{\vee}( t ) \in \ker( \hat{s} )$ for 
all $t \in \mathbb{F}_q^*$.
\end{lem}
\begin{prf}
Let $\gamma \in Y$ such that~$s$ corresponds to $\gamma +(F-1)Y$ 
under~(\ref{TF}). Then $\alpha^\dagger( s ) = 
\exp(2\pi\sqrt{-1} \iota^{-1}(\langle \alpha^\dagger, \gamma \rangle))$ 
by~(\ref{EvaluationOfCharacter}), and thus 
$\langle \alpha^\dagger, \gamma \rangle \in \mathbb{Z}$ by assumption.

Let $\chi \in X$ with $\delta( \chi ) = \gamma$. Then $\chi \otimes 1$ is in the
kernel of $F - 1$ on $X \otimes \mathbb{Q}_{p'}/\mathbb{Z}$ as~$s$ is 
$F$-stable. By definition, $\hat{s}$ is the character of $T_0$ which corresponds 
to $\chi \otimes 1$ under the isomorphism~(\ref{XotimesF}). The inverse image
of $\alpha^\vee( t )$ under~(\ref{TF}) is an element of the form 
$m\alpha^\vee + (F - 1)Y \in Y/(F-1)Y$ for some $m \in \mathbb{Z}$.
The claim now follows from Equation~(\ref{EvaluationOfCharacter}) as  
$\langle \chi, m\alpha^\vee \rangle = 
m\langle \alpha^\dagger, \delta(\chi) \rangle \in \mathbb{Z}$.
\end{prf}

\addtocounter{subsection}{1}
\subsection{Twisting and duality}\label{TwistingAndDuality}
We record a basic fact about duality. Let $g, g^* \in \mathbf{G}$ such that 
$\dot{w} := F(g)g^{-1}$ and $\dot{w}^* := F(g^*){g^*}^{-1}$ 
normalize~$\mathbf{T}_0$, and write~$w$, $w^*$ for the images of $\dot{w}$ and 
$\dot{w}^*$ in~$W$. Put $\mathbf{T} := \mathbf{T}_0^g$ and $\mathbf{T}^* := 
\mathbf{T}_0^{g^*}$. Let $\delta : X( \mathbf{T}_0 ) \rightarrow 
Y( \mathbf{T}_0 )$ denote the duality isomorphism introduced in 
Subsection~\ref{Duality}. We then have the following commutative diagram
\begin{equation}
\label{DualityDiagram}
\begin{xy}
\xymatrix@C+30pt{
X( \mathbf{T} ) \ar[r]^{\delta_{g^*,g}} & Y( \mathbf{T}^* ) \\
X( \mathbf{T}_0 ) \ar[u]^{\CJxxx_g} \ar[r]_{\delta} &
Y( \mathbf{T}_0 ) \ar[u]_{\CJxxx_{g^*}}
}
\end{xy}
\end{equation}
with $\delta_{g^*,g} = \CJxxx_{g^*} \circ \delta \circ \CJxxx_g^{-1}$,
where the maps on the vertical arrows are induced by the conjugation
maps $\CJxxx_g: \mathbf{T}_0 \rightarrow \mathbf{T}$
and $\CJxxx_{g^*}: \mathbf{T}_0 \rightarrow \mathbf{T}^*$.

Suppose further that $w^\dagger = {w^*}^{-1}$. Then $ \delta_{g^*,g}$ is an 
isomorphism between $X( \mathbf{T} )$ and $Y( \mathbf{T}^* )$,
and $F \circ \delta_{g^*,g} = \delta_{g^*,g} \circ F$,
i.e.\ $(\mathbf{T}, F)$ and $(\mathbf{T}^*, F)$ are in duality;
see \cite[Proposition~$4.3.4$]{C2}.
%%%%%%%%%%%%%%%%%%%%%%%%%%%%%%%%%%%%%%%%%%%%%%%%%%%%%%%%%%%%%%%%%%%%%%%%%%%%%%%%
%%%%%%%%%%%%%%%%%%%%%%%%%%%%%%%%%%%%%%%%%%%%%%%%%%%%%%%%%%%%%%%%%%%%%%%%%%%%%%%%
%%
%% See sheet 1 of 16.03.2021.
%%
%%%%%%%%%%%%%%%%%%%%%%%%%%%%%%%%%%%%%%%%%%%%%%%%%%%%%%%%%%%%%%%%%%%%%%%%%%%%%%%%
%%%%%%%%%%%%%%%%%%%%%%%%%%%%%%%%%%%%%%%%%%%%%%%%%%%%%%%%%%%%%%%%%%%%%%%%%%%%%%%%

Assume in addition that $\dot{w}$ and $\dot{w}^*$ are $F$-stable. Then
$$\omega^* := \CJxxx_{g^*} \circ F \circ \CJxxx_{g^*}^{-1}$$
and 
$$\omega := \CJxxx_{g} \circ F \circ \CJxxx_{g}^{-1}$$
are Steinberg morphisms of~$\mathbf{G}$ which commute with~$F$. Moreover,
the following diagram commutes
$$
\begin{xy}
\xymatrix@C+30pt{
X( \mathbf{T} ) \ar[r]^{\delta_{g^*,g}} & Y( \mathbf{T}^* ) \\
	X( \mathbf{T} ) \ar[u]^{\omega} \ar[r]_{\delta_{g^*,g}} &
Y( \mathbf{T}^* ) \ar[u]_{\omega^*}
}
\end{xy}
$$
where $\omega$ and $\omega^*$ denote the induced maps on $X( \mathbf{T} )$
and $\mathbf{Y}( \mathbf{T}^* )$, respectively.
In other words,
$(\mathbf{T}, \omega) $ and $(\mathbf{T}^*, \omega^*)$ are in duality.

We now generalize the above considerations to regular subgroups of~$\mathbf{G}$. 
Let~$\Gamma$ denote a parabolic subsystem of~$\Sigma$. Then $\Gamma^\dagger$ 
also is a parabolic subsystem, and $\mathbf{L}_\Gamma$ and 
$\mathbf{L}_{\Gamma^\dagger}$ are dual regular subgroups of~$\mathbf{G}$. 
Let~$\mathbf{M}$ denote a regular subgroup of~$\mathbf{G}$ of type 
$(\Gamma, [v] )$ for some $v \in \Stab_W( \Gamma_0 )$, where~$\Gamma_0$ is a
base of~$\Gamma$; see 
Subsection~\ref{MaximalRank}. If $\mathbf{M}^*$ is a regular subgroup 
of~$\mathbf{G}$ of type $(\Gamma^\dagger, [v^*] )$ with $v^\dagger = 
{v^*}^{-1}$, then the pairs $(\mathbf{M},F)$ and $(\mathbf{M}^*,F)$ are in 
duality in the sense of \cite[Definition~$1.5.17$]{GeMa}. Indeed, by 
replacing~$v$ with a suitable element $w \in W_{\Gamma}v$, putting $w^* := 
{w^\dagger}^{-1}$, and choosing~$g$ and~$g^*$ as above, we may assume 
that~$\mathbf{T}$ and $\mathbf{T}^*$ are maximally split tori of
$\mathbf{M} := \mathbf{L}_{\Gamma}^g$, respectively $\mathbf{M}^* 
:= \mathbf{L}_{\Gamma^\dagger}^{g^*}$. The duality of $(\mathbf{M},F)$ and
$(\mathbf{M}^*,F)$ follows from this and the diagram~(\ref{DualityDiagram}). 
To indicate that $\mathbf{M}^*$ is a subgroup of~$\mathbf{G}$, we usually write 
$\mathbf{M}^\dagger$ instead of~$\mathbf{M}^*$ for a group dual to~$\mathbf{M}$ 
constructed in this way. 
Thus the identification of~$\mathbf{G}$ with its dual~$\mathbf{G}^*$ induces an 
involutive bijection $\mathbf{M} \mapsto \mathbf{M}^\dagger$ between 
$G$-conjugacy classes of regular subgroups of~$\mathbf{G}$, 
where~$\mathbf{M}^\dagger$ is a regular subgroup of~$\mathbf{G}$ dual 
to~$\mathbf{M}$. Notice that $\mathbf{T}_0^\dagger = 
\mathbf{L}_{\emptyset^\dagger} = \mathbf{L}_{\emptyset} = \mathbf{T}_0$.

\subsection{Some split Levi subgroups} 
\label{SomeSplitLeviSubgroups} In this subsection, we describe some 
split Levi subgroups of~$G$ and their normalizers. 
We begin by introducing two closed subgroups of~$\mathbf{G}$ of maximal rank.

Below,~$\gamma$ denotes the lift of $w_0$ constructed in 
Lemma~\ref{LiftOfLongestElementLem}. We choose $h \in \mathbf{G}$ with 
$F(h)h^{-1} = \gamma$, let $g \in \{ 1, h \}$ and put $\mathbf{T} := 
\mathbf{T}_0^g$. We define the parameter~$\varepsilon$, by $\varepsilon := 1$ 
if $g = 1$, and $\varepsilon := -1$, otherwise. We adopt the notation that 
$\SL_3^\varepsilon( q )$ denotes $\SL_3( q )$ if $\varepsilon = 1$ and 
$\SU_3( q )$ if $\varepsilon = -1$.

\addtocounter{thm}{2}
\begin{prop}
\label{C3C}
Let $J_1 := \{ \alpha_1, \alpha_{23} \}$, $J_2 := \{ \alpha_3, \alpha_4 \}$
and $J := J_1 \cup J_2$, and put $\mathbf{L}_0 := \mathbf{L}_J$
and $\mathbf{L}_0^i := \mathbf{K}_{J_i}$, $i = 1, 2$.
Let $\mathbf{L} := \mathbf{L}_0^g$ and $\mathbf{L}^i := (\mathbf{L}_0^i)^g$
for $i = 1, 2$.

%%%%%%%%%%%%%%%%%%%%%%%%%%%%%%%%%%%%%%%%%%%%%%%%%%%%%%%%%%%%%%%%%%%%%%%%%%%%%%%%
%%%%%%%%%%%%%%%%%%%%%%%%%%%%%%%%%%%%%%%%%%%%%%%%%%%%%%%%%%%%%%%%%%%%%%%%%%%%%%%%
%%
%% $\mathbf{L} = \mathbf{M}_{4,1}$.
%%
%%%%%%%%%%%%%%%%%%%%%%%%%%%%%%%%%%%%%%%%%%%%%%%%%%%%%%%%%%%%%%%%%%%%%%%%%%%%%%%%
%%%%%%%%%%%%%%%%%%%%%%%%%%%%%%%%%%%%%%%%%%%%%%%%%%%%%%%%%%%%%%%%%%%%%%%%%%%%%%%%

{\rm (a)} Let $d = \gcd(3, q^2 - 1)$. Then 
$\mathbf{L} = \mathbf{L}^1 \circ_d \mathbf{L}^2$ with $\mathbf{L}^i \cong 
\SL_3( \mathbb{F})$, $i = 1, 2$. Moreover,
$N_{\mathbf{G}}( \mathbf{L} ) = \langle \mathbf{L}, \gamma \rangle$, 
where~$\gamma$ normalizes~$\mathbf{L}^i$, inducing the inverse-transpose
automorphism in the isomorphic copy~$\SL_3( \mathbb{F})$, $i = 1, 2$.
	
{\rm (b)} Let $d = \gcd(3, q - \varepsilon)$. Then 
$L^i \cong \SL^\varepsilon_3( q )$ for $i = 1, 2$ and
$L = \langle L^1 \circ_d L^2, x \rangle$ for some $x \in T$ satisfying the 
following properties. If $d = 1$, then $x = 1$. If $d = 3$, then $x^3 \in 
L^1 \circ_d L^2$ and~$x$ normalizes~$L^i$, inducing a diagonal automorphism 
in the isomorphic copy~$\SL_3^\varepsilon(q)$, $i = 1, 2$.

Finally, $N_G( L ) = \langle L, \gamma \rangle$. 
If $d = 3$, the element~$\gamma$ inverts~$x$, and thus 
$N_G( L )/(L^1 \circ_3 L^2)$ is isomorphic to the symmetric group on three 
letters.
\end{prop}
\begin{prf}
We only give the proof for $g = 1$. The case $g \neq 1$ can be treated by 
conjugating all relevant structures established for~$\mathbf{L}_0$ with~$g$. 
Alternatively, one can replace the pair~$(\mathbf{L},F)$ by~$(\mathbf{L}_0,F\gamma)$.

(a) We employ the notation and results summarized in~\ref{SetupF4}. Observe that 
$\bar{J}_1$ and $\bar{J}_2$ are closed, disjoint and of type~$A_2$, respectively 
$\tilde{A}_2$. (We use the common notational convention to indicate the 
irreducible closed subsystems of~$\Sigma$ of type~$A$ consisting of short roots 
by a tilde.) Moreover, $\bar{J}_1 \cup \bar{J}_2 = \bar{J}$ is closed and 
thus $[\mathbf{L}^1,\mathbf{L}^2] = \{ 1 \}$. As $\bar{J}_1 \cup \bar{J}_2$ 
has rank~$4$, we also
have $\mathbf{L} = \mathbf{L}^1\mathbf{L}^2$. Finally, $\mathbf{L}^i
\cong \SL_3( \mathbb{F} )$ for $i = 1, 2$. If $d = 3$, then $\mathbf{L}^1$
and $\mathbf{L}^2$ intersect in a group of order~$3$, as $|Z( \mathbf{L} )| = 3$
(the latter assertion can be verified by a computation with CHEVIE using
\cite[$8.1.8$]{springer}). Hence $\mathbf{L} = \mathbf{L}^1 \circ_d 
\mathbf{L}^2$ as claimed. As $\Stab_W(J) = \langle w_0 \rangle$, we obtain
$N_{\mathbf{G}}( \mathbf{L} ) = \langle \mathbf{L}, \gamma \rangle$.
As $\alpha_1$ and $\alpha_{23}$ are long, whereas $\alpha_3$ and $\alpha_4$
are short, $N_{\mathbf{G}}( \mathbf{L} )$ stabilizes~$\mathbf{L}^i$ for
$i = 1, 2$. It follows from Lemma~\ref{LiftOfLongestElementLem},
that~$\gamma$ acts on each $\mathbf{L}^i$ as inverse-transpose automorphism
(in the natural $3$-dimensional representation of~$\mathbf{L}^i$), $i = 1, 2$.
(The corresponding relation for $\pm \alpha_{23}$ is easily verified with
CHEVIE.)

(b) Clearly, $L^i \cong \SL_3(q)$, $i = 1, 2$. Moreover, $L = L^1 \times L^2$ 
if $d = 1$. Suppose that $d = 3$, and let $z \in Z( \mathbf{L} )$ be an element 
of order~$3$. Let $x_i \in \mathbf{T} \cap \mathbf{L}^i$, $i = 1, 2$, such that 
$F( x_1 )^{-1}x_1 = z = F( x_2 )x_2^{-1}$. Then $x := x_1x_2$ is $F$-stable 
and $L = \langle L^1 \circ_3 L^2, x \rangle$. By construction,~$x$ normalizes 
$L^1$ and $L^2$. Moreover,~$x$ acts as a diagonal automorphism on each of~$L^i$,
$i = 1, 2$; see Remark~\ref{DiagonalActionOfx} below. As $x_i^3 \in L^i$ for 
$i = 1, 2$, we have $x^3 \in L^1 \circ_3 L^2$. As $\gamma$ inverts the elements 
of~$\mathbf{T}$, it inverts~$x$. As $x^2 \not\in L^1 \circ_3 L^2$, the group 
$\langle L, \gamma \rangle/(L^1 \circ_3 L^2)$ is not abelian and hence 
isomorphic to the symmetric group of order~$6$.

From~(a) we get $N_G( \mathbf{L} ) = \langle L, \gamma \rangle$. As the
latter is a maximal subgroup of~$G$ by \cite{LiSaSe}, we obtain $N_G( L ) =
N_G( \mathbf{L} )$.
\end{prf}

\medskip
\noindent 
We can be more specific in the choice of the element~$x$ in 
Proposition~\ref{C3C} in case $3 \mid q^2 - 1$. This shows in particular, 
that~$x$ acts as diagonal automorphism on each of~$L^1$ and $L^2$. 

\begin{rem}
\label{DiagonalActionOfx}
{\rm
Keep the notation and assumptions of Proposition \ref{C3C} and assume that 
$3 \nmid q$, i.e.\ $3 \mid q^2 - 1$. Let~$e$ be the minimal positive integer 
such that $\SL_3^\varepsilon(q) \leq \GL_3(q^e)$. Then $e \in \{ 1, 2 \}$.
We may think of $\mathbf{L}^1$ in its natural representation and of~$F$ acting
on~$\mathbf{L}^1$ as standard Frobenius endomorphism, if $\varepsilon = 1$,
and as standard Frobenius endomorphism followed by the inverse-transpose 
automorphism if $\varepsilon = -1$. We have $z =
\diag( \zeta, \zeta, \zeta)$ for $\zeta \in \mathbb{F}_{q^e}$ a third root
of~$1$. Let $\xi \in \mathbb{F}$ be of $3$-power order with 
$\xi^{1 - \varepsilon q} = \zeta$. Then 
$x_1 := \diag( \xi^{-1},\xi^{-1},\xi^2 )$ satisfies
$F( x_1 )^{-1}x_1 = z$. Also, $\xi^3 \in \mathbb{F}_{q^e}$, 
and~$x_1$ acts on $\SL^\varepsilon_3( q )$ as conjugation by 
$\diag( 1, 1, \xi^3) \in \GL_3( q^e )$,
i.e.\ as a diagonal automorphism. If $3^a$ is the $3$-part of $q - \varepsilon$, 
then $|\xi^3| = 3^a$. An analogous choice can be made for~$x_2$.
	
Recall that $x_i^3 \in T \cap L^i$ for $i = 1, 2$. Hence $x_2^3 = x_1^{-3}x 
\in \langle L^1, x \rangle \cap (T \cap L^2)$. Now $L_{\{ 2, 23 \}} = 
\langle L^1, T \cap L^2, x \rangle$ and $Z( L_{\{ 2, 23 \}} ) =  T \cap L^2$. 
Thus $L_{\{ 2, 23 \}} = 
\langle L^1, x \rangle \circ_{[3^a]} Z( L_{\{ 2, 23 \}} )$. (See also
Proposition~\ref{NormalizerM13} below.)
}\hfill{$\Box$}
\end{rem}

\medskip
\noindent
We also note that the group~$\mathbf{L}$ defined in Proposition~\ref{C3C} is
selfdual.

\begin{rem}
\label{SelfDualityOfL}
{\rm
Let $J$ and $\mathbf{L} = \mathbf{L}_J$ be as in Proposition~\ref{C3C}. 
Then~$\mathbf{L}$ is self dual. Indeed, the $\mathbb{Q}$-linear isomorphism 
$X \otimes \mathbb{Q} \rightarrow Y \otimes \mathbb{Q}$, determined by mapping 
the four-tuple $(\alpha_{1}, \alpha_{23}, \alpha_3, \alpha_4)$ of roots to the 
four-tuple $(\alpha_{23}^\vee, \alpha_1^\vee, \alpha_3^\vee, \alpha_4^\vee)$ 
of coroots, restricts to a $\mathbb{Z}$-linear isomorphism $X \rightarrow Y$, 
which defines an isomorphism of the root data of~$\mathbf{L}$ and its dual group.
}\hfill{$\Box$}
\end{rem}

\medskip

\noindent
In the following, if $M \leq G$ and $\chi \in \Irr( M )$, we write 
$N_G( M, \chi )$ for the stabilizer of~$\chi$ in $N_G( M )$.

\begin{prop}
\label{NormalizerM13}
Let $J_1 := \{ \alpha_1, \alpha_{23} \}$, $J_2 := \{ \alpha_3, \alpha_4 \}$,
and let $\mathbf{M}_0 := \mathbf{L}_J$ with $J \in \{ J_1, J_2 \}$. Put 
$\mathbf{M} := \mathbf{M}_0^g$ and $\mathbf{M}' := [\mathbf{M},\mathbf{M}]$.
%%%%%%%%%%%%%%%%%%%%%%%%%%%%%%%%%%%%%%%%%%%%%%%%%%%%%%%%%%%%%%%%%%%%%%%%%%%%%%%%
%%%%%%%%%%%%%%%%%%%%%%%%%%%%%%%%%%%%%%%%%%%%%%%%%%%%%%%%%%%%%%%%%%%%%%%%%%%%%%%%
%%
%% $\mathbf{M} = \mathbf{M}_{13,k}$ with $k = 1$ or~$6$ according as
%% $\varepsilon = 1$ or~$-1$.
%%
%% or $\mathbf{M} = \mathbf{M}_{17,k}$ with $k = 1$ or~$6$ according as
%% $\varepsilon = 1$ or~$-1$.
%%
%%%%%%%%%%%%%%%%%%%%%%%%%%%%%%%%%%%%%%%%%%%%%%%%%%%%%%%%%%%%%%%%%%%%%%%%%%%%%%%%
%%%%%%%%%%%%%%%%%%%%%%%%%%%%%%%%%%%%%%%%%%%%%%%%%%%%%%%%%%%%%%%%%%%%%%%%%%%%%%%%

{\rm (a)} Let $d = \gcd(3, q^2 - 1)$. Then 
$\mathbf{M} = Z(\mathbf{M}) \circ_d \mathbf{M}'$, with 
$\mathbf{M}' \cong \SL_3( \mathbb{F} )$. Moreover,
$$N_{\mathbf{G}}( \mathbf{M} ) = 
\langle \mathbf{N}' \circ_d \mathbf{M}', \gamma \rangle = 
\mathbf{N}\mathbf{M}',$$
with $Z( \mathbf{M} ) \leq \mathbf{N}' \leq \mathbf{N} \leq 
N_{\mathbf{G}}( \mathbf{T} )$, where $\mathbf{N} = 
\langle \mathbf{N'}, \gamma \rangle$. Furthermore, 
$\mathbf{N}' \cap \mathbf{M} =
\mathbf{N} \cap \mathbf{M} = Z( \mathbf{M} )$ and 
$\mathbf{N}'/Z(\mathbf{M}) \cong S_3$ and 
$\mathbf{N}/Z(\mathbf{M}) \cong D_{12}$. 
Finally, $\mathbf{N}' = C_{N_{\mathbf{G}}( \mathbf{M} )}( \mathbf{M}' )$.
In particular,
$\mathbf{N}'$ is normal in $N_{\mathbf{G}}( \mathbf{M} )$.

{\rm (b)} Let $d = \gcd(3, q - \varepsilon)$. Then 
$M = \langle Z( M ) \circ_d M', x \rangle$
with $x = 1$, if $d = 1$, and $x \in T$ as in {\rm Proposition~\ref{C3C}}, if 
$d = 3$. Moreover,
$$N_{G}( \mathbf{M} ) =
\langle N' \circ_d M', x, \gamma \rangle = NM$$
with $N'/Z(M) \cong S_3$ and $N/Z(M) \cong D_{12}$.

{\rm (c)} If $N_G( \mathbf{M} ) = N_G( M )$, every $\chi \in \Irr(M)$ extends to 
$N_{G}( M, \chi )$.
\end{prop}
\begin{prf}
We only consider the case $J = J_2$; the other case is proved analogously. We
also assume $g = 1$. The case $g \neq 1$ can be treated by conjugating all 
relevant structures established for~$\mathbf{M}_0$ with~$g$. Alternatively, 
one can replace the pair $(\mathbf{M},F)$ by~$(\mathbf{M}_0, F\gamma)$.

(a) We use the notation of Proposition~\ref{C3C}. Let $\mathbf{M} = 
\mathbf{L}_{J_2}$ and put $\mathbf{K} := \mathbf{K}_{J_1}$. Then 
$Z( \mathbf{M} ) = \mathbf{T} \cap \mathbf{K}$.  The assertion about the structure 
of~$\mathbf{M}$ follows from Proposition~\ref{C3C}.
To establish the claim about $N_{\mathbf{G}}(\mathbf{M})$, 
write $\mathbf{N'} :=
N_{\mathbf{K}}( Z( \mathbf{M} ) )$. As~$\mathbf{N}'$ normalizes
$Z( \mathbf{M} )$, it also normalizes $\mathbf{M} =
C_{\mathbf{G}}( Z( \mathbf{M} ) )$. 
Moreover, $\mathbf{N}'$ centralizes~$\mathbf{M}'$ as $[\mathbf{M},\mathbf{M}] =
\mathbf{K}_{J_2}$ and $[\mathbf{K},\mathbf{K}_{J_2}] = 1$.
Also, $\mathbf{N}' \cap \mathbf{M} = Z( \mathbf{M} )$ and 
$\mathbf{N}'/Z( \mathbf{M} ) \cong W_{J_1} \cong S_3$.
Now $W_{J_1}$ is a subgroup of $\Stab_W( J_2 )$, which is a
dihedral group of order~$12$; see \cite[p.~$74$]{How}. In particular,
$\mathbf{N}' \circ_d \mathbf{M}'$ is a subgroup of $N_\mathbf{G}( \mathbf{M} )$ 
of index $2$. As $\gamma \in N_\mathbf{G}( \mathbf{M} ) \setminus 
\mathbf{N}' \circ_d \mathbf{M}'$, and~$\gamma$ does not 
centralize~$\mathbf{M}'$, we obtain all our claims.

(b) The structure of~$M$ follows from Proposition~\ref{C3C}.
Now $x \in T$ normalizes $\mathbf{K}$ and $Z( \mathbf{M} )$,
hence~$N'$. In particular, $\langle N' \circ_d M', x \rangle
= N'M$, and $N'$ and~$M$ are invariant under $\gamma$, giving
the structure of $N_G( \mathbf{M} )$. As $F$ acts trivially on~$W$, we 
get the assertions on $N'/Z(M)$ and $N/Z(M)$.

(c) Let $\chi \in \Irr( M )$. In the considerations to follow, we will make use
of the facts summarized in~\ref{CentralProducts} for characters of central
products. As~$M$ is a central product $Z( M ) \circ \langle M', x \rangle$,
we have $\chi = \lambda \psi$ for $\lambda \in \Irr( Z(M) )$ and 
$\psi \in \Irr( \langle M', x \rangle )$. Let~$\chi'$ and~$\psi'$ denote the 
restrictions of~$\chi$ to $Z( M ) \circ_d M'$ and of~$\psi$ to~$M'$, 
respectively. Then $\chi' = \lambda \psi'$. Put 
$H := \langle N' \circ_d M', \gamma \rangle = NM'$. Then 
$H \cap M = Z(M) \circ_d M'$ and $H M = N_G(M)$. 

Let $I := N \cap N_G( M, \chi )$ and $I' := N' \cap N_G( M, \chi )$. Then
$N_G( M, \chi ) = IM$ and $N_G( M, \chi )/M \cong I/Z(M)$. We may 
assume that $4 \mid |I/Z(M)|$, as otherwise the Schur multiplier of 
$I/Z(M)$ is trivial. Clearly,~$I$ stabilizes~$\lambda$.
To determine the action of~$N$ on $Z(M)$, notice that $Z(\mathbf{M})$ is a
maximal $1$-$F$-split torus of $\mathbf{K} \cong \SL_3( \mathbb{F} )$,
that $N'/Z(M) \cong N_{\mathbf{K}}( Z( \mathbf{M} ) )/Z( \mathbf{M} )$,
and that~$\gamma$ acts as inverse-transpose automorphism on~$\mathbf{K}$.
In particular,~$N$ only fixes the trivial element of~$Z(M)$.
By Brauer's permutation lemma,~$N$ only fixes the trivial character of~$Z(M)$.
Thus~$\lambda$ is the trivial character if $I = N$. In this case, we 
define $\hat{\chi} \in \Irr( N'M )$ by $\hat{\chi}( n'm ) = \chi( m )$ for 
$n' \in N', m \in M$; it follows from~(a) that~$M$ normalizes~$N'$, 
so that $\hat{\chi}$ is well defined. 
Then~$\hat{\chi}$ is an extension of~$\chi$ to $N'M$. 
Moreover,~$\hat{\chi}$ is invariant in $N_G(M) = \langle N'M, \gamma \rangle$
as~$\gamma$ stabilizes~$N'$ and~$\chi$, and so $\hat{\chi}$ extends to $N_G(M)$.

Now suppose that $|I/Z(M)| = 4$. Then $I/Z(M) = \langle s, \gamma \rangle$, 
with $s \in \{ s_1, s_{23}, s_{23}^{s_1} \}$. We only treat the case $s = s_1$; 
the other cases are handled in an analogous way (or by symmetry). 
As~$s_1$ stabilizes~$\lambda$, we conclude that $\alpha_1^\vee( z ) \in 
\ker( \lambda )$ for all $z \in \mathbb{F}_q^*$. In particular, $n_1^2 =
\alpha_1^\vee( -1 ) \in \ker( \lambda )$. 
Let~$\lambda'$ denote an extension of~$\lambda$ to~$I'$.
Now $n_1 \in I'$, as $n_1 \in N'$, and thus $\lambda'( n_1 ) = \pm 1$. 
Since~$\gamma$ inverts~$n_1$, it follows that~$\lambda'$ is invariant 
under~$\gamma$. As~$I$ stabilizes~$\psi'$, and either~$\psi'$
is irreducible or has exactly three irreducible constituents, we may
choose an irreducible constituent~$\vartheta$ of~$\psi'$ which is invariant 
under~$I$. Hence $\lambda' \vartheta \in \Irr( I' \circ_d M' )$ is invariant 
in~$IM'$ and thus extends to~$IM'$.

If~$\vartheta = \psi'$, i.e.\ if $\chi'$ is irreducible; then~$\chi$ extends to 
$N_G(M,\chi)$ by \cite[Lemma~$4.1$(a)]{Spaeth10}, as $\lambda \psi'$ extends to
its inertia group~$IM'$ in~$H$. 
Suppose now that $\vartheta \neq \psi'$. Then
$d = 3$ and~$H$ has index~$3$ in $N_G( M )$. 
Now $\lambda' \vartheta$ is not invariant under~$x$,
and thus induces to an irreducible character $\hat{\chi}'$ of
$I'M$ which extends~$\chi$ (notice that~$x$
fixes~$\lambda$ and normalizes~$N'$, so that~$x$ also
normalizes~$I'$). As~$\lambda' \vartheta$ is invariant under~$\gamma$,
the same is true for $\hat{\chi}'$, which thus extends to $IM$.
\end{prf}
%%%%%%%%%%%%%%%%%%%%%%%%%%%%%%%%%%%%%%%%%%%%%%%%%%%%%%%%%%%%%%%%%%%%%%%%%%%%%%%%
%%%%%%%%%%%%%%%%%%%%%%%%%%%%%%%%%%%%%%%%%%%%%%%%%%%%%%%%%%%%%%%%%%%%%%%%%%%%%%%%
%%
%% See sheets 1 - 4 of 30.09.2020
%%
%%%%%%%%%%%%%%%%%%%%%%%%%%%%%%%%%%%%%%%%%%%%%%%%%%%%%%%%%%%%%%%%%%%%%%%%%%%%%%%%
%%%%%%%%%%%%%%%%%%%%%%%%%%%%%%%%%%%%%%%%%%%%%%%%%%%%%%%%%%%%%%%%%%%%%%%%%%%%%%%%
%%%%%%%%%%%%%%%%%%%%%%%%%%%%%%%%%%%%%%%%%%%%%%%%%%%%%%%%%%%%%%%%%%%%%%%%%%%%%%%%
%%%%%%%%%%%%%%%%%%%%%%%%%%%%%%%%%%%%%%%%%%%%%%%%%%%%%%%%%%%%%%%%%%%%%%%%%%%%%%%%
%%
%% See sheets 1, 2 of 03.03.2020
%%
%%%%%%%%%%%%%%%%%%%%%%%%%%%%%%%%%%%%%%%%%%%%%%%%%%%%%%%%%%%%%%%%%%%%%%%%%%%%%%%%
%%%%%%%%%%%%%%%%%%%%%%%%%%%%%%%%%%%%%%%%%%%%%%%%%%%%%%%%%%%%%%%%%%%%%%%%%%%%%%%%

\medskip
\noindent
Notice that $\{ \alpha_1, \alpha_2 \}$ is conjugate in~$W$ to  
$\{ \alpha_1, \alpha_{23} \}$, so that the results established in 
Proposition~\ref{NormalizerM13} hold likewise for 
$\mathbf{L}_{\{ \alpha_1, \alpha_2 \}}$.

\begin{prop}
\label{NormalizerM14}
Let $\mathbf{M}_0 := \mathbf{L}_J$, with $J = \{ \alpha_1, \alpha_4 \}$ and put 
$\mathbf{M} := \mathbf{M}_0^g$.

{\rm (a)} Let $d := \gcd( 2, q - 1 )$.
%%%%%%%%%%%%%%%%%%%%%%%%%%%%%%%%%%%%%%%%%%%%%%%%%%%%%%%%%%%%%%%%%%%%%%%%%%%%%%%%
%%%%%%%%%%%%%%%%%%%%%%%%%%%%%%%%%%%%%%%%%%%%%%%%%%%%%%%%%%%%%%%%%%%%%%%%%%%%%%%%
%%
%% $\mathbf{K} = \mathbf{M}_{14,k} = $ with $k = 1$ or~$4$ according as
%% $\varepsilon = 1$ or~$-1$.
%%
%%%%%%%%%%%%%%%%%%%%%%%%%%%%%%%%%%%%%%%%%%%%%%%%%%%%%%%%%%%%%%%%%%%%%%%%%%%%%%%%
%%%%%%%%%%%%%%%%%%%%%%%%%%%%%%%%%%%%%%%%%%%%%%%%%%%%%%%%%%%%%%%%%%%%%%%%%%%%%%%%
Then $\mathbf{M} = \mathbf{M}_1 \circ_d \mathbf{M}_2$ with connected, reductive,
$F$-stable subgroups~$\mathbf{M}_i$, $i = 1, 2$. Putting $\mathbf{M}_i' := 
[\mathbf{M}_i, \mathbf{M}_i]$, we have 
$\mathbf{M}_i' \cong \SL_2( \mathbb{F} )$, and $Z( \mathbf{M}_i )^\circ$ is a 
torus of rank~$1$, for $i = 1, 2$. Moreover, 
$\mathbf{M}_1 = Z( \mathbf{M}_1 )^\circ \times \mathbf{M}_1'$ and 
$\mathbf{M}_2 \cong \GL_2( \mathbb{F} )$. In particular, $Z( \mathbf{M}_2 )$
is connected and $\mathbf{M}_2 = Z( \mathbf{M}_2 ) \circ_d \mathbf{M}_2'$. 

Furthermore, there are lifts $m_1, m_2 \in \hat{W}^g$ of $s_{22}^g$ and 
$s_{17}^g$ {\rm (see Subsection~\ref{LiftOfLongestElement})}, respectively,
such that~$m_i$ normalizes~$\mathbf{M}_i$ for $i = 1, 2$, and such
that $m_1$ centralizes $\mathbf{M}_1'\mathbf{M}_2$ and $m_2$ centralizes
$\mathbf{M}_1\mathbf{M}_2'$. For any such pair of elements we have
$$
N_\mathbf{G}( \mathbf{M} ) = \mathbf{M}_1.2 \circ_d \mathbf{M}_2.2,
$$
with $\mathbf{M}_i.2 = \langle \mathbf{M}_i, m_i \rangle$ for $i = 1, 2$.

{\rm (b)} With the notation of~{\rm (a)}, we have
%%%%%%%%%%%%%%%%%%%%%%%%%%%%%%%%%%%%%%%%%%%%%%%%%%%%%%%%%%%%%%%%%%%%%%%%%%%%%%%%
%%%%%%%%%%%%%%%%%%%%%%%%%%%%%%%%%%%%%%%%%%%%%%%%%%%%%%%%%%%%%%%%%%%%%%%%%%%%%%%%
%%
%% $\mathbf{K} = \mathbf{M}_{14,k} = $ with $k = 1$ or~$4$ according as
%% $\varepsilon = 1$ or~$-1$.
%%
%%%%%%%%%%%%%%%%%%%%%%%%%%%%%%%%%%%%%%%%%%%%%%%%%%%%%%%%%%%%%%%%%%%%%%%%%%%%%%%%
%%%%%%%%%%%%%%%%%%%%%%%%%%%%%%%%%%%%%%%%%%%%%%%%%%%%%%%%%%%%%%%%%%%%%%%%%%%%%%%%
$M = \langle M_1 \circ_d M_2, x \rangle$, where $x = 1$ if $d = 1$, and 
$x^2 \in M_1 \circ_2 M_2$, otherwise. In the latter case,~$x$ centralizes
$Z(M_1)$ and~$M_2$, and normalizes~$M_1$, inducing a diagonal automorphism
on~$M_1'$. Moreover, 
$M_1 = Z( M_1 ) \circ_d M_1'$ and $M_2 = \langle Z(M_2) \circ_d M_2', y \rangle$ 
with $y = 1$ if $d = 1$, and $y^2 \in Z(M_2) \circ_d M_2'$, otherwise. In the
latter case,~$y$ centralizes $Z( M_2 )$ and induces a diagonal automorphism 
on~$M_2'$. Furthermore, 
$$
N_G( \mathbf{M} ) = \langle M_1.2 \circ_d M_2.2, x \rangle.
$$
with $M_i.2 = \langle M_i, m_i \rangle$, $i = 1, 2$. Finally,~$x$ and~$y$
normalize $\langle Z(M_i), m_i \rangle$ and $\langle M_i, m_i \rangle$ 
for $i = 1, 2$.
%suitable elements $m_i$, 
%$M_i.2 \unlhd N_G(\mathbf{M})$ for $i = 1, 2$. Finally,
%$i = 1, 2$, fixed by~$F_1$. 
%
%$M_1.2 = \langle Z(M_1), m_1 \rangle \times M_1'$ 
%and $M_2.2 = \langle \langle Z( M_2 ), m_2 \rangle \circ_d M_2', y \rangle$.
%
%{\rm (c)} Let~$\mathbf{M} = \mathbf{M}_0$. Assume that~$q$ is odd and 
%adopt the notation of~{\rm (a)} 
%and~{\rm (c)}. Let~$\sigma$ be a power of~$F_1$ such that~$F$ is a power 
%of~$\sigma$. Assume that $\mathbf{T}_0$ is $1$-$\sigma$-split. Choose an 
%inverse image $\dot{w} \in N'$ under the epimorphism 
%$N' \rightarrow \Stab_W(J)$. Then 
%$N_{\mathbf{G}}( \mathbf{M} )^{\sigma\dot{w}}$ contains 
%a $\sigma\dot{w}$-stable subgroup of index~$2$. Also, $\sigma\dot{w}$ 
%normalizes~$\mathbf{M}_i'$ and acts as a field automorphism on~$M_i'$, 
%$i = 1,2$.

{\rm (c)} If $N_G( \mathbf{M} ) = N_G( M )$, every $\chi \in \Irr(M)$ extends to 
$N_G( M, \chi )$. Suppose that $d = 2$ and that $N_G( M, \chi ) = N_G( M )$. 
If $\Res_{M_1 \circ_2 M_2}^M( \chi )$ is irreducible, there is an 
extension $\hat{\chi} \in \Irr( N_G( M ) )$ of~$\chi$ such that
$\Res_{M_1.2 \circ_2 M_2.2}^{N_G(M)}( \hat{\chi} )$ is irreducible.
If $\Res_{M_1 \circ_2 M_2}^M( \chi )$ is reducible, there is an
extension $\hat{\chi} \in \Irr( N_G( M ) )$ of~$\chi$ such that
$\Res_{M_1.2 \circ_2 M_2.2}^{N_G(M)}( \hat{\chi} )$ is reducible.
\end{prop}
\begin{prf}
Again, we only prove the assertions for $g = 1$.

(a) Use the notation and results summarized in~\ref{SetupF4}. 
By \cite[p.~$74$]{How}, we have $\Stab_W( J ) \cong 2^2$. In fact, a computation 
with CHEVIE \cite{chevie} shows that 
$\Stab_W( J ) = \langle s_{22}, s_{17} \rangle$. Let 
$J_1 = \{ \alpha_1, \alpha_{22} \}$ and $J_2 = \{ \alpha_4, \alpha_{14} \}$. 
Then $\alpha_{17} \in \bar{J}_2$, and $\bar{J}_1$ and $\bar{J}_2$ are of type 
$A_1A_1$ and $C_2$, respectively. Moreover, 
$\bar{J}_1 \cap \bar{J}_2 = \emptyset$ and $\Gamma := \bar{J}_1 \cup \bar{J}_2$ 
is closed.
%%%%%%%%%%%%%%%%%%%%%%%%%%%%%%%%%%%%%%%%%%%%%%%%%%%%%%%%%%%%%%%%%%%%%%%%%%%%%%%%
%%%%%%%%%%%%%%%%%%%%%%%%%%%%%%%%%%%%%%%%%%%%%%%%%%%%%%%%%%%%%%%%%%%%%%%%%%%%%%%%
%%
%% Computations with GAP3; see ~hiss/papers/alperinF4ppr/GapComputations/F4.g
%%
%%%%%%%%%%%%%%%%%%%%%%%%%%%%%%%%%%%%%%%%%%%%%%%%%%%%%%%%%%%%%%%%%%%%%%%%%%%%%%%%
%%%%%%%%%%%%%%%%%%%%%%%%%%%%%%%%%%%%%%%%%%%%%%%%%%%%%%%%%%%%%%%%%%%%%%%%%%%%%%%%
Put $\mathbf{K}_i := \mathbf{K}_{J_i}$ for $i = 1, 2$.
Then $[\mathbf{K}_1, \mathbf{K}_2] = 1$ by the remark preceding 
Lemma~\ref{AdjustmentOfInverseImages}. A calculation with CHEVIE, using
\cite[8.1.8, 8.1.9]{springer} shows that $\mathbf{K}_1$ has a center of order 
$d^2$, and that each of $\mathbf{K}_2$ and $\mathbf{K}_1\mathbf{K}_2 =
\mathbf{K}_\Gamma$ has a center of order~$d$; thus $\mathbf{K}_1 \cap
\mathbf{K}_2$ is a group of order~$d$. Moreover, $\mathbf{K}_1$ is a direct
product of two copies of $\SL_2( \mathbb{F} )$, and $\mathbf{K}_2 =
\Sp_4( \mathbb{F} )$. Let $\mathbf{M}_i$ denote the Levi subgroup of
$\mathbf{K}_i$ corresponding to $\alpha_1$ if $i = 1$, and to $\alpha_4$ if
$i = 2$. Then $\mathbf{M}_1 = Z( \mathbf{M}_1 )^\circ \times \mathbf{M}_1'$ and
$\mathbf{M}_2 \cong \GL_2( \mathbb{F} )$. Moreover,
$[\mathbf{M}_1, \mathbf{M}_2] = 1$ and $\mathbf{M} = \mathbf{M}_1\mathbf{M}_2$.
Also, $\mathbf{M}_1$ and $\mathbf{M}_2$ intersect in a group of order~$d$.
This implies the structure of~$\mathbf{M}$ asserted in~(a).
%%%%%%%%%%%%%%%%%%%%%%%%%%%%%%%%%%%%%%%%%%%%%%%%%%%%%%%%%%%%%%%%%%%%%%%%%%%%%%%%
%%%%%%%%%%%%%%%%%%%%%%%%%%%%%%%%%%%%%%%%%%%%%%%%%%%%%%%%%%%%%%%%%%%%%%%%%%%%%%%%
%%
%% $\mathbf{M} = \mathbf{M}_1\mathbf{M}_2$ as $\mathbf{M}$ is connected.
%%
%%%%%%%%%%%%%%%%%%%%%%%%%%%%%%%%%%%%%%%%%%%%%%%%%%%%%%%%%%%%%%%%%%%%%%%%%%%%%%%%
%%%%%%%%%%%%%%%%%%%%%%%%%%%%%%%%%%%%%%%%%%%%%%%%%%%%%%%%%%%%%%%%%%%%%%%%%%%%%%%%

We now prove the claims on $N_G( \mathbf{M} )$. 
Put $m_1 := n_{22}$ and $m_2 := n_{17}n_{14}^2$. Then $m_i \in \mathbf{K}_i$
and $m_i$ centralizes $\mathbf{M}_1'\mathbf{M}_2'$ for $i = 1, 2$. Moreover,
$N_{\mathbf{K}_i}(\mathbf{M}_i) =
\langle \mathbf{M}_i, m_i \rangle =: \mathbf{M}_i.2$ for $i = 1,2$. 
This yields the claims on the structure of $N_G( \mathbf{M} )$.
%%%%%%%%%%%%%%%%%%%%%%%%%%%%%%%%%%%%%%%%%%%%%%%%%%%%%%%%%%%%%%%%%%%%%%%%%%%%%%%%
%%%%%%%%%%%%%%%%%%%%%%%%%%%%%%%%%%%%%%%%%%%%%%%%%%%%%%%%%%%%%%%%%%%%%%%%%%%%%%%%
%%
%% See sheet 3 of 06.02.2020.
%%
%%%%%%%%%%%%%%%%%%%%%%%%%%%%%%%%%%%%%%%%%%%%%%%%%%%%%%%%%%%%%%%%%%%%%%%%%%%%%%%%
%%%%%%%%%%%%%%%%%%%%%%%%%%%%%%%%%%%%%%%%%%%%%%%%%%%%%%%%%%%%%%%%%%%%%%%%%%%%%%%%

(b) Notice that ${M}_1 \circ_d {M}_2$ has index~$d$ in~${M}$. In case~$q$
is odd, we choose a particular $F$-stable element~$x$ of ${M} \setminus
({M}_1 \circ_2 {M}_2)$ as follows. Let $x_1 \in \mathbf{T} \cap \mathbf{M}_1$
and $x_2 \in Z( \mathbf{M}_2 )$ be such that $F(x_1)^{-1} x_1 = F(x_2)x_2^{-1}$ 
is the unique element of order~$2$ in $\mathbf{M}_1 \cap \mathbf{M}_2 =
M_1 \cap M_2$. Then $x := x_1x_2$ is $F$-stable and $x \not\in M_1 \circ_2 M_2$
so that $M = \langle M_1 \circ_2 M_2, x \rangle$. Notice that~$x$ 
centralizes~$M_2$ and normalizes~$M_1$, inducing a diagonal automorhism 
on~$M_1'$. Also,~$x$ centralizes $Z(M_1)$, as $x = x_1x_2$ with 
$x_1 \in T \cap M_1$. With an analogous argument we can find 
$y = y_1y_2 \in M_2$ with $y_1 \in Z( \mathbf{M}_2 )$ and 
$y_2 \in \mathbf{M}_2'$ such that 
$M_2 = \langle Z(M_2) \circ_d M_2', y \rangle$. In particular,~$y$ 
centralizes~$Z(M_2)$ and induces a diagonal automorphism on~$M_2'$ if $d = 2$.
This gives the structure of~$M_2$ as claimed. The structure of $M_1$ is clear 
from~(a). Moreover, $z^{-1}m_i^{-1}zm_i \in Z( M_i )$ for $i = 1, 2$ and 
$z \in \{ x, y \}$, so that~$x$ and~$y$ normalize $\langle Z(M_i), m_i \rangle$ 
and $\langle M_i, m_i \rangle$ for $i = 1, 2$. As 
$(N_{\mathbf{G}}( \mathbf{M} )/\mathbf{M})^F = N_G( \mathbf{M} )/M$, the claims 
on the structure of $N_G( \mathbf{M} )$ in~(b) are established.
%%%%%%%%%%%%%%%%%%%%%%%%%%%%%%%%%%%%%%%%%%%%%%%%%%%%%%%%%%%%%%%%%%%%%%%%%%%%%%%%
%%%%%%%%%%%%%%%%%%%%%%%%%%%%%%%%%%%%%%%%%%%%%%%%%%%%%%%%%%%%%%%%%%%%%%%%%%%%%%%%
%%
%% See sheet 1 of 05.02.2020 and sheet 2 of 06.02.2020.
%%
%%%%%%%%%%%%%%%%%%%%%%%%%%%%%%%%%%%%%%%%%%%%%%%%%%%%%%%%%%%%%%%%%%%%%%%%%%%%%%%%
%%%%%%%%%%%%%%%%%%%%%%%%%%%%%%%%%%%%%%%%%%%%%%%%%%%%%%%%%%%%%%%%%%%%%%%%%%%%%%%%

%(c) We may assume that $\dot{w} \in \{ m_1, m_2' \}$. Then $\sigma \dot{w}$ acts
%as a field automorphism on~$M_i'$, as~$\dot{w}$ centralizes~$M_i'$ for $i = 
%1, 2$. By construction~$m_1$ and $m_2'$ normalize the torus $\mathbf{T}_0'$ 
%that corresponds to 
%$\langle \alpha_1^\vee, \alpha_{22}^\vee, \alpha_{14}^\vee, \alpha_{4}^\vee \rangle
%\otimes_{\mathbb{Z}} \mathbb{Q}_{p'}/\mathbb{Z}$ under the 
%isomorphism~(\ref{YotimesF}). As $\mathbf{T}_0'$ has index~$2$ in $\mathbf{T}_0$,
%the subgroup $\langle N' \circ_Z M', {\mathbf{T}_0'}^F \rangle$ has index~$2$ in
%$N_{G}( \mathbf{ M } )$ since $N_{G}( \mathbf{ M } ) = 
%\langle N' \circ_Z M', T_0 \rangle$ by~(c).

(c) We use the structure of $N_G(\mathbf{M})$ established in~(b) to prove the 
claim. Let $\chi \in \Irr( M )$. If the inertia subgroup of~$\chi$ is strictly 
smaller than~$N_G(M)$, the inertia quotient is cyclic, and $\chi$ extends.
If~$q$ is even, then~$\chi$ clearly extends. We may thus assume that~$\chi$ is 
invariant in $N_G( M )$ and that $q$ is odd. In the following, we will
frequently use the remarks in~\ref{CentralProducts}. Let $\chi'$ denote the 
restriction of~$\chi$ to $M_1 \circ_2 M_2$, and let~$\psi$ be an irreducible 
constituent of~$\chi'$. Then $\psi = \psi_1\psi_2$ with 
$\psi_i \in \Irr( M_i )$, $i = 1, 2$. If $\chi'$ is irreducible, i.e., 
$\chi' = \psi$, then~$\psi_i$ is invariant in $M_i.2$ and thus extends 
to~$M_i.2$, $i = 1, 2$. Hence $\chi' = \psi$ extends to $\hat{\psi} \in 
\Irr( M_1.2 \circ_2 M_2.2 )$. By \cite[Lemma~$4.1$(a)]{Spaeth10}, there is 
an extension $\hat{\chi} \in \Irr( N_G( M ) )$ such that
$\Res^{N_G(M)}_{M_1.2 \circ_2 M_2.2} = \hat{\psi}$.
If~$\chi'$ is reducible, we have $\chi' = \psi_1\psi_2 + \psi_1^x\psi_2^x$, as 
$x^2 \in M_1 \circ_2 M_2$. Moreover, $\psi_2^x = \psi_2$ as~$x$ 
centralizes~$M_2$, and hence $\psi_1^x \neq \psi_1$. Now $N_G( M, \psi )$ has
index~$2$ in $N_G( M )$. Notice that $m_2x \not\in N_G( M, \psi )$, as
$\psi_1^{m_2x} = \psi_1^x \neq \psi_1$. Suppose that $m_1x \in N_G( M, \psi )$.
Then $\psi_1^{m_1} = \psi_1^x$. As $M_1 = Z( M_1 ) \circ_2 M_1'$, we can write
$\psi_1 = \lambda_1 \psi_1'$ with $\lambda_1 \in \Irr( Z( M_1 ) )$ and $\psi_1'
\in \Irr( M_1' )$. However, as~$x$ centralizes~$Z( M_1 )$ and~$m_1$ 
centralizes~$M_1'$, we cannot have $\psi_1^{m_1} = \psi_1^x$ unless 
$\psi_1^x = \psi_1$. This contradiction shows that 
$m_1x \not\in N_G( M, \psi )$. It follows that 
$N_G( M, \psi ) = M_1.2 \circ_2 M_2.2$. Hence $\psi = \psi_1\psi_2$ extends
to $\hat{\psi} \in \Irr( N_G( M, \psi ) )$, and $\hat{\psi}^x \neq \hat{\psi}$.
Then $\hat{\chi} := \Ind_{N_G( M, \psi )}^{N_G(M)}( \hat{\psi} )$ is an 
extension of~$\chi$.
\end{prf}

\begin{prop}
\label{NormalizerM18}
Let~$\mathbf{M}_0 := \mathbf{L}_J$, where~$J$ is one of the sets 
	$\{ \alpha_1 \}$, 
%%%%%%%%%%%%%%%%%%%%%%%%%%%%%%%%%%%%%%%%%%%%%%%%%%%%%%%%%%%%%%%%%%%%%%%%%%%%%%%%
%%%%%%%%%%%%%%%%%%%%%%%%%%%%%%%%%%%%%%%%%%%%%%%%%%%%%%%%%%%%%%%%%%%%%%%%%%%%%%%%
%%
%% $\mathbf{K} = \mathbf{M}_{18,k} = [q-\varepsilon]^3 \circ_d \SL_2(q).d$ 
%% with $k = 1$ or~$10$ according as $\varepsilon = 1$ or~$-1$.
%%
%%%%%%%%%%%%%%%%%%%%%%%%%%%%%%%%%%%%%%%%%%%%%%%%%%%%%%%%%%%%%%%%%%%%%%%%%%%%%%%%
%%%%%%%%%%%%%%%%%%%%%%%%%%%%%%%%%%%%%%%%%%%%%%%%%%%%%%%%%%%%%%%%%%%%%%%%%%%%%%%%
	$\{ \alpha_4 \}$ or 
%%%%%%%%%%%%%%%%%%%%%%%%%%%%%%%%%%%%%%%%%%%%%%%%%%%%%%%%%%%%%%%%%%%%%%%%%%%%%%%%
%%%%%%%%%%%%%%%%%%%%%%%%%%%%%%%%%%%%%%%%%%%%%%%%%%%%%%%%%%%%%%%%%%%%%%%%%%%%%%%%
%%
%% $\mathbf{K} = \mathbf{M}_{18,k} = [q-\varepsilon]^3 \circ_d \SL_2(q).d$
%% with $k = 1$ or~$10$ according as $\varepsilon = 1$ or~$-1$.
%%
%%%%%%%%%%%%%%%%%%%%%%%%%%%%%%%%%%%%%%%%%%%%%%%%%%%%%%%%%%%%%%%%%%%%%%%%%%%%%%%%
%%%%%%%%%%%%%%%%%%%%%%%%%%%%%%%%%%%%%%%%%%%%%%%%%%%%%%%%%%%%%%%%%%%%%%%%%%%%%%%%
	$\{ \alpha_2, \alpha_3 \}$ 
%%%%%%%%%%%%%%%%%%%%%%%%%%%%%%%%%%%%%%%%%%%%%%%%%%%%%%%%%%%%%%%%%%%%%%%%%%%%%%%%
%%%%%%%%%%%%%%%%%%%%%%%%%%%%%%%%%%%%%%%%%%%%%%%%%%%%%%%%%%%%%%%%%%%%%%%%%%%%%%%%
%%
%% $\mathbf{K} = \mathbf{M}_{15,k} = [q-\varepsilon]^2 \circ_d (\Sp_4(q).d)$ 
%% with $k = 1$ or~$3$ according as $\varepsilon = 1$ or~$-1$.
%%
%%%%%%%%%%%%%%%%%%%%%%%%%%%%%%%%%%%%%%%%%%%%%%%%%%%%%%%%%%%%%%%%%%%%%%%%%%%%%%%%
%%%%%%%%%%%%%%%%%%%%%%%%%%%%%%%%%%%%%%%%%%%%%%%%%%%%%%%%%%%%%%%%%%%%%%%%%%%%%%%%
of simple roots of~$\Sigma$. Put $\mathbf{M} := \mathbf{M}_0^g$ and $\mathbf{M}' 
:= [\mathbf{M},\mathbf{M}]$. Let $d := \gcd( 2, q - 1)$.

{\rm (a)} We have $\mathbf{M} = Z(\mathbf{M}) \circ_d \mathbf{M}'$
with $\mathbf{M}' \cong
\SL_2( \mathbb{F} )$ in the first two cases, and $\mathbf{M}' \cong 
\Sp_4( \mathbb{F} )$ in the last case. Furthermore,
$$N_{\mathbf{G}}( \mathbf{M} ) = \mathbf{N}' \circ_d \mathbf{M}'$$
with $Z( \mathbf{M} ) \leq \mathbf{N}' \leq N_{\mathbf{G}}( \mathbf{T}_0 )$ and 
$\mathbf{N}' \cap \mathbf{M} = Z( \mathbf{M} )$. Furthermore,
$\mathbf{N}'/Z(\mathbf{M}) \cong \Stab_W( J )$, where 
$\Stab_W(J) \cong W(C_3)
\cong 2^3.S_3$ in the first two cases, and $\Stab_W(J) \cong W( C_ 2 ) 
\cong D_8$ in the last case. 
Finally, $\mathbf{N}' = C_{N_{\mathbf{G}}( \mathbf{M} )}( \mathbf{M}' )$.

{\rm (b)} We have $M = \langle Z( M ) \circ_d M', x \rangle$, with $x = 1$ 
if~$q$ is even, and $x^2 \in Z( M ) \circ_2 M'$ if $q$ is odd; in the
latter case,~$x$ induces a diagonal automorphism on $M'$.
Moreover,
$$N_G( \mathbf{M} ) = \langle N' \circ_d M', x \rangle.$$

{\rm (c)} If $N_{G}( \mathbf{M} ) = N_G( M )$, then every $\chi \in \Irr( M )$ 
extends to $N_{G}( M, \chi )$.
\end{prop}
\begin{prf}
Once more, we only prove the case $g = 1$.
	
(a) We have 
$\mathbf{M} = Z(\mathbf{M})\mathbf{M}'$ and $\mathbf{M}'$ is a simple, simply 
connected algebraic group of type $A_1$ and $C_2$, respectively. 
Hence~$\mathbf{M}'$ is isomorphic to $\SL_2( \mathbb{F} )$, respectively
$\Sp_4( \mathbb{F} )$. If~$q$ is even, $Z( \mathbf{M}' )$ is trivial and thus 
$\mathbf{M} = Z( \mathbf{M} ) \times \mathbf{M}'$. If~$q$ is odd, 
$Z( \mathbf{M}' )$ has order~$2$ and lies in $Z( \mathbf{M} )$, and thus 
$\mathbf{M} = Z( \mathbf{M} ) \circ_2 \mathbf{M}'$. This proves our first 
assertion on the structure of~$\mathbf{M}$.
	
To investigate the structure of $N_{\mathbf{G}}( \mathbf{M} )$, we will choose 
a closed subsystem $\Delta \subseteq \Sigma$ with the following 
properties: $\Delta \cap \Sigma_J = \emptyset$, $\rk( \Delta ) + |J| = 4$, 
$W_{\Delta} = \Stab_W( J )$ and $W_\Delta$ stabilizes each element of~$J$. 
Moreover, $\tilde{\Delta}$ will denote a closed subsystem of~$\Delta$ (hence of 
$\Sigma$) with $\rk(\Delta) = \rk(\tilde{\Delta})$ and such that $\tilde{\Delta} 
\cup \Sigma_J$ is closed in $\Sigma$. Put $\mathbf{K} := \mathbf{K}_\Delta$, 
$\tilde{\mathbf{K}} := \mathbf{K}_{\tilde{\Delta}}$, and $\mathbf{S} := 
\mathbf{T}_0 \cap \mathbf{K} = \mathbf{T}_0 \cap \tilde{\mathbf{K}}$. 
Then $\mathbf{S} = Z( \mathbf{M} )$.
%%%%%%%%%%%%%%%%%%%%%%%%%%%%%%%%%%%%%%%%%%%%%%%%%%%%%%%%%%%%%%%%%%%%%%%%%%%%%%%%
%%%%%%%%%%%%%%%%%%%%%%%%%%%%%%%%%%%%%%%%%%%%%%%%%%%%%%%%%%%%%%%%%%%%%%%%%%%%%%%%
%%
%% See sheet 1 of 20.10.2020 for the first of the above assertions.
%% The second one follows from C_{\mathbf{G}}( \mathbf{K}_J ) = \mathbf{K}.
%%
%%%%%%%%%%%%%%%%%%%%%%%%%%%%%%%%%%%%%%%%%%%%%%%%%%%%%%%%%%%%%%%%%%%%%%%%%%%%%%%%
%%%%%%%%%%%%%%%%%%%%%%%%%%%%%%%%%%%%%%%%%%%%%%%%%%%%%%%%%%%%%%%%%%%%%%%%%%%%%%%%
Write $\mathbf{N} := N_{\mathbf{K}}( \mathbf{S} )$ and $\tilde{\mathbf{N}} := 
N_{\tilde{\mathbf{K}}}( \mathbf{S} )$. As $\mathbf{N}$ normalizes $\mathbf{S} = 
Z( \mathbf{M} )$, it also normalizes $\mathbf{M} = 
C_{\mathbf{G}}( Z( \mathbf{M} ) )$ and~$\mathbf{M}'$. Furthermore, 
$\mathbf{N} \cap \mathbf{M} = Z( \mathbf{M} ) = \mathbf{S}$ and 
$\mathbf{N}/\mathbf{S} = W_\Delta = \Stab_W( J )$. Thus $N_G( \mathbf{M} ) = 
\mathbf{N}\mathbf{M} = \mathbf{N}\mathbf{M}'$. Moreover, $\mathbf{S} \leq 
\tilde{\mathbf{N}} \leq \mathbf{N}$ and $\tilde{\mathbf{N}}$ 
centralizes~$\mathbf{M}'$ as $\mathbf{M}' = \mathbf{K}_J$ and 
$[\tilde{\mathbf{K}},\mathbf{K}_J] = 1$. 

We now choose $\Delta$ and $\tilde{\Delta}$ in the respective cases. Suppose 
first that $J = \{\alpha_1\}$. Here, we let $\Delta = \tilde{\Delta}$ to be 
the closed subsystem of $\Sigma$ generated by 
$\{ \alpha_3, \alpha_{4}, \alpha_{14} \}$. Then $\Delta$ is of type $C_3$.
%%%%%%%%%%%%%%%%%%%%%%%%%%%%%%%%%%%%%%%%%%%%%%%%%%%%%%%%%%%%%%%%%%%%%%%%%%%%%%%%
%%%%%%%%%%%%%%%%%%%%%%%%%%%%%%%%%%%%%%%%%%%%%%%%%%%%%%%%%%%%%%%%%%%%%%%%%%%%%%%%
%%
%% The positive roots in \Delta are [ 3, 4, 7, 14, 17, 19, 20, 21, 22 ].
%%
%%%%%%%%%%%%%%%%%%%%%%%%%%%%%%%%%%%%%%%%%%%%%%%%%%%%%%%%%%%%%%%%%%%%%%%%%%%%%%%%
%%%%%%%%%%%%%%%%%%%%%%%%%%%%%%%%%%%%%%%%%%%%%%%%%%%%%%%%%%%%%%%%%%%%%%%%%%%%%%%%
If $J = \{ \alpha_4 \}$, we choose $\Delta$ and~$\tilde{\Delta}$ to be the
closed subsystems of~$\Sigma$ generated by $\{ \alpha_1, \alpha_2, \alpha_{13} \}$ and
$\{ \alpha_1, \alpha_{2}, \alpha_{22} \}$, respectively. Then $\Delta$ is of type $B_3$
and $\tilde{\Delta}$ of type~$A_3$. 
%%%%%%%%%%%%%%%%%%%%%%%%%%%%%%%%%%%%%%%%%%%%%%%%%%%%%%%%%%%%%%%%%%%%%%%%%%%%%%%%
%%%%%%%%%%%%%%%%%%%%%%%%%%%%%%%%%%%%%%%%%%%%%%%%%%%%%%%%%%%%%%%%%%%%%%%%%%%%%%%%
%%
%% The positive roots in \Delta are [ 1, 2, 5, 13, 15, 17, 22, 23, 24 ]
%%
%% The positive roots in \tilde{\Delta} are [ 1, 2, 5, 22, 23, 24 ].
%%
%%%%%%%%%%%%%%%%%%%%%%%%%%%%%%%%%%%%%%%%%%%%%%%%%%%%%%%%%%%%%%%%%%%%%%%%%%%%%%%%
%%%%%%%%%%%%%%%%%%%%%%%%%%%%%%%%%%%%%%%%%%%%%%%%%%%%%%%%%%%%%%%%%%%%%%%%%%%%%%%%
If $J = \{ \alpha_2, \alpha_3 \}$, we choose $\Delta$ and $\tilde{\Delta}$ to be 
the closed subsystems of~$\Sigma$ generated by $\{ \alpha_8, \alpha_{16} \}$, 
and $\{ \alpha_{16}, \alpha_{24} \}$, respectively. Then $\Delta$ is of type 
$C_2$ and $\tilde{\Delta}$ of type~$A_1A_1$.
%%%%%%%%%%%%%%%%%%%%%%%%%%%%%%%%%%%%%%%%%%%%%%%%%%%%%%%%%%%%%%%%%%%%%%%%%%%%%%%%
%%%%%%%%%%%%%%%%%%%%%%%%%%%%%%%%%%%%%%%%%%%%%%%%%%%%%%%%%%%%%%%%%%%%%%%%%%%%%%%%
%%
%% The positive roots in \Delta are [ 8, 16, 21, 24 ].
%%
%% The positive roots in \tilde{\Delta} are [ 16, 24 ].
%%
%%%%%%%%%%%%%%%%%%%%%%%%%%%%%%%%%%%%%%%%%%%%%%%%%%%%%%%%%%%%%%%%%%%%%%%%%%%%%%%%
%%%%%%%%%%%%%%%%%%%%%%%%%%%%%%%%%%%%%%%%%%%%%%%%%%%%%%%%%%%%%%%%%%%%%%%%%%%%%%%%

If $J = \{ \alpha_1 \}$, we put $\mathbf{N}' := \mathbf{N}$. Notice that in this
case, $\mathbf{N}' = \mathbf{N} = \tilde{\mathbf{N}}$ centralizes~$\mathbf{M}'$. 
Suppose that $\tilde{\mathbf{N}} \lneq \mathbf{N}$. Then 
$[\mathbf{N}\colon\!\tilde{\mathbf{N}}] = 2$, as 
$\tilde{\mathbf{N}}/\mathbf{S} \cong W_{\tilde{\Delta}}$, and thus 
$\tilde{\mathbf{N}}\mathbf{M}'$ has index~$2$ in $N_{\mathbf{G}}(\mathbf{M})$. 
If $J = \{ \alpha_2, \alpha_3 \}$, let $n := n_8$, and if $J = \{ \alpha_4 \}$,
let $n := n_{13}$. Then $\mathbf{N} = \langle \tilde{\mathbf{N}}, n \rangle$. As 
the image of~$n$ in~$W$ lies in $W_\Delta$, which fixes the elements of~$J$, we 
find that $n^{-1} u_\alpha( s ) n = u_{\alpha}( \pm s )$ for all 
$\alpha \in J \cup (-J)$ and all $s \in \mathbb{F}$;
see \cite[Lemma~$7.2.1$(i)]{C1}. In particular,~$n$ 
centralizes~$\mathbf{M}'$ if $q$ is even, in which case we set $\mathbf{N}' := 
\mathbf{N}$.  Assume now that $q$ is odd. By 
Lemma~\ref{AdjustmentOfInverseImages}, we may replace~$n$ by $n' = tn$ for a 
suitable $t \in {T}_0$ in such a way that $n'$ centralizes~$\mathbf{M}'$. 
Putting $\mathbf{N}' := \langle \tilde{\mathbf{N}}, n' \rangle$, we obtain 
$N_{\mathbf{G}}(\mathbf{M}) = \mathbf{N}'\mathbf{M}'$ and 
$[\mathbf{N}',\mathbf{M}'] = 1$, i.e., $\mathbf{N}'\mathbf{M}' = 
\mathbf{N}' \circ_d \mathbf{M}'$. As $\mathbf{N}'$ normalizes $Z( \mathbf{M} )$ 
and centralizes $\mathbf{M}'$, it normalizes $\mathbf{T}_0 = 
Z( \mathbf{M} )(\mathbf{M}' \cap \mathbf{T}_0)$. Finally, 
$\mathbf{N}' = C_{\mathbf{G}}( \mathbf{M}' ) \cap N_{\mathbf{G}}( \mathbf{M} )
= C_{N_{\mathbf{G}}( \mathbf{M} )}( \mathbf{M}' )$.

(b) If~$q$ is even, our claims are obvious from~(a). Let~$q$ be odd. Then there 
are $x_1 \in Z( \mathbf{M} )$ and $x_2 \in \mathbf{T} \cap \mathbf{M}'$ such 
that $F( x_1 )^{-1} x_1 = F( x_2 ) x_2^{-1}$ is the unique 
element of order~$2$ in $Z( \mathbf{M} ) \cap \mathbf{M}'$. Hence 
$x := x_1x_2$ is $F$-stable and $M = \langle Z( M ) \circ_d M', x \rangle$. 
Notice that~$x$ normalizes~$\mathbf{M}'$, inducing a diagonal 
automorphism on~$M'$.
%%%%%%%%%%%%%%%%%%%%%%%%%%%%%%%%%%%%%%%%%%%%%%%%%%%%%%%%%%%%%%%%%%%%%%%%%%%%%%%%
%%%%%%%%%%%%%%%%%%%%%%%%%%%%%%%%%%%%%%%%%%%%%%%%%%%%%%%%%%%%%%%%%%%%%%%%%%%%%%%%
%%
%% See sheet 1 of 05.02.2020 and sheet 2 of 06.02.2020.
%%
%%%%%%%%%%%%%%%%%%%%%%%%%%%%%%%%%%%%%%%%%%%%%%%%%%%%%%%%%%%%%%%%%%%%%%%%%%%%%%%%
%%%%%%%%%%%%%%%%%%%%%%%%%%%%%%%%%%%%%%%%%%%%%%%%%%%%%%%%%%%%%%%%%%%%%%%%%%%%%%%%
This gives our claim for the structure of $M$. The claim for 
$N_G( \mathbf{M} )$ follows from this and~(a). 

(c) By~(b) we have 
$M = \langle Z(M) \circ_d M', x \rangle$ and 
$N_G( \mathbf{M} ) = \langle N' \circ_d M', x \rangle$; moreover, $Z(M) = 
[q-1]^c$, with $c = 3$ or~$2$. Notice that the 
groups~$\mathbf{S}$,~$\mathbf{K}$ and $\tilde{\mathbf{K}}$ introduced in the 
proof of~(a) are~$F$-invariant. Notice also that $F$ acts trivially on 
$\mathbf{N}/Z( \mathbf{M} )$, and we implicitly assume that the inverse images
of elements of this group used below are $F'$-stable.

Let $\chi \in \Irr( M )$. In the considerations to follow, we will make use
of the facts summarized in~\ref{CentralProducts} for characters of central 
products. As $M = Z(M) \circ \langle M', x \rangle$ is a central product, we
have $\chi = \lambda \psi$ for 
$\lambda \in \Irr( Z(M) )$ and $\psi \in \Irr( \langle M',x \rangle )$.
Let~$\chi'$ and~$\psi'$ denote the restrictions 
of~$\chi$ to $Z( M ) \circ_d M'$ and of~$\psi$ to~$M'$, respectively. Then 
$\chi' = \lambda \psi'$. 

Recall from the proof of~(a) that $N_G( M ) = NM$ with $M \cap N = Z(M)$, and 
that $\tilde{N}$ is a subgroup of~$N$ of index~$2$, containing $Z(M)$ and
centralizing~$M'$.
Put $I := N_G( M, \chi ) \cap N$ and $\tilde{I} := I \cap \tilde{N}$. Then 
$N_G( M, \chi ) = IM$ and thus $N_G( M, \chi )/M \cong I/Z(M)$. 
Also,~$\tilde{I}$ has index at most~$2$ in~$I$. Clearly,~$I$ 
stabilizes~$\lambda$ and $\psi'$.
It follows from \cite[Theorem~$1.1$]{Spaeth10}, applied to 
$N_{\mathbf{K}}( \mathbf{S} )$, that 
$\lambda$ extends to an irreducible character $\hat{\lambda}$ of~$I$. 
(Notice that~$\Delta$ is conjugate in~$W$ to a parabolic subsystem of~$\Sigma$, 
so that~$\mathbf{K}$ is simply connected.) Thus there is an
extension $\tilde{\lambda}$ of~$\lambda$ to~$\tilde{I}$ which is invariant under~$I$.

To continue, suppose first that $\chi' = \lambda \psi'$ is irreducible. Then the 
inertia group of~$\chi'$ in $NM'$ equals $IM'$ which contains $\tilde{I}M'$ as a 
subgroup of index at most~$2$.
%%%%%%%%%%%%%%%%%%%%%%%%%%%%%%%%%%%%%%%%%%%%%%%%%%%%%%%%%%%%%%%%%%%%%%%%%%%%%%%%
%%%%%%%%%%%%%%%%%%%%%%%%%%%%%%%%%%%%%%%%%%%%%%%%%%%%%%%%%%%%%%%%%%%%%%%%%%%%%%%%
%%
%% See sheet 4 of 12.02.2020.
%%
%%%%%%%%%%%%%%%%%%%%%%%%%%%%%%%%%%%%%%%%%%%%%%%%%%%%%%%%%%%%%%%%%%%%%%%%%%%%%%%%
%%%%%%%%%%%%%%%%%%%%%%%%%%%%%%%%%%%%%%%%%%%%%%%%%%%%%%%%%%%%%%%%%%%%%%%%%%%%%%%%
By construction, $\tilde{\lambda} \psi' \in \Irr( \tilde{I}M' )$ is invariant 
under $IM'$, as~$\psi'$ is
invariant under~$I$ and $\tilde{\lambda}$ is invariant under $I$ and~$M'$, the
latter as $M'$ centralizes~$\tilde{N}$. Thus $\tilde{\lambda} \psi'$ extends to
an irreducible character of $IM'$, which is an extension of~$\lambda \psi'$.
The claim follows from \cite[Lemma~$4.1$(a)]{Spaeth10}, applied to the subgroup $NM'$
of $NM = N_G(M)$.

Suppose now that $\psi' = \vartheta + \vartheta^x$ for some
$\vartheta \in \Irr( M' )$ and write $N_G( M, \lambda\vartheta )$ for the 
inertia subgroups of $\lambda\vartheta \in \Irr( Z(M) \circ_d M' )$ in $N_G(M)$.
%%%%%%%%%%%%%%%%%%%%%%%%%%%%%%%%%%%%%%%%%%%%%%%%%%%%%%%%%%%%%%%%%%%%%%%%%%%%%%%%
%%%%%%%%%%%%%%%%%%%%%%%%%%%%%%%%%%%%%%%%%%%%%%%%%%%%%%%%%%%%%%%%%%%%%%%%%%%%%%%%
%%
%% See sheet 5 of 12.02.2020.
%%
%%%%%%%%%%%%%%%%%%%%%%%%%%%%%%%%%%%%%%%%%%%%%%%%%%%%%%%%%%%%%%%%%%%%%%%%%%%%%%%%
%%%%%%%%%%%%%%%%%%%%%%%%%%%%%%%%%%%%%%%%%%%%%%%%%%%%%%%%%%%%%%%%%%%%%%%%%%%%%%%%
Then $\tilde{I}M' \leq N_G( M, \lambda\vartheta ) \leq IM$, where the latter
inclusion has index~$2$. If $N_G( M, \lambda\vartheta ) \in 
\{ \tilde{I}M', IM' \}$, there is an extension 
$\widehat{\lambda\vartheta}$ of
$\lambda\vartheta$ to $N_G( M, \lambda\vartheta )$. Now
$\widehat{\lambda\vartheta}$ is not invariant under~$x$, as 
$\widehat{\lambda\vartheta}^x$ is an extension of $\lambda\vartheta^x$. Thus
$\Ind_{N_G( M, \lambda\vartheta )}^{N_G( M, \chi )}( \widehat{\lambda\vartheta} )$
is an extension of~$\chi$ to $N_G( M, \chi )$.
We may thus assume that $\tilde{I}M' \lneq N_G( M, \lambda\vartheta ) 
\lneq IM$, so that in particular $\tilde{I} \lneq I$.
%%%%%%%%%%%%%%%%%%%%%%%%%%%%%%%%%%%%%%%%%%%%%%%%%%%%%%%%%%%%%%%%%%%%%%%%%%%%%%%%
%%%%%%%%%%%%%%%%%%%%%%%%%%%%%%%%%%%%%%%%%%%%%%%%%%%%%%%%%%%%%%%%%%%%%%%%%%%%%%%%
%%
%% See sheet 6 of 12.02.2020.
%%
%%%%%%%%%%%%%%%%%%%%%%%%%%%%%%%%%%%%%%%%%%%%%%%%%%%%%%%%%%%%%%%%%%%%%%%%%%%%%%%%
%%%%%%%%%%%%%%%%%%%%%%%%%%%%%%%%%%%%%%%%%%%%%%%%%%%%%%%%%%%%%%%%%%%%%%%%%%%%%%%%
%%%%%%%%%%%%%%%%%%%%%%%%%%%%%%%%%%%%%%%%%%%%%%%%%%%%%%%%%%%%%%%%%%%%%%%%%%%%%%%%
%%%%%%%%%%%%%%%%%%%%%%%%%%%%%%%%%%%%%%%%%%%%%%%%%%%%%%%%%%%%%%%%%%%%%%%%%%%%%%%%
%%
%% See sheet 7 of 12.02.2020.
%%
%%%%%%%%%%%%%%%%%%%%%%%%%%%%%%%%%%%%%%%%%%%%%%%%%%%%%%%%%%%%%%%%%%%%%%%%%%%%%%%%
%%%%%%%%%%%%%%%%%%%%%%%%%%%%%%%%%%%%%%%%%%%%%%%%%%%%%%%%%%%%%%%%%%%%%%%%%%%%%%%%
Now $IM/\tilde{I}M'$ is generated by the images of~$x$ and~$n$.
Moreover, $N_G( M, \lambda\vartheta ) \neq \tilde{I}M$ as~$x$ does not
stabilize $\lambda\vartheta$. Thus there is $\tilde{n} \in \tilde{N}$ such that
$N_G( M, \lambda\vartheta ) = \langle \tilde{I}M', xn\tilde{n} \rangle$. 
In particular, $n\tilde{n}$ stabilizes~$\lambda$, and $xn\tilde{n}$ 
stabilizes~$\vartheta$. As $\tilde{N}$ centralizes~$M'$, and $n^2$ induces an
inner automorphism of~$M'$, the latter condition can also be written as 
$\vartheta^x = \vartheta^n$. In particular, this case does not occur if 
$J = \{ \alpha_1 \}$, as then~$N$ centralizes~$M'$. We claim that there is an 
$xn\tilde{n}$-stable extension 
$\tilde{\lambda}$ of~$\lambda$ to~$\tilde{I}$. Provided this claim holds, then
$\tilde{\lambda}\vartheta \in \Irr( \tilde{I}M' )$ is $xn\tilde{n}$-invariant
and thus extends to $\widehat{\lambda\vartheta} \in 
\Irr( \langle \tilde{I}M', xn\tilde{n} \rangle )$, which is not $x$-invariant.
We conclude as in the previous case.

It remains to prove the claim. Suppose first that $J = \{ \alpha_2, \alpha_3 \}$.
In this case $n = n_8$, and $s_8$ swaps the two roots $\alpha_{16}$ and 
$\alpha_{24}$. Thus there are decompositions 
$Z( M ) = [q - 1]^2$ and $\tilde{N} = ([q - 1].2)^2$ 
such that~$n$ permutes the two direct factors of the latter group. Let $\lambda = 
\lambda_1 \boxtimes \lambda_2$, where $\lambda_i \in \Irr( [q - 1] )$,
for $i = 1, 2$. 
As $n\tilde{n}$ stabilizes~$\lambda$, we conclude that $\lambda_2 = \lambda_1$
or $\lambda_2 = \lambda_1^{-1}$. If $\lambda_1 \neq \lambda_1^{-1}$, we also
have $\lambda_2 \neq \lambda_2^{-1}$, and then $\tilde{I} = Z( M )$. In this 
case, the claim is trivially true. Otherwise, $\tilde{I} = \tilde{N}$ and thus
$\langle \tilde{I}M', xn\tilde{n} \rangle = \langle \tilde{N}M', xn \rangle$. 
%%%%%%%%%%%%%%%%%%%%%%%%%%%%%%%%%%%%%%%%%%%%%%%%%%%%%%%%%%%%%%%%%%%%%%%%%%%%%%%%
%%%%%%%%%%%%%%%%%%%%%%%%%%%%%%%%%%%%%%%%%%%%%%%%%%%%%%%%%%%%%%%%%%%%%%%%%%%%%%%%
%%
%% See sheet 10 of 18.02.2020.
%%
%%%%%%%%%%%%%%%%%%%%%%%%%%%%%%%%%%%%%%%%%%%%%%%%%%%%%%%%%%%%%%%%%%%%%%%%%%%%%%%%
%%%%%%%%%%%%%%%%%%%%%%%%%%%%%%%%%%%%%%%%%%%%%%%%%%%%%%%%%%%%%%%%%%%%%%%%%%%%%%%%
Let $\tilde{\lambda}_1$ denote an extension of $\lambda_1$ to 
$[q - 1].2$, the first factor of the above decomposition 
of~$\tilde{N}$. Define $\tilde{\lambda}_2$ by 
$(\tilde{\lambda}_1 \boxtimes 1)^{xn} = 1 \boxtimes \tilde{\lambda}_2$. Then 
$\tilde{\lambda} := \tilde{\lambda}_1 \boxtimes \tilde{\lambda}_2$ is an
$xn\tilde{n}$-stable extension of~$\lambda$ to $\tilde{I} = \tilde{N}$, since 
$(xn)^2$ acts as an inner automorphism on $\tilde{N}$. 
%%%%%%%%%%%%%%%%%%%%%%%%%%%%%%%%%%%%%%%%%%%%%%%%%%%%%%%%%%%%%%%%%%%%%%%%%%%%%%%%
%%%%%%%%%%%%%%%%%%%%%%%%%%%%%%%%%%%%%%%%%%%%%%%%%%%%%%%%%%%%%%%%%%%%%%%%%%%%%%%%
%%
%% See sheet 11 of 18.02.2020.
%%
%%%%%%%%%%%%%%%%%%%%%%%%%%%%%%%%%%%%%%%%%%%%%%%%%%%%%%%%%%%%%%%%%%%%%%%%%%%%%%%%
%%%%%%%%%%%%%%%%%%%%%%%%%%%%%%%%%%%%%%%%%%%%%%%%%%%%%%%%%%%%%%%%%%%%%%%%%%%%%%%%
Suppose finally that $J = \{ \alpha_4 \}$, in which case $n = n_{13}$. As 
$\vartheta^x = \vartheta^n$ and as~$x$ induces a diagonal automorphism on~$M'$, 
the same must be true for~$n$. Now $n^{-1} u_4( t ) n = u_4( - t )$ for all
$t \in \mathbb{F}$. As~$n$ induces a non-inner automorphism on~$M'$, we must 
have $4 \mid q + 1$. It follows that~$x_1$ may be chosen as an element of 
order~$4$ in $Z( \tilde{ \mathbf{K} } )$. (Recall that $\tilde{ \mathbf{K} } 
\cong \SL_4( \mathbb{F} )$.) Hence $x = x_1x_2$ centralizes~$\tilde{N}$,
and thus~$x$ stabilizes any extension of~$\lambda$ to~$\tilde{I}$. As there
is an $n\tilde{n}$-invariant such extension, our claim follows.
\end{prf}

\begin{cor}
\label{MaximalExtendability}
Let~$\mathbf{M}$ denote an $e$-split Levi subgroup of~$\mathbf{G}$ for some 
$e \in \{ 1, 2, 3, 4, 6, 8, 12 \}$. Suppose that $N_G(\mathbf{M}) = N_G( M )$.
Then~$M$ satisfies the maximal extendibility 
condition, i.e.\ every irreducible character of~$M$ extends to its inertia
subgroup in~$N_G(M)$.
\end{cor}
\begin{prf}
If~$\mathbf{M}$ is a maximal torus, the result follows 
from~\cite[Theorem~$1.1$]{Spaeth10}. If the semisimple rank of~$\mathbf{M}$ 
equals~$3$, then $N_G(\mathbf{M})/M$ is cyclic, and thus the claim also holds. 
The remaining cases are exactly those treated in 
Propositions~\ref{NormalizerM13},~\ref{NormalizerM14} and~\ref{NormalizerM18},
as every parabolic subsystem of~$\Sigma$ of rank~$2$ is conjugate in~$W$ to
one with base $\{ \alpha_1, \alpha_2 \}$, $\{ \alpha_3, \alpha_4 \}$,
$\{ \alpha_1, \alpha_4 \}$ or $\{ \alpha_2, \alpha_3 \}$.
\end{prf}

\medskip
\noindent
We end this section by clarifying the condition in 
Corollary~\ref{MaximalExtendability}.

\begin{lem}
\label{NormalizersOfSplitLeviSubgroups}
Let~$\mathbf{M}$ denote an $e$-split Levi subgroup of~$\mathbf{G}$ for some
$e \in \{ 1, 2, 3, 4, 6, 8, 12 \}$. Then $N_G(\mathbf{M}) = N_G( M )$
unless $(e,q) = (1,2)$.
\end{lem}
\begin{prf}
First assume that $e = 1$,  and let $\mathbf{T} = \mathbf{T}_0$ denote the 
standard $1$-$F$-split maximal torus of~$\mathbf{M}$. Thus~$T$ is a complement
to the Sylow $p$-subgroup~$U_M$ in a Borel subgroup~$B_M$ of~$M$. By the 
Schur-Zassenhaus theorem, any two such complements are conjugate in~$B_M$. 
Now~$M$ is a finite group with a split $BN$-pair of characteristic~$p$; 
see~\cite[1.18]{C2}. The Bruhat decomposition implies in particular that 
$N_M( U_M ) = B_M$. Hence~$T^y$ and~$T$ are conjugate in~$M$ for all 
$y \in N_G(M)$. It follows that 
	$$N_G( M ) = (N_G(T) \cap N_G(M))M.$$ 
Now let $y \in N_G( \mathbf{M} )$. Then~$\mathbf{T}^y$ is a $1$-$F$-split 
maximal torus of~$\mathbf{M}$, and thus conjugate to~$\mathbf{T}$ in~$M$.
Hence $N_G( \mathbf{M} ) \leq (N_G( \mathbf{T} ) \cap N_G(M))M$. In fact, 
	$$N_G( \mathbf{M} ) = (N_G( \mathbf{T} ) \cap N_G(M))M.$$
(Although this is well known, we sketch a proof for the readers convenience. 
We may assume that $\mathbf{M} = \langle \mathbf{T}, 
\mathbf{U}_\alpha \mid \alpha \in \Gamma \rangle$ for some closed subsystem
$\Gamma \subseteq \Sigma$. Then $M = \langle T, \mathbf{U}_\alpha^F \mid 
\alpha \in \Gamma \rangle$; see \cite[$1.18$]{C2}. Hence an element 
$n \in N_G( \mathbf{T} ) \cap N_G(M)$ permutes the finite root subgroups 
$\mathbf{U}_\alpha^F$ for $\alpha \in \Gamma$. But then $n$ also permutes 
the $\mathbf{U}_\alpha$ for $\alpha \in \Sigma$ and thus 
$n \in N_G( \mathbf{M} )$.) 
	
It thus suffices to show that $N_G( \mathbf{T} ) 
= N_G(T)$. Applying \cite[Theorem~$3.5.3$(i)]{C2}, we find that 
$C_{\mathbf{G}}( T ) = \mathbf{T}$ unless $q = 2$. Thus, apart from $q = 2$,
we have $C_G(T) = T$ which implies that
$N_G( T ) = T.W = N_G( \mathbf{T} )$ and hence our claim.
%%%%%%%%%%%%%%%%%%%%%%%%%%%%%%%%%%%%%%%%%%%%%%%%%%%%%%%%%%%%%%%%%%%%%%%%%%%%%%%%
%%%%%%%%%%%%%%%%%%%%%%%%%%%%%%%%%%%%%%%%%%%%%%%%%%%%%%%%%%%%%%%%%%%%%%%%%%%%%%%%
%%
%% See sheet of 11.03.2021.
%%
%%%%%%%%%%%%%%%%%%%%%%%%%%%%%%%%%%%%%%%%%%%%%%%%%%%%%%%%%%%%%%%%%%%%%%%%%%%%%%%%
%%%%%%%%%%%%%%%%%%%%%%%%%%%%%%%%%%%%%%%%%%%%%%%%%%%%%%%%%%%%%%%%%%%%%%%%%%%%%%%%
	
Now assume that $e \geq 2$.
If $\mathbf{M} = C_{\mathbf{G}}( Z( M ) )$, then $N_G( \mathbf{M} ) \leq 
N_G ( M ) \leq N_G( Z( M ) ) \leq N_G( C_{\mathbf{G}}( Z( M ) ) ) = 
N_G( \mathbf{M} )$, and hence the claim holds. If there is a semisimple element 
$s \in G$ with $\mathbf{M} = C_{\mathbf{G}}( s )$, then $\mathbf{M} \leq 
C_{\mathbf{G}}( Z( M ) ) \leq C_{\mathbf{G}}( s ) = \mathbf{M}$, since
$Z(M) = Z( \mathbf{M} )^F \leq Z( \mathbf{M} )$. Going through 
the tables in~\cite{LL}, we find that such an element exists unless
$e = 2$ and $q \leq 3$ or $(e,q) \in \{ (3,2), (4,2), (6,3) \}$. 
In the former two cases, we find $C_{\mathbf{G}}( Z(M) ) = \mathbf{M}$
by applying \cite[Theorem~$3.5.3$(i)]{C2} to the elements of order~$3$,
respectively~$4$ of~$Z(M)$.
In the latter three cases, $|Z( M )|$ is divisible by a prime larger than~$3$, 
and the claim follows from \cite[Proposition~$2.3$]{MalleAb}.
\end{prf}

%\markboth{The blocks of $F_4(q)$}
\section{The blocks of $F_4(q)$}
\label{ThreeBlocks}

We keep the notation introduced in Subsection~\ref{SetupF4}. Our aim in this
section is to describe the $\ell$-blocks of $G = F_4( q )$ for primes 
$\ell \nmid q$. If $b$ is such a block, we write $d(b)$ for the defect of~$b$ 
and $l(b)$ for the number of its irreducible Brauer characters.

\subsection{The $\ell$-blocks for good primes}
\label{GoodPrimesBlocks}
Let $\ell > 3$ be a prime dividing $|G|$. Then there is a unique $e \in 
\{ 1, 2, 3, 4, 6, 8, 12 \}$ such that $\ell \mid \Phi_e( q )$, i.e.\ 
Condition~(*) of \cite{MalleAb} is satisfied. Notice that $e = e_\ell(q)$ is
the order of~$q$ in the multiplicative group of the field $\mathbb{F}_\ell$.

We recall the description of the $\ell$-blocks and their defect groups as 
summarized in \cite[Theorem~$3.6$]{MalleAb}. Let $b$ be an $\ell$-block with 
defect group~$D$. Suppose that $b \subseteq \cE_\ell( G, s )$ for some 
semisimple $\ell'$-element $s \in G$. Then there is an $e$-split Levi subgroup 
$\mathbf{L}$ of~$\mathbf{G}$ such that $s \in L^\dagger$ (recall the notion 
$\mathbf{L}^\dagger$ introduced in the last paragraph of 
Subsection~\ref{TwistingAndDuality}), and~$D$ is the Sylow $\ell$-subgroup 
of~$Z( L )$. Moreover, there is $\vartheta \in \cE( L, s )$ which is 
$e$-cuspidal, has~$D$ in its kernel, is of defect~$0$ when viewed as 
a character of $L/D$, and $\Irr( b ) \cap \cE(G,s)$ is the set of constituents 
of $R_{\mathbf{L}}^{\mathbf{G}}( \vartheta )$. This \textit{$e$-cuspidal pair}
$(\mathbf{L},\vartheta)$ is determined by~$b$ up to $G$-conjugacy.
If $D$ is non-cyclic, then $e \neq 8, 12$, and
$\mathbf{L}$ is a maximal torus of~$\mathbf{G}$, unless $e = 1, 2$.
If~$s$ is not quasi-isolated, then $d( b )$ 
and $l( b )$ can be determined from the corresponding numbers in 
$\cE_\ell( M, s )$, where $\mathbf{M}^\dagger$ is a regular subgroup 
of~$\mathbf{G}$ minimal with the property that it contains $C_{\mathbf{G}}( s )$;
see Theorem~\ref{BoRoEtAl}. If $s \neq 1$ is quasi-isolated, the corresponding 
blocks and their invariants are determined in \cite{KeMa,RH}, and for $s = 1$ in 
\cite{BrouMaMi}.
Conversely, every $e$-cuspidal pair $(\mathbf{L}, \vartheta)$ as above gives 
rise to a corresponding $\ell$-block of~$G$.

The $\ell$-blocks and their invariants in those cases, where $C_G( s )$ does not
have a cyclic Sylow $\ell$-subgroup, are described in Tables~\ref{1}--\ref{19}.
Moreover, Table~\ref{21} contains this information for $\ell \in \{ 5, 7 \}$ for
the blocks of the exceptional double cover of $F_4(2)$. The invariants $l(b)$ in
these cases can be found in~\cite{HF42}.

\subsection{The $2$-blocks}
\label{TheTwoBlocks}
Let us assume that~$q$ is odd in this subsection. Let~$s$ be a semisimple 
$2'$-element. Then the $\mathbf{G}$-class type of~$s$ is one of $\{ 1, 4, 6, 7, 
9, 10, 13\text{--}15, 17\text{--}20 \}$. In particular,~$s$ is quasi-isolated if 
and only if $s = 1$ or of $\mathbf{G}$-class type~$4$ (in which case 
$3 \nmid q$). 

Suppose that $s$ is not quasi-isolated. Then $C_{\mathbf{G}}( s )$ 
is a proper regular subgroup of~$\mathbf{G}$, which is a classical group and 
thus has a unique unipotent $2$-block. It follows from 
Theorem~\ref{BoRoEtAl}(a), that $\cE_2( G, s )$ is a single block, whose defect 
groups are isomorphic to, but in general not conjugate to a Sylow $2$-subgroup 
of~$C_G( s )$. 

The case $s = 1$ corresponds to the unipotent blocks of~$G$. By 
\cite[{\sc Th{\'e}or{\`e}me~}A and Table on p.~$349$]{En00}, 
there are three unipotent $2$-blocks of~$G$: the principal block and two blocks 
of defect~$0$, containing the characters $F_4[\theta]$ and $F_4[\theta^2]$, 
respectively. By \cite[Table~$1$]{GeGuide}, we have $l(b) = 28$ for the 
principal $2$-block~$b$ of~$G$.

Now assume that $3 \nmid q$ and let $s$ be of class type $(4,1)$, respectively 
$(4,2)$. The former case occurs if $3 \mid q - 1$, the latter if $3 \mid q + 1$.

\addtocounter{thm}{2}
\begin{prop}
\label{QuasiIsolatedTwoBlock}
Let $s$ be a quasi-isolated element of order~$3$ of $\mathbf{G}$-class 
type~$4$. Then $\cE( G, s )$ is a basic set for $\cE_2( G, s )$. In particular, 
$\cE_2( G, s )$ has exactly~$9$ irreducible Brauer characters. Moreover, if 
$3 \mid q - 1$, the decomposition matrix of $\cE_2( G, s )$ is unitriangular.
\end{prop}
\begin{prf}
In this proof, by the \textit{unipotent decomposition matrix} of~$H$, where~$H$ 
is a finite reductive group, we understand the matrix of scalar products of the 
unipotent characters of~$H$ with the projective indecomposable characters in 
$\cE_2(H,1)$. This is a square matrix whenever $\cE( H, 1 )$ is a basic set for 
$\cE_2(H,1)$. 

We may assume that $C_{\mathbf{G}}( s ) = \mathbf{L}$, the group introduced in 
Proposition~\ref{C3C}, so that $L = C_G( s ) = (L^1 \circ_3 L^2).3$ with 
$L^i \cong \SL_3^\varepsilon( q )$ for $i = 1, 2$. By 
\cite[Theorem~$12$]{CaEn93}, the decomposition matrix of $\cE_2( L^i, 1 )$
and $\cE_2( L, 1 )$, $i = 1, 2$, is the same as the decomposition matrix
of $\cE_2( \GL_3^\varepsilon( q ), 1)$, respectively 
$\cE_2( \GL_3^\varepsilon( q ) \times \GL_3^\varepsilon( q ), 1)$. These facts
will be assumed tacitly in the following.

We first prove the assertions if $3 \mid q - 1$. 
Consider the split Levi subgroup $H := \GL_2( q ) \times [q-1]$ of $\GL_3( q )$.
The unipotent decomposition matrices of~$\GL_3(q)$ and~$H$ are given as follows, 
where we omit entries equal to~$0$; see \cite[Appendix~$1$]{JamesGL10}: 
$$
%\begin{equation}
%\label{SomeUnipotentDecompositionMatrices}
       \begin{array}{ccc} \hline\hline
	       (3) & (2,1) & (1^3) \rule[- 3pt]{0pt}{ 16pt} \\ \hline\hline
               1 & & \rule[ 0pt]{0pt}{ 13pt} \\
                 & 1 & \\
               1 &   & 1 
	       \rule[- 2pt]{0pt}{ 5pt} \\ \hline\hline
       \end{array}
	\quad\quad
	       \begin{array}{cc} \hline\hline
               (2) & (1^2) \rule[- 3pt]{0pt}{ 16pt} \\ \hline\hline
               1 & \rule[ 0pt]{0pt}{ 13pt} \\
               1 & 1
               \rule[- 2pt]{0pt}{ 5pt} \\ \hline\hline
       \end{array}
%\end{equation}
$$
Following \cite{JamesGL10}, we label the columns of these decomposition
matrices by partitions of~$3$, respectively~$2$. The unipotent part of the
projective indecomposable
character corresponding to a column labeled by~$\lambda$ is denoted by 
$\rho_\lambda$. Put $\kappa_{(3)} := \rho_{(3)} + 2\rho_{(2,1)}$ and
$\kappa_{(2,1)} := \rho_{(2,1)} + \rho_{(1^3)}$. Observe that~$\kappa_{(3)}$ 
and~$\kappa_{(2,1)}$ are obtained by Harish-Chandra induction of~$\rho_{(2)}$, 
respectively~$\rho_{(1^2)}$ from~$H$ to~$\GL_3( q )$.

Harish-Chandra inducing the~$12$ projective indecomposable characters of
$\cE_2( \GL_3(q) \times H, 1 )$ and $\cE_2( H \times \GL_3( q ), 1 )$ to 
$\GL_3(q) \times \GL_3(q)$, we obtain~$12$ projective characters in 
$\cE_2( \GL_3(q) \times \GL_3(q), 1 )$, whose unipotent parts are 
$\rho_\lambda \boxtimes \kappa_\mu$ and $\kappa_\mu \boxtimes \rho_\lambda$ 
with $\lambda \in \{ (3), (2,1), (1^3) \}$ and $\mu \in \{ (3), (2,1) \}$. We 
can select~$8$ of these projective characters of~$\GL_3(q) \times \GL_3(q)$ such 
that the matrix of scalar products of the selected characters with the unipotent 
characters of~$\GL_3(q) \times \GL_3(q)$ is unitriangular. Moreover, if we use 
this triangular shape to associate a pair of partitions to each of the~$8$ 
projective characters thus constructed, then only $((1^3),(1^3))$ is not 
associated to any of these. 

Now $\mathbf{L}$ contains the two split Levi subgroups $\mathbf{M}_1 := 
\mathbf{L}_{ \{ 1, 23, 4 \}}$ and $\mathbf{M}_2 
:= \mathbf{L}_{ \{ 1, 3, 4 \}}$ of $\mathbf{G}$, and~$s$ is 
a central element in each of~$M_1$ and~$M_2$. It follows from
\cite[Proposition~$11.4.8$(ii)]{DiMi2} that  $\cE_2( M_i^*, s )$ is Morita 
equivalent to $\cE_2( M_i, 1 )$, $i = 1, 2$. Moreover, Harish-Chandra induction 
from $\cE_2( M_1, 1 )$ to $\cE_2( L, 1 )$ corresponds to Harish-Chandra 
induction from $\cE_2( \GL_3(q) \times H, 1 )$ to 
$\cE_2( \GL_3(q) \times \GL_3(q), 1)$, and likewise for $\cE_2( M_2, 1 )$.

Every irreducible character of $M_i$, $i = 1, 2$ is uniform, as 
$Z( \mathbf{M}_i )$ is connected, each component of 
$[\mathbf{M}_i,\mathbf{M}_i]$ is of type~$A$, and $\mathbf{M}_1^* \cong 
\mathbf{M}_2$ by the remarks in Subsection~\ref{TwistingAndDuality}. Hence 
Harish-Chandra induction $\cE_2( M_i^*, s ) \rightarrow \cE_2( G, s )$ is equal 
to the composition of Harish-Chandra induction 
$\cE_2( M_i, 1 ) \rightarrow \cE_2( L, 1 )$ and Lusztig's Jordan decompositions 
$\cE_2( M_i^*, s ) \rightarrow \cE_2( M_i, 1 )$ and 
$\cE_2( L, 1 ) \rightarrow \cE_2( G, s )$; see \cite[Theorem~$15.8$]{CaEn}. As 
Harish-Chandra induction preserves projective characters, we obtain~$8$ 
projective characters of $\cE_2( G, s )$, whose matrix of scalar products with 
the elements of $\cE( G, s)$ is lower unitriangular. If we also consider the 
restriction of the Gelfand-Graev character of~$G$ to $\cE_2( G, s )$, we 
obtain~$9$ linearly independent projective characters of $\cE_2( G, s )$. By 
Lemma~\ref{GenerationLemma}, there are at most~$9$ irreducible Brauer characters 
in $\cE_2( G, s )$. The unitriangularity of the projective characters thus 
constructed implies the unitriangularity of the decomposition matrix. This 
concludes our proof in case $3 \mid q - 1$.
%%%%%%%%%%%%%%%%%%%%%%%%%%%%%%%%%%%%%%%%%%%%%%%%%%%%%%%%%%%%%%%%%%%%%%%%%%%%%%%%
%%%%%%%%%%%%%%%%%%%%%%%%%%%%%%%%%%%%%%%%%%%%%%%%%%%%%%%%%%%%%%%%%%%%%%%%%%%%%%%%
%%
%% For more details see Sheets 1 -- 3 of 06.04.2022.
%%
%%%%%%%%%%%%%%%%%%%%%%%%%%%%%%%%%%%%%%%%%%%%%%%%%%%%%%%%%%%%%%%%%%%%%%%%%%%%%%%%
%%%%%%%%%%%%%%%%%%%%%%%%%%%%%%%%%%%%%%%%%%%%%%%%%%%%%%%%%%%%%%%%%%%%%%%%%%%%%%%%

The proof in case $3 \mid q + 1$ uses an Ennola duality argument. We replace 
$\mathbf{L}$ and $\mathbf{M}_i$, $i = 1, 2$, by their conjugates with~$g$, 
where~$g$ is as in Proposition~\ref{C3C}. Then $L = (\SU_3( q ) \circ_3 
\SU_3( q )).3$,  and $\mathbf{M}_i$, $i = 1, 2$, is 
a $2$-split Levi subgroup of~$\mathbf{L}$ and of~$\mathbf{G}$. The 
$2$-decomposition matrix of the unipotent characters of $\SU_3( q )$ is a
lower unitriangular matrix, which has been computed by Erdmann in \cite{erd} for
the case $4 \mid q - 1$, and by the second author in \cite[Appendix]{HiHe}
for the case $4 \mid q + 1$. The proof now 
proceeds as above, except that Harish-Chandra induction $\cE( M_i^*, s ) 
\rightarrow \cE( G, s )$ is replaced by Lusztig induction. As the latter 
preserves generalized projective characters, we still obtain a set of~$9$ 
linearly independent generalized projective characters in $\cE_2( G, s )$.
\end{prf}

\addtocounter{subsection}{1}
\subsection{The $3$-blocks}
\label{TheThreeBlocks}
Let us assume that $3 \nmid q$ in this subsection. The $3$-blocks of~$G$ and 
their relevant invariants are described in Tables~\ref{1}--\ref{19} and~\ref{21} 
of the appendix, where the latter table contains the information for the 
exceptional double cover of $F_4(2)$. The derivation of these results is by far
the most laborious part of this work.
We are now going to prove, in a series of lemmas, that the entries in 
Columns~$5$--$7$ of Tables~\ref{1}--\ref{19} are correct. That is, for a 
semisimple $3'$-element $s \in G$ and a block $b \subseteq \cE_3( G, s )$
of positive defect, we determine a label for~$b$, as well as~$d(b)$ and~$l(b)$.
We define $\varepsilon \in \{ 1, -1 \}$ by the condition that 
$3 \mid q - \varepsilon$, and write $e := e_3(q)$ for the order of~$q$ 
in~$\mathbb{F}_3$. Thus $e = 1$ if $\varepsilon = 1$, and $e = 2$, otherwise.

\addtocounter{thm}{1}
\begin{lem}
\label{UnipotentsGenerate}
Let~$s \in G$ be a semisimple $3'$-element. Then every irreducible $3$-modular 
character~$\varphi$ of $\mathcal{E}_3( G, s )$ is a $\mathbb{Z}$-linear 
combination of elements of $\{ \check{\chi} \mid \chi \in \mathcal{E}( G, s ) \}$, 
unless $s = 1$, in which case~$\varphi$ is still a $\mathbb{Q}$-linear 
combination of this set.

In particular, the number of irreducible $3$-modular characters contained in 
$\mathcal{E}_3( G, s )$ is equal to the rank of the $\mathbb{Z}$-span of 
$\{ \check{\chi} \mid \chi \in \mathcal{E}( G, s ) \}$.
\end{lem}
\begin{prf}
Let~$t \in C_G( s )$ be a non-trivial $3$-element. By the list of semisimple 
class types of~$G$ in \cite{LL}, we find that $C_{G}( st )$ satisfies one of the 
hypotheses~(a) or~(b) of Lemma~\ref{GenerationLemma}, unless $s = 1$ and~$t$ 
is a quasi-isolated element of order~$3$ and class type~$4$. In the latter case, 
$C_{G}( st ) = C_{G}( t )$ satisfies hypothesis~(c) of the lemma, as all 
components of $[ C_{\mathbf{G}}( t ), C_{\mathbf{G}}( t ) ]^F$ are of type~$A$
and thus all unipotent characters of $C_G(t)$ are uniform. The result follows. 
For $s = 1$ our assertion follows from \cite[Proposition~$7.14$]{GeGuide} and the
fact that decomposition numbers are integers.

The final statement is a consequence of the fact that every $\check{\chi}$
for $\chi \in \mathcal{E}( G, s )$ is a $\mathbb{Z}$-linear combination of 
irreducible Brauer characters contained in $\mathcal{E}_3( G, s )$.
\end{prf}

\medskip

\noindent 
We now prove our claims for the unipotent blocks. The unipotent $3$-blocks 
of~$G$ have been computed by Enguehard in
\cite[p.~349--351]{En00}.  We follow Carter's book \cite[p.~478f]{C2} for 
the notation of the unipotent characters of~$G$.

\begin{lem}
\label{UnipotentsDefect0}
The unipotent characters of~$G$ of $3$-defect~$0$ are $F_4[i]$, $F_4[-i]$, 
$F_4^I[1]$, $F_4[-1]$ if $\varepsilon = 1$, and $F_4[i]$, $F_4[-i]$, 
$\phi_{4,8}$, $\phi_{16,5}$, if $\varepsilon = -1$.
\end{lem}
\begin{prf}
See {\rm \cite[Table on Page~$349$]{En00}}.
\end{prf}

\begin{lem}
\label{UnipotentsPositiveDefect}
The invariants of the unipotent $3$-blocks of~$G$ of positive defect are as
given in {\rm Columns~$5$--$7$} of {\rm Table~\ref{1}}.
\end{lem}
\begin{prf}
It is known that~$G$ has exactly~$35$ unipotent $3$-modular characters; see 
\cite[Table in~$6.6$]{GeHi}. We sketch an alternative proof for this fact.
By a result of Shoji \cite[$6.2.4$(c) and Proposition~$6.3$]{ShojiI} (see also
the remarks on \cite[p.~$42$]{Ge2019} for the fact that Shoji's results are 
valid without restrictions on~$p$), the two 
almost characters corresponding to the unipotent characters $F_4[\theta]$ and 
$F_4[\theta^2]$ (in the sense of \cite[$(4.24.1)$]{luszbuch}) have values~$0$
except on $3$-singular classes. As~$G$ has exactly~$37$ ordinary unipotent 
characters, and as the almost characters span the same space as the unipotent 
characters; see \cite[Corollary~$4.25$]{luszbuch}, it follows that 
$\{ \check{\chi} \mid \chi \in \mathcal{E}( G, 1 ) \}$ spans a space of dimension 
at most~$35$.  On the other hand, the matrix of values of this set has rank at 
least $35$, as can be checked with the explicit unipotent character table of~$G$ 
computed by K\"ohler \cite{KoePhD}, and, independently, by the third author. For 
this, it suffices to look at the set of unipotent classes and the mixed classes 
on which the other unipotent almost characters have non-zero values, as well as 
some classes of elements of the form $su$, where~$s$ is an involution with 
centralizer of type~$B_4$. Lemma~\ref{UnipotentsGenerate} now implies the 
result.

According to \cite[{\sc Th\'eor\`eme} A and the table on Page~$349$]{En00}, 
there is a unique non-principal block~$b$ of positive defect corresponding to
the $e$-cuspidal pair $([q-\varepsilon]^2.B_2(q), \zeta_e)$, where~$\zeta_e$ is
the $e$-cuspidal unipotent character of the semisimple component~$B_2(q)$ of 
$[q-\varepsilon]^2.B_2(q)$. Moreover,~$b$ contains exactly the unipotent 
characters of the $e$-Harish-Chandra series defined by 
$([q-\varepsilon]^2.B_2(q), \zeta_e)$. These are the characters~$B_{2,1}$, 
$B_{2,\varepsilon}$, $B_{2,r}$, $B_{2,\varepsilon'}$, $B_{2,\varepsilon''}$ if
$\varepsilon = 1$, and $\phi_{4,1}$, $\phi_{4,13}$, $B_{2,r}$, $\phi_{4,7''}$,
$\phi_{4,7'}$, otherwise (for the elements in the $2$-Harish-Chandra series 
see \cite[Table~$2$, case~$1$]{BrouMaMi}).
Using the explicit values of the unipotent characters 
in~$b$, we can check that their restrictions to the $3$-regular elements are 
linearly independent. In view of Lemma~\ref{UnipotentsGenerate}, this shows 
that $l(b) = 5$. This, together with Lemma~\ref{UnipotentsDefect0} implies 
$l(B) = 26$ for the principal $3$-block~$B$ of~$G$.
\end{prf}

\smallskip

\noindent
We next show how to determine the number of irreducible $3$-modular characters 
for the non-unipotent blocks.

\begin{prop}
\label{EllB}
Let $1 \neq s \in G$ be a semisimple $3'$-element. Then $\mathcal{E}( G, s )$ is 
a basic set for $\mathcal{E}_3( G, s )$. In particular, if $b \subseteq 
\mathcal{E}_3( G, s )$ is a $3$-block of~$G$, then 
$l(b) = | \Irr( b ) \cap \mathcal{E}( G, s )|$.
\end{prop}
\begin{prf}
Clearly, the second statement follows from the first.
To prove the first statement, we refine the argument
of \cite[Proposition~$4.1$]{GeHi1}. Let~$\mathcal{S}_3(G)$ denote a set
of representatives for the $G$-conjugacy classes of $3$-elements of~$G$. For 
each $t \in \mathcal{S}_3(G)$, let $\mathcal{C}_{3'}( C_G( t ) )$ be a set of 
representatives for the $C_G(t)$-conjugacy classes of $3'$-elements of $C_G(t)$, 
and denote by $\mathcal{S}_{3'}( C_G( t ) )$ the subset of 
$\mathcal{C}_{3'}( C_G( t ) )$ consisting of semisimple elements. Then 
$\{ ts \mid t \in \mathcal{S}_{3}( G ), s \in 
\mathcal{S}_{3'}( C_G( t ) ) \}$ is a set of representatives for the 
semisimple conjugacy classes of~$G$. 

By $t_0$ we denote the element of $\mathcal{S}_{3}( G )$ of class $(4,1)$ if 
$3 \mid q - 1$, respectively $(4,2)$ if $3 \mid q + 1$. First, consider 
$t \in \mathcal{S}_3(G) \setminus \{ 1, t_0 \}$. We claim that in this case
$\mathbf{M} := C_{\mathbf{G}}( t )$ is a regular subgroup of~$\mathbf{G}$. 
If not, $|( Z( \mathbf{M} )/ Z( \mathbf{M} )^\circ )^F|$ is a $2$-group; see the
tables in \cite{LL}. As~$t$ is a $3$-element, Lemma~\ref{NonRegularCentralizers}
implies that $\mathbf{M} \lneq C_{\mathbf{G}}( t )$, a contradiction.
Thus~$\mathbf{M}$ is regular and there is a regular subgroup 
$\mathbf{M}^\dagger$ of~$\mathbf{G}$ in duality with~$\mathbf{M}$. In 
particular, $Z( \mathbf{M}^\dagger )$ is connected and~$3$ is a good prime 
for~$\mathbf{M}^\dagger$, so that $\mathbf{M}^\dagger$
satisfies the assumptions of \cite[Theorem~$5.1$]{GeHi1}. Putting $M^\dagger
:= (\mathbf{M}^\dagger)^F$ and identifying $\mathbf{M}$ with 
$(\mathbf{M}^\dagger)^*$, we find that 
$|\mathcal{E}(M^\dagger,s)|$ equals the number of irreducible $3$-modular 
characters contained in $\mathcal{E}_3( M^\dagger, s )$ for every $s \in 
\mathcal{S}_{3'}( M )$, and hence
\begin{equation}
\label{Good3Elements}
\sum_{s \in \mathcal{S}_{3'}( M )} |\mathcal{E}( M^\dagger, s )| = 
|\mathcal{C}_{3'}( M^\dagger )|.
\end{equation}

Now consider the case $t = t_0$ and put $\mathbf{L} := C_{\mathbf{G}} (t_0 )$.
If $1 \neq s \in \mathcal{S}_{3'}( L )$, then $\mathbf{M} := C_{\mathbf{L}} ( s ) = 
C_{\mathbf{G}}( t_0s )$ is a regular subgroup of~$\mathbf{G}$ and hence 
of~$\mathbf{L}$. It follows from the Bonnaf{\'e}-Rouquier Morita equivalence
theorem \cite[Th{\'e}or{\`e}me~$11.8$]{BoRo} (see also 
Theorem~\ref{BoRoEtAl}(c)), that $\mathcal{E}_3( L^*, s )$ and 
$\mathcal{E}_3( M, 1 )$ have the same number of irreducible $3$-modular 
characters. As~$\mathbf{M}$ has connected center and~$3$ is a good prime for 
$\mathbf{M}$, the latter number equals
$|\mathcal{E}( M, 1 )|$ by \cite[Theorem~$5.1$]{GeHi1}. By the Jordan 
decomposition of characters, we have $|\mathcal{E}( M, 1 )| = 
|\mathcal{E}( L^*, s )|$.
Thus the number of irreducible $3$-modular characters in
$\mathcal{E}_3( L^*, s )$ equals $|\mathcal{E}( L^*, s )|$.
Let us now compute the number of irreducible $3$-modular characters in
$\mathcal{E}_3( L^*, 1 )$. We have $L = (\SL_3^\varepsilon( q ) \circ_3
\SL_3^\varepsilon( q )).3$, where the outer automorphism of order~$3$ acts as
a simultaneous diagonal automorphism on each of the factors 
$\SL_3^\varepsilon( q )$ of~$L$; see Subsection~\ref{sec3elts} below and 
Remark~\ref{DiagonalActionOfx}. Now~$\SL_3^\varepsilon( q )$ has exactly~$5$
unipotent $3$-modular characters, three of which lie in an orbit under
the outer automorphism of order~$3$. (A reference for the latter two statements
in case of $\SU_3(q)$ is \cite[Theorem~$4.5$]{Geck3U}; the corresponding results
for $\SL_3(q)$ are proved in the same way.) Thus $\SL_3^\varepsilon( q ) \circ_3
\SL_3^\varepsilon( q )$ has $25$ unipotent $3$-modular characters, four
of which are fixed under the outer automorphism of order~$3$, and the other~$21$
lie in~$7$ orbits of length~$3$. It follows that~$L$ has exactly~$11$ unipotent
$3$-modular characters. On the other hand,~$L$ has exactly~$9$ ordinary
unipotent characters. We conclude that
\begin{equation}
\label{Ct0}
\sum_{s \in \mathcal{S}_{3'}( L )} |\mathcal{E}( L^*, s )| = 
-2 + |\mathcal{C}_{3'}( L^* )|.
\end{equation}
Let $c \in \mathbb{Z}$ be such that 
\begin{equation}
\label{G}
\sum_{s \in \mathcal{S}_{3'}( G )} |\mathcal{E}( G, s )| = 
c + |\mathcal{C}_{3'}( G )|.
\end{equation}
Now we use the fact that $|\cE(G, ts )| = |\cE( C_G( t )^*, s )|$ for
all $t \in \mathcal{S}_3(G)$ and all $s \in \mathcal{S}_{3'}( C_G( t ) )$.
This follows from the identification $G = G^*$ and the Jordan decomposition 
of characters: 
\begin{eqnarray*}
|\cE( C_G( t )^*, s )| & = & |\cE( C_{C_G(t)}( s ), 1 )| \\
 & = & |\cE( C_G(ts), 1 )| \\
 & = & |\cE( C_{G^*}(ts), 1 )| \\
 & = & |\cE( G, ts )|.
\end{eqnarray*}
We find, using~(\ref{G}),~(\ref{Ct0}) and~(\ref{Good3Elements}), that
\begin{eqnarray*}
\begin{array}{c} \text{number of conjugacy} \\ \text{classes of\ }G \end{array}
& = & \sum_{t \in \mathcal{S}_3( G )} \sum_{s \in \mathcal{S}_{3'}( C_G( t ) )} |\mathcal{E}( G, ts )| \\
& = & \sum_{t \in \mathcal{S}_3( G )} \sum_{s \in \mathcal{S}_{3'}( C_G( t ) )} |\mathcal{E}( C_G( t )^* , s )| \\
& = & c - 2 +  \sum_{t \in \mathcal{S}_3( G )} |\mathcal{C}_{3'}( C_G( t )^* )|.
\end{eqnarray*}
We claim that
$\sum_{t \in \mathcal{S}_3( G )} |\mathcal{C}_{3'}( C_G( t )^* )|$ equals the 
number of conjugacy classes of~$G$. First, if $t = 1$ or $t = t_0$, then 
$C_{\mathbf{G}}( t )^* \cong C_{\mathbf{G}}( t )$ and $C_G( t )^* \cong 
C_G ( t )$ by Remark~\ref{SelfDualityOfL}. It follows from
Proposition~\ref{Proposition42} that there is a permutation $t \mapsto t'$
of $\mathcal{S}_3(G) \setminus \{ 1, t_0 \}$ such that
$C_{G}( t ) \cong C_{G}( t' )^*$.
Hence $\sum_{t \in \mathcal{S}_3( G )} |\mathcal{C}_{3'}( C_G( t )^* )|
= \sum_{t \in \mathcal{S}_3( G )} |\mathcal{C}_{3'}( C_G( t ) )|$,
giving our claim.

We conclude that $c = 2$. By Table~\ref{1}, the number of unipotent $3$-modular 
characters is two less than $|\mathcal{E}( G, 1 )|$.
By Lemma~\ref{UnipotentsGenerate}, for all $s \in \mathcal{S}_{3'}( G )$, 
the number of irreducible $3$-modular characters in $\mathcal{E}_3( G, s )$ is 
at most equal to $|\mathcal{E}( G, s )|$. We conclude from~(\ref{G}) that 
equality holds for all non-trivial such~$s$. This implies our claim.
\end{prf}

\begin{lem}
\label{NonUnipotentIsolated}
The invariants contained in {\rm Columns~$5$--$7$} of 
{\rm Tables~\ref{2},~\ref{3}} and~{\rm \ref{5}} are correct.
\end{lem}
\begin{prf}
The blocks and their labels are described in \cite[Table~$2$]{KeMa}. The defects 
of these blocks can be derived from \cite[Proposition~$3.2$]{KeMa}. The numbers 
of the irreducible $3$-modular characters in these blocks can be determined with 
Proposition~\ref{EllB}.
\end{prf}

\begin{lem}
\label{NonUnipotentNonIsolated}
The invariants contained in {\rm Columns~$5$--$7$} of 
{\rm Tables~\ref{6}--\ref{19}} and~{\rm \ref{21}} are correct.
\end{lem}
\begin{prf}
Let~$s$ be a semisimple $3'$-element of~$G$ such that $C_{\mathbf{G}}( s ) =
C_{\mathbf{G}^*}(s)$ is contained in a proper regular subgroup~$\mathbf{M}^*$ 
of $\mathbf{G}^* = \mathbf{G}$. We 
choose~$\mathbf{M}^*$ as in Theorem~\ref{BoRoEtAl}. As every component of
such an~$\mathbf{M}^*$ is of classical type, we may apply Theorem~\ref{BoRoEtAl}.
By appealing to~\cite{FS82} and~\cite{Fong}, we easily obtain all the entries of 
Tables~\ref{6}--\ref{19}.

For the invariants in Table~\ref{21} see \cite{HF42}.
\end{prf}

\addtocounter{subsection}{6}
\subsection{The action of outer automorphisms}
\label{DefinitionOfF0}

Recall that~$p$ is a prime and that $q = p^f$ for some positive integer~$f$. 
Recall also the definition of~$\sigma_1$ and $\Aut_1( \mathbf{G} ) = 
\mathbf{G} \rtimes \langle \sigma_1 \rangle$ from 
Subsection~\ref{Automorphisms}.
Let us put $f' := f$ if~$p$ is odd, and $f' := 2f$ if $p = 2$. Then 
$F = \sigma_1^{f'}$. In particular,~$\sigma_1$ commutes with~$F$ and thus~$G$
is $\sigma_1$-invariant. We tacitly use the symbol~$\sigma_1$ to also denote
the restriction of~$\sigma_1$ to~$G$. With this notation,
$\Aut(G) = G \rtimes \langle \sigma_1 \rangle$ and 
$\Out(G) = \langle G\sigma_1 \rangle$ is cyclic of order~$f'$; see
\cite[Theorem~$2.5.12$(a),(e)]{Gor}. In particular, every subgroup of 
$\langle \sigma_1 \rangle \leq \Aut( G )$ is of the form 
$\langle \sigma_1^{m'} \rangle$ for some integer $m'$ with $m' \mid f'$. 
Notice also that every automorphism of~$G$ extends to an element 
of~$\Aut_1(\mathbf{G})$ which commutes with~$F$.

Let $\sigma := \sigma_1^{m'}$ for some positive integer~$m'$ with $m' \mid f'$. 
The restriction of~$\sigma$ to~$G$, also denoted by~$\sigma$, has order~$f'/m'$. 
We define~$m$ by $m := m'/2$ if $m'$ and~$p$ are even and by $m := m'$, 
otherwise. Then $m \mid f$ and $F = (F_1^m)^{f/m}$ in all cases; moreover 
$\sigma = F_1^m$ if either~$m'$ and~$p$ are even or if~$p$ is odd, and
$\sigma^2 = F_1^m$, otherwise. 
The field automorphism of~$G$ of order $f/m'$ is induced 
by~$\sigma$ if~$p$ is odd, and by~$\sigma^2$ if $p = 2$. 

As always, we identify~$\mathbf{G}$ with its dual~$\mathbf{G}^*$ and~$G$ 
with~$G^*$. We let $\ell$ be an odd prime different from~$p$ (we do not 
assume $\ell = 3$ here) and write $e = e_\ell(q)$ for the order of~$q$ in the 
multiplicative group of the field $\mathbb{F}_\ell$.
We will need the following consequence of the Lang-Steinberg theorem.

\addtocounter{thm}{1}
\begin{lem}
\label{SigmaStableConjugate}
Let $s \in G$ be semisimple such that $\sigma(s)$ is $G$-conjugate to~$s$.
Then some $G$-conjugate~$t$ of~$s$ satisfies $\sigma(t) = t$.
\end{lem}
\begin{prf}
Let~${C}$ denote the $\mathbf{G}$-conjugacy classes of~$s$. Our hypothesis 
implies that~$C$ is $\sigma$-stable. An application of the Lang-Steinberg
theorem shows that $C \cap \mathbf{G}^\sigma$ is non-empty and that $C \cap G$
is a conjugacy class of~$G$; see, e.g.\ \cite[Example~$1.4.10$]{GeMa}. 
If $t \in C \cap \mathbf{G}^\sigma$, then $t,s \in C \cap G$ yielding our claim.
\end{prf}

\medskip
\noindent We can now prove the main result of this subsection.

\begin{prop}
\label{ActionAutomorphisms}
Let~$s \in G$ be a semisimple $\ell'$-element. Then 
\begin{equation}
\label{FPrimeAndLusztigSeries}
\sigma( \mathcal{E}( G, s ) ) = \mathcal{E}( G, \sigma^{-1}(s) ) \text{\rm\ and\ }
\sigma( \mathcal{E}_\ell( G, s ) ) = \mathcal{E}_\ell( G, \sigma^{-1}(s) ).
\end{equation}
Suppose that~$\sigma(s)$ is conjugate to~$s$ in~$G$. Then~$\sigma$ stabilizes
$\cE( G, s )$ and $\cE_\ell( G, s )$. 

{\rm (a)} Suppose that $p = 2$ and that~$m'$ is odd. Then either 
$C_{\mathbf{G}}( s )$ is a maximal torus, or the $G$-class type of~$s$ is one of 
$(1,1)$, $(4,1)$, $(4,2)$, $(14,1)$, $(14,4)$, $(15,1)$, $(15,3)$ or $(15,5)$.
In these cases, the non-trivial $\sigma$-orbits on $\cE( G, s )$ have length~$2$, 
and the number of such orbits is as given in the last column of the following 
table.
\setlength{\extrarowheight}{0.5ex}
\begin{equation}
\label{TableIn415a}
\begin{array}{rc|c|c}
i & k & C_G( s ) & \text{\rm no.} \\ \hline\hline
1 & 1 & G & 8 \\
4 & 1,2 & (\SL_3^\varepsilon(q) \circ_3 \SL_3^\varepsilon(q)).3 & 3 \\
14 & 1, 4 & [q-\varepsilon]^2 \times \SL_2(q)^2 & 1 \\
15 & 1, 3 & [q-\varepsilon]^2 \times \Sp_4(q) & 1 \\
15 & 5 & [q^2+1] \times \Sp_4(q) & 1 \\ \hline\hline
\end{array}
\end{equation}
Here, $\varepsilon \in \{ -1, 1 \}$, and the two values for~$k$ in the second 
column, if present, correspond to the cases $\varepsilon = 1$ and $\varepsilon = -1$,
respectively.

{\rm (b)} Suppose that~$p$ is odd or that $p = 2$ and~$m'$ is even, so that
$\sigma = F_1^m$. Then either~$\sigma$ fixes every element of 
$\mathcal{E}( G, s )$, or the $G$-class type of~$s$ is one of $(12,1)$, $(16,1)$ 
or~$(16,9)$. In these cases, the non-trivial $\sigma$-orbits on $\cE( G, s )$ 
have length~$2$, and the number of such orbits is as given in the last column of 
the following table.
\setlength{\extrarowheight}{0.5ex}
\begin{equation}
\label{TableIn415b}
\begin{array}{rc|c|c}
i & k & C_G( s ) & \text{\rm no.} \\ \hline\hline
12 & 1 & ([q-1]\circ_{2} (\SL_2(q)^2 \circ_2 \SL_2( q )).2).2 & 2 \\
16 & 1, 9 & ([q-\varepsilon]^2\circ_{2} \SL_2(q)^2).2 & 1  \\ \hline\hline
\end{array}
\end{equation}
Here, the same conventions regarding~$\varepsilon$ are used as in the table
displayed in~{\rm (\ref{TableIn415a})}. Moreover, the following conditions 
hold.
Let~$t$ be a $\sigma$-stable $G$-conjugate of~$s$. If the 
$G$-class type of~$t$ is $(12,1)$, then the $\mathbf{G}^{\sigma}$-class type 
of~$t$ is $(12,2)$ or $(12,4)$ and $f/m$ is even.
If the $G$-class type of~$t$ is $(16,1)$, the $\mathbf{G}^{\sigma}$-class type 
of~$t$ is $(16,k)$ for $k \in \{ 3, 4, 7, 10 \}$ and $f/m$ is even, or 
the $\mathbf{G}^{\sigma}$-class type of~$t$ is $(16,8)$ and $4 \mid f/m$.
If the $G$-class type of~$t$ is $(16,9)$, the $\mathbf{G}^{\sigma}$-class type 
of~$t$ is $(16,8)$ and $f/m \equiv 2\,(\text{\rm mod\ } 4)$. 

{\rm (c)} Let
$b \subseteq \mathcal{E}_\ell( G, s )$ be a $\sigma$-stable $\ell$-block.
Then $\sigma$ stabilizes $\cE( G, s )$. Moreover, the
permutation actions of~$\langle \sigma \rangle$ on~$\IBr(b)$ and
on $\Irr( b ) \cap \cE( G, s )$ are equivalent, unless $s = 1$ and $\ell = 3$.

{\rm (d)} Suppose that $s = 1$ and $\ell = 3$. If $p$ and $m'$ are as in 
Case~(i), then~$\sigma$ has exactly~$7$ orbits of length~$2$ on $\IBr(b)$ if~$b$ 
is the principal block, and exactly one orbit of length~$2$ on $\IBr(b)$ if~$b$ 
is the non-principal unipotent block of positive defect. The other unipotent
$3$-modular characters are fixed by~$\sigma$. If $p$ and $m'$ are as in 
Case~(ii), then~$\sigma$ fixes every element of $\IBr(b)$.
\end{prop}
\begin{prf}
Notice that~$\sigma$ is a Steinberg morphism of~$\mathbf{G}$, and thus in 
particular an isogeny. The equalities in~(\ref{FPrimeAndLusztigSeries}) follow 
from \cite[Corollary~$9.3$(ii)]{DiMiPara}. Now suppose that~$\sigma(s)$ is 
conjugate in~$G$ to $s$. By Lemma~\ref{SigmaStableConjugate} we may and will 
assume that $\sigma(s) = s$. Then $\mathcal{E}( G, s )$ and hence also 
$\mathcal{E}_\ell( G, s )$ are invariant under~$\sigma$ 
by~(\ref{FPrimeAndLusztigSeries}). 

(a) Suppose that~$p$ is even and that~$m'$ is odd. Since $\sigma(s) = s$, 
we have $C_{\mathbf{G}}( s ) = \sigma( C_{\mathbf{G}}(s) )$. We can also assume 
that $C_{\mathbf{G}}(s)$ is not a torus. 
As~$\sigma_1$ interchanges the long root subgroups
of~$\mathbf{G}$ with its short root subgroups, the tables in~\cite{LL} only
leave the possibilities $(1,1)$, $(4,k)$, $(14,k)$ or
$(15,k)$ for the $G$-class types of~$s$. As~$\sigma$ interchanges the
maximal tori of type $(20,12)$ and $(20,17)$, it also swaps the
classes of type $(14,2)$ and $(14,3)$ and the classes of type $(15,2)$
and~$(15,4)$. Let us now prove the remaining statements.

If~$s$
is of class type $(4,k)$, $k = 1, 2$, then $\ell > 3$, as $s$ is an element of
order~$3$.
Moreover,~$\sigma$ interchanges the two components of $C_G( s )$ of type~$A_2$,
as one is a long root subgroup, and the other one a short root subgroup.
Thus~$\sigma$ has three orbits of length~$2$ on the set of unipotent characters
of~$C_G(s)$ and fixes its other unipotent characters. The unipotent characters
of $C_G( s )$ are uniform functions, hence there is a bijection between
$\mathcal{E}( C_G ( s ), 1 )$ and $\mathcal{E}( G , s )$ which commutes
with~$\sigma$ by \cite[Corollary~$9.2$]{DiMiPara}. It follows that the number of
orbits of~$\sigma$ of length two on $\mathcal{E}( G , s )$ equals~$3$ and
that~$\sigma$ fixes the other elements of $\cE( G, s )$. If~$s$ is of class
type $(14,k)$, $k = 1, 4$, then~$\sigma$ swaps the two components of $C_G( s )$
of type~$A_1$, as one is a short root subgroup and the other one a long root
subgroup. Thus~$\sigma$ interchanges the two unipotent characters of $C_G( s )$
of degree~$q$ and fixes the other unipotent characters. If~$s$ is of class type
$(15,k)$, $k = 1, 3, 5$, the semisimple component of $C_G( s )$ is of type
$C_2$, and $\sigma$ induces the exceptional
graph automorphism of $\Sp_4(q)$ on $C_G( s )$. It follows from
\cite[Theorem~$2.5$(c)]{MalleExt}, that~$\sigma$ interchanges the two unipotent
principal series characters of $C_G( s )$ of the same degree and fixes the
other unipotent characters. In theses cases,
$C_{\mathbf{G}}( s )$ is a regular subgroup of~$\mathbf{G}$, and
the Lusztig induction map corresponding to~$C_{\mathbf{G}}( s )$ yields a
bijection between $\mathcal{E}( C_G ( s ), 1 )$ and $\mathcal{E}( G , s )$.
As this map
commutes with the action of~$\sigma$ (see \cite[Corollary~$9.2$]{DiMiPara}),
this gives the entries of the table displayed in~(\ref{TableIn415a})
in case $s \neq 1$. In case $s = 1$, the entry is determined in
\cite[Theorem~$2.5$(e)]{MalleExt}.

(b) Suppose that~$p$ is odd or that $p = 2$ and~$m'$ is even. 
Then~$\sigma = F_1^m$. Consider the statement:
\begin{equation}
\label{sigmaFixesOrdinaries}
\sigma \text{\ fixes every element of\ } \mathcal{E}(G,s).
\end{equation}
We are going to determine the cases for which~(\ref{sigmaFixesOrdinaries}) 
holds. First of all,~(\ref{sigmaFixesOrdinaries}) is true if $s = 1$ by 
\cite[Proposition~$6.6$]{DiMiPara}. In the following, we will need to compute
the degrees of the unipotent characters of the groups~$C_G( s )$. These are
easily determined using Jean Michel's extension of CHEVIE; see \cite{Michel}.
Suppose the unipotent characters of $C_{G}( s )$ have pairwise distinct
degrees. Then the same is true for the elements of $\cE(G,s)$ by the Jordan
decomposition of characters. In this case~(\ref{sigmaFixesOrdinaries}) trivially
holds. Considering the tables in~\cite{LL}, it remains to investigate the cases
where the $G$-class type of~$s$ is one of $(2,1)$, $(11,1)$, $(11,2)$, 
$(12,1)$, $(12,3)$, $(14,k)$, $1 \leq k \leq 4$, $(15,k)$, $1 \leq  k \leq 5$ or 
$(16,k)$ with $k \in \{ 1, 2, 5, 6, 9 \}$. If~$s$ is of $G$-class type $(2,1)$,
the elements of  $\cE( G, s )$  are distinguished by their degrees, except for
two characters $\chi, \chi' \in \cE( G, s )$ of equal degree, which correspond, 
via the Jordan decomposition, to the unipotent characters of~$C_G( s )$ labeled 
by the bipartitions $(21^2, - )$ and $( -, 31 )$, respectively. There is a
$\sigma$-stable split Levi subgroup~$\mathbf{M}$ of~$\mathbf{G}$ such that~$s$
is contained in a $\sigma$-stable dual $\mathbf{M}^* \leq \mathbf{G}$ (in fact
we may take $\mathbf{M} =_{\mathbf{G}^\sigma} \mathbf{L}_{\{2,3,4\}}$ and 
$\mathbf{M}^* =_{\mathbf{G}^\sigma} \mathbf{L}_{\{1, 2, 3\}}$ in the notation
of Subsection~\ref{SetupF4}), and a 
$\sigma$-stable element $\psi \in \cE( M , s )$ such that the Harish-Chandra 
induced character $R_{\mathbf{M}}^{\mathbf{G}}( \psi )$ contains $\chi$ with 
multiplicity~$1$ and~$\chi'$ with multiplicity~$0$. This follows from the
compatibility of Harish-Chandra induction and the Jordan decomposition
(see \cite[Corollary~$9.2$]{DiMiPara})
with a computation of unipotent characters in $C_G( s ) = \Spin_9( q )$.
Hence~$\sigma$ fixes~$\chi$ and so also~$\chi'$.

In the remaining cases, put 
$\mathbf{M}^* := C_{\mathbf{G}}( Z( C_{\mathbf{G}}( s ) )^\circ )$.
Then~$\mathbf{M}^*$ is a $\sigma$-stable regular subgroup of~$\mathbf{G}$ and we
choose a $\sigma$-stable regular subgroup $\mathbf{M} \leq \mathbf{G}$ dual
to~$\mathbf{M}^*$. In these cases,~$\mathbf{M}$ and~$C_{\mathbf{G}}( s )$ are of 
classical type. Hence the elements of $\cE( C_G( s ), 1 )$ and of~$\cE( M, s )$ 
are uniquely determined by their multiplicities in the Deligne--Lusztig 
characters; see, e.g.\ \cite[Theorem~$15.8$]{CaEn}. The latter theorem, together 
with \cite[Corollary~$9.2$]{DiMiPara} then implies that there is a 
$\sigma$-equivariant bijection $\cE( C_G( s ), 1 ) \rightarrow \cE( M, s )$.
Also, Lusztig induction with respect to~$\mathbf{M}$ induces a bijection 
$\cE( M, s ) \rightarrow \cE( G, s )$, which is $\sigma$-equivariant by 
\cite[Corollary~$9.2$]{DiMiPara}. To show that~(\ref{sigmaFixesOrdinaries}) is
satisfied, it suffices therefore to show that $\sigma$ fixes every element 
of~$\cE( C_G( s ), 1 )$. If~$s$ is of $G$-class type~$(11,1)$,~$(11,2)$, 
or~$(15,k)$, $1 \leq k \leq 5$, then~$\sigma$ fixes the unipotent characters of 
$C_G( s )$ by \cite[Theorem~$2.5$(c)]{MalleExt}. If~$s$ is of $G$-class type 
$(14,k)$, $1 \leq k \leq 4$, then~$\sigma$ fixes every unipotent character 
of~$C_G( s )$, as~$\sigma$ stabilizes the two semisimple factors of type~$A_1$ 
of~$C_G( s )$, one being a short root subgroup, the other one a long root 
subgroup. 

Suppose that the $G$-class type of~$s$ is 
$(12,1)$ or $(12,3)$ and that~$\sigma$ does not fix $\cE( G, s )$ element-wise.
Then~$\sigma$ does not fix every element of $C_G( s )$, again
by \cite[Corollary~$9.2$]{DiMiPara}. This implies that $\sigma$ permutes the two 
long root components of $C_G( s )$ of type~$A_1$, which happens if the 
$\mathbf{G}^\sigma$-class type of~$s$ is one of $(12,2)$ or $(12,4)$.
But then the $G$-class type of~$s$ equals $(12,1)$ by 
Lemma~\ref{TypesOfMaximalRankSubgroupsUnderPowersOfSigma}.
We obtain two $\sigma$-orbits of length~$2$ on $\cE( G, s )$ in this case.
The closed subsystem $\Gamma$ of $\Sigma$ which gives rise to the centralizer
of an element of $\mathbf{G}$-class type~$16$ has stabilizer in~$W$ isomorphic
to $2 \times D_8$, as is easily computed with CHEVIE. This group has ten
conjugacy classes giving rise to the ten $G$-class types $(16,k)$, 
$1 \leq k \leq 10$. Using Lemma~\ref{TypesOfMaximalRankSubgroupsUnderPowersOfSigma}
and Table~\ref{cc}, one finds that the class types $(16,5)$ and $(16,8)$
correspond to the conjugacy classes of elements of order~$4$, whereas the
class type $(16,9)$ corresponds to the squares of these elements. This
information is enough to prove the statements for~$s$ of $G$-class type $(16,k)$.

(c) Suppose that $\ell \neq 3$ or that $s \neq 1$. Then Proposition~\ref{EllB} and 
\cite[Theorem~$5.1$]{GeHi1} imply that $\Irr( b ) \cap \cE( G , s)$ is a basic 
set for $\IBr( b )$. By Lemma~\ref{OrbitCount}, the numbers of fixed points 
of~$\sigma$ on $\IBr( b )$ and on $\Irr( b ) \cap \cE( G , s)$ are the same. 
The analogous statement holds for every power of~$\sigma$, which gives our
claim. 

(d) Suppose that $\ell = 3$ and $s = 1$. If we are in case~(ii),~$\sigma$ 
satisfies~(\ref{sigmaFixesOrdinaries}) by~(a), and Lemma~\ref{UnipotentsGenerate} 
implies our claim. 
Hence assume that we are in case~(i).
Let $M$ denote the $(37 \times 37)$-block diagonal Fourier transform matrix 
for the unipotent characters of~$G$; see \cite[Section~$13.6$]{C2}. 
Furthermore, let~$P$ denote the permutation matrix arising from the
permutation of $\sigma$ on the set of unipotent characters of~$G$
as given by~\cite[Theorem~$2.5$(e)]{MalleExt}. One then checks that
$MP = PM$, i.e.~$\sigma$ permutes the set of almost characters of~$G$ in the
same way as it permutes the unipotent characters. In particular,
the number of fixed points of~$\sigma$ on
the set of unipotent characters is the same as on the set of almost characters.
Consequently,~$\sigma$ has exactly~$8$ orbits of length~$2$ on the set of almost 
characters. 
The two almost characters corresponding to $F_4[\theta]$ and $F_4[\theta^2]$ 
are fixed by $\sigma$ and vanish on the $3$-regular classes of~$G$; see
\cite[$6.2.4$(c) and Proposition~$6.3$]{ShojiI} and \cite[p.~$42$]{Ge2019}.
If we denote by $\mathcal{U}$ the set of almost characters with the latter 
two characters removed, then Lemma~\ref{UnipotentsGenerate} implies that 
$\mathcal{U}$ satisfies the hypotheses of Lemma~\ref{OrbitCount} for $\ell = 3$ 
and the union of unipotent blocks.
Thus~$\sigma$ has~$8$ orbits of length~$2$ on the set of unipotent Brauer 
characters. We may apply the same lemma to the unique non-principal unipotent
$3$-block of positive defect, to see that~$\sigma$ has exactly one orbit of
length~$2$ on this block. As~$\sigma$ fixes the unipotent defect~$0$ 
characters, it follows that $\sigma$ has exactly~$7$ orbits of length~$2$ 
on each of the set of irreducible ordinary and irreducible Brauer characters 
of the principal block. This completes our proof.
\end{prf}

\begin{cor}
\label{OrbitOnPBEvenqLg3}
Suppose that $p = 2$, that $\ell > 3$ and that the Sylow $\ell$-subgroups 
of~$G$ are non-cyclic. Let~$b$ denote the principal $\ell$-block of~$G$.
Then the number of non-trivial orbits of~$\sigma_1$ on $\IBr(b)$ is as
given in the following table.
\setlength{\extrarowheight}{0.5ex}
\begin{equation}
\begin{array}{c||ccccc}\hline\hline
e & 1 & 2 & 3 & 4 & 6 \\ \hline
\text{\rm no.} & 7 & 7 & 6 & 0 & 6 \\ \hline\hline
\end{array}
\end{equation}
(See {\rm Subsection~\ref{GoodPrimesBlocks}} for the significance of~$e$.)
%In particular,~$\sigma_1$ fixes every element of~$\IBr(b)$ if $e = 4$.
\end{cor}
\begin{prf}
This follows from Proposition~\ref{ActionAutomorphisms}(a), in connection with
the action of~$\sigma_1$ on the set of unipotent characters of~$G$; see
\cite[Theorem~$2.5$(e)]{MalleExt}.
\end{prf}

\medskip

\noindent
We finally consider the case of the exceptional covering group of $F_4(2)$.
Observe that the automorphism $\sigma_1$ of $F_4(2)$ lifts to an automorphism
of its double cover, also denoted by $\sigma_1$.

\begin{prop}
\label{ActionAutomorphisms2F42}
Let $\hat{G} = 2.F_4(2)$ denote the exceptional double cover of $F_4(2)$
and let $\ell \in \{  3, 5, 7 \}$. Let $b$ be an $\ell$-block of~$\hat{G}$
of non-cyclic defect. Then $\sigma_1$ fixes~$b$ unless $\ell = 3$ and~$b$ is 
one of the two blocks of defect~$2$, which are swapped by~$\sigma_1$. 

If~$\sigma_1$ fixes~$b$, the non-trivial orbits of~$\sigma_1$ on $\IBr( b )$ 
have length~$2$, and the number of such orbits is as given in the table below.

\setlength{\extrarowheight}{0.5ex}
$$
\begin{array}{c|ccc} \hline\hline
\ell & 3 & 5 & 7 \\ \hline
\text{\rm no.} & 4 & 6 & 6 \\ \hline\hline
\end{array}
$$
\end{prop}
\begin{prf}
This follows by inspecting the character table of~$\hat{G}$ in the 
Atlas~\cite{Atlas} and the decomposition matrices in~\cite{HF42}, in connection 
with Lemma~\ref{OrbitCount}.
\end{prf}

%\markboth{Radical $3$-subgroups and defect groups of $F_4(q)$}
\section{Radical $3$-subgroups and defect groups of $F_4(q)$}\label{secRad3}

\noindent In this section, let $G = F_4(q)$ with $q = p^f$ and $p \neq 3$. We
keep the notation introduced in Subsection~\ref{TheThreeBlocks}. We let
$\varepsilon \in \{ \pm 1 \}$ be defined by $3 \mid q - \varepsilon$, and 
put $e = e_3(q)$, i.e.\ $e = 1$ if $\varepsilon = 1$, and $e = 2$, otherwise. 
Also, the positive integer~$a$ is defined by
$(q - \varepsilon)_3 = 3^a$. In addition, we put $d := \gcd( 2, q - 1)$. The 
radical $3$-subgroups of~$G$ are classified in \cite{AnDF4,AH2} and we state 
the main results here. We also derive some consequences needed later on. For 
our notation for groups and their extensions used below, see 
Subsection~\ref{Groups}.

\subsection{Elements of order~$3$ in $G$}\label{sec3elts} 

\noindent By \cite[Table 4.7.3]{Gor}, the group~$G$ contains exactly three 
conjugacy classes of elements of order~$3$, called $3\A$, $3\B$, and $3\C$, 
with representatives $z_\A$, $z_\B$, and $z_\C$, such that for each 
$X \in \{ \A,\B,\C \}$ the elements
$z_X$ and $z_X^{-1}$ are conjugate. 
Moreover, the local structure is as follows:
\begin{align*}
C_G(z_\A) & = \langle [q-\varepsilon] \circ_{d} \Sp_6(q), x_\A\rangle,& 
N_G(\langle z_\A\rangle) & = \langle C_G(z_\A) , \gamma_\A \rangle, \\
C_G(z_\B) & = \langle [q-\varepsilon] \circ_{d} \Spin_7(q), x_\B\rangle, & 
N_G(\langle z_\B\rangle)& = \langle C_G(z_\B) , \gamma_\B \rangle, \\
C_G(z_\C)&= \langle \SL_3^\varepsilon(q)\circ_3 \SL_3^\varepsilon(q),x_\C\rangle, & 
N_G(\langle z_\C\rangle)& = \langle C_G(z_\C) , \gamma_\C \rangle,
\end{align*}
where the actions of $\gamma_X$ and $x_X$ are determined by:
\begin{longtable}{rcp{9.0cm}}
$\gamma_X=(\iota\mycol 1)$:&& for $X\in \{\A,\B\}$, the action of $\gamma_X$ on 
$[q - \varepsilon]$ is by $x\mapsto x^{-1}$, and trivial on $O^{p'}(C_G(z_X))$;\\
$\gamma_\C=(\gamma\mycol\gamma)$:&& $\gamma_{\C}$ acts as 
the order~$2$ graph automorphism, i.e.\ inverse-transpose, on each 
factor $\SL^\varepsilon_3(q)$ of $O^{p'}(C_G(z_\C))$;\\
$x_X=(1 \mycol d)$: && for $X\in \{\A,\B\}$, the action of $x_X$ on 
$[q-\varepsilon]$ is trivial; on $O^{p'}(C_G(z_X))$ it is trivial, if $d = 1$, 
and by the diagonal automorphism of order~$2$ if $d = 2$;\\
$x_\C = (x_1 \mycol x_2)$:&& each $x_i$ acts on $\SL^\varepsilon_3(q)$ by 
the diagonal automorphism of order~$3$.
\end{longtable}
\noindent Notice that $C_G(z_\C)$ is conjugate in~$G$ to the group~$L$ introduced
in Proposition~\ref{C3C}. Note also that $x_\A$ does not necessarily have 
order~$2$; we only know that $x_\A^2\in [q-\varepsilon] \circ_{d} \Sp_6(q)$; 
similarly for $x_\B, x_\C$ and $\gamma_X$. For example, as discussed in
Remark~\ref{DiagonalActionOfx}, we  have 
$x_\C^3\in \SL_3^{\varepsilon}(q)\circ_3 \SL_3^{\varepsilon}(q)$, and 
there is a choice of $x_\C$ such that each $x_i$ acts on $\SL^\varepsilon_3(q)$ 
in the same way as ${\rm diag}(1,1,\xi^3) \in \GL_3^\varepsilon(q)$ acts on
$\SL^\varepsilon_3(q)$ by conjugation; here, $\xi \in \mathbb{F}^*$ is
an element of order $3^{a+1}$. Then $x_i^3$ acts as conjugation with 
${\rm diag}(1,1,\xi^9)$ and thus as the inner automorphism 
${\rm diag}(\xi^{-3},\xi^{-3},\xi^6)$, $i = 1, 2$. 

\subsection{Radical 3-subgroups of $\SL_3^{\varepsilon}(q)$}\label{secradLe} 
We first recall the classification of the radical $3$-subgroups of $\Le := 
\SL^\varepsilon_3(q)$ and 
fix notation. We view~$\Le$ as a subgroup of $\GL_3^{\varepsilon}( q )$, where, 
in case of $\varepsilon = -1$, the latter group is defined with respect to the 
standard diagonal hermitian form.  In particular, $\GL_3^{\varepsilon}( q )$ 
and $\Le$ are invariant under the inverse-transpose automorphism.

In order to simplify notation, we suppress the subscript~$\varepsilon$ from the
objects associated to~$\Le$, although these may depend on~$\varepsilon$. In 
particular, we put $L := \Le = \SL_3^{\varepsilon}( q )$. We begin by describing 
various elements and subgroups of~$L$. Let $\xi \in \mathbb{F}^*$ have 
order~$3^{a+1}$, and put $\zeta := \xi^{3^a}$. Next, let $y := 
\diag(1, \zeta, \zeta^{-1} ) \in L$. Further, let $c, t \in L$ be
defined by 
$$
c := 
\left( \begin{array}{ccc} 0 & 1 & 0 \\ 0 & 0 & 1 \\ 1 & 0 & 0 \end{array} \right)
\quad\quad 
t := 
\left( \begin{array}{rrr} 0 & -1 & 0 \\ -1 & 0 & 0 \\ 0 & 0 & -1 \end{array} \right).
$$
Finally, let~$x$ denote the diagonal automorphism of~$L$ which acts on~$L$ 
in the same way as $\diag(1,1,\xi^3)$. Notice that~$x$ centralizes~$t$ and 
that~$x^3$ is an inner automorphism of~$L$.

Put $H := C_{L}( \diag(1,1,\zeta) ) \cong \GL_2^{\varepsilon}( q )$, let~$Z$ 
denote the center of~$L$, and~$T$ the subgroup of diagonal matrices in~$L$, 
i.e.\ $T \cong [q - \varepsilon]^2$. We also let $D := 
\langle O_3( T ), c \rangle$; then $D \in \Syl_3( L )$. Finally, we put 
$K := \langle y, c \rangle$. Notice that if $a = 1$, then $K = D$ and 
$K^x = K$. 

The next lemma merges \cite[Lemma 4.4]{AnDF4} and \cite[Lemma 3.1(2)]{AH2}. 

\addtocounter{thm}{2}
\begin{lem}\label{lemRadLi}
Let the notation be as above. Then representatives for $\R_3(L)/L$, the set of 
$L$-conjugacy classes of radical $3$-subgroups of~$L$ is given as follows.

$$\begin{array}{cc}\hline
\text{\rm Conditions} & \text{\rm Representatives for\ } \R_3(L)/L 
\rule[- 7pt]{0pt}{ 20pt} \\ \hline
q = 2 & D  \rule[ 0pt]{0pt}{ 13pt} \\
q = 4 & Z,\; D \\
q = 8 & Z,\; O_3(H),\; K,\; K^x,\; K^{x^2},\; D \\
q \geq 5, a = 1 & Z,\;  O_3(T),\; D \\
q \geq 17, a \geq 2 &  Z,\; O_3(H),\; K,\; K^x,\; K^{x^2},\; O_3(T),\; D
\rule[- 7pt]{0pt}{ 5pt} \\ \hline
\end{array}
$$

\smallskip

\noindent
The structure of $C_{L}( R )$, $N_{L}( R )$ and $\Out_{L}( R )$
for the radical $3$-subgroups~$R$ of $L$ as above is given in the following 
table.
$$\begin{array}{ccccc}\hline
R & \text{\rm Struc.} & C_{L}( R ) & N_{L}( R ) & \Out_{L}( R ) 
\rule[- 7pt]{0pt}{ 20pt} \\ \hline
Z & 3 & L & L & 1 \rule[ 0pt]{0pt}{ 13pt} \\
O_3(H) & 3^a & H & H & 1 \\
K, a \geq 2 &  3_+^{1+2} & Z & K.\SL_2(3) & \SL_2(3) \\
D, a = 1 & 3_+^{1+2} & Z & D.Q_8 & Q_8 \\
D, a \geq 2 & [3^a]^2.3 & Z & D.2 = \langle D, t \rangle & 2 \\
O_3( T ) & [3^a]^2 & T & T.S_3 
= \langle T, c, t \rangle & S_3 \rule[- 7pt]{0pt}{ 5pt} \\ \hline
\end{array}
$$
\end{lem}

\addtocounter{subsection}{1}
\subsection{Radical 3-subgroups of $F_4(q)$}\label{secradF4} 
Let $L^1, L^2 \leq C_G(z_\C)$ denote the subgroups with 
$L^1 \cong \SL_3^\varepsilon( q ) \cong L^2$, such that 
$C_G(z_\C) = \langle L^1 \circ_3 L^2, x_\C\rangle$ and 
$N_G( \langle z_\C \rangle) = \langle C_G( z_\C ), \gamma_\C\rangle$, where 
$x_\C = (x_1{:} x_2)$ and $\gamma_\C = (\gamma_1{:}\gamma_2)$ act as described 
%%%%%%%%%%%%%%%%%%%%%%%%%%%%%%%%%%%%%%%%%%%%%%%%%%%%%%%%%%%%%%%%%%%%%%%%%%%%%%%%
%%%%%%%%%%%%%%%%%%%%%%%%%%%%%%%%%%%%%%%%%%%%%%%%%%%%%%%%%%%%%%%%%%%%%%%%%%%%%%%%
%%
%% For this notational convention see Sheet 1 of 13.03.2018 (GH)
%%
%%%%%%%%%%%%%%%%%%%%%%%%%%%%%%%%%%%%%%%%%%%%%%%%%%%%%%%%%%%%%%%%%%%%%%%%%%%%%%%%
%%%%%%%%%%%%%%%%%%%%%%%%%%%%%%%%%%%%%%%%%%%%%%%%%%%%%%%%%%%%%%%%%%%%%%%%%%%%%%%%
in Subsection~\ref{sec3elts}. Notice that $L^1$ and $L^2$ have the same meaning
as in Proposition~\ref{C3C}. In particular,~$L^1$ is of
Dynkin type $A_2^\varepsilon( q )$, and $L^2$ of Dynkin type 
$\tilde{A}_2^\varepsilon( q )$ so that~$L^1$ is a long root subgroup and~$L^2$ 
is a short root subgroup. Notice that $L^1 \cap L^2 = Z( C_G(z_\C) )$.

The radical $3$-subgroups of~$G$ are classified in~\cite{AH2} and~\cite{AnDF4}.

\addtocounter{thm}{1}
\begin{thm}\label{thmradicalf4} 
The non-trivial radical $3$-subgroups of~$G$, up to 
conjugacy, are as given in {\rm Table~\ref{tabradicalageq2}}. \hfill{$\Box$}
\end{thm}

\medskip
\noindent
It turns out that the Sylow $3$-subgroups of the centralizers of semisimple 
elements of $G$ are radical $3$-subgroups.

\begin{prop}\label{SylowStructure}
Let $\mathbf{M} = \C_{\mathbf{G}}( s )$ denote the centralizer of the 
semisimple element $s \in G$, and let~$P$ be a Sylow $3$-subgroup of~$M$.
Then~$P$ is a radical $3$-subgroup of~$G$, and the conjugacy 
class of~$P$ in $\cR_3(G)$ is as given in the fourth column of
{\rm Table~\ref{SylowsInCentralizersI}}.
\end{prop}
\begin{prf}
We enumerate the possibilities for (the $G$-conjugacy classes) of~$\mathbf{M}$
by the pairs $(i,k)$ according to the tables in~\cite{LL}. 
Assume first that $i = 20$, i.e.\ that~$\mathbf{M}$ is a maximal torus 
of~$\mathbf{G}$. Using Table~\ref{cc} for $i \in \{ 4, 6, 10 \}$, we find that one of 
the following cases occurs:
\begin{itemize}
\item[(a)] $P = 1$,
\item[(b)] $\mathbf{M} \leq C_{\mathbf{G}}( z )$ for some element~$z$ in
the conjugacy class~$3\C$, 
\item[(c)] $P = [3^a]$ and $P \leq Z( C_{\mathbf{G}}( z ) )$ 
for some~$z$ in class~$3\A$ or~$3\B$, or
\item[(d)] $P = [3^a]^2$ and $P \leq Z( C )$, where $\mathbf{C} \leq \mathbf{G}$
is the centralizer of a semisimple element of~$G$ of class type $(15,1)$ if 
$\varepsilon = 1$, and of class type $(15,3)$ if $\varepsilon = -1$.
\end{itemize}
By the construction of the radical $3$-subgroups of~$\SL_3^\varepsilon( q )$ in
Lemma~\ref{lemRadLi}, we find that~$P$ is a radical $3$-subgroup of~$G$ in 
Case~(b). Unless $k \in \{ 14, 15, 19, 20 \}$, the conjugacy class of~$P$ in 
this case is determined from~$|P|$ and the types of the characteristics (see 
Subsection~\ref{secradF4Appendix} for this notion) 
of the abelian radical $3$-subgroups of~$G$,
in connection with the information given in Table~\ref{cc} on
the number of elements in $Z(M)$ belonging to other class
types. For example, if~$M$ is of class type $(20,12)$ and $\varepsilon = 1$, 
then $|P| = 3^{3a}$ and $Z(M)$ contains~$8$ elements of class type $(4,1)$
and~$12$ elements of class type $(6,1)$, as well as~$6$ elements of class
type $(10,1)$. Thus $P =_G R_{16}$. To identify the conjugacy class of~$P$
in $\cR_3(G)$ in the remaining cases, suppose first that $\varepsilon = 1$.
We may then choose notation such that $P =_G R_6$ if $\mathbf{M} =_G 
\mathbf{M}_{20,14}$, and $P =_G R_7$ if $\mathbf{M} =_G \mathbf{M}_{20,19}$. By 
looking at the two classes of maximal tori in each of $\mathbf{M}_{18,3}$ and 
$\mathbf{M}_{19,8}$, which occur as the centralizers of~$P$ in the respective 
cases, we find that $R_6 \leq_G R_{11}$ and $R_7 \leq_G R_{12}$. Now in 
case $\varepsilon = -1$, the two classes of maximal tori in each of
$\mathbf{M}_{18,7}$ and $\mathbf{M}_{19,9}$ show that $P =_G R_6$
if $\mathbf{M} =_G \mathbf{M}_{20,15}$, and $P =_G R_7$ if $\mathbf{M} =_G
\mathbf{M}_{20,20}$.
In Case~(c) we have $C_G( P ) = C_G ( z )$ 
is conjugate to the centralizer of~$R_2$, respectively~$R_3$, showing that~$P$ 
is conjugate to~$R_2$ respectively~$R_3$. In Case~(d), we have $C = C_G( R )$ 
for some radical subgroup $R =_G R_{10}$, and thus $P = O_3( Z( C ) ) 
=_G R_{10}$. 

Now assume that $i \neq 20$. Here, we use the following general argument.
Suppose that $z \in G$ is semisimple and $s \in Z( C_G( z ) )$. Then 
$C_G( z ) \leq C_G( s ) = M$. If, in addition, $|C_G(s)|_3 = |C_G(z)|_3$,
the two centralizers have conjugate Sylow $3$-subgroups. Applying this
argument in case $C_{\mathbf{G}}( z )$ is a maximal torus, yields all the
entries of Table~\ref{SylowsInCentralizersI} in the cases where $P$ is abelian.
We give a sample argument for one of the other cases. Suppose that~$z$ is of 
class type $(13,k)$, where $1 \leq k \leq 6$. We may then assume 
that~$C_{\mathbf{G}}(z)$ equals the regular 
subgroup $\mathbf{M}_{13,k}$ of~$\mathbf{G}$ of type~$A_2$ corresponding to the 
long roots; see Subsections~\ref{ClassTypes} and~\ref{ConstructionOfCentralizers}.
Suppose further that $\varepsilon$ is such that the
Sylow $3$-subgroups of $M_{13,k}$ are non-abelian. Then the center of~$M_{13,k}$ 
contains elements of class type $(4,1)$, if $\varepsilon = 1$, and of class type 
$(4,2)$, if $\varepsilon = -1$, i.e.\ $3\C$-elements. Thus $M_{13,k} \leq 
C_G( z_\C ) = L$. 
If $k = 1$ or $k = 6$, according as $\varepsilon = 1$ or $\varepsilon = -1$, 
respectively, we have $C_G( s ) = (L^1 \circ_3 [q - \varepsilon]^2).3$. Thus a 
Sylow $3$-subgroup of~$C_G( z_\C )$ equals $(D_1 \circ_3 [3^a]^2).3 =_G R_{34}$. 
Now in the respective cases for~$\varepsilon$, we find that $Z(M_{13,k})$ contains 
elements~$s$ of class types $(8,1), (8,4)$ and $(i,k)$ for
$i = 3, 6, 7, 8$ and $k = 1, 2$, as well as $(2,1)$. Since $|M_{13,k}|_3 = 
|C_G( s )|_3$ in all these cases, we obtain the corresponding entries 
in Table~\ref{SylowsInCentralizersI}. (This argument also works in case $q = 4$, 
where there is no element $z \in G$ such that $C_{\mathbf{G}}( z ) =_G
\mathbf{M}_{13,1}$.) The remaining cases are treated similarly.
\end{prf}
%%%%%%%%%%%%%%%%%%%%%%%%%%%%%%%%%%%%%%%%%%%%%%%%%%%%%%%%%%%%%%%%%%%%%%%%%%%%%%%%
%%%%%%%%%%%%%%%%%%%%%%%%%%%%%%%%%%%%%%%%%%%%%%%%%%%%%%%%%%%%%%%%%%%%%%%%%%%%%%%%
%%
%% Checked this proof once more on 29.07.2020 (GRH)
%%
%%%%%%%%%%%%%%%%%%%%%%%%%%%%%%%%%%%%%%%%%%%%%%%%%%%%%%%%%%%%%%%%%%%%%%%%%%%%%%%%
%%%%%%%%%%%%%%%%%%%%%%%%%%%%%%%%%%%%%%%%%%%%%%%%%%%%%%%%%%%%%%%%%%%%%%%%%%%%%%%%

\medskip

\noindent
We collect a few further results on the radical subgroups of~$G$.
The following lemma describes the regular subgroups of~$\mathbf{G}$ that occur
as centralizers of abelian radical $3$-subgroups (not conjugate to one 
of~$R_1$,~$R_8$ or~$R_{15}$).

\begin{lem}
\label{CentralizersAbelianRadicalSubgroups}
Let~$R$ be a non-trivial, abelian radical $3$-subgroup of~$G$, but $R \not\in_G 
\{ R_1, R_8, R_{15} \}$. Then $C_{\mathbf{G}}( R ) =_G \mathbf{M}_{i,k}$, with 
$(i,k)$ as in the following table. The two cases for $k$ correspond to 
$\varepsilon = 1$ and $\varepsilon = - 1$, respectively.
\setlength{\extrarowheight}{0.5ex}
$$\begin{array}{cccl}\\ \hline\hline
R & i & k & \multicolumn{1}{c}{\text{\rm Rem}} \\ \hline\hline
R_2 & 10 & 1, 2 & \\
R_3 & 6 & 1, 2 & \\
R_4 & 9 & 1, 2 & a \geq 2 \\
R_5 & 7 & 1, 2 & a \geq 2 \\
R_6 & 18 & 3, 7 & \\
R_7 & 19 & 8, 9 &  \\
R_9 & 14 & 1, 4 &  \\
R_{10} & 15 & 1, 3 &  \\ \hline\hline
\end{array}\quad\quad
\begin{array}{ccc}\\ \hline\hline
R & i & k \\ \hline\hline
R_{11} & 20 & 7, 8   \\
R_{12} & 20 & 4, 5   \\
R_{13} & 17 & 1, 6   \\
R_{14} & 13 & 1, 6   \\
R_{16} & 18 & 1, 10   \\
R_{17} & 19 & 1, 10   \\
R_{18} & 20 & 1, 2   \\ \hline\hline
\end{array}$$
Notice that $C_{\mathbf{G}}( R ) \leq \mathbf{G}$ is regular, and, unless 
$R \in_G \{ R_6, R_7, R_{11}, R_{12} \}$, also $e$-split. Moreover, $R$ 
is the Sylow $3$-subgroup of $Z( C_G( R ) )$.
\end{lem}
\begin{prf}
Let $(i, k)$ be as in the table corresponding to the row containing~$R$. Then, 
by Proposition~\ref{SylowStructure}, unless $R \not\in_G 
\{ R_3, R_4, R_{13}, R_{14} \}$, there is a maximal torus~$\mathbf{S} \leq 
\mathbf{M}_{i,k}$ such that $R \in \Syl_3( S )$, and $|S|_3 = |Z( M_{i,k})|_3$.
Comparing with $|C_G( R )|$, as given in Table~\ref{tabradicalageq2}, we find 
$C_G( R ) = M_{i,k}$ as claimed. 

Now, suppose that $\varepsilon = 1$ and $(i,k) = (17,1)$. Then 
$\mathbf{M}_{17,1}$ embeds into $\mathbf{L} = C_{\mathbf{G}}( z_\C )$ such that 
$[\mathbf{M}_{17,1},\mathbf{M}_{17,1}] = \mathbf{L}^2$. Moreover, 
$Z( \mathbf{M}_{17,1} )$ is the intersection of maximal tori $G$-conjugate to
$\mathbf{M}_{20,7}$ and $\mathbf{M}_{20,17}$, whose Sylow $3$-subgroups 
are conjugate to~$R_{11}$ and $R_{17}$, respectively; see Tables~\ref{cc} 
and~\ref{SylowsInCentralizersI}. On the other hand, by the 
construction of the radical $3$-subgroups (see Subsection~\ref{secradF4}), the 
intersection of a suitable conjugate of $R_{11}$ and~$R_{17}$ is a radical 
$3$-subgroup of~$G$ conjugate to~$R_{13}$, as can be seen by the type of its 
characteristic; see Subsection~\ref{secradF4Appendix}. The analogous 
argument works for $R =_G R_{14}$ and also for $\varepsilon = - 1$.
Finally suppose that $\varepsilon = 1$, $R =_G R_4$ and $(i,k) = (9,1)$. 
Again, $\mathbf{M}_{17,1}$ embeds into $\mathbf{L}$. Also, 
$Z( \mathbf{M}_{17,1} )$ is contained in the torus of type $(20,14)$,
whose Sylow $3$-subgroup is conjugate to~$R_6$. If $a \geq 2$,
then the Sylow $3$-subgroup of $Z( M_{17,1} )$ is a radical $3$-subgroup
of~$G$, which we denote by~$R_4$. By construction $R_4 \leq_G R_6$.
The remaining cases are proved similarly.

As $\mathbf{M}_{i,k}$ is a regular subgroup of~$\mathbf{G}$ in all cases, and
$\mathbf{M}_{i,k}$ is $e$-split except for $R =_G R_6, R_7$, the penultimate 
assertion follows. Inspection gives the last claim.
\end{prf}

\medskip

\noindent 
The centralizers of most of the non-abelian radical $3$-subgroups of~$G$ are
connected.
\begin{lem}
\label{RegularCentralizers}
Let~$R$ be a non-abelian radical $3$-subgroup of~$G$, but
$$R \not\in_G \{ R_{21}, R_{22}, R_{35}\text{--}R_{38} \}.$$
Then $C_{\mathbf{G}}( R )$ is a closed, connected reductive subgroup 
of~$\mathbf{G}$.
\end{lem}
\begin{prf}
By the construction of~$R$ as a subgroup of $L = (L^1 \circ_3 L^2).3$, we find 
that $C_{\mathbf{G}}( R ) \leq \mathbf{L}$, i.e.\ $C_{\mathbf{G}}( R ) = 
C_{\mathbf{L}}( R )$. As indicated in Table~\ref{tabradicalageq2}, the latter
is equal to a regular subgroup of $\mathbf{L}^i$, for $i \in \{ 1, 2 \}$.
\end{prf}

\medskip
\noindent
It turns out that the abelian, characteristic subgroups of certain non-abelian
radical $3$-subgroups exhibited in the next lemma play a crucial role in the
following.
\begin{lem}
\label{MaximalAbelianSubgroups}
Let $R \in_G \{ R_{29}\text{--}R_{34} \}$. 
Then~$R$ has a unique maximal abelian normal subgroup~$Q$ which is again a radical
subgroup of~$G$. The conjugacy class of~$Q$ is as given in the following table.
\setlength{\extrarowheight}{0.5ex}
$$
\begin{array}{l|cccccc} \hline\hline
R & R_{29} & R_{30} & R_{31} & R_{32} & R_{33} & R_{34} \\ \hline
Q & R_{12} & R_{11} & R_{16} & R_{17} & R_{18} & R_{18} \\ \hline\hline
\end{array}
$$
In particular, the conjugacy class of~$R$ in~$\cR_3(G)/G$ is determined by the 
conjugacy class of~$Q$, unless $R \in_G \{ R_{33}, R_{34} \}$.
\end{lem}
\begin{prf}
Realize~$R$ as a Sylow $3$-subgroup of a suitable $M_{i,k}$ according to 
Table~\ref{SylowsInCentralizersI}. Using~Table~\ref{SylowsInCentralizersI}, 
choose a maximal torus of~$M_{i,k}$, whose Sylow $3$-subgroup is a radical 
$3$-subgroups~$Q$ of the class indicated in the table of the lemma. Now 
$|R\colon\!Q| = 3$, and thus $[R,R] \leq Q \unlhd R$. From the type of the 
characteristic of~$R$ (recall that the characteristic of~$R$ equals 
$\Omega_1( [R,R] )$ in these cases) given in Table~\ref{tabradicalageq2}, 
we find that $[R,R] \not\leq Z( R )$. Hence $Q = C_R( [R,R] )$. Thus~$Q$ is
characteristic in~$R$.

Let~$A$ be any maximal abelian normal subgroup of~$R$ and suppose that~$A \neq Q$.
Then $R = AQ$ and thus $R/A \cong Q/(A \cap Q)$ is abelian. Thus $[R,R] \leq A$
and hence $A \leq C_R( [R,R] ) = Q$, a contradiction.
\end{prf}

\addtocounter{subsection}{5}
\subsection{Duality of radical $3$-subgroups} 
\label{DualRadicalSubgroups}
There is an involutive, inclusion preserving bijection $R \mapsto R^\dagger$ on 
the set $\cR_3(G)/G$ such that $C_{\mathbf{G}}( R^\dagger ) \cong 
C_{\mathbf{G}}( R )^*$, unless 
$R \in_G \{ R_{19}\text{--}R_{22}, R_{27}, R_{28}, R_{35}\text{--}R_{38} \}$. 
In the latter cases, we put $R^\dagger := R$.
The groups $R^\dagger$ are given in the last column of Table~\ref{tabradicalageq2}.

\subsection{Radical subgroups and Brauer pairs} In this subsection, we record a
couple of preliminary results. Let $(R,b_R)$ denote a centric Brauer pair.
Recall the notation $\cW( R, b_R )$ introduced in~(\ref{NumberOfWeights0}).
If $\cW( R, b_R ) \neq 0$, then~$R$ is a radical $3$-subgroup 
of~$G$ and the canonical character~$\theta_R$ of~$b_R$ is of $3$-defect 
zero, viewed as a character of $C_G(R)/Z(R)$; see the results summarized in
Subsection~\ref{BlocksAndWeights}.

We begin with the following result on irreducible characters
of $\GL^\varepsilon_2(q)$ and~$\SL^\varepsilon_3(q)$.

\addtocounter{thm}{2}
\begin{lem}
\label{GL2andSL3Characters}
Let~$H$ be one of~$\GL^\varepsilon_2(q)$ or~$\SL^\varepsilon_3(q)$. Let
$\psi \in \Irr(H)$ be such that $O_3(Z(H))$ is in the kernel of~$\psi$ and
that~$\psi$ is of $3$-defect zero as character of $H/O_3(Z(H))$. Then one of
the following occurs.

{\rm (a)} We have $H = \GL^\varepsilon_2(q)$ and~$\psi$ is Lusztig induced from
a character in general position of the Coxeter torus~$[q^2-1]$ of~$H$.

{\rm (b)} We have $H =  \SL^\varepsilon_3(q)$ and~$\psi$ is Lusztig induced
from a character in general position of the Coxeter torus
$[q^2 + \varepsilon q + 1]$ of~$H$.

{\rm (c)} We have $a = 1$ and $H =  \SL^\varepsilon_3(q)$ and~$\psi$ is Lusztig
induced from a character in general position of the torus $[q^2 - 1]$ of~$H$.

Moreover, in cases~{\rm (b)} and~{\rm (c)}, the character~$\psi$, viewed as a 
character of $H/O_3(Z(H)) = \PSL^\varepsilon_3(q)$, extends to a character 
of~$\PGL_3^\varepsilon(q)$.
\end{lem}
\begin{prf}
This follows from the known character tables of~$H$ and Deligne--Lusztig theory.
For the character degrees of~$H$ one may consult the tables in~\cite{LLCD}.
The last statement follows from Clifford theory.
\end{prf}

\medskip

\noindent
We now investigate the existence of centric Brauer pairs $(R,b_R)$ for 
$R \in_G \{ R_1, R_8, R_{15} \}$.

\begin{lem}
\label{ExceptionalAbelianRadicalSubgroups}
Let $R \in_G \{ R_1, R_8, R_{15} \}$ and let $(R, b_R)$ be a centric Brauer 
pair. Then $R \neq_G R_1$ and one of the following holds.

{\rm (a)} We have $R =_G R_8$ and $\Out_G( R, b_R ) = 1$. In this case $(R,b_R)$ is 
a maximal $b$-Brauer pair for a block~$b$ of~$G$ with $b = \cE_3( G, s )$ for 
a semisimple $3'$-element~$s$ of class type $(20,9)$ or $(20,10)$, according as 
$\varepsilon = 1$ or $\varepsilon = -1$, respectively. In particular, 
$R^\dagger =_G R \leq_G C_G( s )$.

{\rm (b)} We have $R =_G R_8$ and $|\Out_G( R, b_R )| = 3$, so that 
$\cW(R,b_R) = 0$.

{\rm (c)} We have $R =_G R_{15}$, and $b_R$ is the principal block 
of~$C_G( R )$. In particular, $\Out_G( R, b_R ) = \SL_2( 3 )$ and $\cW(R,b_R) = 1$.
Moreover,~$R$ is not a defect group of any $3$-block of~$G$.
\end{lem}
\begin{prf}
Assume that $R =_G R_{1}$. Then $C_G( R ) = (L^1 \circ_3 L^2).3$ with
$L^i \cong \SL_3^\varepsilon( q )$, $i = 1, 2$. Put $L' := L^1 \circ_3 L^2 \leq
C_G( R )$. Then $R \leq L'$. The canonical character~$\theta_R$ of~$b_R$ is a 
$3$-defect zero character of $C_G( R )/R$. Hence the restriction of $\theta_R$ 
to~$L'/R$ splits into a sum of three $3$-defect zero characters of $L'/R \cong 
\PSL^\varepsilon_3( q ) \times \PSL^\varepsilon_3( q )$. On the other hand,
a $3$-defect zero character of $\PSL^\varepsilon_3( q )$ extends to an
irreducible character of $\PGL^\varepsilon_3( q )$; see the last statement of
Lemma~\ref{GL2andSL3Characters}. As an outer element of
order~$3$ of $(L^1 \circ_3 L^2).3$ induces a diagonal automorphism on 
each factor~$L^i$, $i = 1, 2$ (see Subsection~\ref{sec3elts}), every $3$-defect 
zero character of~$L'/R$ extends to~$C_G( R )/R$, a contradiction.

Next, assume that $R =_G R_8$. Then, by Tables~\ref{tabradicalageq2} 
and~\ref{SylowsInCentralizersI}, we find that $C_G(R) = S.3$, where 
$S = \mathbf{S}^F$ is a maximal torus of~$G$ of type $(20,9)$ or $(20,10)$, 
according as $\varepsilon = 1$ or~$-1$, respectively. We also have $N_G( R ) = 
N_G( S ) = S.W_\mathbf{G}( \mathbf{S} )^F $ with 
$W_\mathbf{G}( \mathbf{S} )^F = 3.\SL_2(3)$, and there is an order preserving
$W_\mathbf{G}( \mathbf{S} )^F$-equivariant bijection between~$S$ and 
$\Irr(S)$, as $\mathbf{S} \cong \mathbf{S}^*$; see 
Subsections~\ref{Twisting},~\ref{TwistingAndDuality} and~\ref{Duality}. 
Let $\theta_R$ denote the 
canonical character of~$b_R$. Then $\theta_R(1) = 3$, and thus
$\Res^{S.3}_S( \theta_R )$ is the sum of three irreducible characters of 
$3'$-order. Let~$\xi$ be one of these constituents. If~$\xi$ is in general 
position, i.e.\ $\xi$ lies in a regular orbit of $W_\mathbf{G}( \mathbf{S} )^F$ 
on $\Irr(S)$, then $\Out_G(R, b_R) = 1$, and thus $(R,b_R)$ is a 
maximal Brauer pair for some $3$-block of~$G$. On the other hand, each regular
orbit of $W_\mathbf{G}( \mathbf{S} )^F$ on~$S$ yields a $G$-conjugacy class of
$3'$-elements of class type $(20,9)$ or $(20,10)$, respectively. Let $s \in S$
be a $3'$-element which lies in a regular $W_\mathbf{G}( \mathbf{S} )^F$-orbit,
and let~$D$ be a defect group of the block $b := \cE_3( G, s )$. We claim that
$D =_G R$. If not, $a = 1$ and $D \cong 3^2$ is conjugate to one of
$R_9$,~$R_{10}$,~$R_{13}$ or~$R_{14}$. In all these cases, $C_{\mathbf{G}}( D )$
is a regular subgroup of~$\mathbf{G}$ by
Lemma~\ref{CentralizersAbelianRadicalSubgroups}. Then, up to conjugation in~$G$,
we have $s \in C_G(D)^* = C_G( D^\dagger )$ by
Lemma~\ref{ConjugateSemisimpleLabels}. Hence $D^\dagger \leq C_G( s ) \cong
[q + \varepsilon q + 1]^2$, and thus $D^\dagger \in \Syl_3( C_G( s ) )$. This
contradicts Table~\ref{SylowsInCentralizersI}.

The number of centric Brauer pairs $(R,b_R)$ with $\Out_G(R,b_R) 
= 1$ equals the number of regular orbits of $W_\mathbf{G}( \mathbf{S} )^F$ on
the set of $3'$-elements of $\Irr(S)$. The latter number is also equal to the
number of blocks of the form $\cE_3( G, s )$, where $s \in S$ is a $3'$-element
of class type $(20,9)$, respectively $(20,10)$. As all these blocks have~$R$ as
defect group, we have proved (a). 

Now assume that~$\xi$ does not lie in a regular 
$W_\mathbf{G}( \mathbf{S} )^F$-orbit. Then $\Out_G( R, b_R ) \neq 1$ and, under the 
above bijection,~$\xi$ 
corresponds to a semisimple $3'$-element $t \in S$ such that~$t$ lies in a class 
type different from $(20,9)$, respectively~$(20,10)$. The tables in~\cite{LL} 
show that~$t$ is of class type $(13,5)$ or~$(17,5)$, or $(13,4)$ or $(17,4)$, 
respectively. Notice that $N_G( R )$ stabilizes $\mathbf{S} = 
{C_{\mathbf{G}}( R )^\circ}$ and hence~$S$. Every element of $N_G( R , b_R)$ 
fixes~$\theta_R$ and thus the orbit of~$t$ under $C_G(R) = S.3$. In particular, 
an element of $N_G( R , b_R)$ either permutes all three elements of this orbit 
or else centralizes~$t$. For each possible class type of~$t$, the relative Weyl 
group in $C_G( t )$ of~$S$ is cyclic of order~$3$. Thus $3 \mid |\Out_G( R, b_R)|$ 
and $2 \nmid |\Out_G( R, b_R)|$. It follows that $|\Out_G(R,b_R)| = 3$, as 
$N_G( R )/C_G( R ) \cong \SL_2(3)$. This yields~(b).

Finally, assume that $R =_G R_{15}$. Then $C_G(R) = R$, and thus $b_R$ is the
principal block of~$C_G(R)$. Hence $\Out_G(R,b_R) = N_G(R)/R \cong \SL_3(3)$ by
Table~\ref{tabradicalageq2}. Since this is not a
$3'$-group,~$R$ is not a defect group. Now~$\theta_R$ extends to 
$N_G(R)$, and as~$\SL_3(3)$ has a unique $3$-defect zero character, we obtain
$\cW(R,b_R) = 1$ by Equation~(\ref{NumberOfWeights}). 
\end{prf}

\medskip

\noindent
Now let $(Q, b_Q)$ denote a centric Brauer pair with~$Q$ abelian. Assume that 
$(Q,b_Q)$ is a $b$-Brauer pair for some $3$-block~$b$ of~$G$ with $b \subseteq 
\cE_3( G, s )$, where~$s$ is a semisimple $3'$-element of~$G$.
If we assume that $Q \not\in_G \{R_1, R_8, R_{15} \}$, then 
$\mathbf{M} := C_{\mathbf{G}}( Q )$ is a regular subgroup of~$\mathbf{G}$
by Lemma~\ref{CentralizersAbelianRadicalSubgroups}. In these cases, if we let
$\mathbf{M}^*$ denote a regular subgroup of~$\mathbf{G}$ dual to $\mathbf{M}$, 
and let $t \in \mathbf{M}$ be a semisimple $3'$-element such that  
$b_Q \subseteq \cE_3( M^*, t )$, then~$t$ is conjugate to~$s$ in~$G$; see 
Lemma~\ref{ConjugateSemisimpleLabels}. In the following lemma we list the 
candidates for the pairs $(Q, \theta_Q )$.

\begin{lem}
\label{PossibleCanonicalCharacters}
Let~$Q$ be a non-trivial, abelian radical $3$-subgroup of~$G$, but
$Q \not\in_G \{ R_1, R_2, R_3, R_8, R_{15} \}$. Then $\mathbf{M} :=
C_{\mathbf{G}}( Q )$ is a regular subgroup of~$\mathbf{G}$ by
{\rm Lemma~\ref{CentralizersAbelianRadicalSubgroups}}. Choose a regular
subgroup~$\mathbf{M}^*$ of~$\mathbf{G}$ dual to~$\mathbf{M}$. Let
$\theta \in \cE( M, t )$ for some semisimple $3'$-element $t \in M^*$.
Then $\theta( 1 )_3 = |M/Q|_3$ in exactly one of the following cases.

{\rm (a)} We have $Q =_G R_{10}$, $t \in Z( M^* )$ and $\theta = \hat{t}\psi$,
where $\hat{t}$ is the linear character of~$M$ that corresponds to~$t$ by
duality, and~$\psi$ is the unipotent character of~$M$ of
degree $q(q - \varepsilon)^2/2$.

{\rm (b)} We have $\theta = \pm R_{\mathbf{S}^*}^{\mathbf{M}}( t )$,
where~$\mathbf{S}^*$ is an $F$-stable maximal torus of~$\mathbf{M}^*$, and
$t \in S^*$ is in general position (with respect to~$\mathbf{M}^*$). The cases
for $Q$ and $\mathbf{S}^*$ are as given the following table,
where the conjugacy class of~$\mathbf{S}^*$ in~$\mathbf{G}$ is given by the
integer~$k$, if $\mathbf{S}^* =_G \mathbf{M}_{20,k}$ as in
{\rm Table~\ref{SylowsInCentralizersI}}.
\setlength{\extrarowheight}{0.5ex}
$$
\begin{array}{c|ccccccccc} \hline\hline
Q &  R_6 & R_7 & R_9 & R_{10} & R_{11} & R_{12} & 
R_{16} & R_{17} & R_{18} \\ \hline
\mathbf{S}^* & 14, 15 & 19, 20 & 22,22 & 3,3; 23,24 & 7, 8 & 4, 5 & 
17, 18 & 12, 13 & 1, 2 \\ \hline\hline
\end{array}
$$
Except for $Q =_G R_{10}$, there is a unique such torus~$\mathbf{S}^*$,
whereas there are two for $Q =_G R_{10}$. The two values for~$k$ separated by 
a comma correspond to the cases $\varepsilon = 1$, and $\varepsilon = -1$, 
respectively. We have $t \in Z(M^*)$ if and only if $Q \in_G
\{ R_{11}, R_{12}, R_{18} \}$. In these cases, $\mathbf{S}^* = \mathbf{M}^*$.
\end{lem}
\begin{prf}
Suppose first that $t \in Z( M^* )$. Then $\cE( M, t ) = \hat{t} \cE( M, 1 )$,
so that $\theta = \hat{t} \psi$ for some unipotent character
$\psi \in \cE( M, 1 )$. Now $\theta( 1 )_3 = |M/Q|_3$, if and only if~$\psi$
is a $3$-defect zero character of~$M/Z( M )$. The possibilities for~$M$ can be
determined from Lemma~\ref{CentralizersAbelianRadicalSubgroups} and
Table~\ref{SylowsInCentralizersI}. Looking at the degrees of the unipotent
characters of these groups, we get exactly the cases in~(a) and those in~(b) for
$Q \in_G \{ R_{11}, R_{12}, R_{18} \}$.

Now assume that $t \notin Z( M^* )$ and that $Q =_G R_{10}$. As $\mathbf{M} =
\mathbf{M}^*$ in this case and $Z( \mathbf{M} )$ is connected, the character
degrees of~$M$ and those of~$M/Z( M )$ agree. Now $M/Z( M ) = 
(\mathbf{M}/Z( \mathbf{M} ))^F$, and $\mathbf{M}/Z( \mathbf{M} )$ is a
simple group of adjoint type~$C_2$. The character degrees of $M/Z(M)$ can
be found in the tables~\cite{LLCD}, yielding our claim.

Finally, assume that $Q \not\in_G \{ R_{10} \}$. Then $\mathbf{M}$ is of
type~$A$, and $\theta = \pm 
R_{\mathbf{L}}^{\mathbf{M}}( \hat{t} \nu )$, where~$\mathbf{L}$ is a
regular subgroup of~$\mathbf{M}$ with $\mathbf{L}^* = C_{\mathbf{M}^*} ( t )$,
and $\nu \in \cE( L, 1 )$.
As $Q \leq Z( M ) \leq L$, we have $\theta(1)_3 = |M/Q|_3$ if and only if
$\nu(1)_3 = |L/Q|_3$. Also,~$\mathbf{L}$ is a proper subgroup of~$\mathbf{M}$,
since $t \not\in Z( M^* )$. Now~$L$ has a unipotent character~$\nu$ with
$\nu(1)_3 = |L/Q|_3$, if and only if~$L$ is a maximal torus of~$M$ and
$|L/Q|_3 = 1$. This rules out the possibilities $Q \in_G \{ R_4, R_5, R_{13},
R_{14} \}$, and gives exactly the remaining cases in~(b).
\end{prf}

\begin{cor}
\label{AssumptionsSatisfied}
Let the notation be as in {\rm Lemma~\ref{PossibleCanonicalCharacters}}, and
suppose that $(Q, b_Q)$ is a centric Brauer pair with $b_Q \subseteq
\cE_3( M, t )$. Then $(Q, b_Q)$ satisfies the hypotheses of 
{\rm Proposition~\ref{cormain}(a)}. In particular, $\Out_G( Q, b_Q ) \cong
\Out_{C_G( t )}( Q^\dagger )$ (recall that we have identified~$\mathbf{G}$
with~$\mathbf{G}^*$).
\end{cor}
\begin{prf}
This is clear in case~(a) of Lemma~\ref{PossibleCanonicalCharacters}, as~$\psi$
is invariant under all automorphisms of~$\Sp_4(q)$ (for odd~$q$ this follows
from \cite[Remarks on p.~$159$]{Lu} and for even~$q$ from
\cite[Theorem~$2.5$]{MalleExt}). In case~(b), $b_Q = \cE_3( M, t )$, as~$t$ is
in general position in~$\mathbf{S}^*$.
\end{prf}

\addtocounter{subsection}{4}
\subsection{Defect groups  of $F_4(q)$}\label{secdefF4}
Let~$b$ be a $3$-block of~$G$ with defect group~$D$. Then $b \subseteq 
\cE_3(G, s)$ for a semisimple $3'$-element $s \in G^*$.
As usual, we identify~$\mathbf{G}$ with~$\mathbf{G}^*$ and $G$
with~$G^*$. By Tables~\ref{1}--\ref{19}, we know that~$b$ is uniquely determined
by the conjugacy class of~$s$ in~$G$ and its defect group~$D$ (even by the
conjugacy class of~$s$ and $|D|$). We label the block~$b$ by $(s, D)$.
We begin with a lemma. 

\addtocounter{thm}{1}
\begin{lem}
\label{DualityInT}
Let $\mathbf{L}, \mathbf{L}^1, \mathbf{L}^2$ and~$\mathbf{T}$ be defined as in 
{\rm Proposition~\ref{C3C}}, and put $T_i := T \cap L^i$ for $i = 1, 2$. The 
identification of $\mathbf{G}$ with $\mathbf{G}^*$ induces a bijection
$T \rightarrow \Irr(T), s \mapsto \hat{s}$ such that the following holds for all
$3'$-elements $s \in T$.

We have $s \in T_1$ if and only if $T_1$ is contained in the kernel 
of~$\hat{s}$.
\end{lem}
\begin{prf}
We only prove the claim for $\varepsilon = 1$. The other case is proved by
twisting with $w_0$. Recall that $\mathbf{L}^1 =
[\mathbf{L}_{\{1, {23}\}},\mathbf{L}_{\{1, {23}\}}]$
and $\mathbf{L}^2 =
[\mathbf{L}_{\{3, {4}\}},\mathbf{L}_{\{3, {4}\}}]$.
Also, $\mathbf{T}_i = \mathbf{T} \cap \mathbf{L}^i$, for $i = 1, 2$.

Let $s \in T$ be a $3'$-element. Suppose that 
$s \in T_1 = \langle \alpha_1^\vee( t ), \alpha_{23}^\vee( t )
\mid t \in \mathbb{F}_q^* \rangle$. Then $s$ is in the kernel of~$\hat{s}$ by 
Lemma~\ref{KernelsOfDualCharacters}. Similarly, $s$ is in the kernel 
of~$\hat{s}$ if $s \in T_2$.
Suppose that $T_1 \leq \ker( \hat{s} )$. As 
$T_1T_2$ has index $3$ in~$T$, we have $s \in T_1T_2$. Write $s = s_1 s_2$ with 
$s_i \in T_i$ for $i = 1, 2$. By the above, $T_1 \leq \ker(\hat{s}_1^{-1})$. In 
turn, $T_1 \leq \ker( \hat{s}_1^{-1}\hat{s} )$. As $\hat{s}_1^{-1}\hat{s} = 
\hat{s}_2$, it follows that $\hat{s}_2 = 1_T$, and thus $s_2 = 1$.
\end{prf}

\begin{prop}\label{AllDefectGroups}
Let~$b$ be a $3$-block of~$G$ labeled by $(s, D)$. Then the following 
statements hold.
	
{\rm (a)} We have $D^\dagger \leq C_G ( s )$ for some $D^\dagger$ dual to~$D$ 
(in the sense of {\rm Subsection~\ref{DualRadicalSubgroups}}). 

{\rm (b)} The conjugacy class of~$D$ in $\cR_3(G)$ is as given in 
{\rm Column~$9$} of the row corresponding to~$b$ in {\rm Table}~{i}, where~$i$ 
is the $\mathbf{G}$-class type of~$s$.
\end{prop}
\begin{prf}
To prove~(a), we begin with the case that~$D$ is abelian and let $(D,b_D)$ 
denote a maximal $b$-Brauer pair. If $D =_G R_8$, the claim follows from
Lemma~\ref{ExceptionalAbelianRadicalSubgroups}(a). Hence assume that 
$D \neq_G R_8$ in the following. Then $C_{\mathbf{G}}( D )$ is a regular 
subgroup of~$\mathbf{G}$ by 
Lemmas~\ref{PossibleCanonicalCharacters},~\ref{ExceptionalAbelianRadicalSubgroups} 
and~\ref{CentralizersAbelianRadicalSubgroups}, and~(a)
follows from Lemma~\ref{ConjugateSemisimpleLabels}.

Assume next that~$s$ is not quasi-isolated and let $\mathbf{M}^\dagger$ denote 
the regular subgroup of~$\mathbf{G}$ defined in Theorem~\ref{BoRoEtAl}, which is 
minimal with the property that $C_{\mathbf{G}}( s ) \leq \mathbf{M}^\dagger$. 
Let $\mathbf{M}$ denote a regular subgroup of~$\mathbf{G}$ dual 
to~$\mathbf{M}^\dagger$. By Theorem~\ref{BoRoEtAl}(a), we may assume that 
$D \leq M$, and hence $D^\dagger \leq M^\dagger$. Now if $|D| = |M|_3$, we may 
choose $D^\dagger \leq C_G( s )$, which proves~(a) in these cases.
If $|D| < |M|_3$, which is the case if~$s$ has $\mathbf{G}$-class 
type~$11$ or~$12$, then the Sylow $3$-subgroups of~$C_G( s )$ are abelian;
see Tables~\ref{SylowsInCentralizersI} and~\ref{tabradicalageq2}. 
Then~$D$ is abelian by Theorem~\ref{BoRoEtAl}(c), a case we have already
settled.

If $s = 1$, the claim is trivially satisfied. We are left with the case that~$D$ 
is non-abelian and that $s \neq 1 $ is isolated, 
i.e.\ that the $\mathbf{G}$-class type of~$s$ is one of~$2$,~$3$ or~$5$. 
Then $|D| = 3^{4a+1}$ by Tables~\ref{2}--\ref{5}, and~$D$ has a normal subgroup 
conjugate to $R_{18}$ by \cite[Proposition~$3.5$]{KeMa}. Let~$P$ denote a
Sylow $3$-subgroup of~$C_G( s )$. Table~\ref{tabradicalageq2}
implies that $D \in_G \{ P, P^\dagger \}$, and thus
$D \in_G \{ R_{33}, R_{34} \}$. 
Choose $Q^\dagger = Q =_G R_{18}$ such that $Q^\dagger \leq P$ is the maximal 
abelian normal subgroup of~$P$ described in Lemma~\ref{MaximalAbelianSubgroups}. 
Then $s \in C_G(Q^\dagger) \cong C_G( Q )^*$. As $C_G(Q) = T$ is a maximal 
torus, $b_Q := \cE_3(C_G(Q), s)$ is a block
of~$C_G(Q)$. Let~$b'$ be the unique block of~$G$ such that that $(1,b') \leq 
(Q, b_{Q})$; see \cite[Theorem~$2.9$]{Kessar}. Choose a maximal Brauer pair
$(D',b_{D'})$ containing $(Q, b_{Q})$. Then~$D'$ is a defect group of~$b'$.
By Lemma~\ref{ConjugateSemisimpleLabels}, applied to the $b'$-Brauer pair
$(Q,b_Q)$, we may assume that $b' \subseteq 
\cE_3(G, s)$. From Columns~$6$ of Tables~\ref{2}--\ref{5} we conclude that
${Q} \lneq D'$ and thus $|D| = |D'|$. As $b$ is uniquely determined by $(s,|D|)$,
we conclude that $b' = b$ and $D =_G D' = R$. Thus assume that $D = D'$ and so
$Q \leq D$ in the following. Then $(Q,b_Q) \leq (D,b_D)$. We have $C_G( Q ) 
= T$, and $C_G( D ) \in \{ T_1, T_2 \}$; see Table~\ref{tabradicalageq2}.
In particular, $\theta_Q = \hat{s}$, where $\hat{s} \in \Irr( T )$ corresponds
to $s \in T$ by duality, and $\theta_D$ is the restriction of~$\hat{s}$ to
$C_G ( D )$. As~$s$ is a $3'$-element, $s \in C_G(D) \cup C_G(D^\dagger) =
T_1 \cup T_2$ and $\theta_D \neq 1$, we obtain $s \in C_G ( D^\dagger )$ by 
Lemma~\ref{DualityInT}. Hence $D^\dagger =_G P$ as claimed. The proof
of~(a) is complete.

To prove~(b), assume first that~$b$ corresponds to the 
principal block of $C_G( s )$ via Theorem~\ref{BoRoEtAl} (i.e.\ the label in 
Columns~$5$ of Tables~\ref{2}--\ref{19} equals~$1$). Then $|D| = |C_G(s)|_3$, so
that~$D$ is conjugate in~$G$ to a Sylow $3$-subgroup of~$C_G(s)$ by~(a).
In this case, the claim follows from Table~\ref{SylowsInCentralizersI}.
Assume that~$b$ does not correspond to the principal block of~$C_G( s )$.
Then~$D$ is abelian by \cite[Theorem~$1.2$]{KeMa} for quasi-isolated~$s$, and by 
Theorem~\ref{BoRoEtAl} in the other cases. Using the knowledge of $d(b)$ and
Table~\ref{tabradicalageq2}, we find that $D \in_G \{ [3^a]^2, [3^a], 1 \}$. 
Suppose first that $D = [3^a]^2$. If $D =_G R_8$, the claim holds by 
Lemma~\ref{ExceptionalAbelianRadicalSubgroups}. We are left with the cases 
$D \in_G \{ R_9, R_{10} \}$ by Lemma~\ref{PossibleCanonicalCharacters}.
In particular, $D =_G D^\dagger$.
The possible class types for~$s$ can be found in Tables~\ref{1}--\ref{19}.
It follows from this and Table~\ref{SylowsInCentralizersI}, that there is a
conjugate~$Q$ of~$R_9$ and a conjugate~$R$ of~$R_{10}$ such that $C_G( Q )
\leq C_G( s )$ and $C_G( R ) \leq C_G( s )$, unless~$s$ is of class type
$(15,1)$ or $(15,3)$. Thus in the former cases, we may assume that
$s \in Z( C_G( D^\dagger ) ) = Z( C_G ( D )^* )$, and hence $D =_G R_{10}$ by
Lemma~\ref{PossibleCanonicalCharacters}. In the latter cases, assume that
$D =_G R_9$. Lemmas~\ref{ConjugateSemisimpleLabels} 
and~\ref{PossibleCanonicalCharacters} imply that~$s$ lies in a
maximal torus of $C_{\mathbf{G}}( D )$ of class type~$(20,22)$. However,
$C_{\mathbf{G}} ( s )$ does not have a maximal torus of this type, as
Table~\ref{cc} shows. This contradiction shows that $D =_G R_{10}$
in these cases as well.

Finally, suppose that $D = [3^a]$. Then $D \in_G \{ R_2, R_3 \}$ by
Lemmas~\ref{ExceptionalAbelianRadicalSubgroups} 
and~\ref{PossibleCanonicalCharacters}. Recall that $R_2^\dagger =_G R_3$ and 
$R_3^\dagger =_G R_2$. The possible class types of~$s$ can be found in 
Tables~\ref{1}--\ref{19}. In particular,~$s$ is not quasi-isolated. 
By~(a), we may assume that $D^\dagger \leq C_G( s )$. Put $\mathbf{K}^\dagger := 
C_{\mathbf{G}}( D^\dagger ) \cap \mathbf{M}^\dagger$, where~$\mathbf{M}^\dagger$ 
is as in the proof of~(a). Let $t \in G$ with $D^\dagger = \langle t \rangle$. 
Then $C_{\mathbf{G}}( st ) =
C_{\mathbf{G}}( D^\dagger ) \cap C_{\mathbf{G}}( s ) \leq \mathbf{K}^\dagger$,
with equality if $\mathbf{M}^\dagger = C_{\mathbf{G}}( s )$. The latter holds
unless $s$ is of class type $(11,k)$. If $D^\dagger =_G R_2$, 
respectively~$R_3$, then $C_{\mathbf{G}}( D^\dagger)$ is $G$-conjugate to 
$\mathbf{M}_{10,k}$, respectively $\mathbf{M}_{6,k}$ with $k = 1, 2$; see 
Lemma~\ref{CentralizersAbelianRadicalSubgroups}.
By running through the possible class types for~$s$, we can determine 
$\mathbf{K}^\dagger$ in each case by comparing the maximal tori fusing into
$C_{\mathbf{G}}( D^\dagger)$ respectively $\mathbf{M}^\dagger$; see 
Table~\ref{cc}. We give more 
details in case $\varepsilon = 1$. Then $C_{\mathbf{G}}( D^\dagger ) \in_G 
\{ \mathbf{M}_{6,1}, \mathbf{M}_{10,1} \}$. It turns out that, unless the class 
type of~$s$ equals~$(11,2)$, there is a unique $i \in \{ 6, 10 \}$ such that the 
intersection of a $G$-conjugate of $\mathbf{M}_{i,1}$ with $\mathbf{M}^\dagger$ 
can be the centralizer of~$st$. A variant of this argument also works in the 
remaining case, where $\mathbf{M}^\dagger =_G \mathbf{M}_{10,2}$. We have thus 
determined~$D^\dagger$. The same proof works for $\varepsilon = -1$.
%%%%%%%%%%%%%%%%%%%%%%%%%%%%%%%%%%%%%%%%%%%%%%%%%%%%%%%%%%%%%%%%%%%%%%%%%%%%%%%%
%%%%%%%%%%%%%%%%%%%%%%%%%%%%%%%%%%%%%%%%%%%%%%%%%%%%%%%%%%%%%%%%%%%%%%%%%%%%%%%%
%%
%% For details see Sheets 1 -- 4 of my notes of 24.03.2020.
%%
%%%%%%%%%%%%%%%%%%%%%%%%%%%%%%%%%%%%%%%%%%%%%%%%%%%%%%%%%%%%%%%%%%%%%%%%%%%%%%%%
%%%%%%%%%%%%%%%%%%%%%%%%%%%%%%%%%%%%%%%%%%%%%%%%%%%%%%%%%%%%%%%%%%%%%%%%%%%%%%%%
\end{prf}

\medskip

\noindent
The classification of the defect groups for the abelian $3$-blocks of~$G$ 
reveals an analogous behavior as in the case of $\ell$-blocks for $\ell > 3$.
\begin{cor}\label{Abelian3DefectGroups}
Let $s \in G$ be a semisimple $3'$-element and let $b \subseteq \cE_3( G, s )$ be a 
$3$-block of~$G$ with a non-cyclic abelian defect group~$D$. Then there is an 
$e$-split Levi subgroup~$\mathbf{M}$ of~$\mathbf{G}$ such that~$D$ is the Sylow 
$3$-subgroup of~$Z(M)$, unless the $G$-class type of~$s$ is one of $(18,3)$ or 
$(19,8)$ if $e = 1$, respectively $(18,7)$ or $(19,9)$ if $e = 2$.
\end{cor}
\begin{prf}
This follows from Proposition~\ref{AllDefectGroups} in connection with 
Lem\-ma~\ref{CentralizersAbelianRadicalSubgroups}.
\end{prf}
%%%%%%%%%%%%%%%%%%%%%%%%%%%%%%%%%%%%%%%%%%%%%%%%%%%%%%%%%%%%%%%%%%%%%%%%%%%%%%%%
%%%%%%%%%%%%%%%%%%%%%%%%%%%%%%%%%%%%%%%%%%%%%%%%%%%%%%%%%%%%%%%%%%%%%%%%%%%%%%%%
%%
%% \marginpar{Perhaps~\ref{Abelian3DefectGroups} -- \ref{Lemma622} are not needed anymore.}
%%
%% Die naechsten beiden Ergebnisse am 14.08.2020 auskommentiert.
%% Corollary \ref{Abelian3DefectGroups} habe ich erstmal noch behalten.
%%
%%%%%%%%%%%%%%%%%%%%%%%%%%%%%%%%%%%%%%%%%%%%%%%%%%%%%%%%%%%%%%%%%%%%%%%%%%%%%%%%
%%%%%%%%%%%%%%%%%%%%%%%%%%%%%%%%%%%%%%%%%%%%%%%%%%%%%%%%%%%%%%%%%%%%%%%%%%%%%%%%

\addtocounter{subsection}{3}
\subsection{Candidates for weight subgroups}
Let~$b$ be a $3$-block of~$G$ with defect group~$D$. If~$D$ 
is non-abelian and not a Sylow $3$-subgroup of~$G$, it is conjugate to one of 
$R_{29}\text{--}R_{34}$ by Tables~\ref{1}--\ref{19}. Recall from the summary 
given in~\ref{BlocksAndWeights}, that if $(Q,b_Q)$ is a $b$-Brauer pair giving 
rise to a $b$-weight, then $(Q,b_Q)$ is centric and~$Q$ is a radical 
$3$-subgroup of~$G$. By Equation~(\ref{CenterContainment}), we have $Z(D) \leq 
Z(R) \leq R \leq D$ for some conjugate~$R$ of~$Q$. Using these conditions, the 
following lemma restricts the set of candidates for~$Q$ in case~$D$ is 
non-abelian.

\addtocounter{thm}{1}
\begin{lem}
\label{AllCandidates}
Let~$b$ be a $3$-block of $G = F_4(q)$ and let $(D, b_D)$ denote a maximal
$b$-Brauer pair, where $D \in_G \{ R_{29}\text{--}R_{34} \}$. Let $Q \lneq D$ be 
such that the Brauer pair $(Q,b_Q) \leq (D,b_D)$ is centric and $\cW(Q, b_Q) 
\neq 0$. Then $Q$ or $Q^\dagger$ is one of the groups contained in the
following table.

\setlength{\extrarowheight}{0.5ex}
$$\begin{array}{c|l} \hline\hline
D & \multicolumn{1}{c}{Q/Q^\dagger} \\ \hline\hline
R_{33} & R_{16}, R_{18}, R_{25} \\ \hline
R_{31} & R_{16}, R_{19} (a = 1), R_{23} (a \geq 2)\\ \hline
R_{29} & R_{19} \\ \hline\hline
\end{array}
$$
\end{lem}
\begin{prf}
As $(Q,b_Q)$ is a centric $b$-Brauer pair but not maximal, $|\Out_G(Q,b_Q)|$ is 
divisible by~$3$ (this is a consequence of Brauer's First Main Theorem; see 
\cite[Theorem~$6.7.6$(v)]{MLBooks}). Moreover, $\Out_G(Q,b_Q)$ has an irreducible 
projective character of degree divisible by~$3$, since $\cW(Q, b_Q) \neq 0$. 
These conditions exclude the cases $Q \in_G \{ R_1\mbox{--}R_{7} \}$, 
$Q \in_G \{ R_9\mbox{--}R_{14} \}$ and $Q \in_G \{ R_{27}\mbox{--}R_{32} \}$. 
Indeed, by Table~\ref{tabradicalageq2}, in these cases $\Out_G(Q,b_Q)$ is a 
$2$-group or a subgroup of $[6] \times S_3$. No such group of order divisible 
by~$3$ has a projective character of degree divisible by~$3$. By 
Lemma~\ref{ExceptionalAbelianRadicalSubgroups}, the radical subgroup~$R_8$
is also excluded.

Now suppose that $Q$ is one of the other radical $3$-subgroups of~$G$.
From $Q \leq D$ we obtain $C_G ( D ) \leq C_G ( Q )$. Together with 
$Z( D ) \leq Z( Q )$ our claim follows from Table~\ref{tabradicalageq2}. For 
example, if $D =_G R_{33}$, we have $Z(D) = [3^a]^2 =_G R_{13}$ and the type of 
the characteristic of $Z( D )$ equals $3\A_6\C_2$. This yields~$R_{25}$ as only 
non-abelian candidate for~$Q$, and $R_{16}$, $R_{18}$ as the only abelian
candidates.
\end{prf}

%\markboth{Alperin's blockwise weight conjecture for~$F_4(q)$}
\section{Alperin's blockwise weight conjecture for~$F_4(q)$}
\label{SectionWeights}

Let $G = F_4(q)$ with~$q$ a power of the prime $p$. We are now ready to prove 
the first main result of this article.

\begin{thm}\label{thmAWC} 
Let $\ell$ be an odd prime not dividing~$q$. Let~$b$ be an $\ell$-block of 
$G = F_4(q)$ and let $(D, b_D)$ denote a maximal $b$-Brauer pair. Then the
following two statements hold.

{\rm (a)} The number of $G$-conjugacy classes of $b$-weights
equals the number of irreducible Brauer characters of~$b$, i.e.\ 
the Alperin weight conjecture is satisfied for~$b$.

{\rm (b)} Let $(R,b_R)$ denote a centric $b$-Brauer pair contained in $(D,b_D)$, 
i.e.\ $(1,b) \leq (R,b_R) \leq (D,b_D)$ and~$b_R$ has defect group~$Z(R)$. If 
$\cW( R, b_R ) \neq 0$, the 
canonical character~$\theta_R$ of~$b_R$ extends to $N_G( R, b_R )$. In 
particular, $\cW(R, b_R) = |\Irr^0( \Out_G( R, b_R ))|$ by
{\rm Equation~(\ref{NumberOfWeights})}.
\end{thm}
\begin{prf}
Suppose first that $\ell > 3$. Then~(a) follows from~(b) by 
\cite[Corollary~$3.7$]{MalleAb}. As~(b) holds by 
\cite[Proposition~$2.3$(a)]{MalleAb} and
Corollary~\ref{MaximalExtendability}, we are done in this case.

For the remainder of the proof we assume $\ell = 3$.
Let $s \in G$ be a semisimple $3'$-element with $b \subseteq \cE_3(G,s)$.
If $C_{\mathbf{G}}( s )$ is a maximal torus of~$\mathbf{G}$, then~$D$ is 
abelian, $l(b) = 1$ and $|\Out_G( D, b_D )| = 1$, so that both claims hold.
To prove statement~(a), it thus suffices to show that the entries in 
Columns~$9$--$11$ 
of each of the Tables~\ref{1}--\ref{19} are correct. Let $(R,b_R)$ denote a 
centric $b$-Brauer pair contained in $(D,b_D)$ as in~(b). The entries in 
Column~$9$ corresponding to~$b$ list those~$R$, for which $\cW := \cW(R, b_R)$ 
is non-zero, where the first entry in this column is the defect group~$D$. The 
latter has already been determined in Proposition~\ref{AllDefectGroups}. For 
each~$R$ in Column~$9$, the values of~$\Out_G( R, b_R)$ and~$\cW$ are given in 
Columns~$10$ and~$11$, respectively.
In view of the last statement of~(b),~$\cW$ can be determined from the entries 
in Column~$10$. Notice that for~$R$ to appear in Column~$9$, we must have 
$Z(D) \leq Z(R)$ by Equation~(\ref{CenterContainment}). In 
Lemma~\ref{AllCandidates} we have listed the candidates for these~$R$ in 
case~$D$ is non-abelian. If~$D$ is abelian, then $R = D$ and $\cW(R, b_R) = 
\cW(D, b_D)$. If~$D$ is cyclic, statement~(a) is known to hold by 
\cite[Theorem 1.1]{KS16}. 
If~$R$ is cyclic, so is $\Out_G( R,  b_R )$ by Table~\ref{tabradicalageq2}, and 
thus~$\theta_R$ extends to~$N_G( R, b_R )$.
In what follows we will therefore always assume that~$D$ and~$R$ are
non-cyclic.

We first prove our assertions for the principal block.
Let~$b$ be the principal block of~$G$, so that $D \in \Syl_3(G)$ and~$b_D$ is 
the principal block of~$C_G(D)$. Moreover,~$b_R$ is the principal block 
of~$C_G(R)$ and~$\theta_R$ is the trivial character. Thus $|C_G(R)/Z(R)|_3 = 
1 = \theta_R(1)_3$ and $C_G(R)$ is abelian by Table~\ref{tabradicalageq2}. 
Also, $N_G(R, b_R) = N_G(R)$ and $\Out_G(R,b_R) = \Out_G( R )$. In particular,
$\cW(R,b_R) = |\Irr^0( \Out_G(R) )|$. The set of irreducible characters of 
$\Out_G(R)$ is easy to compute with the help of a computer algebra system such 
as MAGMA~\cite{magma} or GAP~\cite{GAP04}.
If $C_G(R) = Z(R)$, we may suppose that 
$$R \in_G \{ R_{15}, R_{21}, R_{22}, R_{35}, R_{36}, R_{37}, R_{38} \}$$ 
by Table~\ref{tabradicalageq2}. If 
$C_G(R) \neq Z(R)$, then $\Irr^0( \Out_G(R) ) = \emptyset$ by 
Table~\ref{tabradicalageq2}, except when $R =_G R_{18}$, in which case 
$\Out_G(R) = W(F_4)$. This completes the proof for the principal block.

From now on we assume that~$b$ is not the principal block, and continue with a 
preliminary consideration. Suppose that~$R$ is abelian. Then 
$R \not\in_G \{ R_1, R_8, R_{15} \}$ by Tables~\ref{1}--\ref{19}, and thus
$C_{\mathbf{G}}( R )$ is a regular subgroup of~$\mathbf{G}$ by
Lemma~\ref{CentralizersAbelianRadicalSubgroups}. By
Lemma~\ref{ConjugateSemisimpleLabels}, we may assume that
$s \in C_G( R^\dagger ) \cong C_G( R )^*$ and that $b_R \subseteq
\cE_3( C_G( R ), s )$. By Lemma~\ref{PossibleCanonicalCharacters} and
	Corollary~\ref{AssumptionsSatisfied}, we have $\Out_G( R, b_R ) \cong
\Out_{C_G( s )}( R^\dagger )$. Put
$$\mathbf{K}^\dagger := C_{\mathbf{G}}( R^\dagger ) \cap C_{\mathbf{G}} ( s ).$$
Lemma~\ref{CentralizersAreRegular} implies that $\mathbf{K}^\dagger =
C_{C_{\mathbf{G}}(s)}( R^\dagger )$ is a regular subgroup of
$C_{\mathbf{G}} ( s )$, as $3 \nmid
|Z(C_{\mathbf{G}}( s )/Z(C_{\mathbf{G}}( s )^\circ|$ and~$3$ is a good prime
for $C_{\mathbf{G}}( s )$. Also, $\mathbf{K}^\dagger =
C_{C_{\mathbf{G}}( R^\dagger )}(s)$ is a regular subgroup
of~$C_{\mathbf{G}}( R^\dagger )$ by~\cite[Corollary~$2.6$]{GeHi1}.
Proposition~\ref{cormain}(b) implies that
$\Out_{C_G( s )}( R^\dagger ) = W_{C_{\mathbf{G}}(s)}( \mathbf{K}^\dagger )^F$.

We proceed to prove our assertions in case of abelian defect groups.
Assume that $D$ is abelian. By Tables~\ref{1}--\ref{19}, we have
$$D \in_G \{ R_9, R_{10}, R_{11}, R_{12}, R_{16}, R_{17}, R_{18} \}.$$
In each case we find $C_{\mathbf{G}}( D^\dagger )$ from 
Lemma~\ref{CentralizersAbelianRadicalSubgroups}.
We now run through all the possibilities for~$s$ and~$D$, 
using the information contained in Table~\ref{cc} on the inclusion of
maximal tori. Except if $D =_G R_{10}$, it turns out that $C_{\mathbf{G}}( s )$ 
and $C_{\mathbf{G}}( D^\dagger )$ have a unique common $F$-stable maximal 
torus~$\mathbf{S}^\dagger$, hence $\mathbf{K}^\dagger = 
C_{\mathbf{G}}( D^\dagger ) \cap C_{\mathbf{G}} ( s ) = \mathbf{S}^\dagger$. In 
the other cases, we have $s \in Z( C_{\mathbf{G}}( D^\dagger ) )$ so that 
$\mathbf{K}^\dagger = C_{\mathbf{G}}( D^\dagger )$. In any case, the relative 
Weyl group $W_{C_{\mathbf{G}}(s)}( \mathbf{K}^\dagger )^F$ and hence the 
isomorphism type of $\Out_G( R, b_R )$ is easy to determine using standard 
methods.
%%%%%%%%%%%%%%%%%%%%%%%%%%%%%%%%%%%%%%%%%%%%%%%%%%%%%%%%%%%%%%%%%%%%%%%%%%%%%%%%
%%%%%%%%%%%%%%%%%%%%%%%%%%%%%%%%%%%%%%%%%%%%%%%%%%%%%%%%%%%%%%%%%%%%%%%%%%%%%%%%
%%
%% For details see sheets 3 -- 8 of my notes of 09.05.2018.
%%
%% For more details see sheets 6 -- 15 of my notes of 24./27.07.2020.
%%
%%%%%%%%%%%%%%%%%%%%%%%%%%%%%%%%%%%%%%%%%%%%%%%%%%%%%%%%%%%%%%%%%%%%%%%%%%%%%%%%
%%%%%%%%%%%%%%%%%%%%%%%%%%%%%%%%%%%%%%%%%%%%%%%%%%%%%%%%%%%%%%%%%%%%%%%%%%%%%%%%

We now show that~$\theta_D$ extends to $N_G( D, b_D )$. If 
$D \in_G \{ R_{11}, R_{12} \}$, this is clear as $\Out_G( D, b_D )$ is cyclic.
If $D =_G R_{18}$, the assertion follows from \cite[Theorem~$1.1$]{Spaeth10} 
and, in the remaining cases, from Propositions~\ref{NormalizerM14} 
and~\ref{NormalizerM18}. This completes the proof for abelian defect groups.

We now proceed to prove our claims for non-abelian defect groups.
Suppose that~$b$ is not the principal block and that~$D$ is non-abelian. Then
$D \in_G \{ R_{29}\text{--}R_{34} \}$ by Tables~\ref{1}--\ref{19}. It suffices
to consider the cases $D \in_G \{ R_{29}, R_{31}, R_{33} \}$, as the other cases
are proved in an analogous way. For each such~$D$, we discuss the possibilities 
for~$R$ determined in Lemma~\ref{AllCandidates} case by case. If 
$C_{\mathbf{G}}(s)$ is a regular subgroup of~$\mathbf{G}$, then, by 
Tables~\ref{1}--\ref{19}, the block~$b'$ of $M := C_G( s )^\dagger$ 
corresponding to~$b$ by Theorem~\ref{BoRoEtAl}(a) is of the form 
$\hat{s} \otimes b_0$, where~$b_0$ is the principal block of~$M$, so that we may 
apply Lemma~\ref{NewOutLemma}, which gives the extendibility of $\theta_R$ to 
$N_G( R, b_R )$.

Suppose first that $D =_G R_{29}$. Then~$s$ is of class type $(13,4)$ or 
$(13,5)$, according as $\varepsilon = -1$ or $\varepsilon = 1$, respectively.
By Lemma~\ref{AllCandidates} we have $\cW(R,b_R) = 0$ or $R \in_G 
\{ R_{19}, D \}$. We apply Lemma~\ref{NewOutLemma} with the subgroup $\mathbf{M} 
=_G \mathbf{M}_{17,4}$ respectively $\mathbf{M} =_G \mathbf{M}_{17,5}$, so that
$\mathbf{M}^\dagger$ is $G$-conjugate to~$C_{\mathbf{G}}( s )$. 
If $R =_G R_{19}$, Lemma~\ref{NewOutLemma} yields $\Out_G( R, b_R ) \cong 
\Out_M( R ) \cong \SL_2(3)$. Similarly, the case $R = D$ yields 
$\Out_G( R, b_R ) \cong \Out_M( R ) = 2$.

Next, assume that $D =_G R_{31}$. Then the class type of~$s$ is one of 
$$\{ (6,1), (6,2), (8,2), (8,3), (13,2), (13,3)\}.$$ 
By Lemma~\ref{AllCandidates}, we are left with the possibilities
\begin{equation}
\label{CandidatesAtR31}
R \in_G \{ R_{16}, R_{19}\,(a = 1), R_{23}\,(a \geq 2), D \}.
\end{equation}
Assume that~$s$ is of class type $(13,2)$ or $(13,3)$, according as
$\varepsilon = 1$ or~$-1$, respectively. To apply Lemma~\ref{NewOutLemma}, 
we take $\mathbf{M} =_G \mathbf{M}_{17,3}$, respectively $\mathbf{M} =_G 
\mathbf{M}_{17,2}$. In case $R =_G R_{16}$ we get $\Out_G( R, b_R ) \cong 
\Out_M( R ) \cong S_3$, and thus $\cW(R,b_R) = 0$. In 
the other cases, Lemmas~\ref{NewOutLemma} and~\ref{lemRadLi}, and the 
considerations on central products in Subsection~\ref{CentralProducts} yield 
the corresponding entries in Table~\ref{13}. 
Now suppose that~$s$ is of class type~$(6,1)$ or~$(6,2)$, according as
$\varepsilon = -1$ or $\varepsilon = 1$, respectively. We then choose
$\mathbf{M}$ as $\mathbf{M} =_G \mathbf{M}_{10,1}$, respectively 
$\mathbf{M} =_G \mathbf{M}_{10,2}$. 
Let $Q \leq D$ denote the abelian normal subgroup of~$D$ of
index~$3$ considered in Lemma~\ref{MaximalAbelianSubgroups}.
Then $Q =_G R_{16}$, and $\theta_Q$ is an extension of~$\theta_D$ to
$C_M( Q ) = S$, where $\mathbf{S}$ is a maximal $e$-split torus of~$\mathbf{M}$. 
We obtain $\Out_G( Q, b_Q ) = \Out_M( Q ) \cong W(C_3)$. 
If~$R$ is non-abelian, we proceed as follows. Notice that~$D$ is conjugate to
a Sylow $3$-subgroup of~$K$, where $\mathbf{K} \leq \mathbf{M}$ is a regular
subgroup conjugate to $\mathbf{M}_{17,3}$ if $\varepsilon = 1$, respectively 
$\mathbf{M}_{17,2}$, if $\varepsilon = -1$; see Table~\ref{cc}. Then 
$RC_G(R) \leq K \leq M$ by Table~\ref{tabradicalageq2}.  As~$M$ contains an 
element which induces the graph automorphism of~$K$ and normalizes~$R$, we also 
have $N_G( R ) \leq M$, hence $\Out_M( R ) = \Out_G( R )$, and thus the
required entries in Table~\ref{6}.
%%%%%%%%%%%%%%%%%%%%%%%%%%%%%%%%%%%%%%%%%%%%%%%%%%%%%%%%%%%%%%%%%%%%%%%%%%%%%%%%%
%%%%%%%%%%%%%%%%%%%%%%%%%%%%%%%%%%%%%%%%%%%%%%%%%%%%%%%%%%%%%%%%%%%%%%%%%%%%%%%%%
%%%
%%% For more details see sheets 1 - 3 of my notes of 14./15.08.2019.
%%%
%%%%%%%%%%%%%%%%%%%%%%%%%%%%%%%%%%%%%%%%%%%%%%%%%%%%%%%%%%%%%%%%%%%%%%%%%%%%%%%%%
%%%%%%%%%%%%%%%%%%%%%%%%%%%%%%%%%%%%%%%%%%%%%%%%%%%%%%%%%%%%%%%%%%%%%%%%%%%%%%%%

Now suppose that~$s$ is of class type $(8,2)$ if $\varepsilon = 1$,
respectively~$(8,3)$, if $\varepsilon = -1$. In these cases we choose
$\mathbf{M} =_G \mathbf{M}_{10,2}$, respectively $\mathbf{M} =_G
\mathbf{M}_{10,1}$, as the regular subgroup of~$\mathbf{G}$ minimal with the 
property that $\mathbf{M}^\dagger$ contains $C_{\mathbf{G}}( s )$.
Let $b' \subseteq \cE_3( M, s )$ denote the block of~$M$ corresponding to~$b$ 
via Theorem~\ref{BoRoEtAl}(a). By Theorem~\ref{BoRoEtAl}(b) we have 
$\Out_G( R, b_R ) = \Out_M( R, b'_R ) \leq \Out_M( R )$. The latter has already 
been determined 
in the previous paragraph. We let $\mathbf{K} \leq \mathbf{M}$ be as in that
paragraph and assume that $D \leq K$. 
Suppose that $R_{16} =_G Q \leq D$. Using Table~\ref{cc}, we find that
$C_{\mathbf{K}}( Q ) = C_{\mathbf{M}}( Q ) = \mathbf{S}$, a maximal torus of 
type $(20,12)$ if $\varepsilon = 1$, and $(20,13)$, if $\varepsilon = - 1$.
By the description of the centralizers in Table~\ref{tabradicalageq2},
this implies that $C_{\mathbf{K}}( R ) \leq \mathbf{S}$ and hence $C_M( R ) 
\leq S$ for all $R \leq D$
with $R$ occurring in~(\ref{CandidatesAtR31}). In turn, $\theta'_R$ is the 
restriction of~$\hat{s}$ to $C_M( R )$, where $\hat{s} \in
\Irr(S)$ corresponds to $s \in S^*$ via duality. Once more by Table~\ref{cc},
we may assume that~$s$ is a central element of $\mathbf{K}^\dagger$, 
where~$\mathbf{K}^\dagger$ is a suitable regular subgroup dual to~$\mathbf{K}$.
Now~$\hat{s}$ extends to a linear character of~$K$, as $s \in Z(K^\dagger)$;
see \cite[Proposition~$2.5.20$]{GeMa}. In particular,
$\hat{s}$ is invariant in $\Out_K( R )$. Observe that $\Out_K( R )$ has
index~$2$ in $\Out_M( R )$ for all $R$ occurring in~(\ref{CandidatesAtR31}).
We may assume that~$R$ is invariant under the graph automorphism of~$K$ 
induced by~$M$, and we claim that $\hat{s}$ is not invariant under this
automorphism. This will then give the desired values for $\Out_M( R, \hat{s} )$. 
By Corollary~\ref{AssumptionsSatisfied}, we have $\Out_M( Q, \hat{s} ) = W(A_3)$ 
for $Q =_G R_{16}$, by Corollary~\ref{AssumptionsSatisfied}. This implies the 
claim.

If $R =_G R_{16}$, the extendibility of $\theta_R$ follows from 
Proposition~\ref{NormalizerM18}. In the other cases, 
the Schur multiplier of $\Out_G( R, b_R )$ is trivial, which gives the 
extendibility.

Finally, assume that $D =_G R_{33}$. Then the class type of~$s$ is one of
$$\{ (2,1), (3,1), (3,2), (6,1), (6,2), (7,1), (7,2), (8,1), (8,4), (13,1), 
(13,6)\}.$$
By Lemma~\ref{ExceptionalAbelianRadicalSubgroups}, we are left with the
possibilities
$$R \in_G \{ R_{16}, R_{18}, R_{25}, D \}.$$
Aiming at a contradiction, assume that $R =_G R_{16}$. Considering the various 
possibilities for~$s$ and using Table~\ref{cc}, we find that
$s \in Z( C_G( R^\dagger ) )$, unless~$s$
is of type~$(8,1)$ or~$(13,1)$ if $\varepsilon = 1$, or of type~$(8,4)$
or~$(13,6)$ if $\varepsilon = -1$. Now $\theta_R \in \cE( C_G( R ), s )$ 
by Lemma~\ref{ConjugateSemisimpleLabels} and 
satisfies $\theta_R(1)_3 = |C_G(R)/R|_3$. By 
Lemma~\ref{PossibleCanonicalCharacters}, this excludes the cases with
$s \in Z( C_G( R^\dagger ) )$, and shows that~$s$ lies in a torus of type
$(20,17)$ if $\varepsilon = 1$, and of type $(20,18)$ if $\varepsilon = -1$.
But then~$s$ cannot be of class type $(8,1)$,~$(13,1)$,~$(8,4)$ nor~$(13,6)$,
once more by Table~\ref{cc}. This contradiction shows that $R \neq_G R_{16}$.

Let $Q \leq D$ denote the abelian normal subgroup of~$D$ of
index~$3$ considered in Lemma~\ref{MaximalAbelianSubgroups}.
Then $Q =_G R_{18}$, and $\theta_Q$ is an extension of~$\theta_D$ to
$C_G( Q ) = T$, the maximal $e$-split torus of~$G$. If $R = Q$, we obtain
$\Out_G( R, b_R )$ from Corollary~\ref{AssumptionsSatisfied}, and the
extendibility of $\theta_R$ from \cite[Theorem~$1.1$]{Spaeth10}. We are
thus left with the cases $R \in_G \{ R_{25}, D \}$.

If the class type of~$s$
is one of $(6,1)$, $(6,2)$, $(8,1)$, $(8,4)$, $(13,1)$ or $(13,6)$, we proceed 
exactly as in the corresponding cases for $R =_G R_{31}$. The case that $s$ 
is of class type $(7,1)$ or $(7,2)$ can also be dealt with analogously.
%%%%%%%%%%%%%%%%%%%%%%%%%%%%%%%%%%%%%%%%%%%%%%%%%%%%%%%%%%%%%%%%%%%%%%%%%%%%%%%%
%%%%%%%%%%%%%%%%%%%%%%%%%%%%%%%%%%%%%%%%%%%%%%%%%%%%%%%%%%%%%%%%%%%%%%%%%%%%%%%%
%%
%% For more details see sheets 1 - 3 of my notes of 16.08.2019.
%%
%%%%%%%%%%%%%%%%%%%%%%%%%%%%%%%%%%%%%%%%%%%%%%%%%%%%%%%%%%%%%%%%%%%%%%%%%%%%%%%%
%%%%%%%%%%%%%%%%%%%%%%%%%%%%%%%%%%%%%%%%%%%%%%%%%%%%%%%%%%%%%%%%%%%%%%%%%%%%%%%

Next, assume that~$s$ is of class type $(2,1)$. Then~$q$ is odd. We choose~$s$ 
in the center of the standard $e$-split Levi subgroup $\mathbf{K}^\dagger$ of 
type $\mathbf{M}_{13,1}$ respectively $\mathbf{M}_{13,6}$ according as 
$\varepsilon = 1$ or $-1$, and we let~$\mathbf{K}$ denote the standard 
$e$-split Levi subgroup of type $\mathbf{M}_{17,1}$, respectively 
$\mathbf{M}_{17,6}$, dual to~$\mathbf{K}^\dagger$. 
Notice that $\mathbf{K}$ and~$\mathbf{K}^\dagger$ are subgroups of $\mathbf{L} = 
C_{\mathbf{G}}( z_\C ) = \mathbf{L}^1 \circ_3 \mathbf{L}^2$; see 
Proposition~\ref{C3C}. In fact, $\mathbf{L}^1 = [\mathbf{K}^\dagger,\mathbf{K}^\dagger]$
and $\mathbf{L}^2 = [\mathbf{K},\mathbf{K}]$. Let~$\mathbf{T}$ denote the
standard $e$-split maximal torus of~$\mathbf{L}$, and put $\mathbf{T}_i = 
\mathbf{T} \cap \mathbf{L}^i$, for $i = 1, 2$.
Then $T_1 = Z( K )$ and $T_2 = Z( K^\dagger )$. As $s \in Z( K^\dagger)$, we have 
$T_2 \in \ker(\hat{s})$ by Lemma~\ref{DualityInT}, and as $T_1T_2$ has index~$3$
in~$T$ we may view $\hat{s}$ as an element of $\Irr( T_1 )$. As $\hat{s}$ has
order~$2$, it extends to a linear character of 
$K = T_1 \circ \langle L^2, x_{\C} \rangle$, with $\langle L^2, x_{\C} \rangle$
in its kernel; see Subsection~\ref{CentralProducts}.
Also $s \in T_2 = Z( K^\dagger ) \cong [q - \varepsilon]^2$ may be taken to be the 
element corresponding to $(-1,-1)$ under the latter isomorphism. By duality, 
$\hat{s} \in \Irr( T_1 ) = [q - \varepsilon]^2$ corresponds to the 
linear character sending a generator of each cyclic factor of 
$[q - \varepsilon]^2$ to~$-1$.

We adopt the notation of Subsections~\ref{sec3elts} and~\ref{secradLe}, as well 
as that of Lemma~\ref{lemRadLi}; the subgroups and elements of 
$\SL_3^\varepsilon(q)$ introduced in Subsection~\ref{secradLe} are 
provided with an index~$i$ here, indicating their association 
to~$L^i$, $i = 1, 2$; see also Subsection~\ref{secradF4Appendix}.
By \cite[Tables~$6$--$8$]{AH2}, we have 
$N_G( R ) = N_{L.2}( R )$ with $L.2 = \langle L, \gamma_{\C} \rangle$. 
(Recall that $R = D$ or $R =_G R_{25}$ at this stage of the proof.)
If $R = D$, then
\begin{equation}
\label{NGR33}
N_{L.2}( D ) = (\langle T_1, t_1 \rangle \circ_3 
\langle D_2, t_2 \rangle).\langle x_C, \gamma_C \rangle,
\end{equation}
and thus $N_G( D )$ stabilizes~$\hat{s}$ (recall that~$D_2$ denotes a particular
Sylow $3$-subgroup of $L^2$). Moreover,~$\hat{s}$ extends to 
$N_G( D ) \leq N_G( Q )$, as~$\hat{s}$ extends to $N_G( Q )$.
Now assume that $R =_G R_{25}$ and that $a \geq 2$. Then
\begin{equation}
\label{NGR25ag1}
N_{L.2}( R ) = (T_1.\langle t_1, c_1 \rangle \circ_3 3_+^{1+2}.\SL_2(3)).
\langle \gamma_C \rangle,
\end{equation}
and thus 
$$N_{G}( R, b_R ) = (T_1.\langle t_1 \rangle \circ_3 3_+^{1+2}.\SL_2(3)).
\langle \gamma_C \rangle,$$
as $c_1$ does not stabilize~$\hat{s}$. It follows that $\Out_G( R, b_R ) = 
(2 \times \SL_2(3)).2$. Clearly, $\hat{s}$ extends to 
$T_1.\langle t_1 \rangle \circ_3 3_+^{1+2}.\SL_2(3)$ such that 
$3_+^{1+2}.\SL_2(3)$ is in the kernel of this extension and 
$\langle \gamma_C \rangle$ stabilizes this extension. It follows that~$\theta_R = 
\hat{s}$ extends to $N_{G}( R, b_R )$.
If $R =_G R_{25}$ and $a = 1$, we have
\begin{equation}
\label{NGR25a1}
N_{L.2}( R ) = (T_1.\langle t_1 \rangle \circ_3 3_+^{1+2}.Q_8).
\langle x_C, \gamma_C \rangle,
\end{equation}
and thus $N_{G}( R, b_R ) = N_{G}( R )$, and $\Out_G( R, b_R ) = 
(2 \times \SL_2(3)).2$. As $\langle x_C, \gamma_C \rangle \cong S_3$ has trivial
Schur multiplier, $\theta_R = \hat{s}$ extends to $N_{G}( R, b_R )$.
%%%%%%%%%%%%%%%%%%%%%%%%%%%%%%%%%%%%%%%%%%%%%%%%%%%%%%%%%%%%%%%%%%%%%%%%%%%%%%%%
%%%%%%%%%%%%%%%%%%%%%%%%%%%%%%%%%%%%%%%%%%%%%%%%%%%%%%%%%%%%%%%%%%%%%%%%%%%%%%%%
%%
%% For more details see sheets 4 - 5 of my notes of 19.08.2019.
%%
%%%%%%%%%%%%%%%%%%%%%%%%%%%%%%%%%%%%%%%%%%%%%%%%%%%%%%%%%%%%%%%%%%%%%%%%%%%%%%%%
%%%%%%%%%%%%%%%%%%%%%%%%%%%%%%%%%%%%%%%%%%%%%%%%%%%%%%%%%%%%%%%%%%%%%%%%%%%%%%%
Finally assume that $s$ is of type $(3,1)$ if $\varepsilon = 1$ or of type 
$(3,2)$, if $\varepsilon = -1$. In the former case, $4 \mid q - 1$ and 
$s$ is contained in the center of a group of type $\mathbf{M}_{7,1}$. In the
latter case, $4 \mid q + 1$ and $s$ is contained in the center of a group of 
type $\mathbf{M}_{7,2}$. As above, we may view $\hat{s}$ as a linear character of 
$M_{9,1}$, respectively $M_{9,2}$, i.e.\ of a central product 
$\GL_2^\varepsilon(q) \circ \langle L^2, x_{\C} \rangle$, such that
$\langle L^2, x_{\C} \rangle$ is in the kernel of~$\hat{s}$. As $|s| = 4$, the 
linear character $\hat{s}$ of $\GL_2^\varepsilon(q) \leq L^1$ has order~$4$, and is
fixed by~$t_1$. On the other hand,~$\hat{s}$ is not fixed by $\gamma_C$, 
which induces the inverse-transpose automorphism on~$\GL^\varepsilon_2(q)$.
The claims for $R = D$ and $R =_G R_{25}$ now follow
from~(\ref{NGR33}),~(\ref{NGR25ag1}) and~(\ref{NGR25a1}).
%%%%%%%%%%%%%%%%%%%%%%%%%%%%%%%%%%%%%%%%%%%%%%%%%%%%%%%%%%%%%%%%%%%%%%%%%%%%%%%%
%%%%%%%%%%%%%%%%%%%%%%%%%%%%%%%%%%%%%%%%%%%%%%%%%%%%%%%%%%%%%%%%%%%%%%%%%%%%%%%%
%%
%% For more details see sheet 6 of my notes of 20.08.2019.
%%
%%%%%%%%%%%%%%%%%%%%%%%%%%%%%%%%%%%%%%%%%%%%%%%%%%%%%%%%%%%%%%%%%%%%%%%%%%%%%%%%
%%%%%%%%%%%%%%%%%%%%%%%%%%%%%%%%%%%%%%%%%%%%%%%%%%%%%%%%%%%%%%%%%%%%%%%%%%%%%%%
\end{prf}

\begin{cor}
\label{WeightSubgroups}
The non-cyclic weight subgroups of~$G$ are exactly the following (up to 
conjugation):
$$\{ R_8\text{--}R_{12}, R_{15}\text{--}R_{26}, R_{29}\text{--}R_{38} \}.$$
Of these, $R_{15}$, $R_{21}$, $R_{22}$ and $R_{35}$--$R_{38}$, occur only as 
weight subgroups of the principal block.
\end{cor}

\noindent
We finally deal with the exceptional double cover of $F_4(2)$.

\begin{prop}
Let $G := F_4( 2 )$ and $\hat{G} := 2.G$, the exceptional double cover of~$G$.
Let $\ell \in \{ 3, 5, 7 \}$, and let~$b$ be an $\ell$-block of~$\hat{G}$
with a non-cyclic defect group containing faithful characters. Then the invariants 
of~$b$ given in {\rm Columns}~$8$--$11$ of {\rm Table~\ref{21}} are correct.
\end{prop}
\begin{prf}
Assume first that $\ell = 3$. We identify the $3$-subgroups of~$\hat{G}$ and~$G$ 
via the canonical map $\hat{G} \rightarrow G$. A $3$-subgroup~$R$ of~$\hat{G}$ 
is radical, if and only if its image in~$G$ is radical. The radical 
$3$-subgroups of~$G$ are given in Table~\ref{tabradicalageq2}.

We use the permutation representation of~$\hat{G}$ on $139\,776$ points from
Wilson's Atlas of Finite Group Representations \cite{WWWW} to compute
in~$\hat{G}$ with GAP, and thus determine $N_{\hat{G}}( R )/R$ for 
each radical $3$-subgroup~$R$ of $\hat{G}$. We find $N_{\hat{G}}( R )/R = 2^4$, 
$2 \times \GL_2(3)$, $2 \times \SL_3(3)$, $2 \times W(F_4)$, 
$2\dot{\ }[(Q_8\times Q_8)\colon\!S_3]$ and $2 \times D_8$ for $R =_G R_{38}$, 
$R_{37}$, $R_{15}$, $R_{18}$, $R_{21}$ and $R_{10}$, respectively.
As $2\dot{\ }[(Q_8\times Q_8)\colon\!S_3]$ has exactly two faithful
defect zero character, we obtain the numbers in Column~$11$.

In case~$\ell > 3$, the defect group of~$b$ is a Sylow $\ell$-subgroup of~$\hat{G}$,
and the invariants are easily computed with GAP.
\end{prf}

\section{The action of the automorphism group on the set of weights}
\label{ActionOnWeights}

Let~$p$,~$q = p^f$, $\mathbf{G}$ and~$F$ be as in Subsection~\ref{SetupF4}, so 
that $\mathbf{G}^F = G = F_4(q)$ is the simple Chevalley group of type~$F_4$. 
Let~$\ell$ be an odd prime with $\ell \neq p$. We aim to prove the inductive 
blockwise Alperin weight condition for~$G$ at~$\ell$; see 
Subsection~\ref{InductiveBlockwiseAlperin}.

Recall that $F = F_1^f$ for a Steinberg morphism of~$\mathbf{G}$ with
$\mathbf{G}^{F_1} = F_4( p )$. Let~$\sigma_1$ be the Steinberg morphism 
of~$\mathbf{G}$ defined in Subsection~\ref{Automorphisms} 
(i.e.\ $\sigma_1 = F_1$ if~$p$ is odd, and $\sigma_1^2 = F_1$ if $p = 2$). 
In this section, our reference torus~$\mathbf{T}_0$ is assumed to be
$\sigma_1$-stable and contained in a $\sigma_1$-stable Borel subgroup 
of~$\mathbf{G}$.

\subsection{Preliminaries}
\label{PreliminariesFinalSection}
As in Subsection~\ref{DefinitionOfF0}, we
put $f' = f$ if~$p$ is odd, and $f' = 2f$ if $p = 2$. Then 
$F = \sigma_1^{f'}$ and $\Out(G) = \langle G\sigma_1 \rangle$ is cyclic of 
order~$f'$. Unless $q = 2$, the Schur multiplier of~$G$ is trivial, 
and hence any $\ell'$-covering group of~$G$ is equal to~$G$. Thus, by 
Remark~\ref{remBWC}, it suffices to verify the conditions of 
Hypothesis~\ref{defBWC} for every $\ell$-block of~$G$. An equivalent set of
conditions is formulated in Hypothesis~\ref{defBWC_Simplified}. The case of the
exceptional covering group of $F_4(2)$ will be dealt with separately.

Let~$b$ denote an $\ell$-block of~$G$. We may and will assume that~$b$ has
non-cyclic defect groups, as the inductive blockwise Alperin weight 
condition is known to hold for blocks with cyclic defect groups by the work of 
Koshitani and Sp{\"a}th~\cite{KS16}. In Theorem~\ref{thmAWC},
we have established the blockwise Alperin weight condition, i.e.\ we have 
already verified Hypothesis~\ref{defBWC_Simplified}(1) for~$b$. It remains
to settle the equivariance condition Hypothesis~\ref{defBWC_Simplified}(2).
As this is trivial if~$b$ has a unique irreducible Brauer character, we will
also assume $|l(b)| > 1$.

Put $A := N_{\Aut(G)}(b)$. Then $A = G \rtimes \langle \sigma \rangle$
for some $\sigma = \sigma_1^{m'}$, where~$m'$ is a positive integer with 
$m' \mid f'$. In fact $\langle \sigma \rangle = 
\Stab_{\langle \sigma_1 \rangle}( b )$. Even more can be said. Suppose that
$s \in G$ is a semisimple $\ell'$-element such that 
$b \subseteq \cE_\ell( G, s )$. Then, except in case $\ell \mid q^2 + 1$
the block~$b$ is uniquely determined in $\cE_\ell( G, s )$ by its defect group.
Thus, unless $\ell \mid q^2 + 1$, as a consequence 
of~(\ref{FPrimeAndLusztigSeries}) in Proposition~\ref{ActionAutomorphisms}, 
we have $\Stab_{\langle \sigma_1 \rangle}( b ) = 
\Stab_{\langle \sigma_1 \rangle}( [s] )$, where $[s]$ denotes the 
$G$-conjugacy class of~$s$.

Notice that every element of $A \setminus G$ naturally extends to a Steinberg 
morphism of~$\mathbf{G}$ which commutes with~$F$. Indeed, an element of 
$A \setminus G$ is of the form $\CJxxx_g \circ \sigma^j$ for some $g \in G$ and 
some positive integer~$j$, where $\sigma$ is tacitly assumed to be restricted 
to~$G$. We may now view $\CJxxx_g \circ \sigma^j$ as a Steinberg morphism 
of~$\mathbf{G}$. In particular, we may view~$A$ as a subgroup of 
$\Aut_1( \mathbf{G} )$.

For easy reference, we list a set of recurring assumptions.
\addtocounter{thm}{1}
\begin{hypo}
\label{HypoInvariance}
{\rm 
Let $b$ be an $\ell$-block of~$G$ with non-cyclic defect groups and with
$l(b) > 1$.

Put $A = N_{\Aut( G )}( b )$ and define $\sigma \in \langle \sigma_1 \rangle$ by 
$A = G \rtimes \langle \sigma \rangle$. Let~$m'$ be the positive integer with
$m' \mid f'$ and $\sigma = \sigma_1^{m'}$. Following the conventions of
Subsection~\ref{DefinitionOfF0}, we put $m := m'/2$ if $m'$ and~$p$ are 
even and $m := m'$, otherwise. Then $\sigma = F_1^m$ if either~$m'$ and~$p$ are
even or if~$p$ is odd. If $p = 2$ and $m'$ is odd, we have $m = m'$ and $m \mid f$.
In this case $\sigma^2 = F_1^m$; see Subsection~\ref{DefinitionOfF0}.

Let~$s$ be a $\sigma$-stable semisimple $\ell'$-element such that
$b \subseteq \cE_\ell( G, s )$. (Such an~$s$ exists by
{\rm Proposition~\ref{ActionAutomorphisms}} and
{\rm Lemma~\ref{SigmaStableConjugate}}.)\hfill{$\Box$}
}
\end{hypo}

\medskip
\noindent
Let~$D$ be a defect group of~$b$.
According to Hypothesis~\ref{defBWC_Simplified}, we have to consider $N_A( D )$. 
If~$D$ is $\sigma$-invariant then $N_A( D ) = 
N_G ( D ) \rtimes \langle \sigma \rangle$, a case which is somewhat easier to 
handle. To deal with the general case, observe that $N_A( D )/N_G( D )$ embeds 
into $\langle \sigma \rangle$, and hence $N_A( D ) = 
\langle N_G( D ), \omega \rangle$ for some $\omega \in N_A( D )$. In particular,
$C_{\mathbf{G}}( D )$ is $\omega$-invariant. The next lemma explores this 
situation.

\begin{lem}
\label{TwistingAndNormalizers}
Assume {\rm Hypothesis~\ref{HypoInvariance}} and let~$D$ be a defect group 
of~$b$. Let $\sigma' \in A$ be such that 
$A = G \rtimes \langle \sigma' \rangle$.
Suppose that~$\mathbf{M}$ is a regular subgroup of~$\mathbf{G}$ such that~$D$ 
is the Sylow $\ell$-subgroup of~$Z( M )$. In particular,~$D$ is abelian. 
Then $N_G( D ) = N_G( \mathbf{M} ) = N_G( M )$ and $N_A( D ) = N_A( \mathbf{M} ) 
= N_A( M )$.
	
Let~$\mathbf{M}_0$ denote a $\sigma_1$-stable regular subgroup of~$\mathbf{G}$ 
containing~$\mathbf{T}_0$. Suppose that $\mathbf{M} = \mathbf{M}_0^g$ for some 
$g \in \mathbf{G}$ such that $\dot{w} := F(g)g^{-1}$ is $\sigma'$-stable and 
$\dot{w} \in N_{\mathbf{G}}( \mathbf{T}_0 ) \cap N_{\mathbf{G}}( \mathbf{M}_0 )$.
	
Then $\CJxxx_g$ induces an isomorphism between $N_A( M )$ and 
$N_{\mathbf{G}}( \mathbf{M}_0 )^{F\dot{w}} \rtimes \langle \sigma' \rangle$.
In particular
$$N_A( M ) = N_G( M ) \rtimes \langle \omega \rangle$$ 
with $\omega = \CJxxx_g \circ \sigma' \circ \CJxxx_g^{-1}$. 
Hence $N_A(D) = \langle N_G(D), \omega \rangle$.
\end{lem}
\begin{prf}
As $C_{\mathbf{G}}( D )$ is a regular subgroup of~$\mathbf{G}$ by 
Lemma~\ref{CentralizersAbelianRadicalSubgroups} and 
\cite[Proposition~$2.3$(a)]{MalleAb}, we have $C_{\mathbf{G}}( D ) 
= \mathbf{M}$ and thus obtain $N_G( D ) = N_G( \mathbf{M} ) = N_G( M )$. 
As~$\mathbf{G}$ is normalized by~$A$ and the elements of~$A$ commute with~$F$, 
we also get $N_A( M ) \leq N_A( D ) \leq N_A( C_{\mathbf{G}}( D ) ) = 
N_A( \mathbf{M} ) \leq N_A( M )$. 

Put $\mathbf{A} := \mathbf{G} \rtimes \langle \sigma \rangle =
\mathbf{G} \rtimes \langle \sigma' \rangle \leq \Aut_1( \mathbf{G} )$. 
Observe that~$\mathbf{A}$ is normal in $\Aut_1(\mathbf{G})$ and that~$A$ 
consists exactly of those elements of~$\mathbf{A}$ which commute with~$F$. Hence
\begin{eqnarray*}
N_A( \mathbf{M} ) 
& = & \{ \omega \in N_{\mathbf{A}}( \mathbf{M} ) \mid \omega \text{\ commutes with\ } F \}\\
& = & \{ \CJxxx_g^{-1} \circ \omega_0 \circ \CJxxx_g \mid \omega_0 \in N_{\mathbf{A}}( \mathbf{M}_0 ),
\omega_0 \text{\ commutes with\ } F \dot{w}\}.
\end{eqnarray*}
As $\sigma'$ stabilizes $\mathbf{M}_0$, we have $N_{\mathbf{A}}( \mathbf{M}_0 ) = 
N_{\mathbf{G}}( \mathbf{M}_0 ) \rtimes \langle \sigma' \rangle$. As~$\sigma'$ 
commutes with $F \dot{w}$, an element $\CJxxx_h \circ (\sigma')^j \in 
N_{\mathbf{G}}( \mathbf{M}_0 ) \rtimes \langle \sigma' \rangle$ with 
$h \in N_{\mathbf{G}}( \mathbf{M}_0 )$ and $j \in \mathbb{Z}$ commutes 
with~$F\dot{w}$, if and only if 
$h \in N_{\mathbf{G}}( \mathbf{M}_0 )^{F\dot{w}}$.
%%%%%%%%%%%%%%%%%%%%%%%%%%%%%%%%%%%%%%%%%%%%%%%%%%%%%%%%%%%%%%%%%%%%%%%%%%%%%%%%
%%%%%%%%%%%%%%%%%%%%%%%%%%%%%%%%%%%%%%%%%%%%%%%%%%%%%%%%%%%%%%%%%%%%%%%%%%%%%%%%
%%
%% See sheets 2, 3 of 03.08.2020.
%%
%%%%%%%%%%%%%%%%%%%%%%%%%%%%%%%%%%%%%%%%%%%%%%%%%%%%%%%%%%%%%%%%%%%%%%%%%%%%%%%%
%%%%%%%%%%%%%%%%%%%%%%%%%%%%%%%%%%%%%%%%%%%%%%%%%%%%%%%%%%%%%%%%%%%%%%%%%%%%%%%%
\end{prf}

\medskip
\noindent
In case~$A$ acts trivially on $\IBr( b )$, we have to show that $N_A( D )$ 
fixes all the conjugacy classes of weights associated to~$b$. In the next 
subsection, we prepare for the necessary arguments to establish these results.

\addtocounter{subsection}{2}
\subsection{Preparations}
\label{PreparationsFinalSection}
Assume Hypothesis~\ref{HypoInvariance}. Below, we will make use of the notation 
introduced in Subsection~\ref{BlocksAndWeights}. A centric $b$-Brauer 
pair~$(R,b_R)$ is called \textit{relevant}, if $\cW( R, b_R ) \neq 0$. In this 
case we say that the radical $\ell$-subgroup~$R$ of~$G$ is 
$b$-\textit{relevant}. This notion is also used for $G$-conjugacy classes of 
$b$-Brauer pairs and radical $\ell$-subgroups of~$G$. We first have the 
following easy observation, which is the basis of our proof. 

\addtocounter{thm}{1}
\begin{lem}\label{remAWC} 
Let $(R,b_R)$ be a relevant $b$-Brauer pair and suppose that $\omega \in A$ 
stabilizes $(R, b_R)$. If, in addition,~$\omega$ fixes every element in 
$\Irr^0( N_G( R, \theta_R ) \mid \theta_R )$, then~$\omega$ stabilizes each 
$b$-weight $(R, \varphi)$ with $\varphi$ lying above~$\theta_R$.

Recall that by {\rm Theorem~\ref{thmAWC}(b)} the canonical character~$\theta_R$
of~$b_R$ extends to~$N_G(R, b_R)$. If~$\omega$ stabilizes one such extension, as
well as each $\xi \in \Irr^0(\Out_G( R, b_R))$, then~$\omega$ fixes every 
element in $\Irr^0( N_G( R, \theta_R ) \mid \theta_R )$.
\end{lem}
\begin{prf}
This is straightforward. Since $\omega$ stabilizes $(R, b_R)$, it also 
stabilizes $C_G(R)$, $N_G(R)$, $\theta_R$ and $N_G( R, \theta_R ) = 
N_G( R, b_R )$. If 
$(R, \varphi)$ is a $b$-weight with $\varphi$ lying above~$\theta_R$, then 
$\varphi = \Ind_{N_G( R, \theta_R )}^{N_G( R )}( \zeta )$ for some $\zeta \in 
\Irr^0( N_G( R, \theta_R ) \mid \theta_R )$. By assumption,~$\zeta$, and 
hence~$\varphi$ is fixed by $\omega$.

Let $\hat{\theta}_R$ denote an extension of~$\theta_R$ to $N_G( R, \theta_R )$.
The elements of $\Irr^0( N_G( R, \theta_R ) \mid \theta_R )$ are of the form
$\hat{\theta}_R\xi$ for $\xi \in \Irr^0( \Out(  R, b_R ) )$, and thus fixed 
by~$\omega$, if $\hat{\theta}_R$ and $\xi$ are $\omega$-stable.
\end{prf}

\medskip
\noindent
In the following lemma we use the classification of the defect groups and their 
centralizers; see Tables~\ref{1}--\ref{19} and~\ref{tabradicalageq2}. 

\begin{lem}
\label{lemBAWsubpairPreparation}
Let~$D$ be a defect group of~$b$. Assume that $|D|$ equals the 
$\ell$-part of $|C_G( s )|$. Let $\omega \in A$ be such that $N_A( D ) = 
\langle N_G( D ), \omega \rangle$. If~$D$ is abelian, put~$Q := D$, otherwise 
let~$Q$ denote the maximal abelian normal subgroup of~$D$ exhibited in 
{\rm Lemma~\ref{MaximalAbelianSubgroups}}.
	
Let $\mathbf{M} := C_{\mathbf{G}}( Q )$. Then~$\mathbf{M}$ is an
$\omega$-stable regular subgroup of~$\mathbf{G}$.
Suppose that $\mathbf{M}^\dagger \leq \mathbf{G}$ is a
regular subgroup with $s \in \mathbf{M}^\dagger$ and such that 
$(\mathbf{M},F)$ and $(\mathbf{M}^\dagger,F)$ are in duality.
	
View~$\omega$ as a Steinberg morphism of~$\mathbf{G}$ in the natural way; 
then $\omega$ commutes with~$F$. Suppose further that $\omega^\dagger$ is a 
Steinberg morphism of~$\mathbf{G}$ which commutes with~$F$, 
fixes~$\mathbf{M}^\dagger$ and~$s$, and that $(\mathbf{M}, \omega)$
and $(\mathbf{M}^\dagger, \omega^\dagger)$ are in duality.

Then there is a maximal $b$-Brauer pair $(D,b_D)$ with $\omega(b_D) = b_D$;
moreover, $b_D \subseteq \cE_\ell( C_G( D ), s )$ if~$D$ is abelian.
\end{lem}
\begin{prf}
As $(\mathbf{M}, \omega)$ and $(\mathbf{M}^\dagger, \omega^\dagger)$ are in 
duality, $\omega|_{\mathbf{M}}$ and $\omega^\dagger|_{\mathbf{M}^\dagger}$ are
dual isogenies in the sense of~\cite[Proposition~$11.1.11$]{DiMi2}.
As~$s$ is $\omega^\dagger$-stable and $\mathbf{M}$ is $\omega$-stable, 
\cite[Corollary~$9.3$(ii)]{DiMiPara} implies that $\cE( M, s )$ is 
$\omega$-stable, too. 

By the description of the defect groups in Tables~\ref{1}--\ref{19},
either~$\mathbf{M}$ is a maximal torus, or else $e \in \{ 1, 2 \}$ 
and~$\mathbf{M}$ is an $e$-split Levi subgroup of~$\mathbf{G}$ which is 
$G$-conjugate to one of $\mathbf{M}_{i,k}$ with 
$k \in \{ 13, 14, 15, 17, 18, 19 \}$ (and $k$ is such that $\mathbf{M}_{i,k}$ 
is $e$-split). Notice that the cases $i = 13, 17$ only occur for $\ell > 3$.

Let $(D,b_D')$ denote a maximal $b$-Brauer pair, and let $(Q,b_Q')$ be the 
$b$-Brauer pair with $(Q,b_Q') \leq (D,b_D')$. If $t \in M^\dagger$ such 
that $b_Q' \subseteq \cE_\ell( M, t )$, then~$s$ and~$t$ are conjugate in~$G$
by Lemma~\ref{ConjugateSemisimpleLabels}. Now consider $\mathbf{K}^\dagger 
:= C_{\mathbf{M}^\dagger}( t )$.

By going through the various cases for~$s$ and~$\mathbf{M}^\dagger$, and using
Table~\ref{cc}, we find that any two $F$-stable maximal tori 
of~$\mathbf{M}^\dagger$ which have some $G$-conjugate in 
$C_{\mathbf{G}}( t )$ are conjugate in~$M^\dagger$. In particular, 
$\mathbf{K}^\dagger$ is a maximal 
torus of $\mathbf{M}^\dagger$, and $t \in K^\dagger$ is in general position
with respect to~$M^\dagger$. Notice that the class type of~$t$ determines the 
type of~$\mathbf{K}^\dagger$. As~$s$ and~$t$ are conjugate in~$G$, it follows 
that $\mathbf{S}^\dagger := C_{\mathbf{M}^\dagger}( s )$ is an $F$-stable maximal
torus of~$\mathbf{M}^\dagger$, and $\mathbf{K}^\dagger$ is conjugate in~$M^\dagger$
to~$\mathbf{S}^\dagger$. Conjugation inside~$M^\dagger$ does not affect
$\cE_\ell( M , t )$, so that we may assume
that $t \in \mathbf{S}^\dagger$. Choose an $F$-stable maximal torus~$\mathbf{S}$ 
of~$\mathbf{M}$ dual to~$\mathbf{S}^\dagger$. Now $\theta_t := 
\pm R_{\mathbf{S}}^{\mathbf{M}}( \hat{t} ) \in \cE( M, t )$ is a character of~$M$
of central defect, and $\Res^M_{C_G(D)}( \theta_t )$ is the canonical
character of~$b_D'$. This follows from Lemma~\ref{PossibleCanonicalCharacters} 
in case $\ell = 3$, but the case $\ell > 3$ is proved analogously. Put
$\theta_s := \pm R_{\mathbf{S}}^{\mathbf{M}}( \hat{s} ) \in \cE( M, s )$,
let $b_Q \subseteq \cE_\ell( M, s )$ denote the $\ell$-block of~$M$ containing~$\theta_s$ 
and let~$b_D$ denote the $\ell$-block of $C_G( D )$ containing 
$\Res_{C_G(D)}^M( \theta_s )$. Then $(D,b_D)$ is a centric Brauer pair with 
$(Q,b_Q) \leq (D,b_D)$. As $\theta_s$ is the unique character in~$\cE( M, s )$,
it follows that $(Q,b_Q)$ is $\omega$-invariant. In turn $(D, b_D)$ is
$\omega$-invariant.
	
Notice that $\Out_G( D, b_D' ) = \Out_G( Q, b_Q' )$ and $\ell \nmid 
|\Out_G( D, b_D' )|$, as $(D,b_D')$ is a maximal $b$-Brauer pair. We also have 
$\Out_G( D, b_D ) = \Out_G( Q, b_Q ) \leq \Out_G( M )$. 
To prove that $(D,b_D)$ is a 
$b$-Brauer pair, we claim that $\ell \nmid |\Out_G( D, b_D )|$. This is clear if 
$\mathbf{M}$ is not a torus or if $\ell > 3$, as $\ell \nmid |\Out_G( M )|$ in 
these cases. Suppose that $\ell = 3$ and that $\mathbf{M} = \mathbf{S}$ is a 
torus. The argument in the proof of Proposition~\ref{cormain}(a) gives
$\Out_G( Q, b_Q ) \cong W_{C_{\mathbf{G}}( s )}( \mathbf{S}^\dagger)$ and
$\Out_G( Q, b_Q' ) \cong W_{C_{\mathbf{G}}( t )}( \mathbf{S}^\dagger)$.
Now~$s$ and~$t$ are conjugate in $W_{\mathbf{G}}( \mathbf{S}^\dagger )$;
see \cite[Proposition~$3.7.1$]{C2}. If $\mathbf{S}^\dagger$ is the
maximally $e$-split torus of~$\mathbf{G}$, then $W_{\mathbf{G}}( \mathbf{S}^\dagger )
= W_{\mathbf{G}}( \mathbf{S}^\dagger )^F$, and thus 
$W_{C_{\mathbf{G}}( s )}( \mathbf{S}^\dagger) = 
\{ w \in W_{\mathbf{G}}( \mathbf{S}^\dagger ) \mid s^w = s \}$ and
$W_{C_{\mathbf{G}}( t )}( \mathbf{S}^\dagger)$ are conjugate in~$G$. The only 
case remaining is when
$Q \in_G \{ R_{11}, R_{12} \}$, and hence~$s$ is of class type
$(13,4)$, $(13,5)$, $(17,4)$, $(17,5)$, $(18,3)$, $(18,7)$, $(19,8)$ or 
$(19,9)$. In these cases, $W_{C_{\mathbf{G}}( s )}( \mathbf{S}^\dagger)$ has 
order~$2$, as $S^\dagger = M^\dagger \cong 
[q - \varepsilon] \times [q^3 - \varepsilon]$ in these cases. This proves our 
claim.

Now let~$b'$ denote the $\ell$-block of~$G$ such that $(1,b') \leq (Q,b_Q)$.
By Lemma~\ref{ConjugateSemisimpleLabels}, we have 
$b' \subseteq \cE_\ell( G, s )$. As~$b'$ is the Brauer correspondent of~$b_D$,
a defect group of~$b'$ equals~$D$. As~$b$ is the unique block in 
$\cE_\ell( G, s )$ with defect group~$D$, it follows that $b = b'$.
\end{prf}

\addtocounter{subsection}{2}
\subsection{Proofs for the non-unipotent blocks}
\label{ProofsForNonUnipotentBlocks}

Let~$b$,~$\sigma$ and~$s$ be as in Hypothesis~\ref{HypoInvariance}. In addition, 
suppose that~$b$ is non-unipotent, and that $\sigma = F_1^m$. In this subsection
we also put $G_m := \mathbf{G}^\sigma$. 

%%%%%%%%%%%%%%%%%%%%%%%%%%%%%%%%%%%%%%%%%%%%%%%%%%%%%%%%%%%%%%%%%%%%%%%%%%%%%%%%
%%%%%%%%%%%%%%%%%%%%%%%%%%%%%%%%%%%%%%%%%%%%%%%%%%%%%%%%%%%%%%%%%%%%%%%%%%%%%%%%
%%
%% Outline of strategy
%%
Let us outline the strategy of proof in the generic case, namely when the order 
of a defect group of~$b$ equals the $\ell$-part of $|C_G( s )|$. Recall that 
this condition holds exactly when the entry corresponding to~$b$ in Column~$5$ 
of Tables~\ref{1}--\ref{19} equals~$1$. Starting from $s$ and $\sigma$, we 
construct a pair of regular subgroups~$\mathbf{M}$ and $\mathbf{M}^\dagger$ 
of~$\mathbf{G}$ such that the Sylow $\ell$-subgroup~$D$ of $Z(M)$ is a defect 
group of~$b$. Simultaneously, we construct a pair~$\omega$ and $\omega^\dagger$
of Steinberg morphisms of $\mathbf{G}$ such that $N_A( D ) = 
\langle N_G( D ), \omega \rangle$, and such that
the hypotheses of Lemma~\ref{lemBAWsubpairPreparation} are satisfied. This 
yields an $\omega$-stable maximal $b$-Brauer pair $(D,b_D)$. We then have to
investigate the action of~$\omega$ on~$N_G( D, b_D )$ and~$\Out_G( D, b_D )$
in order to apply Lemma~\ref{remAWC}. This is done by describing 
$\Out_G( D, b_D )$ in terms of suitable structures inside the Weyl group~$W$.
It turns out that, in most of the cases, $\omega$ acts trivially on 
$\Out_G( D, b_D )$.

To construct~$\mathbf{M}$ and~$\mathbf{M}^\dagger$, we exhibit a pair $y, y^*
\in \mathbf{G}$ such that $\dot{u} := F(y)y^{-1}$ and $\dot{u}^* := 
F(y^*){y^*}^{-1}$ normalize~$\mathbf{T}_0$, and that the images~$u$ and~$u^*$ of
$\dot{u}$, respectively~$\dot{u}^*$ in $W$ are related by $u^* = {u^\dagger}^{-1}$. 
We also choose a suitable element $v \in W$ which determines the $G_m$-class 
type of~$s$. Putting $\mathbf{S} := \mathbf{T}_0^y$ and $\mathbf{S}^\dagger := 
\mathbf{T}_0^{y^*}$, and defining $\omega$ and $\omega^\dagger$ accordingly,
we obtain a pair of commutative diagrams
$$
\begin{xy}
\xymatrix@C+25pt{
\mathbf{S} \ar[r]^{\omega} & \mathbf{S} \\
	\mathbf{T}_0 \ar[u]^{\CJxxx_y} \ar[r]_{\sigma \dot{v}} &
\mathbf{T}_0 \ar[u]_{\CJxxx_{y}}
}
\end{xy}\quad\quad
\begin{xy}
\xymatrix@C+25pt{
\mathbf{S}^\dagger \ar[r]^{\omega^\dagger} & \mathbf{S}^\dagger \\
\mathbf{T}_0 \ar[u]^{\CJxxx_{y^*}} \ar[r]_{\sigma \dot{v}} &
\mathbf{T}_0 \ar[u]_{\CJxxx_{y^*}}
}
\end{xy}
$$
where~$\dot{v} \in N_{\mathbf{G}}( \mathbf{T}_0 )$ is a suitable lift
of~$v$ (whose choice only becomes relevant in the diagrams below).
Notice that the $F$-stable maximal tori~$\mathbf{S}$ and $\mathbf{S}^\dagger$
are in duality by the facts summarized in Subsection~\ref{TwistingAndDuality}.
We choose our data such that the Sylow $\ell$-subgroup~$D$ of~$S$ is a defect 
group of~$b$, and that $s \in S^\dagger$ is $\omega^\dagger$-stable. Let 
$D^\dagger$ denote the Sylow $\ell$-subgroup 
of~$S^\dagger$ and put $\mathbf{M} := C_{\mathbf{G}}( D )$ and 
$\mathbf{M}^\dagger := C_{\mathbf{G}}( D^\dagger )$. We obtain the following
pair of commutative diagrams
$$
\begin{xy}
\xymatrix@C+25pt{
\mathbf{M} \ar[r]^{\omega} & \mathbf{M} \\
	\mathbf{M}_0 \ar[u]^{\CJxxx_y} \ar[r]_{\sigma \dot{v}} &
\mathbf{M}_0 \ar[u]_{\CJxxx_{y}}
}
\end{xy}\quad\quad
\begin{xy}
\xymatrix@C+25pt{
\mathbf{M}^\dagger \ar[r]^{\omega^\dagger} & \mathbf{M}^\dagger \\
	\mathbf{M}_0 \ar[u]^{\CJxxx_{y^*}} \ar[r]_{\sigma \dot{v}} &
\mathbf{M}_0 \ar[u]_{\CJxxx_{y^*}}
}
\end{xy}
$$
for suitable standard Levi subgroups $\mathbf{M}_0$ and $\mathbf{M}_0^\dagger$
such that $(\mathbf{M}_0,F)$ and $(\mathbf{M}_0^\dagger,F)$ are in duality. We 
may now apply Lemmas~\ref{TwistingAndNormalizers} 
and~\ref{lemBAWsubpairPreparation} to find $N_A( D ) = 
\langle N_G( D ), \omega \rangle$ and to get an $\omega$-stable $b$-Brauer 
pair~$(D,b_D)$.

To describe $N_G( D, b_D )$ and the action of~$\omega$ on $\Out_G( D, b_D )$,
let $\mathbf{L}_0$ denote a standard Levi subgroup such that
$\mathbf{L}^* := \mathbf{L}_0^{y^*} = C_{\mathbf{G}}( s )$. Then 
$\mathbf{M}^\dagger \cap \mathbf{L}^* = \mathbf{S}^\dagger$ in our situation, 
and we have the following commutative diagram
$$
\begin{xy}
\xymatrix@C+1pt{
	N_{\mathbf{G}}( \mathbf{M} )^F \ar[r]^{\omega} &
N_{\mathbf{G}}( \mathbf{M} )^F \ar[r]^{\kappa} & W_{\mathbf{G}}( \mathbf{M} )^F \ar[r]^{\alpha}
	& W_{\mathbf{G}}( \mathbf{M}^\dagger )^F & W_{\mathbf{L}^*}( \mathbf{S}^\dagger )^F \ar[l]_{\iota} \\
	N_{\mathbf{G}}( \mathbf{M}_0 )^{F\dot{u}} \ar[u]^{\CJxxx_{y}} \ar[r]_{\sigma \dot{v}} &
	N_{\mathbf{G}}( \mathbf{M}_0 )^{F\dot{u}} \ar[u]^{\CJxxx_{y}} \ar[r]_{\kappa} &
	W_{\mathbf{G}}( \mathbf{M}_0 )^{F\dot{u}} \ar[u]^{\CJxxx_{y}} \ar[r]_{\dagger} &
	W_{\mathbf{G}}( \mathbf{M}_0^\dagger )^{F\dot{u}} \ar[u]_{\CJxxx_{y^*}} &
	W_{\mathbf{L}_0}( \mathbf{T}_0 )^{F\dot{u}} \ar[u]_{\CJxxx_{y^*}} \ar[l]^\iota %&
}
\end{xy}
$$
where $\alpha$ is an isomorphism, $\kappa$ denotes a canonical epimorphism 
and~$\iota$ an embedding. By Lemma~\ref{lemBAWsubpairPreparation}, we have
$N_G( D ) = N_G( \mathbf{M} )$. Moreover, by Proposition~\ref{cormain}, the 
image in $W_{\mathbf{G}}( \mathbf{M}^\dagger )^F$ of $N_G( D, b_D )$ under 
$\alpha \circ \kappa$ equals the image of 
$W_{\mathbf{L}^*}( \mathbf{S}^\dagger )^F$ under~$\iota$. Thus $N_G( D, b_D ) = 
\kappa^{-1}(\alpha^{-1}( \iota( W_{\mathbf{L}^*}( \mathbf{S}^\dagger )^F ) ) )$.
We may now transfer the computations to the bottom line of the above
diagram. The inverse image of $W_{\mathbf{L}^*}( \mathbf{S}^\dagger )^F$ 
under $\CJxxx_{y^*}$ equals 
$C_W( u^* ) \cap W_{\mathbf{L}_0}( \mathbf{T}_0 ) \leq W$. The inverse
image under~$\dagger$ of the latter group equals 
$C_W( u ) \cap W_{\mathbf{L}_0}( \mathbf{T}_0 )^\dagger$. Thus 
$\Out_G( D, b_D ) \cong C_W( u ) \cap W_{\mathbf{L}_0}( \mathbf{T}_0 )^\dagger 
\leq W_{\mathbf{G}}( \mathbf{M}_0 )^{F\dot{u}}$. We now have to investigate 
the action of $\sigma \dot{v}$ on 
$C_W( u ) \cap W_{\mathbf{L}_0}( \mathbf{T}_0 )^\dagger$
and on $\kappa^{-1}( C_W( u ) \cap W_{\mathbf{L}_0}( \mathbf{T}_0 )^\dagger )$.

We basically follow this outline in Proposition~\ref{AbelianD} below, where many 
more details are worked out. The proof of Proposition~\ref{lemBAWsubpair} also
proceeds along these lines, but with fewer details. In 
Proposition~\ref{AbelianD}, the 
element~$u$ of the above considerations has the form $v^{f/m}w$, where~$v$ 
and~$w$ are suitable elements from Table~\ref{CT}. Our approach requires the 
derivation of further properties of these elements,
which is done in Lemma~\ref{ProofCT2} below. The crucial group 
$C_W( u ) \cap W_{\mathbf{L}_0}( \mathbf{T}_0 )^\dagger$ corresponds to the
group $\bar{C}_0$ of Definition~\ref{DefineM} below. This finishes the outline 
of our proofs.
%% End of outline
%%
%%%%%%%%%%%%%%%%%%%%%%%%%%%%%%%%%%%%%%%%%%%%%%%%%%%%%%%%%%%%%%%%%%%%%%%%%%%%%%%%
%%%%%%%%%%%%%%%%%%%%%%%%%%%%%%%%%%%%%%%%%%%%%%%%%%%%%%%%%%%%%%%%%%%%%%%%%%%%%%%%

Recall that~$e$ denotes the order of~$q$ modulo~$\ell$. In accordance with our 
previous usage, we put $\varepsilon := 1$ if $e = 1$, and $\varepsilon := -1$ 
if $e = 2$. Recall that $\dagger$ stands for the ``duality'' isomorphism of~$W$; 
see Subsection~\ref{Duality}.

\addtocounter{thm}{1}
\begin{dfn}
\label{DefineM}
{\rm 
Let $\Gamma \subseteq \Sigma$, $i, k \in \mathbb{Z}$, $e \in \{ 1, 2 \}$ and 
$v, w \in W$ such that some row of {\rm Table~\ref{CT}} has the values 
$(\Gamma^\dagger, (i,k) , e, -, v^\dagger, -, w^\dagger, - , -)$. Put 
	$$\bar{C}_0 := C_{W_\Gamma}( {v^{f/m}}w ).$$
If $e > 2$, put $\mathbf{M}_0^\dagger := \mathbf{T}_0$.
If $e \leq 2$, let $\mathbf{T}_0'$ denote the subtorus of $\mathbf{T}_0$ on which 
$F{v^\dagger}^{f/m}w^\dagger$ acts as $t \mapsto t^{\varepsilon q}$ and and put 
$\mathbf{M}_0^\dagger := C_{\mathbf{G}}( \mathbf{T}_0' )$.

Then $\mathbf{M}_0^\dagger = \mathbf{L}_{\Delta^\dagger}$ for a
subset $\Delta \subseteq \Sigma^+$ which is $W$-conjugate to a subset
of $\{ \alpha_1, \ldots , \alpha_4 \}$. Assume that $|\Delta| \leq 2$
and put $\mathbf{M}_0 := \mathbf{L}_\Delta$.\hfill{$\Box$}
}
\end{dfn}

\medskip

\noindent
We have excluded the case $|\Delta| = 3$ in the above definition, as this
corresponds to cyclic defect groups. Recall the definition of $W_{\Gamma}$
from Subsection~\ref{SetupF4}, namely $W_{\Gamma} = 
\langle s_{\alpha} \mid \alpha \in \Gamma \rangle$. In particular, 
$W_{\Gamma}^\dagger = W_{\Gamma^\dagger}$. Notice that $\Gamma^\dagger$ is a 
base of the closed subsystem $\overline{\Gamma^\dagger}$ of~$\Sigma$ by the
conventions adopted in Table~\ref{CT} and thus $W_{\Gamma^\dagger} = 
W_{\overline{\Gamma^\dagger}}$. However, $\Gamma$ need not be a base of 
$\bar{\Gamma}$, in which cases we have $W_{\Gamma} \lneq W_{\bar{\Gamma}}$.

\begin{lem}
\label{ProofCT2}
Assume the notation of {\rm Definition~\ref{DefineM}}. Recall that $\hat{W}$ 
denotes the subgroup of $N_{\mathbf{G}}( \mathbf{T}_0 )$ generated by 
$n_1, \ldots , n_4$; see {\rm Subsection~\ref{LiftOfLongestElementLem}}. 
Then there is a subgroup $C_0 \leq \hat{W}$ which maps to $\bar{C}_0$ under
the natural epimorphism $\hat{W} \rightarrow W$, and the following statements 
hold. 
	
{\rm (a)} The group $C_0$ normalizes~$\mathbf{M}_0$. 
	
{\rm (b)} Suppose that $f/m$ is even or that the commutator factor group of 
$\bar{C}_0$ does not have order~$2$ or~$4$. Then~$C_0$ centralizes 
$[\mathbf{M}_0,\mathbf{M}_0]$ and there are inverse 
images $\dot{v}, \dot{w} \in \hat{W}$ of~$v$ and~$w$ respectively, such that 
the following conditions are satisfied.
	
\begin{itemize}
\item[{\rm (i)}] If $\Delta \neq \emptyset$, then 
$\dot{v} \in [\mathbf{M}_0,\mathbf{M}_0]$
unless $w \neq 1$ and $(i,k)$ is one of $(18,6), (18,4), (19,6), (19,7)$.

\item[{\rm (ii)}] Unless 
$(i,k) \in \{ (12,2), (12,4), (16,k) \mid k \in \{ 3, 4, 7, 8, 10 \} \}$,
the element $\dot{v}$ centralizes~$C_0$.

\item[{\rm (iii)}] If $w \neq 1$ then $\dot{w} \in Z(C_0)$. 
\end{itemize}

Suppose that $w \neq 1$ and that $\Delta \neq \emptyset$.
Then $e > 1$. In these cases,
the first and second column of {\rm Table~\ref{CentralizerLiftsI}} give the 
parameters~$i$ and~$\Delta$, the third and 
fourth column display generators for~$\bar{C}_0$ and~$C_0$, respectively.

\begin{table}[h]
\caption{\label{CentralizerLiftsI} Lifts of some centralizers}
$$
\begin{array}{rccc} \hline\hline
i & \Delta & \bar{C}_0 & C_0 \rule[- 6pt]{0pt}{ 19pt} \\ \hline \hline
3 & 22 & s_1, s_4, s_4s_{17}, (s_4s_{17})^{s_3} &
        n_1, n_4n_3^2, n_4n_{17}, (n_4n_{17})^{n_3}n_{22}n_3^2  \rule[- 3pt]{0pt}{ 16pt}  \\
4 & 18, 10 & s_5, s_{19} & n_5, n_{19}n_3^2 \\
6 & 24 & s_2, s_3, s_4 & n_2, n_3, n_4 \\
7 & 20, 19 & s_1, s_7 & n_1, n_7n_3^2 \\
8 & 1, 13 & s_4, s_4s_{17}, (s_4s_{17})^{s_3} & n_4n_3^2, n_4n_{17}, (n_4n_{17})^{n_3}n_{22}n_3^2 \\
9 & 23, 15 & s_5, s_4 & n_5, n_4n_3^2 \\
10 & 21 & s_1, s_2, s_3 & n_1, n_2 , n_3n_4^2 \\
11 & 8 & s_2, s_3, s_{21} & n_2, n_3n_4^2, n_{21}n_4^2 \\
12 & 14 & s_1, s_3, s_{21} & n_1, n_3n_4^2, n_{21} \\
14 & 22, 17 & s_1, s_4 & n_1, n_4n_3^2 \\
15 & 8, 16 & s_2, s_3 & n_2, n_3n_4^2 \\
16 & 1, 6 & s_3, s_{21} & n_3n_4^2, n_{21}n_1^2 \\
18 & 17 & s_4 & n_4n_3^2 \\
19 & 22 & s_1 & n_1
\rule[- 2pt]{0pt}{ 5pt} \\ \hline\hline
\end{array}
$$
\end{table}

Suppose that $w = 1$ and that $\Delta \neq \emptyset$. Then $e = 1$. In these
cases, 
the first three columns of {\rm Table~\ref{CentralizerLiftsII}} give the 
parameters $i, k, \Delta$, the fourth and fifth column display generators 
for~$\bar{C}_0$ and~$C_0$, respectively.

\begin{table}[h]
\caption{\label{CentralizerLiftsII} Lifts of some centralizers (cont.)}
$$
\begin{array}{rrccc} \hline\hline
	i & k & \Delta & \bar{C}_0 & 
	C_0 \rule[- 6pt]{0pt}{ 19pt} \\ \hline \hline
3 & 2 & 22 & s_1, s_4, s_4s_{17}, (s_4s_{17})^{s_3} &
        n_1, n_4n_3^2, n_4n_{17}, (n_4n_{17})^{n_3}n_{22}n_3^2  \rule[- 3pt]{0pt}{ 16pt}  \\
11 & 2 & 8 & s_2, s_3, s_{21} & n_2, n_3n_4^2, n_{21}n_4^2 \\
12 & \neq 1 & 14 & s_1, s_3, s_{21} & n_1, n_3n_4^2, n_{21} \\
13 & 4, 5 & 1, 23 &  s_3, s_4 & n_3, n_4 \\
16 & \neq 1 & 1, 6 & s_3, s_{21} & n_3n_4^2, n_{21}n_1^2 \\
17 & 4, 5 & 4, 21 & s_1, s_2 & n_1, n_2 \\
18 & 4, 6 & 1, 13 & s_4 & n_4n_3^2 \\
18 & 3, 7 & 1, 2 & s_4 & n_4 \\
19 & 6, 7 & 4, 14 & s_1 & n_1 \\
19 & 8, 9 & 3, 4 & s_1 & n_1
\rule[- 2pt]{0pt}{ 5pt} \\ \hline\hline
\end{array}
$$
\end{table}
{\rm (c)} The elements of~$C_0$ as well as the elements $\dot{v}$ and~$\dot{w}$
constructed in {\rm (b)} are $\sigma$-stable.
\end{lem}
\begin{prf}
The explicit computations described below are performed with CHEVIE.

By construction, $C_W( {v^\dagger}^{f/m} w^\dagger )$ normalizes
$W_{\mathbf{M}_0^\dagger}( \mathbf{T}_0 )$. Hence $\bar{C}_0
\leq C_W( {v^\dagger}^{f/m} w^\dagger )^\dagger$ normalizes 
$W_{\mathbf{M}_0^\dagger}( \mathbf{T}_0 )^\dagger = 
W_{\mathbf{M}_0}( \mathbf{T}_0 )$. Thus, in any case, there is 
$C_0 \leq \hat{W}$ normalizing~$\mathbf{M}_0$ and mapping to~$\bar{C}_0$.
However, to meet the conditions in~(b), we need specific lifts, in 
particular those indicated in Tables~\ref{CentralizerLiftsI} 
and~\ref{CentralizerLiftsII}, so that we have to check~(a) for these cases.
CHEVIE can be used to show that the elements given in the respective 
columns of Tables~\ref{CentralizerLiftsI} and~\ref{CentralizerLiftsII} do
indeed generate $\bar{C}_0$. The corresponding entries in the columns
headed by $C_0$ are lifts of these generators in~$\hat{W}$
and normalize $\mathbf{M}_0 = 
\mathbf{L}_{\Delta}$ (a trivial fact if $\Delta = \emptyset$,
in which case $\mathbf{M}_0 = \mathbf{T}_0$).

(b) If~$\Delta \neq \emptyset$, we choose~$C_0$ as in 
Tables~\ref{CentralizerLiftsI} and~\ref{CentralizerLiftsII}. 
If $w = 1$, put $\dot{w} := 1$.
Suppose that $\Delta = \emptyset$ and that $w \neq 1$. Then $e \geq 2$
by Table~\ref{CT}. If $e = 2$, either $w^\dagger = w_0$ or
${v^\dagger}^{f/m} \neq 1$. In the former case, which occurs exactly 
for $i \in \{ 2, 5 \}$, we also have
$w_0 \in \bar{C}_0$, and we choose $C_0$ such that it
contains the lift $\gamma$ of~$w_0$. In the latter case, which only occurs for
$i \in \{ 15, 16 \}$, we take $C_0$ as in the row corresponding to~$i$ of
Table~\ref{CentralizerLiftsI}. In case $e \in \{ 3, 6 \}$, the group
$\bar{C}_0$ has odd
order, so that we may lift it to an isomorphic group~$C_0$ in the extended Weyl
group. In case $e = 4$, we put $C_0 := \langle n_3n_2, n_4n_3n_8n_4 \rangle$.
If $\Delta = \emptyset$ and $w = 1$, we define the groups $C_0$ in our 
proof of~(ii). In any case, $C_0 \leq \hat{W}$.

Recall that $\bar{C}_0  = C_{W_{\Gamma}}( v^{f/m}w ) = 
C_{\Gamma^\dagger}( {v^\dagger}^{f/m}w^\dagger )^\dagger$ in the notation of 
Remark~\ref{ProofCT}. Notice also that if $e \in \{ 1, 2 \}$ then $\Delta = 
\emptyset$ if and only if the corresponding entry in the ``cl''-column of 
Table~\ref{CT} equals~$1$ or~$2$. Indeed, these are exactly the cases in which
$\mathbf{T}_0' = \mathbf{T}_0$ in the notation of Definition~\ref{DefineM}.
In particular, the conditions in~(b) on the commutator factor group 
of~$\bar{C}_0$ and the non-emptiness of~$\Delta$ are easily checked from 
the entries of Table~\ref{CT}. Assume that~$\Delta \neq \emptyset$.
To find $\Delta^\dagger$ and hence~$\Delta$, we determine 
the $\varepsilon$-eigenspace~$Y'$ of ${v^\dagger}^{f/m}w^\dagger$ on 
$Y = Y( \mathbf{T}_0 )$. We then choose~$\Delta^\dagger$ as a base of the 
subsystem of~$\Sigma$ of roots perpendicular to~$Y'$; then $\mathbf{M}_0^\dagger 
= \mathbf{L}_{\Delta^\dagger}$. This yields the corresponding entries for~$\Delta$
in Tables~\ref{CentralizerLiftsI} and~\ref{CentralizerLiftsII}.

To show that~$C_0$ commutes with 
$[\mathbf{M}_0,\mathbf{M}_0]$, we may assume $\Delta \neq \emptyset$.
Using the formulas in \cite[Lemmas~$7.2.1$(i), $6.4.4$(i)]{C1}, one checks that 
the generators of $C_0$ given in Tables~\ref{CentralizerLiftsI} 
and~\ref{CentralizerLiftsII} indeed
centralize the root subgroups corresponding to the roots in $\pm\Delta$.

Observe that~$v$ is given in Table~\ref{CT} as a product of~$s_j$'s and~$w_0$. We 
lift each $s_j$ to~$n_j$ and $w_0$ to~$\gamma$, and we lift~$v$ to the 
corresponding product~$\dot{v}$ of the $n_j$s and~$\gamma$. Then $\dot{v} \in 
\hat{W}$.

(i) With the above choice of~$\dot{v}$, we check $\dot{v} \in 
[\mathbf{M}_0,\mathbf{M}_0]$ from Tables~\ref{CT},~\ref{CentralizerLiftsI} 
and~\ref{CentralizerLiftsII}, except in the 
cases listed in the statement.

(ii) Suppose first that $\Delta \neq \emptyset$. Then the claim follows from (i) 
and the fact that~$C_0$ centralizes $[\mathbf{M}_0,\mathbf{M}_0]$ which we have 
already proved, except in the cases 
excluded in~(i). If $(i,k) = (18,6)$, we have $v = s_{13}s_1$, and one checks 
that $n_{13}n_1$ commutes with $n_4n_3^2$. The same argument applies for
$(i,k) = (18,4)$ if we lift~$w_0$ to~$\gamma$. If $(i,k) = (19,6)$, we have
$v = s_{14}s_4$, and one checks that~$n_{14}n_4$ commutes with~$n_1$. The same
argument applies for $(i,k) = (19,7)$ if we lift~$w_0$ to~$\gamma$.

Suppose that $\Delta = \emptyset$ and $w \neq 1$. Then $v \in \{ 1, w_0 \}$ 
unless $i \in \{ 15, 16 \}$. In the former case, our claim clearly holds. In the
latter case, it follows from our choice of $C_0$ and from what we have already 
proved above. Suppose finally that $\Delta = \emptyset$ and $w = 1$. Then 
$e = 1$. If $v \in \{ 1, w_0 \}$, there is nothing to prove. Suppose that there 
is an element $w \neq 1$ such that 
$(\Gamma^\dagger, (i,k) , -, -, v^\dagger, -, w^\dagger, - , -)$ also is a line 
in Table~\ref{CT} with the same value of~$f/m$ and with corresponding $\Delta$ 
satisfying $|\Delta| \leq 2$ and such that 
$C_{\Gamma^\dagger}( {v^\dagger}^{f/m}w^\dagger )^\dagger = 
W_{\Gamma}$. Then the claim follows from what we 
have already proved, as $C_{\Gamma^\dagger}( {v^\dagger}^{f/m} )^\dagger \leq 
W_{\Gamma}$. We now discuss the remaining cases. If $(i,k) = (8,2)$ and 
$v = s_1$, we lift $\bar{C}_0$ to the group generated by $n_3$ and the group
given in the row corresponding to $i = 8$ in Table~\ref{CentralizerLiftsI}. As 
$n_1$ commutes with~$n_3$, we obtain our claim. The analogous argument works
for $(i,k) = (8,3)$ and $v = w_0s_4$. If $(i,k) \in \{ (13,2), (13,3), (13,4), (13,5) \}$,
we have $\bar{C}_0 = \langle s_3, s_4 \rangle$, which we lift to $C_0 := 
\langle n_3, n_4 \rangle$. As $n_1$ and $n_{23}$ commute with each of~$n_3$ and~$n_4$,
our claim follows. The analogous argument works in case $i = 17$. In the cases
occurring for $i = 18$, we lift $\bar{C}_0$ to $\langle n_4n_3^2 \rangle$. Then~$C_0$
centralizes $n_1$, $n_2$ and $n_{13}$, giving our claim. For $i = 19$,
we lift $\bar{C}_0$ to $\langle n_1 \rangle$. As this commutes with $n_3$, $n_4$ and 
$n_{14}$, we are done.

(iii) If $i \in \{ 2, 5 \}$, we have $w_0 \in \bar{C}_0$, 
and our claim follows from Lemma~\ref{LiftOfLongestElementLem}. In the other 
cases,  we determine the order of 
$Z( C_0 )$ as well as of the intersection of $Z( C_0 )$ with the normal subgroup 
of the extended Weyl group spanned by $n_j^2$, $j = 1, \ldots , 4$. This shows 
that in each case $Z( C_0 )$ projects onto $Z( \bar{C}_0 )$;
as the latter group contains~$w$, this proves our claim. 

(c) This is clear as~$\hat{W}$ consists of~$\sigma$-stable elements.
\end{prf}

\medskip
\noindent
We are now ready to obtain the main structural results about centralizers of
abelian defect groups and the corresponding inertia groups.

\begin{prop}
\label{AbelianD}
Let~$b$,~$\sigma$ and~$s$ be as described at the beginning of this subsection 
and assume in addition that~$b$ has abelian defect groups. Then there is a 
maximal $b$-Brauer pair $(D,b_D)$ and an element
$\omega \in N_A(D)$ with $N_A(D) = \langle N_G( D ), \omega \rangle$ such
that~$\omega$ fixes~$b_D$. Moreover,~$\omega$ centralizes $\Out_G( D, b_D )$
except in the cases listed in~{\rm (a)} and~{\rm (b)} below.

Suppose from now on that $f/m$ is even or that the commutator quotient of 
$\Out_G( D, b_D )$ does not have order~$2$ or~$4$.
Put $\mathbf{M} := C_{\mathbf{G}}( D )$. Then there is $N'' \leq N_G( D, b_D )$
with $Z(M) \leq N'' \leq C_G( [\mathbf{M}, \mathbf{M}] )$ and 
$N'' \cap M = Z(M)$ such that $N_G( D, b_D ) = N''M$. In particular, 
$\Out_G( D, b_D ) \cong N''/Z(M)$. Moreover, every coset of 
$N''/Z(M)$ contains an $\omega$-stable element, except in cases~{\rm (a)} 
and~{\rm (b)} below. Finally, if $\mathbf{M}$ is not a torus,~$\omega$ acts as a 
field automorphism on $M' = [\mathbf{M},\mathbf{M}]^F$. 

{\rm (a)} The $G_m$-class type of~$s$ is one of $(12,2)$ or $(12,4)$ and~$f/m$ 
is even. In this case, $\Out_G( D, b_D ) \cong 2^3$, and the group of 
$\omega$-fixed points on $\Out_G( D, b_D )$ has order~$4$.

{\rm (b)} The $G_m$-class type of~$s$ is $(16,k)$ for 
$k \in \{ 3, 4, 7, 8, 10 \}$ and $f/m$ is even. In all these cases, 
$\Out_G( D, b_D ) \cong 2^2$, and the group of $\omega$-fixed points
on $\Out_G( D, b_D )$ has order~$2$.

In each of the cases listed in~{\rm (a)} and~{\rm (b)}, let $\lambda \in 
\Irr(Z( G ))$ be defined by $\Res_{Z(M)}^M( \theta_D )
= \lambda \theta_D( 1 )$, where $\theta_D$ denotes the canonical character 
of~$b_D$. Then there is  an $\omega$-stable extension of~$\lambda$ to~$N''$. 
\end{prop}
\begin{prf}
First, we deal with the case that the order of a defect group of~$b$ is smaller 
than the $\ell$-part of $|C_G( s )|$. In this case, $e \in \{ 1, 2 \}$ and the 
defect groups of~$b$ are $G$-conjugate to $R_{10,\ell}$; see Table 
\ref{1}--\ref{19}. In turn, the centralizers in~$\mathbf{G}$ of the defect 
groups are $G$-conjugate to $\mathbf{M}_{15,1}$ if $e = 1$ and to 
$\mathbf{M}_{15,3}$ if $e = 2$. Moreover, the $\mathbf{G}$-class type of~$s$ 
equals~$i$ with $i \in \{ 2, 5, 6, 10, 11 \}$. Let $g \in \mathbf{G}$ be such 
that $\sigma(g){g}^{-1} = \gamma$, with 
$\gamma \in N_{\mathbf{G}}( \mathbf{T}_0)$ as in 
Lemma~\ref{LiftOfLongestElementLem}. 
Put $\mathbf{L}^* := \mathbf{M}_{i,1}$ if the $G_m$-class type of~$s$ equals 
$(i,1)$ and $e = 1$, and $\mathbf{L}^* := \mathbf{M}_{i,1}^{g}$, otherwise.
Then $\mathbf{M}_{i,1}^{g} =_{G_m} \mathbf{M}_{i,1}$ for $i \in \{ 2, 5 \}$.
If $i \in \{ 6, 10, 11 \}$ and $e = 2$, the $G$-class type of~$s$ equals 
$(i,2)$, as otherwise~$b$ would have cyclic defect. But then the $G_m$-class
type of~$s$ must be $(i,2)$ as well.
We may thus assume that $\mathbf{L}^* = C_{\mathbf{G}}( s )$. 
Now $\mathbf{M}_{i,1}$ contains 
$\mathbf{M}_{15,1}$ as a subgroup of maximal rank. Put $\mathbf{M} := 
\mathbf{M}_{15,1}$ if the $G_m$-class type of~$s$ equals $(i,1)$ and $e = 1$, 
and $\mathbf{M} := \mathbf{M}_{15,1}^{g}$, otherwise. Then~$\mathbf{M}$ is 
$\sigma$-stable, $s\in Z( M )$ and $D := O_\ell( Z( M ) )$ 
is a defect group of~$b$. Notice that we can identify~$\mathbf{M}$ 
with $\mathbf{M}^\dagger$ and~$D$ with~$D^\dagger$ in this case. As~$D$
is $\sigma$-stable, we obtain $N_A( D ) = \langle N_G( D ), \sigma \rangle$.
Let $b_D \subseteq \cE_\ell( M , s )$ denote the $\ell$-block of~$M$
which covers the $\ell$-block of $M' = \Sp_4(q)$ containing the unipotent 
$\ell$-defect zero character~$\psi$ of $M'$. Now $\cE_\ell( M , s )$ is 
$\sigma$-stable by~(\ref{FPrimeAndLusztigSeries}), and thus~$b_D$ is 
$\sigma$-stable, as $b_D$ is uniquely determined by~$s$ and~$\psi$.
Thus $(D, b_D)$ is a centric, 
$\sigma$-stable Brauer pair with $\ell \nmid |\Out_G( D, b_D )|$.
Let $b'$ denote the $\ell$-block of~$G$ with $(1, b') \leq (D,b_D)$.
Then~$b'$ has defect group~$D$, and lies in $\cE_{\ell}( G, s )$
by Lemma~\ref{ConjugateSemisimpleLabels}, and thus $b' = b$ by uniqueness.
Let $N' \leq N_G( \mathbf{M} )$ be as in Proposition~\ref{NormalizerM18}. 
As~$\mathbf{T}_0$ is $1$-$\sigma$-split, this proposition, applied with~$\sigma$ 
instead of~$F$, shows that every coset of $N'/Z( M )$ contains a $\sigma$-stable 
element. Putting $N'' := N' \cap N_G( D, b_D )$, we obtain the claims about
the structure of~$N_G( D, b_D )$. Clearly, $\omega = \sigma$ acts
as a field automorphism on~$M'$. We have now proved all assertions in case
the order of a defect group of~$b$ is smaller than the $\ell$-part of $|C_G(s)|$.

Suppose now that the order of a defect group of~$b$ is equal to the $\ell$-part 
of~$|C_G( s )|$. We begin by choosing a pair $\mathbf{S}$, $\mathbf{S}^\dagger$ 
of $F$-stable maximal tori of~$\mathbf{G}$ and a pair of elements 
$h^*, h \in \mathbf{G}$ with the following properties.

(i) The pairs $(\mathbf{S},F)$ and $(\mathbf{S}^\dagger,F)$ are in duality. 

(ii) The Steinberg morphisms $\omega := \CJxxx_h^{-1} \circ \sigma \circ \CJxxx_h$
and $\omega^\dagger := \CJxxx_{h^*}^{-1} \circ \sigma \circ \CJxxx_{h^*}$
commute with~$F$ and stabilize~$\mathbf{S}$ and $\mathbf{S}^\dagger$, respectively.
Moreover, the pairs $(\mathbf{S}, \omega)$ and $(\mathbf{S}^\dagger, \omega^\dagger)$ 
are in duality.

(iii) The element~$s$ lies in $S^\dagger := {\mathbf{S}^\dagger}^F$ and is 
fixed by~$\omega^\dagger$.

(iv) Let~$D$ and $D^\dagger$ denote the Sylow $\ell$-subgroup of~$S$, 
respectively $S^\dagger$. Then~$D$ is a defect group of~$b$.

(v) Put $\mathbf{M} := C_{\mathbf{G}}( D )$ and $\mathbf{M}^\dagger := 
C_{\mathbf{G}}( D^\dagger )$. Then $\mathbf{M}$ and $\mathbf{M}^\dagger$ are 
$\omega$-stable, respectively $\omega^\dagger$-stable regular subgroups of 
$\mathbf{G}$, such that $(\mathbf{M}, F)$ and $(\mathbf{M}^\dagger, F)$ as well as
$(\mathbf{M}, \omega)$ and $(\mathbf{M}^\dagger, \omega^\dagger)$  are in duality. The 
corresponding duality isomorphism
$W_{\mathbf{G}}( \mathbf{M} ) \stackrel{\alpha}{\longrightarrow}
W_{\mathbf{G}}( \mathbf{M}^\dagger )$ satisfies 
$\alpha \circ \omega = \omega^\dagger  \circ \alpha$ and $\alpha \circ F = 
F \circ \alpha$. 

To choose the pairs $\mathbf{S}$, $\mathbf{S}^\dagger$ and $h^*, h$, we proceed 
as follows.  First, we write $u^* := {u^\dagger}^{-1}$ for $u \in W$.
Suppose that the $G_m$-class type of~$s$ equals $(i,k)$ and let $\mathbf{L}_0 
:= \mathbf{M}_{i,1}$, with $\mathbf{M}_{i,1} = \mathbf{L}_{\Gamma_i}$ as in 
Subsection~\ref{ClassTypes}. For each~$e$, choose a pair $v, w \in W$ such that 
$(v^*, w^*)$ is as given in the row of Table~\ref{CT} corresponding to the tuple 
$(i,k,e,f/m)$.  In particular, $v^*$ and~$w^*$ satisfy the conditions described 
in Remark~\ref{ProofCT}, so that $w^* \in W_{\Gamma_i} = 
W_{\mathbf{L}_0}( \mathbf{T}_0 )$. 
Choose inverse images $\dot{v}, \dot{v}^*, \dot{w} \in 
N_{\mathbf{G}}( \mathbf{T}_0 )$ and 
$\dot{w}^* \in N_{\mathbf{L}_0}( \mathbf{T}_0 )$ of $v, v^*, w, w^*$, 
respectively, with the following properties. The elements $\dot{v}$ and 
$\dot{v}^*$ are $\sigma$-stable, $\dot{w}$ is $\sigma \dot{v}$-stable and 
$\dot{w}^*$ is $\sigma \dot{v}^*$-stable. The latter two conditions hold 
automatically if $\dot{w}$ and $\dot{w}^*$ are $\sigma$-stable and commute with 
$\dot{v}$, respectively $\dot{v}^*$. In any case, these conditions can be met
as~$v^*$ and~$w^*$ commute. We tacitly assume that the lifts 
$\dot{v}$ and $\dot{w}$ satisfy, whenever applicable, the stronger properties 
listed in Lemma~\ref{ProofCT2}. Choose elements $g, g^* \in \mathbf{G}$ 
with $\sigma(g)g^{-1} = \dot{v}$ and $\sigma(g^*){g^*}^{-1} = \dot{v}^*$.
Next, choose elements $h \in \mathbf{G}$ and
$h^* \in \mathbf{L}_0^{g^*}$ with $F(h)h^{-1} = g^{-1}\dot{w}g$ and
$F(h^*){h^*}^{-1} = {g^*}^{-1}\dot{w}^*g^*$. Take each of the
elements $g, h, g^*, h^*$ as~$1$, if the corresponding elements in~$W$ equal~$1$.

Putting $\mathbf{L}^* := \mathbf{L}_0^{g^*}$, we may and will now assume that 
$C_{\mathbf{G}}( s ) = \mathbf{L}^*$ (recall that $\sigma$-twisting 
$\mathbf{L}_0 = \mathbf{M}_{i,1}$ with $v^*$ yields a group in the 
$G_m$-conjugacy class of $\mathbf{M}_{i,k}$). Let $\mathbf{S}^\dagger := 
\mathbf{T}_0^{g^*{h^*}}$ and $\mathbf{S} := \mathbf{T}_0^{gh}$. 
Put $u := v^{f/m}w$, $u^* := {v^*}^{f/m}w^*$, $\dot{u} := 
{\dot{v}}^{f/m}\dot{w}$ and $\dot{u}^* := {(\dot{v}^*)}^{f/m}\dot{w}^*$.
Notice that $u^* = {u^\dagger}^{-1}$ as~$v$ and~$w$ commute.
Now $F(gh)(gh)^{-1} = F(g) {g}^{-1}\dot{w} = {\dot{v}}^{f/m}\dot{w} = \dot{u}$, 
and thus~$\mathbf{S}$ is an $F$-stable maximal torus of $F$-type 
$(\emptyset,[u])$. Similarly, $\mathbf{S}^\dagger$ is an 
$F$-stable maximal torus of $F$-type $(\emptyset,[u^*])$.
This implies~(i) by \cite[Proposition~$4.3.4$]{C2}. By 
construction,~$\omega$ stabilizes~$\mathbf{S}$, and $\omega$ commutes with~$F$ 
since $\sigma$ fixes $F(h)h^{-1} = g^{-1}\dot{w}g$. The analogous argument 
applies for $\omega^\dagger$. An elementary computation shows the duality
of $(\mathbf{S}, \omega) $ and $(\mathbf{S}^\dagger, \omega^\dagger)$;
see Subsection~\ref{TwistingAndDuality}. Thus~(ii) holds.
Property~(iii) is clear by the choice of $h^* \in C_{\mathbf{G}}( s )$.
The fact that~$D$ is a defect group of~$b$ follows from Tables~\ref{1}--\ref{19},
in conjunction with Table~\ref{SylowsInCentralizersI}. To prove~(v), first
notice that $\mathbf{M}$ is $\omega$-stable, since~$D$ is. The same argument 
applies to $\mathbf{M}^\dagger$. Next, observe 
that $\mathbf{M} = \mathbf{S}$, unless $e \in \{ 1, 2\}$. In these 
cases,~$\mathbf{M}$ is the centralizer of the subtorus of~$\mathbf{S}$,
on which~$F$ acts as $t \mapsto t^{\varepsilon q}$. Thus $\mathbf{M}
= \mathbf{M}_0^{gh}$, where $\mathbf{M}_0$ is the centralizer of the subtorus
of~$\mathbf{T}_0$, on which $F u$ acts as $t \mapsto t^{\varepsilon q}$
(recall that $\varepsilon = 1$ if $e = 1$, and $\varepsilon = -1$ if $e = 2$).
The analogous construction applies to~$\mathbf{M}^\dagger$, where $F u$
is replaced by $F u^*$. We obtain the following commutative
diagram of relative Weyl groups,
$$
\begin{xy}
\xymatrix@C+25pt{
W_{\mathbf{G}}( \mathbf{M} ) \ar[r]^{\omega} & W_{\mathbf{G}}( \mathbf{M} ) \ar[r]^{\alpha}
	& W_{\mathbf{G}}( \mathbf{M}^\dagger ) \ar[r]^{\omega^\dagger} & W_{\mathbf{G}}( \mathbf{M}^\dagger ) \\
       W_{\mathbf{G}}( \mathbf{M}_0 ) \ar[u]^{\CJxxx_{gh}} \ar[r]_{\sigma v} & 
W_{\mathbf{G}}( \mathbf{M}_0 ) \ar[u]_{\CJxxx_{gh}} \ar[r]_{\dagger}
        & W_{\mathbf{G}}( \mathbf{M}_0^\dagger ) \ar[u]_{\CJxxx_{g^* h^*}}
          \ar[r]_{\sigma {v^*}}
        & W_{\mathbf{G}}( \mathbf{M}_0^\dagger ) \ar[u]_{\CJxxx_{g^* h^*}}
}
\end{xy}
$$
where~$\alpha$ is defined such that the middle square commutes. The duality of 
$(\mathbf{M}_0,F)$ and $(\mathbf{M}_0^\dagger, F)$ gives rise to the 
isomorphism~$\dagger$ and implies the duality of $(\mathbf{M},F)$ and
$(\mathbf{M}^\dagger,F)$; see Subsection~\ref{TwistingAndDuality}. Similarly, the 
duality of $(\mathbf{M}_0,\sigma)$ and 
$(\mathbf{M}_0^\dagger, \sigma)$ implies the duality of $(\mathbf{M}_0, \sigma \dot{v})$
and $(\mathbf{M}_0^\dagger, \sigma \dot{v}^*)$ and thus of $(\mathbf{M},\omega)$ and
$(\mathbf{M}^\dagger, \omega^\dagger)$. This gives all the claims of~(v).

By Lemma~\ref{ProofCT2}(b), the element $\dot{w}$ is $\sigma$-stable 
if~$f/m$ is even. If $f/m$ is odd, and $v = \dot{v} = 1$, then~$\dot{w}$ is
$\sigma$-stable as it is $\sigma \dot{v}$-stable by our choice. As~$b$ has 
non-cyclic defect, no other cases occur. Thus $F(gh)(gh)^{-1} = 
\dot{v}^{f/m}\dot{w}$ is $\sigma$-stable and hence $N_A( D ) = 
\langle N_G( D ), \omega \rangle$ by Lemma~\ref{TwistingAndNormalizers}.
Lemma~\ref{lemBAWsubpairPreparation} implies the existence of 
an $\omega$-stable maximal $b$-Brauer pair $(D,b_D)$ with
$b_D \subseteq \cE_\ell( M, s )$.
Consider the following commutative diagram:
$$
\begin{xy}
\xymatrix@C+30pt{
N_{\mathbf{L}^* }( \mathbf{S}^\dagger )\ar[r]^{\omega^\dagger } &
N_{\mathbf{L}^* }( \mathbf{S}^\dagger ) \\
N_{\mathbf{L}_0}( \mathbf{T}_0 ) \ar[u]^{\CJxxx_{g^* h^* }} \ar[r]_{\sigma \dot{v}^*} &
N_{\mathbf{L}_0}( \mathbf{T}_0 ) \ar[u]_{\CJxxx_{g^* h^*}}
}
\end{xy}
$$
%%%%%%%%%%%%%%%%%%%%%%%%%%%%%%%%%%%%%%%%%%%%%%%%%%%%%%%%%%%%%%%%%%%%%%%%%%%%%%%%%%
%%%%%%%%%%%%%%%%%%%%%%%%%%%%%%%%%%%%%%%%%%%%%%%%%%%%%%%%%%%%%%%%%%%%%%%%%%%%%%%%%%
%%%
%%% See sheets 1 - 2 of 20.11.2020 and 3 - 4 of 02.12.2020.
%%%
%%%%%%%%%%%%%%%%%%%%%%%%%%%%%%%%%%%%%%%%%%%%%%%%%%%%%%%%%%%%%%%%%%%%%%%%%%%%%%%%%%
%%%%%%%%%%%%%%%%%%%%%%%%%%%%%%%%%%%%%%%%%%%%%%%%%%%%%%%%%%%%%%%%%%%%%%%%%%%%%%%%%%
%%%%%%%%%%%%%%%%%%%%%%%%%%%%%%%%%%%%%%%%%%%%%%%%%%%%%%%%%%%%%%%%%%%%%%%%%%%%%%%%%
%%%%%%%%%%%%%%%%%%%%%%%%%%%%%%%%%%%%%%%%%%%%%%%%%%%%%%%%%%%%%%%%%%%%%%%%%%%%%%%%%
%%%
%%% See sheets 1 - 3 of 23.11.2020.
%%%
%%%%%%%%%%%%%%%%%%%%%%%%%%%%%%%%%%%%%%%%%%%%%%%%%%%%%%%%%%%%%%%%%%%%%%%%%%%%%%%%%
%%%%%%%%%%%%%%%%%%%%%%%%%%%%%%%%%%%%%%%%%%%%%%%%%%%%%%%%%%%%%%%%%%%%%%%%%%%%%%%%%
We have 
$$\CJxxx_{g^* h^*}^{-1}( N_{\mathbf{L}^*}( \mathbf{S}^\dagger )^F ) =
N_{\mathbf{L}_0}( \mathbf{T}_0 )^{F\dot{u}^*}.$$
In particular, $(g^* h^*)^{-1}$ conjugates 
$W_{\mathbf{L}^*}( \mathbf{S}^\dagger )^F$ to 
$$N_{\mathbf{L}_0}( \mathbf{T}_0 )^{F\dot{u}^*}/\mathbf{T}_0^{F\dot{u}^*} = 
C_W( u^* ) \cap W_{\mathbf{L}_0}( \mathbf{T}_0 ).$$
By the definition of~$\mathbf{M}^\dagger$ we have 
$W_{\mathbf{L}^*}( \mathbf{S}^\dagger )^F \leq 
W_{\mathbf{G}}( {\mathbf{M}^\dagger} )^F$, so that 
$\alpha^{-1}(W_{\mathbf{L}^*}( \mathbf{S}^\dagger )^F) \leq 
W_{\mathbf{G}}( \mathbf{M} )^F$, as~$\alpha$ is $F$-equivariant.
The first commutative diagram above implies that $(gh)^{-1}$ 
conjugates $\alpha^{-1}(W_{\mathbf{L}^*}( \mathbf{S}^\dagger )^F )$ to 
$$\bar{C}_0 := C_W( u ) \cap W_{\mathbf{L}_0}( \mathbf{T}_0 )^\dagger.$$
Proposition~\ref{cormain}(a) implies that $\Out_G( D, b_D ) =
\alpha^{-1}(W_{\mathbf{L}^*}( \mathbf{S}^\dagger )^F)$.

Notice that $\bar{C}_0$ is exactly the group associated in Lemma~\ref{ProofCT2}
to the tuple $(i,k,e,f/m)$, as $W_{\mathbf{L}_0}( \mathbf{T}_0 ) = 
W_{\Gamma^\dagger}$ if $\mathbf{L}_0 = \mathbf{L}_{\Gamma^\dagger}$. The action 
of $\omega$ on $\Out_G( D, b_D )$ corresponds to the action of $\sigma v$ on 
$\bar{C}_0$. By Lemma~\ref{ProofCT2}(a), the latter action is trivial unless we 
are in one of the cases described in~(a) or~(b). This gives the second statement 
of our proposition, namely that~$\omega$ centralizes $\Out_G( D, b_D )$. In the 
situations of (a) and~(b) we have to compute the number of orbits of~$v$ 
on~$\bar{C}_0$. This is done with a CHEVIE computation confirming the 
assertions.

Now assume that $f/m$ is even or that the quotient of $\Out_G( D, b_D )$ by its 
commutator subgroup does not have order~$2$ or~$4$ (this includes the cases 
of~(a) and~(b)). Let~$C_0$ denote  the inverse image of~$\bar{C}_0$
exhibited in Lemma~\ref{ProofCT2}. In particular, $C_0$ consists of
$\sigma$-stable elements and centralizes $\dot{w}$ and 
$[\mathbf{M}_0,\mathbf{M}_0]$. Moreover, $C_0$ centralizes~$\dot{v}$ except in 
the cases listed in (a) and~(b). In the latter cases, $C_0$ still 
centralizes~$\dot{u}$. Thus, in all cases, 
$C_0 \leq N_{\mathbf{G}}( \mathbf{M}_0 )^{F\dot{u}}$.
By the considerations above, $N_G( D, b_D )$ is conjugate, under $(gh)^{-1}$,
to $C_0\mathbf{M}_0^{F\dot{u}}$. 
Now define $N_0'' := Z( M_0 ) C_0$ and
$N'' := (N_0'')^{gh}$. Then $N_G( D, b_D ) = N''M$, and $Z( M ) \leq N''
\leq C_{\mathbf{G}}( [ \mathbf{M}, \mathbf{M} ] )$. As $N''$ centralizes
$[\mathbf{M}, \mathbf{M}]$, we also have $N'' \cap M = Z(M)$
by Propositions~\ref{NormalizerM13}, \ref{NormalizerM14} and~\ref{NormalizerM18}.
Unless we are in either the situation~(a) or~(b), the elements of~$C_0$
are fixed by $\sigma \dot{v}$, and thus every coset of $N''/Z( M )$ contains an
$\omega$-stable element. Suppose that $\mathbf{M}_0$ is not a torus. 
Then $\dot{v} \in [\mathbf{M}_0,\mathbf{M}_0]^{F\dot{u}}$ by 
Lemma~\ref{ProofCT2}(b), unless 
the $G_m$-class type of~$s$ is one of $(18,6)$, $(18,4)$, $(19,6)$ or $(19,7)$.
In the latter cases, $[\mathbf{M}_0,\mathbf{M}_0]$ has semisimple rank~$1$,
and if $\alpha \in \Sigma^+$ is the positive root of $[\mathbf{M}_0,\mathbf{M}_0]$,
we have $u_{\pm \alpha}( t )^{\dot{v}} = u_{\pm \alpha}( -t )$ for all 
$t \in \mathbb{F}$; this follows from the corresponding entries in 
Table~\ref{CT} and \cite[Lemma~$7.2.1$]{C1}. As $f/m$ is even, either~$q$ is 
even or $4 \mid q - 1$, and thus~$\dot{v}$ acts
as an inner automorphism on $[\mathbf{M}_0,\mathbf{M}_0]^{F\dot{u}}$. 
In each case, $\sigma \dot{v}$ acts as a field automorphism on 
$[\mathbf{M}_0,\mathbf{M}_0]^{F\dot{u}}$, and hence~$\omega$
acts as a field automorphism on~$M'$. 

We finally prove the last assertion. Suppose first that~$s$ is as described in~(a).
We have $N'' = \langle Z( M ), n_{1}^{gh}, (n_{3}n_4^2)^{gh}, n_{21}^{gh} \rangle$
and that the elements $n_{1}^{gh}$, $(n_{3}n_4^2)^{gh}$, $n_{21}^{gh} $ are $\omega$-stable. 
One checks that $(n_3n_4^2)^2 = 1$, and that the kernel of the natural epimorphism
$\langle n_{1}, n_{3}n_4^2, n_{21} \rangle \rightarrow \langle s_1, s_3, s_{21} \rangle$
is generated by $n_1^2$ and $n_{21}^2$.
By \cite[Lemma~$6.4.4$]{C1} we have $n_i^2 = \alpha_i^\vee(-1)$ for all
$1 \leq i \leq 24$. 
By Lemma~\ref{KernelsOfDualCharacters}, the elements of~$S$ of the form
$\alpha_i^\vee( t )^{gh}$ for $t \in \mathbb{F}_q^*$ and
$i \in \{ 1, 3, 21 \}$, are contained in the kernel of~$\hat{s}$.
As $\lambda$ is the restriction of~$\hat{s}$ to $Z(M)$, and as $N''/Z(M)$
is elementary abelian, $\langle (n_1^2)^{gh}, (n_{21}^2)^{gh} \rangle \leq \ker( \lambda )$.
Thus $Z(M)/\ker( \lambda )$ has an $\omega$-stable complement
in $N''/\ker( \lambda )$. View~$\lambda$ as a character of $Z(M)/\ker( \lambda )$.
As such, it has a trivial extension to $N''/\ker( \lambda )$, and this is
$\omega$-invariant. Exactly the same argument applies for elements~$s$
as in~(b).
%%%%%%%%%%%%%%%%%%%%%%%%%%%%%%%%%%%%%%%%%%%%%%%%%%%%%%%%%%%%%%%%%%%%%%%%%%%%%%%%%%
%%%%%%%%%%%%%%%%%%%%%%%%%%%%%%%%%%%%%%%%%%%%%%%%%%%%%%%%%%%%%%%%%%%%%%%%%%%%%%%%%%
%%%
%%% See sheets 1 - 2 of 20.11.2020 and 3 - 4 of 02.12.2020.
%%%
%%%%%%%%%%%%%%%%%%%%%%%%%%%%%%%%%%%%%%%%%%%%%%%%%%%%%%%%%%%%%%%%%%%%%%%%%%%%%%%%%%
%%%%%%%%%%%%%%%%%%%%%%%%%%%%%%%%%%%%%%%%%%%%%%%%%%%%%%%%%%%%%%%%%%%%%%%%%%%%%%%%%%
%%%%%%%%%%%%%%%%%%%%%%%%%%%%%%%%%%%%%%%%%%%%%%%%%%%%%%%%%%%%%%%%%%%%%%%%%%%%%%%%%
%%%%%%%%%%%%%%%%%%%%%%%%%%%%%%%%%%%%%%%%%%%%%%%%%%%%%%%%%%%%%%%%%%%%%%%%%%%%%%%%%
%%%
%%% See sheets 1 - 3 of 23.11.2020.
%%%
%%%%%%%%%%%%%%%%%%%%%%%%%%%%%%%%%%%%%%%%%%%%%%%%%%%%%%%%%%%%%%%%%%%%%%%%%%%%%%%%%
%%%%%%%%%%%%%%%%%%%%%%%%%%%%%%%%%%%%%%%%%%%%%%%%%%%%%%%%%%%%%%%%%%%%%%%%%%%%%%%%%
%%%%%%%%%%%%%%%%%%%%%%%%%%%%%%%%%%%%%%%%%%%%%%%%%%%%%%%%%%%%%%%%%%%%%%%%%%%%%%%%
%%%%%%%%%%%%%%%%%%%%%%%%%%%%%%%%%%%%%%%%%%%%%%%%%%%%%%%%%%%%%%%%%%%%%%%%%%%%%%%%
%%
%% For the CHEVIE (GAP 3) computations used in this proof see
%% ~hiss/papers/alperinF4ppr/GapComputations/F4Lemma814.g
%%
%%%%%%%%%%%%%%%%%%%%%%%%%%%%%%%%%%%%%%%%%%%%%%%%%%%%%%%%%%%%%%%%%%%%%%%%%%%%%%%%
%%%%%%%%%%%%%%%%%%%%%%%%%%%%%%%%%%%%%%%%%%%%%%%%%%%%%%%%%%%%%%%%%%%%%%%%%%%%%%%%
\end{prf}

\medskip

\noindent
We also need the following result on the action of field automorphisms on certain
characters.

\begin{lem}
\label{FieldAutomorphism}
Let~$M'$ be one of the groups $\SL_2( q )$ or $\Sp_4( q )$ with~$q$ odd, 
or $\SL_3( q )$ with $3 \mid q - 1$. Let $d = 2$ in the first two cases, and
$d = 3$ in the last case.

Let $\nu \in \Irr( M' )$ which is not invariant under a diagonal 
automorphism~$\kappa$ of~$M'$, and let $\nu_1, \ldots, \nu_d$ denote the 
distinct conjugates of~$\nu$ under~$\kappa$.
In case of $M' = \SL_3( q )$, assume that
$\nu(1) = (q-1)^2(q+1)/3$, and in case of $M' = \Sp_4( q )$ assume that
$\nu(1) = (q \pm 1 )^2(q^2 + 1)/2$.

Let~$\omega$ be a field automorphism of $M'$ which permutes 
$\nu_1, \ldots , \nu_d$.
Then~$\omega$ fixes at least one of $\nu_1, \ldots , \nu_d$.
\end{lem}
\begin{prf}
Suppose first that $M' = \SL_2(q )$. Then $d = 2$ and $\nu_1$ and $\nu_2$
differ on a unipotent element which is fixed by~$\omega$; see the character
table of~$M'$ given in~\cite{Frob}.
If $M' = \SL_3( q )$, the three characters $\nu_1, \nu_2 , \nu_3$ take two
distinct values on a unipotent conjugacy class of~$M'$ which is fixed by~$\omega$;
see the character table given in \cite{SiFr}. Thus $\omega$ fixes at least
one of these characters.
If $M' = \Sp_4( q )$, we use the character table computed by Srinivasan 
in~\cite{BSSp4}. If $\nu(1) = (q-1)^2(q^2+1)/2$, then $\{ \nu_1, \nu_2 \} =
\{ {\xi'_{21}}{(k)}, {\xi'_{22}}{(k)} \}$ for some value of~$k$. This follows from
the description of the conjugacy classes of $\Sp_4(q)$ in 
\cite[p.~$489$-$491$]{BSSp4} and its character table in \cite[p.~$516$]{BSSp4}.
Next, $\nu_1$ and $\nu_2$ differ on the class $A_{41}$, which is fixed by~$\omega$.
This gives the result. The same argument works for $\nu_1(1) = (q+1)^2(q^2+1)/2$.
(The article \cite{Przygocki} of Przygocki contains a few corrections to Srinivasan's
character table of $\Sp_4(q)$, but these changes do not affect our argument.)
\end{prf}

\medskip
\noindent
We next consider the blocks with non-abelian defect groups. These only occur 
for $\ell = 3$.

\begin{prop}\label{lemBAWsubpair}
Let $\ell = 3$. Assume {\rm Hypothesis~\ref{HypoInvariance}}, and in addition
that~$b$ is a non-unipotent block with a non-abelian defect group. Assume
furthermore that $\sigma = F_1^{m}$. 

Then there is a maximal $b$-Brauer pair $(D,b_D)$ fixed by~$\sigma$.
Moreover,~$D$, and the $b$-Brauer pairs $(R, b_R) \leq (D, b_D)$ with 
$\cW( R, b_R ) > 1$ can be chosen such that $\sigma(R) = R$ and 
$N_G( R, b_R ) \leq RC_G(R)N_{\mathbf{G}}(R)^\sigma$, except that
the latter condition does not always hold in the following situations:

The $G_m$-class type of~$s$ is one of $(7,2)$ or $(9,2)$ and~$f/m$ is even,
and~$R$ 
is non-abelian and properly contained in a defect group of~$b$. In this 
case,~$\sigma$ fixes the two central factors of $N_G( R, b_R )$ isomorphic
to $[q - 1]^2.2$ and $3_+^{1+2}.\SL_2(3)$, respectively, and centralizes an
element in $[q - 1]^2.2 \setminus [q - 1]^2$.
\end{prop}
\begin{prf}
The possible class types of~$s$ and the $b$-relevant radical $3$-subgroups can 
be read off from Tables~\ref{1}--\ref{19}. In particular, the defect groups 
of~$b$ are $G$-conjugate to one of $R_{29}$--$R_{34}$ and thus their 
centralizers in~$\mathbf{G}$ are connected reductive groups by 
Lemma~\ref{RegularCentralizers}. 
Let $(i,k)$ and $(i,k')$ denote the $G$-class
type, respectively the $G_m$-class type of~$s$.
Define $\varepsilon' \in \{ -1, 1 \}$ by the condition that
$3 \mid p^m - \varepsilon'$. Notice that $\varepsilon = -1$ implies
$\varepsilon' = -1$. Also, if $f/m$ is even, $\varepsilon = 1$.

Suppose first that $(\varepsilon, k) = (\varepsilon',k')$ or that the $G$-class 
type of~$s$ is one of $(3,1)$ or $(3,2)$. Let~$\mathbf{L}$ denote the 
centralizer of a $3\C$-element $z_{\C} \in G_m$ as described in 
Proposition~\ref{C3C}. We let~$D$ and~$D^\dagger$ denote a pair of 
$\sigma$-stable radical subgroups of~$L$ such that~$D$ is a defect group of~$b$, 
such that $(C_{\mathbf{L}}( D ),F)$ and $(C_{\mathbf{L}}( D^\dagger ),F)$ are in 
duality, and such that $N_{G}( D ) = DC_{G}( D )N_{G_m}( D )$. Such a pair 
exists by the construction of $R_{29}$--$R_{34}$; see 
Subsections~\ref{secradLe},~\ref{secradF4} and Lemma~\ref{lemRadLi}. Moreover, 
in every $G$-conjugacy class of $b$-relevant radical subgroups we can choose a 
representative~$R$ with $R \leq D$ and $N_{G}( R ) = RC_{G}( R )N_{G_m}( R )$.
Observe that $D_m := D \cap G_m$ and $D^\dagger_m := D^\dagger \cap G_m$ are 
radical $3$-subgroups of $G_m$ of the same type as~$D$, 
respectively~$D^\dagger$. It suffices to show that $s \in C_G( D^\dagger )$, as 
this implies the existence of a $\sigma$-stable $b$-Brauer pair $(D,b_D)$ by 
Lemma~\ref{lemBAWsubpairPreparation}. 
Now~$D_m$ is a defect group of a $3$-block of $\cE_3( G_m, s )$ by our 
assumption on~$s$ and $(\varepsilon,k)$. In particular, 
$D_m^\dagger$ is conjugate in~$G_m$ to a subgroup of~$C_{G_m}( s )$ by 
Proposition~\ref{AllDefectGroups}(a). By replacing~$s$ by a $G_m$-conjugate, we 
may assume that $s \in C_{G_m}( D^\dagger_m )$. As $C_{\mathbf{G}}( D^\dagger_m ) = 
C_{\mathbf{G}}( D^\dagger )$ by Lemma~\ref{RegularCentralizers}, we have 
$s \in C_G( D^\dagger )$, and we are done.

We are left with the cases 
$i \in \{ 6, 7, 8, 9, 10, 13, 17 \}$ and either $k = k'$ and $\varepsilon \neq 
\varepsilon'$ or $k \neq k'$. If $\varepsilon \neq \varepsilon'$, then 
$(\varepsilon,\varepsilon') = (1, - 1)$. In this case, $f/m$ is even.
If $k \neq k'$, then $f/m$ is even by Table~\ref{CT},
unless $i \in \{ 13, 17 \}$ and $k' \in \{ 4, 5 \}$. If $(i,k') 
\in \{ (13,5), (17,5 ) \}$ we have $(\varepsilon, \varepsilon') = (1,-1)$,
and $f/m$ is even and not divisible by~$3$. In particular, $k = k' = 5$.
If $(i,k') \in \{ (13,4), (17,4 ) \}$ we have $(\varepsilon, \varepsilon') 
= (-1,-1)$, and~$f/m$ is odd and divisible by~$3$. In this case, $k = 6$.
	
We put $\mathbf{L}_0^* := \mathbf{M}_{i,1}$, as constructed in~Table~\ref{CT}. 
If $k \neq k'$, choose $v \in W$ according to the following table.
$$
\begin{array}{rrrrrc} \hline\hline
i & k & k' & \varepsilon & \varepsilon' & v \rule[- 6pt]{0pt}{ 19pt} \\ \hline \hline
6 & 1 & 2 & 1 & 1 & s_{24} \\
6 & 1 & 2 & 1 & -1 & w_0 \\
7 & 1 & 2 & 1 & \pm 1 & w_0 \\
9 & 1 & 2 & 1 & \pm 1 & w_0 \\
8 & 1 & 4 & 1 & \pm 1 & w_0 \\
8 & 1 & 2 & 1 & \pm 1 & s_1 \\
8 & 1 & 3 & 1 & \pm 1 & w_0s_1 \\
10 & 1 & 2 & 1 & 1 & s_{21} \\
10 & 1 & 2 & 1 & -1 & w_0 
\rule[- 2pt]{0pt}{ 5pt} \\ \hline\hline
\end{array}
%%$
\quad\quad
%%$
\begin{array}{rrrrrc} \hline\hline
i & k & k' & \varepsilon & \varepsilon' & v \rule[- 6pt]{0pt}{ 19pt} \\ \hline \hline
13 & 1 & 6 & 1 & \pm 1 & w_0 \\
13 & 1 & 2 & 1 & \pm 1 & s_1 \\
13 & 1 & 3 & 1 & \pm 1 & w_0s_1 \\
13 & 6 & 4 & -1 & -1 & w_0s_1s_{23} \\
13 & 5 & 5 & 1 & \pm 1 & s_1s_{23}\\
17 & 1 & 6 & 1 & \pm 1 & w_0 \\
17 & 1 & 2 & 1 & \pm 1 & w_0s_4 \\
17 & 1 & 3 & 1 & \pm 1 & s_4 \\
17 & 6 & 4 & -1 & -1 & w_0s_4s_{19} \\
17 & 5 & 5 & 1 & \pm 1 & s_4s_{19} 
\rule[- 2pt]{0pt}{ 5pt} \\ \hline\hline
\end{array}
$$

\medskip
\noindent
We first deal with the cases $i \neq 8$. Then $\mathbf{L}_0^*$ is a regular 
subgroup of~$\mathbf{G}$ and the set of subgroups thus defined contains a dual 
of each of its members; we let~$\mathbf{L}_0$ denote this dual.
If $k = k'$, let $\mathbf{L} := \mathbf{L}_0$. Otherwise, choose a 
$\sigma$-stable lift $\dot{v} \in N_{\mathbf{G}}( \mathbf{T}_0 )$ of~$v$ with 
$\dot{v} = \gamma$ if $v = w_0$. Let $g \in \mathbf{G}$ with $\sigma(g) g^{-1} 
= \dot{v}$ and put $\mathbf{L} := \mathbf{L}_0^g$ and $\mathbf{T} := 
\mathbf{T}_0^g$. By Table~\ref{CT}, the $G_m$-conjugacy class 
of~$\mathbf{L}$ contains a group dual to $\mathbf{M}_{i,k'}$, and 
we may assume that $\mathbf{L} = C_{\mathbf{G}} ( s )^*$.
As $\dot{v}$ is $\sigma$-stable, the $F$-type of $\mathbf{L}$ equals
$(\bar{\Gamma}^\dagger,[v^{f/m}])$, where $\bar{\Gamma}$ is the closed subsystem 
of~$\Sigma$ defining the $G$-class type~$i$, and the $F$-type
of~$\mathbf{T}$ equals~$(\emptyset,[v^{f/m}])$. 
The action of~$\sigma$ on~$\mathbf{L}$ corresponds to the action 
of~$\sigma \dot{v}$ on $\mathbf{L}_0$. In particular, an element of 
$N_{\mathbf{L}}( \mathbf{T} )$ is $\sigma$-stable, if and only if it
is the $g$-conjugate of a $\sigma \dot{v}$-stable element of 
$N_{\mathbf{L}_0}( \mathbf{T}_0 )$. One checks that~$v$ centralizes
$W_{\mathbf{L}_0}( \mathbf{T}_0 )$, so that $\sigma$ fixes every 
element of $W_{\mathbf{L}}( \mathbf{T} )$. Thus the elements of
$W_{\mathbf{L}}( \mathbf{T} )$ can be lifted to $\sigma$-stable
elements of $N_{\mathbf{L}}( \mathbf{T} )$.

Let~$Q$ denote the Sylow $3$-subgroup of~$T$. Then~$Q$ is 
$\sigma$-stable and $Q$ is $G$-conjugate to one of $R_{11}$, $R_{12}$ or 
$R_{18}$. As $\Out_L(Q) = W_{\mathbf{L}}( \mathbf{T} )$, we find that 
$N_L( Q ) = C_L(Q)N_{\mathbf{L}}( Q )^\sigma$. There is an element~$u$ of 
order~$3$ in $W_{\mathbf{L}}( \mathbf{T} )$, which lifts to a $\sigma$-stable 
element $\dot{u} \in N_{\mathbf{L}}( \mathbf{T} )$ such that 
$D := \langle Q, \dot{u} \rangle$ is a Sylow $3$-subgroup of~$L$.
Moreover, $\Out_L( D )$ is generated by elements of order~$2$ in 
$W_{\mathbf{L}}( \mathbf{T} )$. It follows that 
$N_L( D ) = DC_L(D)N_{\mathbf{L}}( D )^\sigma$.

As~$D$ is a $\sigma$-stable defect group of~$b$, and $\sigma(s) = s$, we obtain 
a $\sigma$-stable $b$-Brauer pair $(D,b_D)$ by 
Lemma~\ref{lemBAWsubpairPreparation}, applied with $\omega = \sigma$ and 
a radical $3$-subgroup $Q^\dagger$ such that $C_{\mathbf{G}}( Q^\dagger )$
is dual to $C_{\mathbf{G}}( Q )$ and $s \in C_G( Q^\dagger )$; see 
Proposition~\ref{AllDefectGroups}. Let $(R, b_R) \leq (D, b_D)$ be a 
$b$-Brauer pair with $\cW( R, b_R ) > 0$. If $R \lneq D$ and~$R$ is abelian, we 
get $R = Q$, so that~$R$ is $\sigma$-stable. In turn, $(R, b_R)$ is 
$\sigma$-stable. Now suppose that $R \lneq D$ is non-abelian. Then 
$i \in \{ 6, 7, 9, 10 \}$ and we let $\mathbf{M}_0$ denote the standard Levi 
subgroup of~$\mathbf{L}_0$ of type~$A_2$, respectively $\tilde{A}_2$, 
and put $\mathbf{M} := \mathbf{M}_0^g$. Then $\mathbf{M} =_{G} 
\mathbf{M}_{13,1}$ if $i \in \{ 6, 7 \}$, and $\mathbf{M}_{17,1}$ if 
$i \in \{ 9, 10 \}$. Thus $M = \langle Z(M) \circ_3 M', x \rangle$
with~$x$ as in Proposition~\ref{NormalizerM13}, $Z(M) \cong [q-1]^2$ and 
$M' = \SL_3( q )$. Now
$R = Q_0 \circ_3 Q_1$, where $Q_0$ is the Sylow $3$-subgroup of $Z( M )$, and 
$Q_1 \leq M'$ is isomorphic to~$3_+^{1+2}$. By the construction described in 
Subsection~\ref{secradLe}, we may choose~$Q_1$ to be $\sigma$-stable.
In particular,~$R$ and thus $(R,b_R)$ are $\sigma$-stable.
If $i \in \{ 7, 9 \}$, there is a $\sigma$-stable element $n \in L$ of order~$4$
normalizing~$Z(M)$ and centralizing $[\mathbf{M}, \mathbf{M}]$; see
Proposition~\ref{NormalizerM13}(b). Thus $N_L( R ) 
= \langle Z(M), n \rangle \circ_3 N_{M'}( Q_1 )$ if $a \geq 2$. In case $a = 1$, 
there is a $\sigma$-stable element~$x' \in M \setminus M'$ which induces a 
diagonal automorphism on~$M'$, normalizes~$Q_1$ and with ${x'}^3 \in Z( M' )$;
see Proposition~\ref{NormalizerM13} with~$F$ repaced by~$\sigma \dot{v}$. 
As $N_{M'}( Q_1 ) \cong 3_+^{1+2}.Q_8$ in this case, we obtain 
$N_L( R ) = \langle Z(M), n \rangle \circ_3 \langle N_{M'}( Q_1 ), x' \rangle$.
In either case, $N_L( R ) = RC_L( R )N_{\mathbf{L}}( R )^\sigma \cong 
[q-1]^2.2 \circ_3 3_+^{1+2}.\SL_2(3)$.
If $i \in \{ 6, 10 \}$, we may assume that $N_{M'}( Q_1 )$
respectively $\langle N_{M'}( Q_1 ), x' \rangle$ consists of $\sigma$-stable elements.
Indeed, $[\mathbf{M},\mathbf{M}]^\sigma = \SL_3^{\varepsilon'}( p^m )$, and 
thus~$Q_1$ can be chosen such that $N_{M'}( Q_1 ) \leq \mathbf{M}^\sigma$.
As there is a $\sigma$-stable element in~$[\mathbf{L},\mathbf{L}]$ which induces the
inverse transpose automorphism on $[\mathbf{M},\mathbf{M}]$, 
we obtain $N_L( R ) = RC_L( R )N_{\mathbf{L}}( R )^\sigma$.

For arbitrary $R \leq D$, as $RC_{G}( R ) \leq L$, we get $N_G( R, b_R ) = 
N_{L}( R )$ by 
Lemma~\ref{NewOutLemma}(c); to apply this lemma notice that $N_{L}( R )$ fixes 
$b_R$, as the canonical character of~$b_R$ equals $\hat{s} \in \Irr(C_G(R))$. 
From $N_L( R ) = RC_L( R )N_{\mathbf{L}}( R )^\sigma$ we conclude 
$N_G( R, b_R ) \leq RC_G(R)N_{\mathbf{G}}( R )^\sigma$.

It remains to consider the case $i = 8$. If $k' = 1$, put $\mathbf{L}^* := 
\mathbf{L}_0^*$. If $k' \neq 1$, let~$\mathbf{L}^*$ and $\mathbf{T}^\dagger$ 
denote the groups obtained by $\sigma$-twisting~$\mathbf{L}_0^*$ 
respectively~$\mathbf{T}_0$ with~$v^\dagger$. We may then assume that 
$\mathbf{L}^* = C_{\mathbf{G}}( s )$ by Table~\ref{CT}. Notice that 
$W_{\mathbf{L}_0^*}( \mathbf{T}_0 ) = \langle s_1, s_2, s_{24} \rangle$. Let 
$\mathbf{M}_0^\dagger$ denote the parabolic subgroup of~$\mathbf{L}_0^*$ of 
type~$A_2$ corresponding to the roots $\alpha_1$ and~$\alpha_2$.
Let~$\mathbf{M}$ and~$\mathbf{T}$ denote the $\sigma$-stable subgroups 
of~$\mathbf{G}$ obtained from~$\mathbf{M}_0^\dagger$ respectively~$\mathbf{T}_0$ 
by $\sigma$-twisting with~$v$. Then~$\mathbf{M}$ is $G_m$-conjugate to
$\mathbf{M}_{17, k'}$ if $k' \in \{ 1, 2, 3 \}$, and to $\mathbf{M}_{17, 6}$ if 
$k' = 4$; see Table~\ref{CT}. Let~$Q$ denote the Sylow $3$-subgroup of~$T$ 
contained in a $\sigma$-stable Sylow $3$-subgroup~$D$ of~$M$.
By Lemma~\ref{lemBAWsubpairPreparation}, there is a $\sigma$-stable $b$-Brauer 
pair~$(D, b_D)$. This satisfies our assertion by what we have already proved 
for $i = 17$. Now let $(R, b_R) \leq (D,b_D)$ with $R := Q$. Then $(R,b_R)$
is $\sigma$-stable. By Proposition~\ref{cormain}, $\Out_G( R, b_R )$ is sent to 
$W_{\mathbf{L}^*}( \mathbf{T}^\dagger )$ under duality. Thus $N_G( R, b_R )$ is 
$\mathbf{G}$-conjugate to an extension of~$T_0$ by 
$\langle s_{4}, s_{3}, s_{17} \rangle \cong W(A_3)$. One checks that~$n_1$ 
centralizes $\langle n_4, n_3, n_{17} \rangle$, and thus $N_G( R, b_R ) =
C_G( R )N_{\mathbf{G}}( R )^\sigma$.
\end{prf}

\medskip
\noindent
We remark that the exceptional cases do occur. For example, let $p = 13$, 
$f = 2$ and $m = 1$. Then $\varepsilon = \varepsilon' = 1$. Now assume that~$s$
is of $G_m$-class type $(7,2)$. Adopt the notation of Proposition~\ref{lemBAWsubpair}.
By Theorem~\ref{BoRoEtAl}(a), there is a $b$-Brauer pair $(R,b_R)$ with $R =_G R_{25}$
such that $R \leq L$. Now $\mathbf{L}^\sigma \cong \SL_2( 13 ) \times 7 \times \SU_3( 13 )$.
As this group does not have any subgroup isomorphic to the quaternion group~$Q_8$,
we cannot have $N_L( R ) = RC_G(R)N_{\mathbf{L}}( R )^\sigma$.

The following proposition contains the essential arguments for the invariance of 
the $b$-weights in case of non-unipotent blocks.
%%%%%%%%%%%%%%%%%%%%%%%%%%%%%%%%%%%%%%%%%%%%%%%%%%%%%%%%%%%%%%%%%%%%%%%%%%%%%%%%
%%%%%%%%%%%%%%%%%%%%%%%%%%%%%%%%%%%%%%%%%%%%%%%%%%%%%%%%%%%%%%%%%%%%%%%%%%%%%%%%
%%
%%  The proposition lemma was included on January 3, 2020, revised and extended
%%  on January 10, 2020, on June 24, 2020, on June 25, 2020, on July 2
%%  and on July 7, 2020.
%%
%%  Further revisions on January 8, 2021 and on January 22, 2021.
%%
%%%%%%%%%%%%%%%%%%%%%%%%%%%%%%%%%%%%%%%%%%%%%%%%%%%%%%%%%%%%%%%%%%%%%%%%%%%%%%%%
%%%%%%%%%%%%%%%%%%%%%%%%%%%%%%%%%%%%%%%%%%%%%%%%%%%%%%%%%%%%%%%%%%%%%%%%%%%%%%%%

\begin{prop}
\label{InvariantExtensions}
Let~$b$,~$\sigma$ and~$s$ be as in {\rm Hypothesis~\ref{HypoInvariance}}, and 
assume in addition that~$b$ is a non-unipotent block and that $\sigma = F_1^m$. 
Then there is a maximal 
$b$-Brauer pair $(D,b_D)$ and an element $\omega \in N_A(D)$ with $N_A(D) = 
\langle N_G( D ), \omega \rangle$ such that~$\omega$ fixes~$b_D$. Moreover, the
following holds.
	
Let $(R, b_R) \leq (D,b_D)$ be a relevant $b$-Brauer pair fixed by~$\omega$ and
let~$\theta_R$ denote the canonical character of~$b_R$. Then some extension 
of~$\theta_R$ to~$N_G( R, b_R )$ is $\omega$-stable. In particular,~$\omega$ 
stabilizes every element of $\Irr^0( N_G( R, \theta_R ) \mid \theta_R )$, 
unless~$b$ is one of the blocks described in 
{\rm Proposition~\ref{AbelianD}(a),(b)}.
\end{prop}
\begin{prf}
If the defect groups of~$b$ are abelian, choose $(D,b_D)$ and~$\omega$ as in 
Proposition~\ref{AbelianD}. Otherwise, let $\omega := \sigma$, and 
choose~$(D,b_D)$ and the subpairs $(R,b_R) \leq (D,b_D)$ as in 
Proposition~\ref{lemBAWsubpair}. Then $(D,b_D)$ and $(R,b_R)$ are 
$\omega$-stable, so that our first statement holds. Moreover,~$\omega$
centralizes~$\Out_G( R, b_R )$, up to the exceptional cases listed in 
Propositions~\ref{AbelianD}(a),(b) and~\ref{lemBAWsubpair}. Thus, in view of 
Lemma~\ref{remAWC}, the last statement follows from the penultimate one in the 
non-exceptional cases.
Write $\theta := \theta_R$. Then $N_G( R, b_R ) = N_G( R, \theta )$. 
In particular,~$\omega$ fixes~$\theta$.

Suppose first that the defect groups of~$b$ are non-abelian. Then $\ell = 3$ and
$\omega = \sigma$. If $\theta(1) = 1$, the assertion follows from 
Proposition~\ref{lemBAWsubpair} and Lemma~\ref{StableInLinearCase}. If 
$(R, b_R)$ is as in one of the exceptional cases listed in 
Proposition~\ref{lemBAWsubpair}, the two elements of 
$\Irr^0( N_G( R, \theta) \mid \theta )$ are $\sigma$-stable by
Lemma~\ref{StableInLinearCase} and the considerations on central products in
Subsection~\ref{CentralProducts}. We have now proved our claim in case the
$G$-class type of~$s$ is one of  $(13,k)$ or $(17,k)$ with $k \in \{ 4, 5 \}$,
as then either $|\Irr^0( N_G( R, \theta) \mid \theta )| = 1$ or $\theta(1) = 1$.
Thus we exclude these cases in the discussion to follow. If $f/m$ is even, then 
$C_G( s )$ is $G$-conjugate to one of $M_{i,1}$ with $i \in 
\{ 2, 3, 5, \ldots , 10, 13 , 17 \}$; see Table~\ref{CT}. In these instances,
$C_G( R )$ is abelian, hence $\theta(1) = 1$, a case we have already considered.
We may thus assume that~$f/m$ is odd and that $\theta(1) \neq 1$.
Then~$R$ is $G$-conjugate to one of $R_{19}$, $R_{20}$, $R_{23}$, $R_{24}$,
$R_{31}$ or $R_{32}$. Here, $|\Irr^0( N_G( R, \theta) \mid \theta )|$ is even
and at most equal to~$4$, which implies that~$\sigma$ fixes at least one element
of $\Irr^0( N_G( R, \theta) \mid \theta )$.

Suppose now that~$b$ has abelian defect groups. If~$f/m$ is odd and the 
commutator quotient of $\Out_G( D, \theta )$ has order~$2$ or~$4$, there is an 
invariant extension of~$\theta$, and we are done. Suppose from now on that
$f/m$ is even or that the commutator quotient of $\Out_G( D, \theta )$ 
does not have order~$2$ or~$4$. Put $\mathbf{M} := C_{\mathbf{G}}( D )$. By 
Proposition~\ref{AbelianD}, we have $N_G( D, \theta ) = N''C_G(D)$ with
$Z(M) \leq N'' \leq C_G( [\mathbf{M}, \mathbf{M}] )$, and every coset of
$N''/Z(M)$ contains an $\omega$-stable element. Thus, 
if $\theta( 1 ) = 1$, every extension of~$\theta$ to $N_G( D, \theta )$
is $\omega$-invariant by Lemma~\ref{StableInLinearCase}.
Suppose then that $\theta(1) \neq 1$. Then $e \in \{1,2\}$ and~$D$ is 
$G$-conjugate to one of $R_{9,\ell}$, $R_{10,\ell}$, $R_{11,\ell}$, 
$R_{12,\ell}$, $R_{16,\ell}$ or 
$R_{17,\ell}$, where $R_{11,\ell}$ and $R_{12,\ell}$ only occur for $\ell > 3$. 
Also,~$\mathbf{M}$ is an $e$-split Levi subgroup of~$\mathbf{G}$ which 
is~$G$-conjugate to one of the standard regular subgroups described in 
Propositions~\ref{NormalizerM13},~\ref{NormalizerM14} and~\ref{NormalizerM18}.  
As $C_G( D ) = M$ and~$D$ is the Sylow $\ell$-subgroup of~$Z(M)$, we also have 
$N_G( D ) = N_G( M )$.

Recall that $d = \gcd( 3, q - 1 )$ if $D \in_G 
\{ R_{11,\ell}, R_{12,\ell} \}$ and that $d = \gcd( 2, q - 1 )$ in all other 
cases. Suppose first that $R \neq_G R_{9, \ell}$. 
Then~$M$ is of the form $M = \langle Z(M) \circ_d M', x' \rangle$ and 
$M' = [\mathbf{M},\mathbf{M}]^F$. Moreover, $x' = 1$ if $d = 1$, and~$x'$ is an 
arbitrary element of $M \setminus (Z(M) \circ_d M')$ if $d \neq 1$. By the above 
description of $N_G( D, \theta )$, we have  $N_G( D, \theta ) = N_G( M, \theta ) 
= \langle N'' \circ_d M', x' \rangle$. Now~$\omega$ stabilizes~$\mathbf{M}$ and 
commutes with~$F$, and thus $Z(M)$, $M' = [\mathbf{M},\mathbf{M}]^F$ 
and $N'' = C_{N_G( D, \theta )}( M' )$ are also stabilized by~$\omega$. 
As~$M$ is a central product of $Z(M)$ with $\langle M', x' \rangle$, we can 
write $\theta = \lambda \psi$ for irreducible characters~$\lambda$ of $Z(M)$
and $\psi$ of $\langle M', x' \rangle$, respectively.
Hence the restriction of $\theta$ to $Z(M) \circ_d M'$ is of the form 
$\theta' = \lambda \psi'$, where $\psi'$ denotes the restriction of~$\psi$ 
to~$M'$. It follows that~$\omega$ fixes~$\lambda$ and~$\psi'$, as~$\omega$ 
fixes~$\theta$ and stabilizes $Z(M)$ and $M'$. If $d = 1$, then $M' = M$ and
$\psi' = \psi$. Otherwise,~$\omega$ acts as a field automorphism on~$M'$ by 
Proposition~\ref{AbelianD}, and thus~$\omega$ fixes at least one 
irreducible constituent of~$\psi'$ by Lemma~\ref{FieldAutomorphism}.

Let~$\hat{\theta}$ denote an extension of~$\theta$ to $N_G( R, \theta )$,
write~$\hat{\theta}'$ for its restriction to $N''M' = N'' \circ_d M'$, and let 
$\hat{\lambda} \psi''$ be an irreducible constituent of~$\hat{\theta}'$. Then 
$\hat{\lambda}$ is an extension of $\lambda$ to $N''$ and $\psi''$ is an 
irreducible constituent of~$\psi'$. By what we have said above, we may assume 
that $\psi''$ is $\omega$-invariant. The distinct extensions of $\theta$ 
to~$N''M$ and the distinct extensions of $\lambda$ to $N''$ are of the form 
$\hat{\theta}\xi$, respectively $\hat{\lambda}\xi$, where~$\xi$ runs through 
the irreducible characters of $N''/Z(M) = N_G( D, \theta )/M$. If we are in 
one of the cases (a) or (b) of Proposition~\ref{AbelianD}, the last statement 
of this proposition shows that we may choose $\hat{\theta}$ such that 
$\hat{\lambda}$ is $\omega$-invariant. In the other cases, 
Lemma~\ref{StableInLinearCase} implies that~$\omega$ stabilizes~$\hat{\lambda}$.
It follows that $\omega$ fixes $\hat{\lambda}\psi''$. If $\hat{\theta}'$ is 
reducible, i.e.\ if $\psi'' \neq \psi'$, then $\hat{\theta}$ is 
$\omega$-invariant, as $\hat{\theta}$ is induced from $\hat{\lambda}\psi''$. 
We may thus assume that $\hat{\theta}' = \hat{\lambda}\psi'$ is irreducible 
and fixed by $\omega$.  Then $\theta' = \lambda \psi'$ is irreducible, and
Corollary~\ref{FixingReductionCor}, applied with $N = N''M'$, implies 
that~$\omega$ fixes $\hat{\theta}$.

%%%%%%%%%%%%%%%%%%%%%%%%%%%%%%%%%%%%%%%%%%%%%%%%%%%%%%%%%%%%%%%%%%%%%%%%%%%%%%%%
%%%%%%%%%%%%%%%%%%%%%%%%%%%%%%%%%%%%%%%%%%%%%%%%%%%%%%%%%%%%%%%%%%%%%%%%%%%%%%%%
%%
%% Included on 04.12.2020; expanded and adjusted on 09.12.2020.
%%
%% Further adjusted on 08.01.2021, on 15.02.2021 and on 16.02.2021.
%%
In case of $R_{9,\ell}$ we use a variant of this argument. Let~$s$ denote an
$\omega$-stable semisimple $\ell'$-element such that $b \subseteq 
\cE_\ell( G, s )$. An inspection of Tables~\ref{1}--\ref{19} gives the 
possible $G$-class types of~$s$ and also shows that $\Out_G( D, \theta )$ has 
order~$2$ or~$4$. We may thus assume that~$f/m$ is even. The case that~$s$ has
$G$-class type~$(16,7)$ is one of the exceptional cases of 
Proposition~\ref{AbelianD}(b) and thus excluded. Table~\ref{CT} then
shows that the $G$-class type of~$s$ is one of 
$\{ (4,1), (7,1), (9,1), (14,1) \}$; for example, every $\sigma$-stable
semisimple element of $\mathbf{G}$-class type~$14$ has $G$-class type~$(14,1)$
and no $\sigma$-stable element of $\mathbf{G}$-class type~$13$ can have 
%\marginpar{Read the proofs for $R_{9,\ell}$ once more. Check compatibility with
%Proposition~\ref{NormalizerM14}.}
$G$-class type $(13,2)$ or $(13,3)$. From Tables~\ref{4},~\ref{7},~\ref{9}
and~\ref{14} we then get that $\ell \mid q + 1$, that $\Out_G( D, \theta ) 
\cong 2^2$ and thus that $N_G( D, \theta ) = N_G( M )$. 

In Proposition~\ref{AbelianD}, the group $\mathbf{M} = C_{\mathbf{G}}( D )$
was constructed by twisting $\mathbf{L}_{\Delta}$ for 
$\Delta = \{ \alpha_{22}, \alpha_{17} \}$ with $s_1s_4$; see also
Definition~\ref{DefineM} and Lemma~\ref{ProofCT2}. In 
Proposition~\ref{NormalizerM14}, the standard $2$-split Levi subgroup 
of~$\mathbf{G}$ in the class containing $\mathbf{M}_{14,4}$ was constructed
by twisting $\mathbf{L}_\Gamma$ for $\Gamma = \{ \alpha_1, \alpha_4 \}$ 
with~$w_0$, the longest element of~$W$.
Now there is $z \in W$ that maps the four-tuple 
$(\alpha_{22},\alpha_{17},\alpha_1,\alpha_4)$ of roots 
%\marginpar{The action of $\omega$ on $\mathbf{M}$ corresponds to the action
%of $\sigma\dot{v}$ on~$\mathbf{M}_0$. How are these actions transformed
%under the conjugation $\mathbf{L}_{ \alpha_{22}, \alpha_{17} } 
%\rightarrow \mathbf{L}_{ \alpha_{1}, \alpha_{4} }$? Study the proof of
%Proposition~\ref{AbelianD} on Page 84.}
to $(\alpha_1,\alpha_4,\alpha_{22},\alpha_{17})$. Choose a $\sigma$-stable
lift $\dot{z} \in \hat{W}$. Then $\mathbf{L}_{\Delta}^{\dot{z}} = 
\mathbf{L}_\Gamma$. Conjugation by $\dot{z}$ sends the elements $n_1$, 
$n_4n_3^2$ of Lemma~\ref{ProofCT2} to lifts $m_1$ and $m_2$ of $s_{22}$ 
and $s_{17}$, respectively, 
satisfying the properties of Proposition~\ref{NormalizerM14}(a). 
The Steinberg morphism $\sigma \dot{v}$ with $v = s_{22}s_{17}$ on 
$\mathbf{L}_{\Delta}$ is transformed to the Steinberg morphism 
$\sigma \dot{v}'$ with $v' = s_1s_4$. Also, $\mathbf{M}^{\dot{z}}$ is
obtained by twisting $\mathbf{L}_{\Gamma}$ with $s_{22}s_{17}$. As 
$w_0 = s_{22}s_{17}s_1s_4$, the group $\mathbf{M}^{\dot{z}}$ is a
representative of the $G$-conjugacy class $\mathbf{M}_{14,4}$ of
$2$-split Levi subgroups corresponding to the closed subsystem with
base $\{ \alpha_1, \alpha_4 \}$; see the discussion in 
Subsection~\ref{MaximalRank}. We may therefore replace~$\mathbf{M}$
by this standard copy and use the notation of Proposition~\ref{NormalizerM14}. 
By Lemma~\ref{ProofCT2}, the automorphism $\omega$ then acts on~$\mathbf{M}$ as 
$\sigma \dot{v}'$ for some $\dot{v}' \in [\mathbf{M},\mathbf{M}]^F = M'$. In 
particular,~$\omega$ acts as a field automorphism on~$M'$, normalizes~$M_1$ 
and~$M_2$ and fixes~$m_1$ and $m_2$. 

We have $N_G( M, \theta ) = N_G( M ) 
= \langle M_1.2 \circ_d M_2.2, x \rangle$, where $x \in M$ 
centralizes~$Z(M_1)$ and~$M_2$ and acts as a diagonal automorphism 
on~$M_1$ if $d = 2$. Moreover, $M_i.2$ is $\omega$-invariant for $i = 1, 2$, as
$\dot{v}' \in M' = M_1' \circ_d M_2'$. Put $N := M_1.2 \circ_d M_2.2$. Then
$[N_G(M)\colon\!N] = d$ and $M \cap N = M_1 \circ_d M_2$.
Let $\hat{\theta}$ denote an extension of~$\theta$ to $N_G( M )$, satisfying the
additional properties exhibited in Proposition~\ref{NormalizerM14}(c) in case
$d = 2$. 
Let $\hat{\theta}_1\hat{\theta}_2$ with $\hat{\theta}_i \in \Irr( M_i.2 )$ 
for $i = 1, 2$, denote an 
irreducible constituent of $\Res_{N}^{N_G(M)}( \hat{\theta} )$, and put
$\theta_i = \Res_{M_i}^{M_i.2}( \hat{\theta_i} )$, $i = 1, 2$.
If $\Res_{N}^{N_G(M)}( \hat{\theta} )$ is reducible, then $d = 2$ and 
$\hat{\theta}$ is induced from $\hat{\theta}_1\hat{\theta}_2$ by 
Proposition~\ref{NormalizerM14}(c). If $\Res_{N}^{N_G(M)}( \hat{\theta} )$ is 
irreducible, so is $\Res_{M \cap N}^{N_G(M)}( \hat{\theta} ) = 
\Res_{M\cap N}^M( \theta )$, once more by Proposition~\ref{NormalizerM14}(c).
In either case, $\theta_i \in \Irr(M_i)$ for $i = 1, 2$ and $\theta_1\theta_2$ 
is an irreducible constituent of $\Res^M_{M\cap N}( \theta )$. Thus it suffices 
to show that $\hat{\theta}_1$ and $\hat{\theta}_2$ are $\omega$-invariant, using 
Corollary~\ref{FixingReductionCor} in case $\Res_{N}^{N_G(M)}( \hat{\theta} )$ 
is irreducible and $d = 2$.

%The restriction of $\theta$ to
%$M_1 \circ_d M_2$ is irreducible, hence of the form $\theta_1\theta_2$ with
%$\theta_i \in \Irr( M_i )$ for $i = 1, 2$. Similarly, the restriction of 
%$\hat{\theta}$ to $M_1.2 \circ_d M_2.2$ is irreducible, hence of the form 
%$\hat{\theta}_1 \hat{\theta}_2$, where $\hat{\theta}_i$ is an extension 
%of~$\theta_i$ to $M_i.2$ for $i = 1, 2$. Applying Corollary~\ref{FixingReductionCor}
%with $N =  M_1.2 \circ_2 M_2.2$ if $d = 2$, it is enough to show that 
%$\hat{\theta}_1 \hat{\theta}_2$ is $\omega$-invariant.
%This will follow if $\hat{\theta}_1$ and $\hat{\theta}_2$
%are $\omega$-invariant. It suffices to prove the latter statement for $i = 2$.

Recall that $M_1 = Z(M_1) \circ_d M_1'$, so that we can write $\theta_1 =
\lambda_1 \theta_1'$ for $\lambda_1 \in \Irr( Z( M_ 1) )$ and $\theta_1'
\in \Irr( M_1' )$. If $\Res_{M \cap N}^M( \theta )$ is irreducible,~$\theta_1$
is $\omega$-invariant. Otherwise,
$\Res_{M \cap N}^M( \theta ) = \theta_1\theta_2 + \theta_1^x\theta_2$
with $\theta_1^x \neq \theta_1$. Now $\theta_1^x = \lambda_1^x {(\theta_1')}^x =
\lambda_1 {(\theta_1')}^x$, as~$x$ centralizes~$Z(M_1)$. Thus ${(\theta_1')}^x
\neq \theta_1'$. As $\omega$ acts as a field automorphism on~$M_1'$ by
Proposition~\ref{AbelianD}, it fixes~$\theta_1'$ by 
Lemma~\ref{FieldAutomorphism}. In particular,~$\omega$ does not 
interchange~$\theta_1$ and~$\theta_1^x$.
As~$\omega$ fixes $\theta_1\theta_2 + \theta_1^x\theta_2$, it must fix~$\theta_1$.
It follows that in either case,~$\omega$ fixes~$\lambda_1$ and~$\theta_1'$.
Now any extension of~$\theta_1$ to $M_1.2 = \langle Z(M_1) , m_1 \rangle
\circ_d M_1'$ is of the form $\hat{\lambda}_1 \theta_1'$ for some extension
$\hat{\lambda}_1$ of~$\lambda$ to $\langle Z(M_1) , m_1 \rangle$. 
As~$\omega$ fixes~$m_1$, Lemma~\ref{StableInLinearCase} implies that 
$\hat{\lambda}_1$ and thus $\hat{\lambda}_1 \theta_1'$ are $\omega$-invariant.

Write $\hat{\theta}_2'$ for the restriction of~$\hat{\theta}_2$ to
$\langle Z(M_2) , m_2 \rangle \circ_d M_2'$, and $\theta_2'$ for the
restriction of~$\theta$ to $Z(M_2) \circ_d M_2'$.
If $\hat{\theta}_2'$ is irreducible, so is $\theta_2'$, as otherwise the two 
constituents of $\theta_2'$ would be $m_2$-conjugate as well 
as $y$-conjugate, which is impossible; see Proposition~\ref{NormalizerM14}(b). 
An argument as in the previous paragraph, with~$x$ replaced by~$y$, shows 
that~$\hat{\theta}_2'$ is $\omega$-invariant. As~$\theta_2'$ is 
$\omega$-invariant, Corollary~\ref{FixingReductionCor}, applied to the 
normal subgroups $\langle Z(M_2) , m_2 \rangle \circ_d M_2'$ and~$M_2$ 
of~$M_2.2$ implies that $\hat{\theta}_2$ is $\omega$-invariant. This completes
our proof.
\end{prf}

\addtocounter{subsection}{6}
\subsection{Proof for the unipotent blocks}

\label{ProofUnipotentBlocks}
Here, we assume that~$b$ is a unipotent $\ell$-block of~$G$ of positive defect,
and we write $G_1 := \mathbf{G}^{F_1} = F_4(p)$.
If~$b$ is not the principal block, then $e \in \{ 1, 2, 4 \}$. If $e = 4$, there 
are two such blocks, which are swapped by~$\sigma_1$ if $p = 2$. In the other cases,
there is a unique such block; see the references given in 
Subsection~\ref{GoodPrimesBlocks} and Table~\ref{1}. 

Under our hypothesis, $A = N_{\Aut(G)}( b ) =  \Aut( G ) = 
G \rtimes \langle \sigma_1 \rangle$,
unless $p = 2$, $e = 4$ and $b$ is a non-principal block;
in the latter case, $A = N_{\Aut(G)}( b ) =  G \rtimes \langle F_1 \rangle$.
For the notation refer to Subsection~\ref{PreliminariesFinalSection}. 

We begin with the proof for the principal blocks in case $\ell > 3$.
\addtocounter{thm}{1}
\begin{prop}
\label{EquivariancePBOddq}
Let $\ell > 3$ and let~$D$ denote a Sylow $\ell$-subgroup of~$G$. Assume 
that~$D$ is non-cyclic and let~$b$ be the principal $\ell$-block of~$G$. Then 
$N_A( D )$ fixes every element of $\Irr( N_G( D )/C_G( D ))$, unless~$q$ is even 
and $e \in \{ 1, 2, 3, 6 \}$. In the latter cases, the non-trivial orbits of 
$N_A( D )$ on $\Irr( N_G( D )/C_G( D ))$ have length two, and the number of 
non-trivial orbits is as given in the following table.
\setlength{\extrarowheight}{0.5ex}
$$
\begin{array}{c||cccc}\hline\hline
e & 1 & 2 & 3 & 6 \\ \hline
\text{\rm no.} & 7 & 7 & 6 & 6 \\ \hline\hline
\end{array}
$$
\end{prop}
\begin{prf}
Let $\mathbf{S} := C_{\mathbf{G}}( D )$. Then~$\mathbf{S}$ is an $F$-stable
maximal torus of~$\mathbf{G}$ as $\ell > 3$. In particular, $N_G( D ) = 
N_G ( \mathbf{S} )$, and thus $N_G( D )/C_G( D ) = 
\left( N_{\mathbf{G}} ( \mathbf{S} )/\mathbf{S} \right)^F$. 
	
Suppose that~$\mathbf{S}$ is obtained from~$\mathbf{T}_0$ by twisting with 
$w \in W$. We will find an element $\sigma'_1 \in A$ with 
$A = G \rtimes \langle \sigma_1' \rangle$, and a $\sigma_1'$-stable inverse 
image $\dot{w} \in N_{\mathbf{G}}( \mathbf{T}_0 )$ of~$w$, and choose 
$g \in \mathbf{G}$ with $F(g)g^{-1} = \dot{w}$. Then~$\sigma_1'$ commutes with 
$F\dot{w}$, and we obtain $N_A( D ) = N_G( D ) \rtimes \langle \omega \rangle$ 
by Lemma~\ref{TwistingAndNormalizers}. From the considerations of 
Subsection~\ref{Twisting} we obtain the following commutative diagram:
$$
\begin{xy}
\xymatrix@C+30pt{
\left( N_{\mathbf{G}}( \mathbf{S} )/\mathbf{S} \right)^F \ar[r]^{\omega} &
\left( N_{\mathbf{G}}( \mathbf{S} )/\mathbf{S} \right)^F \\
\left( N_{\mathbf{G}}( \mathbf{T}_0 )/\mathbf{T}_0 \right)^{F\dot{w}} \ar[u]^{\CJxxx_g} \ar[r]_{\sigma_1'} &
\left( N_{\mathbf{G}}( \mathbf{T}_0 )/\mathbf{T}_0 \right)^{F\dot{w}} \ar[u]_{\CJxxx_g}
}
\end{xy}
$$
We may thus replace the action of~$\omega$ on 
$\left( N_{\mathbf{G}}( \mathbf{S} )/\mathbf{S} \right)^F$ by the action 
of~$\sigma_1'$ on 
$\left( N_{\mathbf{G}}( \mathbf{T}_0 )/\mathbf{T}_0 \right)^{F\dot{w}} = 
C_W( w )$. 

Suppose first that $w \in W$ can be chosen to be $\sigma_1$-stable. This is the 
case if~$q$ is odd, as then $\sigma_1 = F_1$ acts trivially on~$W$. If 
$e \in \{ 1, 2 \}$, then $w \in Z(W)$ and thus is fixed by~$\sigma_1$. If $e = 4$, 
a computation with CHEVIE \cite{chevie} shows that there is a $\sigma_1$-stable
element in the conjugacy class $D_4(a_1)$ of~$W$, which gives rise to the
maximal torus~$\mathbf{S}$ of~$\mathbf{G}$ with $S \cong [q^2 + 1] \times [q^2 + 1]$ 
(cf.\ Table~\ref{MaximalTori}). As~$\sigma_1$
is a Steinberg morphism of~$\mathbf{G}$, we may choose $\dot{w}$ to be
$\sigma_1$-stable, and let $\sigma_1' := \sigma_1$ in these cases.
The action of~$\sigma_1$ on the set 
of conjugacy classes of $C_W( w )$ is trivial if~$q$ is odd, and if~$q$ is 
even it can be computed with CHEVIE. If $e = 4$, every conjugacy class of
$C_W( w )$ is fixed by $\sigma_1$, and if $e \in \{ 1, 2 \}$, there are
exactly~$7$ pairs of $\sigma_1$-conjugate conjugacy classes of $C_W(w) = W$.
This yields our claim. 
	
It remains to consider the case that~$q$ is even and $e \in \{ 3, 6 \}$.
We give the proof in case $e = 3$; the other case is analogous. A CHEVIE
computation shows that the $W$-conjugacy class $A_2 + \tilde{A}_2$, which gives 
rise to the maximal torus~$\mathbf{S}$ of~$\mathbf{G}$ with 
$S = [q^2+q+1] \times [q^2 + q + 1]$ contains an element~$w$
such that $C_W( w )$ is $\sigma_1$-stable, but that this class does not 
contain any $\sigma_1$-stable element. 
Now let~$w$ be an element in the conjugacy class $A_2 + \tilde{A}_2$ such
that $n^{-1} \sigma_1( w ) n =  w$ for some $n \in N_W( C_W( w ) )$.
Choose an inverse image $\dot{n} \in N_{\mathbf{G}}( \mathbf{T}_0 )$ of~$n$
and put $\sigma_1' := \sigma_1 \dot{n}$. Then $\sigma_1'$ centralizes~$w$, 
and we may assume that it also centralizes~$\dot{w}$, as $\sigma_1'$ is a 
Steinberg morphism of~$\mathbf{G}$. 
A computation with CHEVIE shows that $\sigma_1'$
has exactly~$6$ orbits of length two on the set of conjugacy classes of 
$C_W( w )$ and fixes the remaining conjugacy classes. This completes the proof.
\end{prf}

\begin{table}
\caption{\label{ActiobOnWeights} Action on weights for the principal block}
$$
\begin{array}{ccccc}  \hline\hline
R & \Out_G( R, b_R) & \mathcal{W} & \text{Rem} & \#
\rule[- 7pt]{0pt}{ 20pt} \\ \hline\hline
 R_{38} & 2^3 & 8 & \null & 2 \rule[ 0pt]{0pt}{ 13pt} \\
 R_{37} & \GL_2(3) & 2 & \null \\
 R_{15} & \SL_3(3) & 1 & \null \\
 R_{18} & W(F_4) & 4 & \null & 1 \\ %\hline
 R_{21} & (Q_8\times Q_8).S_3 & 11 & a = 1 & 4 \\ %\hline
 R_{21} & (\SL_2(3)\times \SL_2(3)).2 & 2 & a \geq 2 \\
 R_{22} & \SL_2(3)\times \SL_2(3) & 1 & a \geq 2 \\
 R_{35}, R_{36} & (\SL_2(3)\times 2).2, (2\times \SL_2(3)).2 & 4, 4 & a \geq 2 & 4 \\ \hline\hline
\end{array}
$$
\end{table}

\medskip
\noindent
\noindent The following prepares the proof for the principal block in case 
$\ell = 3$.

\begin{lem}\label{lemBAWnormalPB}
Let~$b$ be the principal $3$-block of~$G$. Then every $G$-conjugacy class of 
$b$-relevant radical $3$-subgroups of~$G$ contains an $F_1$-stable 
representative~$R$ such that $N_G( R ) = RC_G( R )N_{G_1}( R )$. If $p = 2$,
such a representative can be chosen to be $\sigma_1$-invariant, unless
$R \in_G \{ R_{35}, R_{36} \}$, in which case $\sigma_1$ maps~$R_{35}$ 
to a conjugate of $R_{36}$.
\end{lem}
\begin{prf}
Fix a $\sigma_1$-stable $3\C$-element $z_{\C}$. The existence of $z_{\C}$ is 
clear if~$p$ is odd. For $p = 2$, we can use GAP~\cite{GAP04} to compute the 
class fusion of ${^2\!F}_4(2) = \mathbf{G}^{\sigma_1}$ into 
$F_4(2) = G_1 = \mathbf{G}^{F_1}$; this shows that the elements of order~$3$ 
of~${^2\!F}_4(2)$ fuse into the $3\C$-class  of~$F_4(2)$. Put 
$\mathbf{L} := C_{\mathbf{G}}( z_{\C} )$. Then $\mathbf{L} = \mathbf{L}^1 
\circ_3 \mathbf{L}^2$, and $\sigma_1$ swaps the two simple components 
$\mathbf{L}^i \cong \SL_3( \mathbb{F} )$, $i = 1, 2$, if $p = 2$, and fixes 
them if $p$ is odd; for the notation see Proposition~\ref{C3C}(a). In 
particular, $\mathbf{L}^i$ is 
$F_1$-stable for $i = 1, 2$. Let $e' \in \{ 1, 2 \}$ denote the order of~$3$ 
modulo~$p$. Then there is a $\sigma_1$-stable $e'$-$F_1$-split maximal 
torus~$\mathbf{T}$ of~$\mathbf{L}$. Again, this is clear if~$p$ is odd. If 
$p = 2$, choose a maximal $F_1$-stable torus~$\mathbf{T}_1 \leq \mathbf{L}^1$
with $\mathbf{T}_1^{F_1} \cong 3^2$, put $\mathbf{T}_2 := 
\sigma_1(\mathbf{T}_1)$ and $\mathbf{T} := \mathbf{T}_1\mathbf{T}_2$.
Clearly, $N_{G}( \mathbf{T} ) = TN_{G_1}( \mathbf{T} )$. 
By Proposition~\ref{C3C}(b) and in the notation of 
Subsection~\ref{sec3elts}, we have
$L = \langle L^1 \circ_3 L^2, x_{\C} \rangle$ with $x_{\C} \in T$. 

Let~$Q$ denote the Sylow $3$-subgroup of~$T$. Then $N_L( Q )$ is $\sigma_1$-stable,
and as $N_L( Q ) = T.(S_3 \times S_3)$, we obtain a $\sigma_1$-stable
Sylow $3$-subgroup~$P$ of~$G$ with $Q \leq P$.
%%%%%%%%%%%%%%%%%%%%%%%%%%%%%%%%%%%%%%%%%%%%%%%%%%%%%%%%%%%%%%%%%%%%%%%%%%%%%%%%
%%%%%%%%%%%%%%%%%%%%%%%%%%%%%%%%%%%%%%%%%%%%%%%%%%%%%%%%%%%%%%%%%%%%%%%%%%%%%%%%
%%
%% For more details see sheet 1 - 2 of my notes of 19.06.2020.
%%
%% This is probably not needed anymore.
%%
%%%%%%%%%%%%%%%%%%%%%%%%%%%%%%%%%%%%%%%%%%%%%%%%%%%%%%%%%%%%%%%%%%%%%%%%%%%%%%%%
%%%%%%%%%%%%%%%%%%%%%%%%%%%%%%%%%%%%%%%%%%%%%%%%%%%%%%%%%%%%%%%%%%%%%%%%%%%%%%%%

Suppose for the moment that $p = 2$, and consider $Q_0 := 
\mathbf{T}^{\sigma_1}$. This is an elementary abelian subgroup of 
$\mathbf{G}^{\sigma_1} \leq G$ of type $3\C^2$. By \cite[Lemma~$2.7$(b)]{AnDF4}, 
the $G$-conjugacy class of~$Q_0$ equals~$(3\C^2)_1$. Put $C := C_{G_1}( Q_0 )$.
Then $C$ is $\sigma_1$-stable, $C = \mathbf{T}^{F_1}.3 = 3^4.3$ by 
\cite[Table~$2$]{AnDF4}, and $C =_{G_1} R_{37}$ by
\cite[Case ($M_9$), p.~$567$]{AnDF4}. One can check that~$C$ contains a 
$\sigma_1$-stable elementary abelian subgroup of type $(3\C^3)_1$. 
 Indeed, using the permutation representation of~$G_1$ on $69\,888$ points
 given in \cite{WWWW} and GAP, we construct a copy of $\mathbf{T}^{F_1} \leq C$ 
 in~$G_1$. An explicit computation shows that
 $Z( C ) = Q_0 = [C,C]$, and that the exponent of~$C$ equals~$3$. In 
 particular, $C/Q_0$ is elementary abelian of order~$27$. Thus $C$ has
 exactly $13$ subgroups of order~$27$ containing~$Q_0$, and~$4$ of these
 lie in~$\mathbf{T}^{F_1}$. As~$F_1 = \sigma_1^2$ fixes the elements of~$C$,
 the orbits of~$\sigma_1$ on the subgroups of~$C$ have lengths at most~$2$.
 In particular,~$\sigma_1$ fixes an elementary abelian subgroup~$R$ of~$C$
 of order~$27$ with $Q_0 \leq R \not\leq \mathbf{T}^{F_1}$. Let $y \in C 
 \setminus \mathbf{T}^{F_1}$ such that $R = \langle Q_0, y \rangle$. Then 
 $C_{G_1}( R ) = C_C( y )$. A GAP computation shows that $|C_C( y' )| = 3^3$
 for every $y' \in C \setminus \mathbf{T}^{F_1}$, and thus $C_{G_1}( R ) = R$. 
 This implies that $R =_{G_1} (3\C^3)_1 =_{G_1} R_{15}$ by \cite[Lemma~$2.7$, 
 and proof of Proposition~$3.11$]{AnDF4}.

Let us return to the general case. Let~$R$ be a $b$-relevant radical 
$3$-subgroup of~$G$. We now consider the possibilities for~$R$ as listed in 
Table~\ref{1}. If $R \in_G \{ R_{15}, R_{21} \}$, we may assume that 
$R \leq G_1$. As $|\Out_{G_1}( R )| = |\Out_G( R )|$ we obtain our claim.
Suppose that $p = 2$. Then we may even choose~$R$ to be $\sigma_1$-stable.
Namely, if $R =_G R_{21}$, choose~$R$ as $R := D \circ_3 \sigma_1( D )$, 
where~$D$ is a Sylow $3$-subgroup of $\SU_3( 2 ) \leq (L^1 \cap F_4( 2 ) )$;
see Lemma~\ref{lemRadLi}. If $R =_G R_{15}$, then~$R$ is of type $(3\C^3)_1$,
and the claim follows from what we have noted above.
Now suppose that $R =_G R_{37}$. It follows from \cite{AnDF4, AH2} that 
there is an elementary abelian subgroup $E \leq \mathbf{T}^{\sigma_1}$ of type 
$(3\C^2)_1$ such that $R := C_P( E ) = Q.3$ is a radical subgroup of~$G$ of 
type~$R_{37}$. To be more specific, if~$p$ is odd, refer to 
\cite[Proofs of Cases (3.1), (3.2) and Table~$9$]{AH2}. If $p = 2$, use
\cite[Case ($M_9$), p.~$567$]{AnDF4}, together with the considerations
in the previous paragraph. Now~$R$ is $\sigma_1$-stable, as~$E$ is. Thus 
if~$R$ is conjugate to 
one of $R_{18}$, $R_{37}$ or $R_{38}$, we may choose a $\sigma_1$-stable
representative~$R$ in $N_{G}( \mathbf{T} )$. In each case, $R \cap T = Q
= C_R( [R,R] )$ is characteristic in~$R$ and in~$T$. In particular,~$R$ 
is $F_1$-stable and $N_G( R ) \leq N_G( \mathbf{T} )$. This gives our result
in these cases.

We are left with the case that $R \in_G \{ R_{35}, R_{36} \}$. It suffices to 
prove one of these cases, say $R =_G R_{35}$. Here, we may choose a 
representative~$R$ such that 
$Z(R) = \langle z_{\C} \rangle$ with $z_{\C}$ as above. Hence 
$R \leq C_G( z_\C ) = \langle L^1 \circ_3 L^1, x_\C \rangle$ and 
$N_G( \langle z_\C \rangle ) = \langle C_G( z_\C ), \gamma_C \rangle$, where 
$\gamma_\C \in G_1$; see Subsection \ref{sec3elts}. By the construction of~$R$ 
indicated in Subsections~\ref{secradLe} and~\ref{secradF4}, we may assume 
$R = K \circ_3 D \leq L^1 \circ_3 L^2$, where $K \leq L^1$ and $D \leq L^2$ are
$F_1$-stable and as in Lemma~\ref{lemRadLi}. In particular,~$R$ is $F_1$-stable 
and $\gamma_\C \in N_G( R )$. As $N_G( R ) = 
\langle N_{L^1}( K ) \circ_3 N_{L^2}( D ), \gamma_\C \rangle$, the claim 
$N_G(R) = RC_G(R)N_{G_1}(R)$ is reduced to an analogous assertion in~$L^1$ 
and~$L^2$, which clearly holds. The final statement in case $p = 2$ is also
clear, as~$\sigma_1$ swaps~$L^1$ and~$L^2$.
\end{prf}

\medskip
\noindent
We can now give the proof for the principal block  in case $\ell = 3$.

\begin{prop}
\label{EquivariancePBEvenqL3}
Let $\ell = 3$, and let~$b$ be the principal $3$-block of~$G$. If~$q$ is 
odd,~$\sigma_1$ fixes every conjugacy class of $b$-weights. If~$q$ is 
even,~$\sigma_1$ has exactly seven orbits of length two on the conjugacy classes
of $b$-weights, and fixes the remaining $G$-conjugacy classes of $b$-weights.
\end{prop}
\begin{prf}
Let~$R$ be a $b$-relevant radical $3$-subgroup of~$G$.
	
Suppose first that~$q$ is odd. By Lemma~\ref{lemBAWnormalPB} we may assume
that~$\sigma_1$ fixes~$R$ and hence $N_A( R ) = N_G( R ) \rtimes
\langle \sigma_1 \rangle$. Let $(R,b_R)$ be a $b$-Brauer pair. Then~$b_R$ is the
principal $\ell$-block of~$R$ and thus $N_G( R, b_R ) = N_G( R )$. In particular,
$\theta_R$ is the trivial character of~$RC_G(R)$ and extends to the trivial
character of~$N_G(R)$. As we may also assume that $N_G( R ) =
RC_G( R )N_{G_1}( R )$ by Lemma~\ref{lemBAWnormalPB}, every element
of~$N_A( R )$ fixes the irreducible characters of $N_G(R)/RC_G(R) =
\Out_G(  R, b_R )$. Our assertion now follows from Lemma~\ref{remAWC}.

Suppose now that $q$ is even. By Lemma~\ref{lemBAWnormalPB}, we may 
assume that~$R$ is $\sigma_1$-stable, unless $R =_G R_{35}$ or $R =_G R_{36}$.
In the latter case, we choose $\sigma_1$-conjugate representatives in the
classes containing $R_{35}$, respectively~$R_{36}$. If $a \geq 2$, this gives
rise to~$4$ non-trivial $\sigma_1$-orbits on the set of conjugacy classes 
of $b$-weights.

If~$R$ is $\sigma_1$-stable, we have to count the number of non-trivial 
$\sigma_1$-orbits on $\Out(  R, b_R )$. These numbers are given in the
last column of Table~\ref{ActiobOnWeights}. This table proves the 
proposition.
\end{prf}

\medskip

\noindent
We conclude this subsection with the proof for the non-principal unipotent blocks.

\begin{prop}
\label{NonPrincipalUnipotentBlocks}
Let $e \in \{ 1, 2 \}$ and let~$b$ be the unique non-principal unipotent 
$\ell$-block of $G$ of positive defect and let~$D$ be a defect group of~$b$. 
If~$q$ is odd, $N_A( D )$ fixes every $G$-conjugacy class of $b$-weights.
If~$q$ is even, $N_A( D )$ has a unique orbit of length two on the 
$G$-conjugacy classes of $b$-weights and fixes the remaining
$G$-conjugacy classes of $b$-weights.
\end{prop}
\begin{prf}
Put $\mathbf{M}_0 := \mathbf{L}_{\{2,3\}}$ in the notation of 
Subsection~\ref{SetupF4}. By twisting~$\mathbf{M}$ suitably, we obtain an 
$F$-stable regular subgroup~$\mathbf{M}$ of~$\mathbf{G}$ such that $D := Z( M ) 
=_G R_{10,\ell}$ is a defect group of~$b$. 

Let us now specify our choices. In case $e = 1$, let $g = 1 = \dot{w}$. In case 
$e = 2$, choose a $\sigma_1$-stable inverse image $\dot{w} \in 
N_{\mathbf{G}}( \mathbf{T}_0 )$ of $w_0$, 
and an element $g \in \mathbf{G}$ with $F(g)g^{-1} = \dot{w}$, put $\mathbf{M} 
:= \mathbf{M}_0^g$, and $\omega := \CJxxx_g \circ \sigma_1 \circ \CJxxx_g^{-1}
\in \Aut_1(\mathbf{G})$. As $\sigma_1$ stabilizes~$\mathbf{M}_0$ and commutes with 
$F\dot{w}$, Lemma~\ref{TwistingAndNormalizers} yields
$N_A( D ) = \langle N_G( D ), \omega \rangle$ with 
$\omega = \CJxxx_g \circ \sigma_1 \circ \CJxxx_g^{-1}$.
We obtain the commutative diagram
$$
\begin{xy}
\xymatrix@C+30pt{
 N_{\mathbf{G}}( \mathbf{M} )^F \ar[r]^{\omega} &
 N_{\mathbf{G}}( \mathbf{M} )^F \\
 N_{\mathbf{G}}( \mathbf{M}_0 )^{F\dot{w}} \ar[u]^{\CJxxx_g} \ar[r]_{\sigma_1} &
 N_{\mathbf{G}}( \mathbf{M}_0 )^{F\dot{w}} \ar[u]_{\CJxxx_g}
}
\end{xy}
$$
so that we can replace $N_G( D ) = N_{\mathbf{G}}( \mathbf{M} )^F$ and the action of 
$\langle \omega \rangle$ by $N_{\mathbf{G}}( \mathbf{M}_0 )^{F\dot{w}}$ and the 
action of $\langle \sigma_1 \rangle$. Now $\mathbf{M}$ is the standard $e$-split Levi 
subgroup of~$\mathbf{G}$ corresponding to $\{ \alpha_2, \alpha_3 \}$ 
considered in Proposition~\ref{NormalizerM18} and 
$N_{\mathbf{G}}( \mathbf{M}_0 )^{F\dot{w}} \cong  
N_{\mathbf{G}}( \mathbf{M} )^F$.

The canonical character of~$b$ is unipotent, so has $Z(M)$ in its kernel. It 
thus suffices to look at the central quotient
$N_{\mathbf{G}}( \mathbf{M}_0 )^{F\dot{w}}/Z( \mathbf{M}_0^{F\dot{w}} )$, whose
structure can be determined with Proposition~\ref{NormalizerM18}. 
This central quotient is of the form $\bar{N}' \times \PSp_4(q).d$ with
$d = \gcd( 2, q - 1)$, where $\bar{N}' \cong D_8$ arises from the stabilizer 
$\langle s_8, s_{16} \rangle$ of $\{ \alpha_2, \alpha_3 \}$ in~$W$. 
If~$q$ is odd, $\sigma_1$ acts trivially on $\Irr( \bar{N}' )$, and if $q$ is 
even,~$\sigma_1$ has a unique orbit of length two on $\Irr( \bar{N}' )$, as 
$\sigma_1$ swaps the two roots $\alpha_8$ and $\alpha_{16}$; see 
Table~\ref{PositiveRoots}. Since the $b$-weights are in bijection with 
$\Irr(\bar{N}')$, our claim follows.
\end{prf}

\addtocounter{subsection}{4}
\subsection{The proof for $2.F_4(2)$}
Here, we let $G := F_4( 2 )$ and $\hat{G} := 2.G$, the exceptional double cover 
of~$G$. The automorphism $\sigma_1$ of~$G$ lifts to an automorphism 
of~$\hat{G}$, also denoted by~$\sigma_1$.

\addtocounter{thm}{1}
\begin{prop}
\label{2f42mod3}
Let $\ell \in \{ 3, 5, 7 \}$ and let~$b$ be the non-principal $\ell$-block 
of~$\hat{G}$ of maximal defect.

Then the non-trivial $\sigma_1$-orbits on the set of $\hat{G}$-conjugacy classes 
of $b$-weights have length~$2$, and the number of non-trivial orbits is as 
given in the table below. For $\ell = 3$, the four orbits of length~$2$ are
distributed among the weight subgroups as indicated in the last column of 
{\rm Table~\ref{ActiobOnWeightsFaithful}}.
\setlength{\extrarowheight}{0.5ex}
$$
\begin{array}{c|ccc} \hline\hline
\ell & 3 & 5 & 7 \\ \hline
\text{\rm no.} & 4 & 6 & 6 \\ \hline\hline
\end{array}
$$
\end{prop}
\begin{prf}
Wilson's Atlas \cite{WWWW} contains a representation of the group 
$G.2$. Thomas Breuer noticed that its construction also works for the
group $\hat{G}.2$, providing a permutation representation on $279\,552$ points.
Using this representation and GAP \cite{GAP04}, it is straightforward to verify 
the claimed multiplicities.
\end{prf}

\begin{table}
\caption{\label{ActiobOnWeightsFaithful} Action on weights for the faithful $3$-block of $2.F_4(2)$}
$$
\begin{array}{cccc}  \hline\hline
R & N_{\hat{G}}(R)/R & \mathcal{W} & \#
\rule[- 7pt]{0pt}{ 20pt} \\ \hline\hline
 R_{38} & 2^4 & 8 & 2 \rule[ 0pt]{0pt}{ 13pt} \\
 R_{37} & 2 \times \GL_2(3) & 2  \\
 R_{15} & 2 \times \SL_3(3) & 1  \\
 R_{18} & 2 \times W(F_4) & 4 & 1 \\ 
	R_{21} & 2\dot{\ }[(Q_8\times Q_8)\colon\!S_3] & 2 & 1 \\ \hline\hline
\end{array}
$$
\end{table}

\addtocounter{subsection}{1}
\subsection{Proof in case $p = 2$ and $m'$ odd}

Let $b$,~$\sigma$ and~$s$ be as in Hypothesis~\ref{HypoInvariance}. Assume 
further that $p = 2$, that $m'$ is odd and that~$b$ is a non-unipotent block. 
By Proposition~\ref{ActionAutomorphisms}(a), this can only occur 
in the following cases: the $G$-class type of~$s$ is one of $(14,1)$,
$(14,4)$, $(15,1)$, $(15,3)$ or $(15,5)$, or $\ell > 3$ and the $G$-class type 
of~$s$ is one of $(4,1)$ or~$(4,2)$. In the latter two cases, 
$e \in \{ 1, 2, 3 \}$ and $e \in \{ 1, 2, 6 \}$, respectively.

\addtocounter{thm}{1}
\begin{prop}
\label{mprimeOdd}
Assume the hypothesis and notation introduced at the beginning of this 
subsection. 

Then the non-trivial orbits of~$N_A( D )$ on the set of conjugacy classes of 
$b$-weights have length two, and the number of such orbits equals~$3$ in 
case~$s$ has class type $(4,1)$ and $e \in \{ 1, 3 \}$ or $(4,2)$ and 
$e \in \{ 2, 6 \}$, and~$1$ in all other cases.
\end{prop}
\begin{prf}
Suppose first that the $G$-class type of~$s$ is $(4,1)$ or $(4,2)$. 
We may assume that 
$\sigma_1( s ) = s$, so that $\sigma = \sigma_1$. Moreover,~$\sigma$ fixes 
$C_{\mathbf{G}}( s )$ swapping its two simple components~$\mathbf{L}^1$ 
and~$\mathbf{L}^2$. By the description of the defect groups in Table~\ref{4}, 
this implies that there is a $\sigma_1$-stable maximal torus $\mathbf{S}$ 
of~$C_{\mathbf{G}}( s )$, such that the Sylow $\ell$-subgroup~$D$ of~$S$ is a 
defect group of~$b$. In fact $\mathbf{S} = \mathbf{S}_1\mathbf{S}_2$,
where~$\mathbf{S}_i$ is a maximal torus of~$\mathbf{L}^i$ and $\sigma$ swaps 
$\mathbf{S}_1$ with $\mathbf{S}_2$.
%%%%%%%%%%%%%%%%%%%%%%%%%%%%%%%%%%%%%%%%%%%%%%%%%%%%%%%%%%%%%%%%%%%%%%%%%%%%%%%%
%%%%%%%%%%%%%%%%%%%%%%%%%%%%%%%%%%%%%%%%%%%%%%%%%%%%%%%%%%%%%%%%%%%%%%%%%%%%%%%%
%%
%% To be specific, $\mathbf{S} =_G \mathbf{S}_{1}, \mathbf{S}_{22}, \mathbf{S}_9$,
%% if~$s$ is of class type $(4,1)$ and $e = 1, 2, 3$, respectively, and
%% $\mathbf{S} =_G \mathbf{S}_{2}, \mathbf{S}_{22}, \mathbf{S}_{10}$,
%% if~$s$ is of class type $(4,2)$ and $e = 1, 2, 6$, respectively.
%	
%%%%%%%%%%%%%%%%%%%%%%%%%%%%%%%%%%%%%%%%%%%%%%%%%%%%%%%%%%%%%%%%%%%%%%%%%%%%%%%%
%%%%%%%%%%%%%%%%%%%%%%%%%%%%%%%%%%%%%%%%%%%%%%%%%%%%%%%%%%%%%%%%%%%%%%%%%%%%%%%%
Moreover,~$\mathbf{S}$ is its own dual with respect to~$\sigma$. Hence there is 
a $\sigma$-stable $b$-Brauer pair $(D,b_D)$ by 
Lemma~\ref{lemBAWsubpairPreparation}.
Now $\Out_G( D, b_D )$ may be identified with 
$W_{C_{\mathbf{G}}( s )}( \mathbf{S} )^F$ by Proposition~\ref{cormain},
and $W_{C_{\mathbf{G}}( s )}( \mathbf{S} )^F = 
W_{\mathbf{L}^1}( \mathbf{S}_1 )^F \times W_{\mathbf{L}^2}( \mathbf{S}_2 )^F$,
where $\sigma$ swaps the two factors 
$W_{\mathbf{L}^1}( \mathbf{S}_1 )^F$ and $W_{\mathbf{L}^2}( \mathbf{S}_2 )^F$.
This gives our claim.

The proof in case of class types $(14,1)$ and $(14,4)$ is analogous.

It remains to consider the case that the $G$-class type of $s$ is one of 
$(15,1)$, $(15,3)$ or $(15,5)$. By the facts summarized in 
Hypothesis~\ref{HypoInvariance}, we have that $m = m'$, that~$m$ divides~$f$ 
and that $\sigma^2 = F_1^m$. 
We proceed as in Proposition~\ref{AbelianD}, but only sketch the arguments.
As $p = 2$, the elements $n_j \in N_G( \mathbf{T}_0 )$, $j = 1, \ldots , 4$,
have order $2$. We may and will thus identify 
$W = N_{\mathbf{G}}( \mathbf{T}_0 )/\mathbf{T}_0$ with the subgroup 
of $N_{\mathbf{G}}( \mathbf{T}_0 )$ generated by $n_j$, $ j = 1, \ldots , 4$. 
Let $\mathbf{L}_0 := \mathbf{L}_{ \{ \alpha_2,\alpha_3 \} } 
=_G \mathbf{M}_{15,1}$. If the $G$-class type of~$s$ is $(15,1)$,
take $v = 1$. If the $G$-class type of~$s$ is $(15,3)$ or $(15,5)$, let 
$v = s_{8}$. Then~$v$ centralizes $s_2$ and $s_3$. 
respectively. Notice that $\sigma^2 = F_1^m$ fixes $\sigma(v)v$, and if 
$g \in \mathbf{G}$ with $\sigma(g)g^{-1} = v$, we have $F(g)g^{-1} = 
(\sigma(v)v)^{f/m}$. By $\sigma$-twisting $\mathbf{L}_0$ with $v^\dagger$, 
we obtain 
$\sigma$-stable regular subgroups $\mathbf{L}^\dagger$ of $F$-type 
$( \Gamma^\dagger, [(\sigma(v^\dagger)v^\dagger)^{f/m}])$, where~$\Gamma$
is the parabolic subsystem of~$\Sigma$ spanned by $\{ \alpha_2, \alpha_3 \}$.
If $v = s_{8}$, then $v^\dagger = s_{16}$, and $(\sigma(v^\dagger)v^\dagger) =
s_8s_{16}$ lies in conjugacy class~$23$ of~$W$. Using 
Lemma~\ref{TypesOfMaximalRankSubgroupsUnderPowersOfSigma}, the $2$-power map 
on~$W$ given in
Table~\ref{MaximalTori}, and the fusion of the maximal tori of Table~\ref{cc},
we conclude that $\mathbf{L}^\dagger =_G \mathbf{M}_{15,5}$ if $f/m$ is odd,
and that $\mathbf{L}^\dagger =_G \mathbf{M}_{15,3}$ if $f/m \equiv 0(4)$.
By the results summarized in Subsection~\ref{MaximalRank}, the centralizer of 
every semisimple $\sigma$-stable element of these $G$-class types arises in 
this way. We may thus identify $\mathbf{L}^\dagger$ with $C_{\mathbf{G}}( s )$. 

By $\sigma$-twisting $\mathbf{L}_0$ with~$v$, we obtain a $\sigma$-stable 
regular subgroup $\mathbf{L}$ such that $(\mathbf{L}^\dagger, \sigma)$
and $(\mathbf{L},\sigma)$ are in duality. If $\ell \mid q - 1$, let $w = 1$.
Otherwise, let~$w$ be the longest element of 
$W_{\mathbf{L}_0}( \mathbf{T}_0 ) = W_{\{ \alpha_2,\alpha_3 \}}$. By first 
$\sigma$-twisting $\mathbf{T}_0$ with~$v^\dagger$, respectively~$v$, and then 
$F$-twisting with~$w^\dagger$, respectively~$w$, as in 
the proof of Proposition~\ref{AbelianD}, we obtain a pair of $F$-stable maximal 
tori $\mathbf{S}^\dagger \leq \mathbf{L}^\dagger$ and
$\mathbf{S} \leq \mathbf{L}$, and a pair of Steinberg morphisms $\omega^\dagger$ 
and $\omega$ such that $\omega^\dagger$ fixes~$\mathbf{S}^\dagger$ and~$\omega$
fixes~$\mathbf{S}$ and 
such that $(\mathbf{S}^\dagger, \omega^\dagger)$ and $(\mathbf{S},\omega)$ as
well as $(\mathbf{S}^\dagger,F)$ and $(\mathbf{S},F)$ are
in duality. Let $D^\dagger$ and~$D$ denote the Sylow $\ell$-subgroups of~$S^\dagger$
and~$S$, respectively, and put $\mathbf{M}^\dagger := C_{\mathbf{G}}( D^\dagger )$ and
$\mathbf{M} := C_{\mathbf{G}}( D )$. Then~$D$ is a defect group of~$b$, and we 
obtain an $\omega$-stable $b$-Brauer pair $(D,b_D)$ by 
Lemma~\ref{lemBAWsubpairPreparation}. Moreover, $N_A( D ) = 
\langle N_G( D ), \omega \rangle$ by Lemma~\ref{TwistingAndNormalizers}. 
Finally,~$\omega$ acts as the exceptional
isogeny on $M' = [\mathbf{M},\mathbf{M}]^F \cong \Sp_4( q )$.

By Proposition~\ref{NormalizerM18}, we have $M = Z( M ) \times M'$.
Notice that $S \leq M$ and that $N_G( D, b_D ) = N_M( S ) = S.W_M$,
a split extension with an $\omega$-stable subgroup $W_M \cong \Out_M( S ) 
\cong D_8$. In particular, the canonical character of~$b_D$ has an $\omega$-invariant
extension to $N_G( D, b_D )$. As~$\omega$ interchanges two generators of~$W_M$, 
our claim follows.
\end{prf}

\addtocounter{subsection}{1}
\subsection{The equivariance condition}
\label{BAWCEquivariance}
We can now prove the second main result of this article.

\addtocounter{thm}{1}
\begin{thm}\label{thmmain}
Let $G = F_4(q)$ and let $\ell$ be an odd prime with $\ell \nmid q$. Then~$G$
satisfies the inductive blockwise Alperin weight condition at the prime~$\ell$. 
\end{thm}
\begin{prf}
We have to verify the conditions of Hypothesis~\ref{defBWC_Simplified} for all
$\ell$-blocks of~$G$, and an analogous set of conditions for the $\ell$-blocks of 
the double cover $2.F_4(2)$ of $F_4(2)$ containing faithful characters. For 
blocks with cyclic defect groups the inductive Alperin weight condition is known 
to hold by~\cite{KS16}, so that it holds for all blocks of $2.F_4(2)$ containing 
faithful characters by Propositions~\ref{ActionAutomorphisms2F42} 
and~\ref{2f42mod3}.

Let~$b$ denote an $\ell$-block of~$G$ with a non-cyclic defect group. In 
Theorem~\ref{thmAWC} we have already verified Condition~(1) of 
Hypothesis~\ref{defBWC_Simplified}. Thus it remains to verify the equivariance 
Condition~(2) of this hypothesis. If $l(b) = 1$, equivariance trivially 
holds. We will thus assume that $l(b) > 1$ in the following. The orbits of 
$N_{\Aut(G)}( b )$ are described in Propositions~\ref{ActionAutomorphisms}, and 
we will assume these results in what follows.

Suppose that~$b$ is the principal $\ell$-block. For $\ell = 3$, equivariance
follows from Proposition~\ref{EquivariancePBEvenqL3}, and for $\ell > 3$ from 
Corollary~\ref{OrbitOnPBEvenqLg3} and Proposition~\ref{EquivariancePBOddq}.
Now let~$b$ be a unipotent non-principal block. Then $e \in \{ 1, 2 \}$ and
equivariance follows from Proposition~\ref{NonPrincipalUnipotentBlocks} in 
view of the action of~$\sigma_1$ on the set of unipotent characters 
of~$G$ described in \cite[Theorem~$2.5$]{MalleExt}).

From now on we assume that $b$ is a non-unipotent block.

In case $p = 2$ and~$b$ is stabilized by some odd power of $\sigma_1$,
Proposition~\ref{mprimeOdd} establishes our assertion. We may thus assume 
$\sigma = F_1^m$. In this case,
Proposition~\ref{InvariantExtensions} reduces to the instances listed in 
Proposition~\ref{AbelianD}(a),(b). In the latter cases, there is
at least one~$\omega$-invariant extension of~$\theta$ to $N_G(D,b_D)$.
Thus the $\omega$-orbits on the set
of $b$-weights correspond to the $\omega$-orbits on $\Irr( \Out_G( D, b_D ) )$.
By Brauer's permutation lemma, the number of such orbits equals the number of
orbits of~$\omega$ on $\Out_G( D, b_D )$. As there are only two orbit lengths,
the total number of orbits determines the number of orbits of length two.
\end{prf}

\section*{Acknowledgements}
It is our pleasure to thank the many colleagues whose advice and support made 
this work possible. In particular, we thank (in alphabetical order), Robert
Boltje, Cedric Bonnaf{\'e}, Thomas Breuer, Meinolf Geck, Radha Kessar, Markus 
Linckelmann, Gunter Malle and Britta Sp{\"a}th. We are particularly grateful
to Gunter Malle and an anonymous referee for their detailed comments and
suggestions on a previous version, which lead to a considerable improvement of 
the exposition.

The first author was supported by New Zealand Marsden Fund, via award numbers 
UOA 1323 and UOA 1626.
The second author also thanks the Hausdorff Institute of Mathematics, Bonn;
some of this work was done during his stay in October/November 2018
within the Trimester Program ``Logic and Algorithms in Group Theory''.

This work was initiated during two visits of the second author at the University of 
Auckland in 2014 and 2015. This author would like to express his sincere gratitude
to the Department of Mathematics of the University of Auckland for its hospitality.
The visits were financed by New Zealand Marsden Fund (Grant no.\ UOA1323) and the 
author would like to thank this institution for their generous support.

%\backmatter

\renewcommand\thetable{\arabic{table}}
\cleardoublepage

\section{Appendix}

\subsection{The $\ell$-blocks of $F_4(q)$ and their invariants}
\label{ellBlocksAnTheirInvariants}

This first subsection of the appendix gives information on the $\ell$-blocks of
the group $G = F_4(q)$ for odd primes~$\ell$ not dividing~$q$.
There is one table, which may consist of several parts, for $19$ out of the~$20$ 
geometric class types of semisimple $\ell'$-elements, where the table number 
agrees with the number of the geometric class type in~\cite{LL}. 
There is no table for the geometric class type~$20$, as here the centralizers 
are maximal tori, and thus the corresponding blocks have only one irreducible 
Brauer character. As the $\mathbf{G}$-class  
type~$4$ corresponds to elements of order~$3$, Table~\ref{4} is only valid
for primes $\ell > 3$.
There is an additional Table~\ref{21}, which gives the relevant information for 
those $\ell$-blocks of $2.F_4(2)$ with non-cyclic defect groups which
contain faithful characters.

Recall that~$G$ and its dual group~$G^*$ are identified.
Let $s \in G$ be a semisimple $\ell'$-element of $G$-class type~$(i,k)$
and let $b$ be an $\ell$-block inside $\cE_\ell( G, s )$ of positive defect.
Then Table~$i$ describes~$b$ and some of its invariants, unless the 
Sylow~$\ell$-subgroups of~$C_G( s )$ are cyclic.

We use the following notation. First of all, $e = e_\ell(q)$ denotes the order 
of~$q$ in the multiplicative group of $\mathbb{F}_\ell$. Thus, only the cases
$e \in \{ 1, 2, 3, 4, 6 \}$ are relevant, as otherwise the Sylow 
$\ell$-subgroups of~$G$ are cyclic. Moreover, $e > 2$ only occurs in
Tables~\ref{1},~\ref{2} and~\ref{4}, as otherwise the Sylow~$\ell$-subgroups 
of~$C_G( s )$ are cyclic. In all other cases, $e \in \{ 1, 2 \}$,
and then we let the parameter $\varepsilon \in \{ -1, 1 \}$ be such that 
$\ell \mid q - \varepsilon$. Also, if~$\ell$ is understood from the
context,~$a$ denotes the positive integer 
such that~$\ell^a$ is the highest power of~$\ell$ dividing $\Phi_e( q )$
(recall that $\Phi_e$ denotes the $e$th cyclotomic polynomial).
For example, if $\ell = 3$, the $3$-part $|G|_3$ of $|G|$ equals $3^{4a+2}$.

Let us now explain the structure and the contents of the tables.
These have up to $13$ columns, numbered $1, \ldots , 13$ in the first row 
of the table. The column number determines the content type of the column; 
if this content is not relevant for a table, the corresponding column is omitted.
The second row of the tables contains the column headings, and the remaining
rows between the second and last set of double rules constitute the body of
the tables, containing the desired information on the blocks.

The contents of the table is separated by horizontal rules according to the
following scheme: Rules beginning in Column~$6$ separate the information for
different ranges of~$\ell$, if necessary. Rules beginning in Column~$5$ separate 
the information for different blocks; it turns out that $\cE_\ell( G, s )$ 
contains at most two blocks of positive defect, except in the case $e = 4$.
Rules beginning in Column~$3$ separate the information for the two values 
of~$\varepsilon$ in case $e \in \{ 1, 2 \}$. If such a rule is present, 
Column~$2$ usually contains two values for~$k$, say 
$k_1, k_2$, separated by a comma. The first set of rows between
this horizontal rule beginning in Column~$3$ corresponds to $k_1$,
if~$\varepsilon = 1$, and to $k_2$, otherwise. Likewise, the second set of rows
corresponds to $k_2$, if $\varepsilon = 1$, and to $k_1$, otherwise. 
An analogous convention is used if there are two entries in Column~$9$.
If the entries in Column~$2$ are separated by a slash, there is a 
corresponding and matching pair of entries in Column~$9$.
If their are four entries in Column~$2$, they are grouped in pairs, separated 
by a slash, and the conventions above apply to each of these pairs. 
Finally, a rule beginning in Column~$2$ separates the information for the 
different values of~$k$. The values of~$k$ are contained in the second column.
An entry in a row of the body of the table is effective for all rows below it
up to the next entry or to a separating horizontal rule.

Columns~$3$ and~$4$ describe the structure of
the centralizer $C := C_G( s )$ by printing the order of $Z := Z(C)$, as well 
as the components of $[ C_{\mathbf{G}}( s ), C_{\mathbf{G}}( s ) ]^F$, under 
the heading $[C,C]$. These components are described by their Dynkin type. 
We use the convention that $\tilde{A}_j$ designates a component of type~$A_j$ 
corresponding to short root elements. Finally, we write $A_j^\varepsilon(q)$
and $\tilde{A}_j^\varepsilon(q)$ to denote a linear group if $\varepsilon = 1$,
and a unitary group if $\varepsilon = -1$.

Columns~$5$,~$6$ and~$7$ give a
label for~$b$, its defect $d(b)$ and the number $l(b)$ of its irreducible
Brauer characters, respectively. The label of~$b$ describes an $e$-cuspidal pair
associated to~$b$ in the following way. Suppose first that~$s$ is not quasi-isolated, 
i.e.\ the centralizer $C_{\mathbf{G}}( s )$ is 
contained in a proper Levi subgroup of~$\mathbf{G}$. This is the case if and 
only if~$s$ belongs to a class type $i \geq 6$. As all proper regular subgroups 
of~$\mathbf{G}$ are of classical type, we may apply Theorem~\ref{BoRoEtAl}. 
Thus~$b$ corresponds to a unipotent block~$b'$ of 
$\mathcal{E}_\ell( C_G( s ), 1 )$, and the label of~$b$ is taken to be the 
$e$-cuspidal pair $(L,\zeta)$ associated to~$b'$; see \cite[Theorem]{cab}.
This label is given as $L, \zeta$ in Column~$5$ of Tables~\ref{1}--\ref{19},
unless~$L$ is a torus (and thus $\zeta$ the trivial character), where we just 
write~$1$ and suppress~$L$. If $L = B_2(q)$, we write $\zeta_1$ for the cuspidal 
unipotent character,~$\zeta_2$ for the $2$-cuspidal unipotent character labelled 
by the bipartition~$(1,1)$, and~$\zeta_4$,~$\zeta_4'$ for the two $4$-cuspial 
unipotent characters labelled by the bipartitions~$(1^2,-)$ and $(-,2)$, 
respectively. In case of
$L = A_2^\varepsilon( q )$, we write $(2,1)$ for the $e$-cuspidal unipotent
character labelled by the partition $(2,1)$. 
If~$s$ is quasi-isolated, i.e.\ $s \in \{ 1, 2, 3, 5 \}$, the label of~$b$ describes
the $e$-cuspidal pair of~$G$ defining~$b$, following \cite[p.~$349$]{En00} 
(for $i = 1$) and 
\cite[Table~$2$]{KeMa} (for $i \in \{ 2, 3, 5 \}$), respectively.

Column~$8$ describes the isomorphism type of a defect group $D(b)$ of~$b$, in 
case $D(b)$ is abelian. If $D(b)$ is non-abelian, then $\ell = 3$ and $D(b)$ is
$G$-conjugate to a Sylow $3$-subgroup of a group $M_{i',k'}$ dual to $C_G( s )$;
in this case $D(b)$ is identified by the pair $(i',k')$. (In case $i = 13$ or $i = 17$
and $k \in \{ 2, 3 \}$, we use the convention that $(i',k')$ is such that
$\{ i, i' \} = \{ 13, 17 \}$ and $\{ k, k' \} = \{ 2, 3 \}$ but $i \neq i'$ and
$k \neq k'$.)

Column~$9$
contains the weight subgroups of~$b$, described by their conjugacy class in the set
of radical $\ell$-subgroups. The first row corresponding to~$b$ gives its defect group.
The radical $3$-subgroups occurring as weight subgroups in the tables below are 
named according to Table~\ref{tabradicalageq2}. In order to be consistent with 
this labelling, the defect groups in case $\ell > 3$ and $e \in \{ 1, 2 \}$ are 
named as $R_{j,\ell}$, with $j \in \{ 2, 3, 9, 10, 11, 12, 16, 17, 18 \}$. These 
groups are Sylow $\ell$-subgroups of suitable maximal tori $M_{20,k}$, 
with~$j$ and~$k$ as in the following table, where the two cases for~$k$ in 
some of the columns correspond to the cases $\varepsilon = 1$ and 
$\varepsilon = -1$, respectively.
\setlength{\extrarowheight}{0.5ex}
$$
\begin{array}{l|ccccccccc} \hline\hline
	j & 2 & 3 & 9 & 10 & 11 & 12 & 16 & 17 & 18 \\ \hline
	k & 21 & 16 & 22 & 3 & 7, 8 & 4, 5 & 12, 13 & 17, 18 & 1, 2 \\ \hline\hline
\end{array}
$$
If the description of a block is restricted to the case $\ell = 3$, the second 
index in $R_{j,\ell}$ is omitted. In the tables for $e > 2$, we give the defect 
groups as Sylow $\ell$-subgroups of a maximal torus $M_{20,k}$, denoted 
by~$S_{k,\ell}$.

Column~$10$ gives the structure of $\Out_G( R, b_R ) = N_G( R, b_R )/RC_G(R)$ 
for the weight subgroups $R$ of~$b$, and Column~$11$, headed by~$\cW$ contains 
the numbers
$\cW( R, b_R )$; see (\ref{NumberOfWeights0}) in Subsection~\ref{BlocksAndWeights} 
for this notation. In Columns~$12$ and~$13$, if present, we remark a case
distinction between $a = 1$ and $a \geq 2$, respectively between $\ell = 3$,
$\ell > 3$ and $\ell \geq 3$.

Our naming of the groups follows standard conventions.
For example, $S_3$, $D_8$ and $Q_8$ designate the symmetric group on three letter, 
the dihedral and the quaternion group, respectively, of order~$8$. Further, 
$W( X_j )$ denotes the Weyl group of the root system $X_j$. 

\setlength{\extrarowheight}{0.5ex}
\setcounter{table}{0}

\begin{landscape}
\begin{table}[h]
\caption{\label{1} The unipotent $\ell$-blocks of $F_4(q)$ of positive defect}
$$
\begin{array}{cccccccccccc} 
\multicolumn{12}{c}{\ell = 3} \\ \hline\hline
1 & 2 & 3 & 4 & 5 & 6 & 7 & 8 & 9 & 10 & 11 & 12 \\ 
 i & k & |Z| & [C,C] & b & d(b) & l(b) & D(b) & R & \Out(b_R) & \cW & \text{Rem}
\\ \hline\hline
	1 & 1 & 1 & F_4(q) & {1} & 4a+2 & 26 & (1,1) & R_{38} & 2^3 & 8 & \null \\
 & & & & & & & & R_{37} & \SL_2(3).2 & 2 & \null \\
 & & & & & & & & R_{15} & \SL_3(3) & 1 & \null \\
 & & & & & & & & R_{18} & W(F_4) & 4 & \null \\
 & & & & & & & & R_{21} & (Q_8\times Q_8).S_3 & 11 & a = 1 \\
 & & & & & & & & R_{21} & (\SL_2(3)\times \SL_2(3)).2 & 2 & a \geq 2 \\
 & & & & & & & & R_{22} & \SL_2(3)\times \SL_2(3) & 1 & a \geq 2 \\
 & & & & & & & & R_{35} & (\SL_2(3)\times 2).2 & 4 & a \geq 2 \\
 & & & & & & & & R_{36} & (2\times \SL_2(3)).2 & 4 & a \geq 2 \\ \cline{5-12}
 & & & & B_2(q), \zeta_e & 2a & 5 & [3^a]^2 & R_{10} & D_8 & 5 & \null \\ \hline\hline
\end{array}
$$
\end{table}
\end{landscape}

\clearpage

\begin{landscape}
\setcounter{table}{0}
\begin{table}[h]
\caption{\label{a1} The unipotent $\ell$-blocks of $F_4(q)$ of positive defect (cont.)}
$$
\begin{array}{ccccccccccc} 
\multicolumn{11}{c}{e = 1, 2,\quad \ell > 3} \\ \hline\hline
1 & 2 & 3 & 4 & 5 & 6 & 7 & 8 & 9 & 10 & 11 \\
 i & k & |Z| & [C,C] & b & d(b) & l(b) & D(b) & R & \Out(b_R) & \cW 
\\ \hline\hline
1 & 1 & 1 & F_4(q) & {1} & 4a & 26 & [\ell^a]^4 & R_{18,\ell} & W(F_4) & 26 \\ \cline{5-11}
 & & & & B_2(q), \zeta_e & 2a & 5 & [\ell^a]^2 & R_{10,\ell} & D_8 & 5 \\ \hline\hline
\end{array}
$$

\vspace*{0.6cm}

$$
\begin{array}{ccccccccccc}
\multicolumn{11}{c}{e = 3, 6} \\ \hline\hline
1 & 2 & 3 & 4 & 5 & 6 & 7 & 8 & 9 & 10 & 11 \\
 i & k & |Z| & [C,C] & b & d(b) & l(b) & D(b) & R & \Out(b_R) & \cW
\\ \hline\hline
1 & 1 & 1 & F_4(q) & {1} & 2a & 21 & [\ell^a]^2 & S_{9,\ell}, S_{10,\ell} & \SL_2(3) \times 3 & 21 \\ \hline\hline
\end{array}
$$

\vspace*{0.6cm}

$$
\begin{array}{ccccccccccc}
\multicolumn{11}{c}{e = 4} \\ \hline\hline
1 & 2 & 3 & 4 & 5 & 6 & 7 & 8 & 9 & 10 & 11 \\
 i & k & |Z| & [C,C] & b & d(b) & l(b) & D(b) & R & \Out(b_R) & \cW
\\ \hline\hline
1 & 1 & 1 & F_4(q) & {1} & 2a & 16 & [\ell^a]^2 & S_{6,\ell} & \SL_2(3).[4] & 16 \\ \cline{5-11}
& & & & B_2(q), \zeta_4 & a & 4 & [\ell^a] & R_{10,\ell} & [4] & 4 \\ \cline{5-11}
& & & & B_2(q), \zeta_4' & a & 4 & [\ell^a] & R_{10,\ell} & [4] & 4 \\ \hline\hline
\end{array}
$$

\end{table}
\end{landscape}

\clearpage

\begin{landscape}
\begin{table}[h]
\caption{\label{2} The $\ell$-blocks of $F_4(q)$ of geometric type~$2$
(exist only when~$q$ is odd)}

$$
\begin{array}{cccccccccccc} 
\multicolumn{12}{c}{e = 1, 2} \\ \hline\hline
1 & 2 & 3 & 4 & 5 & 6 & 7 & 8 & 9 & 10 & 11 & 13 \\ 
i & k & |Z| & [C,C] & b & d(b) & l(b) & D(b) & R & \Out(b_R) & \cW & \ell
\\ \hline\hline
2 & 1 & 2 & B_4(q) & {1} & 4a+1 & 20 & (5,1) & R_{33} & 2^3 & 8 & 3 \\
& & & & & & & & R_{25} & (2 \times \SL_2(3)).2 & 4 & \\ 
& & & & & & & & R_{18} & W(B_4) & 8 & \\ \cline{6-12}
& & & & & 4a & 20 & [\ell^a]^4 & R_{18,\ell} & W(B_4) & 20 & > 3 \\ \cline{5-12}
& & & & B_2(q), \zeta_e & 2a & 5 & [\ell^a]^2 & R_{10,\ell} & D_8 & 5 & \geq 3 \\ \hline\hline
\end{array}
$$

\vspace*{1cm}

$$
\begin{array}{ccccccccccc}
\multicolumn{11}{c}{e = 4} \\ \hline\hline
1 & 2 & 3 & 4 & 5 & 6 & 7 & 8 & 9 & 10 & 11 \\
 i & k & |Z| & [C,C] & b & d(b) & l(b) & D(b) & R & \Out(b_R) & \cW
\\ \hline\hline
2 & 1 & 2 & B_4(q) & {1} & 2a & 14 & [\ell^a]^2 & S_{6,\ell} & [4]^2.2 & 14 \\ \cline{5-11}
& & & & B_2(q), \zeta_4 & a & 4 & [\ell^a] & R_{10,\ell} & [4] & 4  \\ \cline{5-11}
& & & & B_2(q), \zeta_4' & a & 4 & [\ell^a] & R_{10,\ell} & [4] & 4  \\ \hline\hline
\end{array}
$$

\end{table}
\end{landscape}

\clearpage

\begin{landscape}
\begin{table}[h]
\caption{\label{3} The $\ell$-blocks of $F_4(q)$ of geometric type~$3$
(exist only when~$q$ is odd)}

$$
\begin{array}{cccccccccccc} 
\multicolumn{12}{c}{e = 1, 2} \\ \hline\hline
1 & 2 & 3 & 4 & 5 & 6 & 7 & 8 & 9 & 10 & 11 & 13 \\ 
i & k & |Z| & [C,C] & b & d(b) & l(b) & D(b) & R & \Out(b_R) & \cW & \ell 
\\ \hline\hline
3 & 1, 2 & 4 & A_3^{\varepsilon}(q)\tilde{A}_1(q) & {1} & 4a+1 & 10 & (5,1) & 
R_{33} & 2^2 & 4 & 3 \\
& & & & & & & & R_{25} & 2 \times \SL_2(3) & 2 & \\ 
& & & & & & & & R_{18} & W(A_3) \times 2 & 4 & \\ \cline{6-12}
& & & & & 4a & 10 & [\ell^a]^4 & R_{18,\ell} & W(A_3) \times 2 & 10 & > 3 \\ \cline{4-12}
& & & A_3^{-\varepsilon}(q)\tilde{A}_1(q) & {1} & 3a & 10 & [\ell^a]^3 & R_{16,\ell} & D_8 \times 2 & 10 & \geq 3
\\ \hline\hline
\end{array}
$$
\end{table}
\end{landscape}

\clearpage

\begin{landscape}
\begin{table}[h]
\caption{\label{4} The $\ell$-blocks of $F_4(q)$ of geometric type~$4$
(exist only when $3 \nmid q$ and $\ell > 3$)}

$$
\begin{array}{ccccccccccc} 
\multicolumn{11}{c}{e = 1, 2} \\ \hline\hline
1 & 2 & 3 & 4 & 5 & 6 & 7 & 8 & 9 & 10 & 11 \\
 i & k & |Z| & [C,C] & b & d(b) & l(b) & D(b) & R & \Out(b_R) & \cW
\\ \hline\hline
4 & 1, 2 & 3 & A_2^{\varepsilon}(q)\tilde{A}^{\varepsilon}_2(q) & {1} & 4a & 9 & [\ell^a]^4 &
R_{18,\ell} & S_3 \times S_3 & 9 \\ \cline{4-11}
& & & A_2^{-\varepsilon}(q)\tilde{A}^{-\varepsilon}_2(q) & {1} & 2a & 4 & [\ell^a]^2 &
R_{9,\ell} & 2^2 & 4 
\\ \hline\hline
\end{array}
$$

\vspace*{1.0cm}

$$
\begin{array}{ccccccccccc}  
\multicolumn{11}{c}{e = 3} \\ \hline\hline
1 & 2 & 3 & 4 & 5 & 6 & 7 & 8 & 9 & 10 & 11 \\
 i & k & |Z| & [C,C] & b & d(b) & l(b) & D(b) & R & \Out(b_R) & \cW
\\ \hline\hline
4 & 1 & 3 & A_2(q)\tilde{A}_2(q) & {1} & 2a & 9 & [\ell^a]^2 & S_{9,\ell} & 3^2 & 9
\\ \hline\hline
\end{array}
$$

\vspace*{1.0cm}

$$
\begin{array}{ccccccccccc}  
\multicolumn{11}{c}{e = 6} \\ \hline\hline
1 & 2 & 3 & 4 & 5 & 6 & 7 & 8 & 9 & 10 & 11 \\
 i & k & |Z| & [C,C] & b & d(b) & l(b) & D(b) & R & \Out(b_R) & \cW
\\ \hline\hline
4 & 2 & 3 & A_2^{-1}(q)\tilde{A}^{-1}_2(q) & {1} & 2a & 9 & [\ell^a]^2 & S_{10,\ell} & 3^2 & 9
\\ \hline\hline
\end{array}
$$
\end{table}
\end{landscape}

\clearpage

\begin{landscape}
%\addtocounter{table}{1}
\begin{table}[h]
\caption{\label{5} The $\ell$-blocks of $F_4(q)$ of geometric type~$5$
(exist only when~$q$ is odd)}

$$
\begin{array}{cccccccccccc}
\multicolumn{12}{c}{e = 1, 2} \\ \hline\hline
1 & 2 & 3 & 4 & 5 & 6 & 7 & 8 & 9 & 10 & 11 & 13 \\ 
i & k & |Z| & [C,C] & b & d(b) & l(b) & D(b) & R & \Out(b_R) & \cW & \ell
\\ \hline\hline
5 & 1 & 2 & C_3(q){A}_1(q) & {1} & 4a+1 & 20 & (2,1) & R_{34} & 2^3 & 8 & 3 \\
& & & & & & & & R_{26} & (\SL_2(3) \times 2).2 & 4 & \\
& & & & & & & & R_{18} & W(C_3) \times 2 & 8 & \\ \cline{6-12}
& & & & & 4a & 20 & [\ell^a]^4 & R_{18,\ell} & W(C_3) \times 2 & 20 & > 3 \\ \cline{5-12}
& & & & B_2(q), \zeta_e & 2a & 4 & [\ell^a]^2 & R_{10,\ell} & 2^2 & 4 & \geq 3 \\ \hline\hline
\end{array}
$$
\end{table}
\end{landscape}

\clearpage

\begin{landscape}
\begin{table}[h]
\caption{\label{6} The $\ell$-blocks of $F_4(q)$ of geometric type~$6$}
$$
\begin{array}{ccccccccccccc} 
\multicolumn{13}{c}{e = 1, 2} \\ \hline\hline
1 & 2 & 3 & 4 & 5 & 6 & 7 & 8 & 9 & 10 & 11 & 12 & 13 \\ 
i & k & |Z| & [C,C] & b & d(b) & l(b) & D(b) & R & \Out(b_R) & \cW & \text{Rem} & \ell
\\ \hline\hline
6 & 1, 2 & q - \varepsilon & B_3(q) & {1} & 4a+1 & 10 & (10,k) & 
R_{33} & 2^2 & 4 & & 3 \\ 
& & & & & & & & R_{25} & \SL_2(3).2 & 2 & \\ 
& & & & & & & & R_{18} & W(C_3) & 4 & \\ \cline{6-13}
& & & & & 4a & 10 & [\ell^a]^4 & R_{18,\ell} & W(C_3) & 10 & & > 3 \\ \cline{5-13}
& & & & B_2(q), \zeta_e & 2a & 2 & [\ell^a]^2 & R_{10,\ell} & 2 & 2 & & \geq 3 \\ \cline{3-13}
& & q + \varepsilon & B_3(q) & {1} & 3a+1 & 10 & (10,k) & R_{31} & 2^2 & 4 & & 3 \\
& & & & & & & & R_{23} & \SL_2(3).2 & 2 & a \geq 2 & \\
& & & & & & & & R_{19} & \SL_2(3).2 & 2 & a = 1 & \\
& & & & & & & & R_{16} & W(C_3) & 4 & \\ \cline{6-13}
& & & & & 3a & 10 & [\ell^a]^3 & R_{16,\ell} & W(C_3) & 10 & & > 3 \\ \cline{5-13}
& & & & B_2(q), \zeta_e & a & 2 & [\ell^a] & R_{2,\ell} & 2 & 2 & &  \geq 3 \\ \hline\hline
\end{array}
$$
\end{table}
\end{landscape}

\clearpage

\begin{landscape}
\begin{table}[h]
\caption{\label{7} The $\ell$-blocks of $F_4(q)$ of geometric type~$7$}
$$
\begin{array}{cccccccccccc} 
\multicolumn{12}{c}{e = 1, 2} \\ \hline\hline
1 & 2 & 3 & 4 & 5 & 6 & 7 & 8 & 9 & 10 & 11 & 13 \\ 
i & k & |Z| & [C,C] & b & d(b) & l(b) & D(b) & R & \Out(b_R) & \cW & \ell
\\ \hline\hline
 7 & 1, 2 & q - \varepsilon & A_2^\varepsilon(q)\tilde{A}_1(q) & {1} & 
4a+1 & 6 & (9,k) & R_{33} & 2^2 & 4 & 3 \\
 & & & & & & & & R_{25} & 2 \times \SL_2(3) & 2 \\ \cline{6-12}
& & & & & 4a & 6 & [\ell^a]^4 & R_{18,\ell} & 2 \times S_3 & 6 & > 3 \\ \cline{3-12}
 & & q + \varepsilon &  A_2^{-\varepsilon}(q)\tilde{A}_1(q) & {1} & 
2a & 4 & [\ell^a]^2 & R_{9,\ell} & 2^2 & 4 & \geq 3 \\ \cline{5-12}
& & & & A_2^{-\varepsilon}(q), (2,1) & a & 2 & [\ell^a] & R_{2,\ell} & 2 & 2 & \geq 3 \\ \hline\hline
\end{array}
$$
\end{table}
\end{landscape}

\clearpage

\begin{landscape}
\begin{table}[h]
\caption{\label{8} The $\ell$-blocks of $F_4(q)$ of geometric type~$8$ 
(exist only when~$q$ is odd)}
$$
\begin{array}{ccccccccccccc} 
\multicolumn{13}{c}{e = 1, 2} \\ \hline\hline
1 & 2 & 3 & 4 & 5 & 6 & 7 & 8 & 9 & 10 & 11 & 12 & 13 \\ 
i & k & |Z| & [C,C] & b & d(b) & l(b) & D(b) & R & \Out(b_R) & \cW & \text{Rem} & \ell 
\\ \hline\hline
8 & 1, 4 & 2(q - \varepsilon) & A_3^{\varepsilon}(q) & {1} & 4a + 1 & 5 & (5,1) & 
R_{33} & 2 & 2 & & 3 \\
 & & & & & & & & R_{25} & \SL_2(3) & 1 \\
 & & & & & & & & R_{18} & W(A_3) & 2 \\ \cline{6-13}
& & & & & 4a & 5 & [\ell^a]^4 & R_{18,\ell} & W(A_3) & 2 & & > 3 \\ \cline{3-13}
 & & 2(q + \varepsilon) & A_3^{-\varepsilon}(q) & {1} & 2a & 5 & [\ell^a]^2 & 
R_{10,\ell} & D_8 & 5 & & \geq 3 \\ \cline{2-13}
	& 2, 3 & 2(q + \varepsilon) & A_3^{\varepsilon}(q) & {1} & 3a + 1 & 5 & (10,2), (10,1) & 
R_{31} & 2 & 2 & & 3 \\
 & & & & & & & & R_{23} & \SL_2(3) & 1 & a \geq 2 \\
 & & & & & & & & R_{19} & \SL_2(3) & 1 & a = 1 \\
 & & & & & & & & R_{16} & W(A_3) & 2 \\ \cline{6-13}
& & & & & 3a & 5 & [\ell^a]^3 & R_{16,\ell} & W(A_3) & 2 & & > 3 \\ \cline{3-13}
 & & 2(q - \varepsilon) & A_3^{-\varepsilon}(q) & {1} & 3a & 5 &  [\ell^a]^3 & 
R_{16,\ell} & D_8 & 5 & & \geq 3 \\ \hline\hline
\end{array}
$$
\end{table}
\end{landscape}

\clearpage

\begin{landscape}
\begin{table}[h]
\caption{\label{9} The $\ell$-blocks of $F_4(q)$ of geometric type~$9$}
$$
\begin{array}{cccccccccccc} 
\multicolumn{12}{c}{e = 1, 2} \\ \hline\hline
1 & 2 & 3 & 4 & 5 & 6 & 7 & 8 & 9 & 10 & 11 & 13 \\ 
i & k & |Z| & [C,C] & b & d(b) & l(b) & D(b) & R & \Out(b_R) & \cW & \ell
\\ \hline\hline
9 & 1, 2 & q - \varepsilon & \tilde{A}_2^{\varepsilon}(q)A_1(q) & {1} & 4a+1 & 6 & (7,k) & 
R_{34} & 2^2 & 4 & 3 \\
& & & & & & & & R_{26} & \SL_2(3) \times 2 & 2 & \\ \cline{6-12}
& & & & & 4a & 6 & [\ell^a]^4 & R_{18,\ell} & S_3 \times 2 & 6 & > 3 \\ \cline{3-12}
& & q + \varepsilon & \tilde{A}_2^{-\varepsilon}(q)A_1(q) & {1} & 2a & 4 & [\ell^a]^2 & 
R_{9,\ell} & 2^2 & 4 & \geq 3 \\ \cline{5-12}
& & & & \tilde{A}_2^{-\varepsilon}(q), (2,1) & a & 2 & [\ell^a] & R_{3,\ell} & 2 & 2 & \geq 3 \\ \hline\hline 
\end{array}
$$
\end{table}
\end{landscape}

\clearpage

\begin{landscape}
\begin{table}[h]
\caption{\label{10} The $\ell$-blocks of $F_4(q)$ of geometric type~$10$}
$$
\begin{array}{ccccccccccccc} 
\multicolumn{13}{c}{e = 1, 2} \\ \hline\hline
1 & 2 & 3 & 4 & 5 & 6 & 7 & 8 & 9 & 10 & 11 & 12 & 13 \\ 
i & k & |Z| & [C,C] & b & d(b) & l(b) & D(b) & R & \Out(b_R) & \cW & \text{Rem} & \ell
\\ \hline\hline
10 & 1, 2 & q - \varepsilon & C_3(q) & {1} & 4a+1 & 10 & (6,k) &
R_{34} & 2^2 & 4 & & 3 \\
 & & & & & & & & R_{26} & \SL_2(3).2 & 2 \\ 
 & & & & & & & & R_{18} & W(C_3) & 4 \\ \cline{6-13}
& & & & & 4a & 10 & [\ell^a]^4 & R_{18,\ell} & W(C_3) & 10 & & > 3 \\ \cline{5-13}
& & & & B_2(q), \zeta_e & 2a & 2 & [\ell^a]^2 & R_{10,\ell} & 2 & 2 & & \geq 3 \\ \cline{3-13}
& & q + \varepsilon & C_3(q) & {1} & 3a+1 & 10 & (6,k) & 
R_{32} & 2^2 & 4 & & 3 \\
 & & & & & & & & R_{24} & \SL_2(3).2 & 2 & a \geq 2 \\
 & & & & & & & & R_{20} & \SL_2(3).2 & 2 & a = 1 \\
 & & & & & & & & R_{17} & W(C_3) & 4 \\ \cline{6-13}
 & & & & & 3a & 10 & [\ell^a]^3 & R_{17,\ell} & W(C_3) & 10 & & > 3 \\ \cline{5-13}
& & & & B_2(q), \zeta_e & a & 2 & [\ell^a] & R_{3,\ell} & 2 & 2 & & \geq 3 \\ \hline\hline
\end{array}
$$
\end{table}
\end{landscape}

\clearpage

\begin{landscape}
\begin{table}[h]
\caption{\label{11} The $\ell$-blocks of $F_4(q)$ of geometric type~$11$ 
(exist only when~$q$ is odd)}
$$
\begin{array}{cccccccccccc} 
\multicolumn{12}{c}{e = 1, 2} \\ \hline\hline
1 & 2 & 3 & 4 & 5 & 6 & 7 & 8 & 9 & 10 & 11 & 13 \\ 
i & k & |Z| & [C,C] & b & d(b) & l(b) & D(b) & R & \Out(b_R) & \cW & \ell 
\\ \hline\hline
 11 & 1, 2 & 2(q - \varepsilon) & C_2(q)A_1(q) & {1} & 4a & 10 & [\ell^a]^4 & 
R_{18,\ell} & D_8 \times 2 & 10 & \geq 3 \\ \cline{5-12}
& & &  & B_2(q), \zeta_e & 2a & 2 & [\ell^a]^2 & R_{10,\ell} & 2 & 2 & \geq 3 \\ \cline{3-12}
 & & 2(q + \varepsilon) & C_2(q)A_1(q) & {1} & 3a & 10 & [\ell^a]^3 & 
R_{17,\ell} & D_8 \times 2 & 10 & \geq 3 \\ \cline{5-12}
& & & & B_2(q), \zeta_e & a & 2 & [\ell^a] & R_{3,\ell} & 2 & 2 & \geq 3 \\ \hline\hline
\end{array}
$$
\end{table}
\end{landscape}

\clearpage

\begin{landscape}
\begin{table}[h]
\caption{\label{12} The $\ell$-blocks of $F_4(q)$ of geometric type~$12$ 
(exist only when~$q$ is odd)}
$$
\begin{array}{cccccccccccc}
\multicolumn{12}{c}{e = 1, 2} \\ \hline\hline
1 & 2 & 3 & 4 & 5 & 6 & 7 & 8 & 9 & 10 & 11 & 13 \\ 
i & k & |Z| & [C,C] & b & d(b) & l(b) & D(b) & R & \Out(b_R) & \cW & \ell
\\ \hline\hline
 12 & 1, 3 & 2(q - \varepsilon) & A_1(q)A_1(q)\tilde{A}_1(q) & {1} & 4a & 8 & [\ell^a]^4 & 
R_{18,\ell} & 2^3 & 8 & \geq 3 \\ \cline{3-12}
 & & 2(q + \varepsilon) & A_1(q)A_1(q)\tilde{A}_1(q) & {1} & 3a & 8 & [\ell^a]^3 & 
R_{16,\ell} & 2^3 & 8 & \geq 3 \\ \cline{2-12}
 & 2, 4 & 2(q - \varepsilon) & A_1(q^2)\tilde{A}_1(q) & {1} & 3a & 4 & [\ell^a]^3 & 
R_{16,\ell} & 2^2 & 4 & \geq 3 \\ \cline{3-12}
 & & 2(q + \varepsilon) & A_1(q^2)\tilde{A}_1(q) & {1} & 2a & 4 &  [\ell^a]^2 & 
R_{10,\ell} & 2^2 & 4 & \geq 3 \\ \hline\hline
\end{array}
$$
\end{table}
\end{landscape}

\clearpage

\begin{landscape}
\begin{table}[h]
\caption{\label{13} The $\ell$-blocks of $F_4(q)$ of geometric type~$13$}
$$
\begin{array}{ccccccccccccc} 
\multicolumn{13}{c}{e = 1, 2} \\ \hline\hline
1 & 2 & 3 & 4 & 5 & 6 & 7 & 8 & 9 & 10 & 11 & 12 & 13 \\ 
i & k & |Z| & [C,C] & b & d(b) & l(b) & D(b) & R & \Out(b_R) & \cW & \text{Rem} & \ell 
\\ \hline\hline
13 & 1, 6 & (q - \varepsilon)^2 & A_2^{\varepsilon}(q) & {1} & 4a+1 & 3 & (17,k) & 
R_{33} & 2 & 2 & & 3 \\
 & & & & & & & & R_{25} & \SL_2(3) & 1 \\ \cline{6-13}
& & & & & 4a & 3 & [\ell^a]^4 & R_{18,\ell} & S_3 & 3 & & > 3 \\ \cline{3-13}
& & (q + \varepsilon)^2 & A_2^{-\varepsilon}(q) & {1} & a & 2 & [\ell^a] & R_{3,\ell} & 2 & 2 & & \geq 3 \\ \cline{2-13}
%& & & & A_2^{-\varepsilon}(q), (2,1) & 0 & 1 & 1 & 1 & 1 & 1 & & \geq 3 \\ \cline{2-13} 
& 2, 3 & q^2 - 1 & A_2^{\varepsilon}(q) & {1} & 3a+1 & 3 & (17,k') & 
R_{31} & 2 & 2 & & 3 \\
 & & & & & & & & R_{23} & \SL_2(3) & 1 & a \geq 2 \\
 & & & & & & & & R_{19} & \SL_2(3) & 1 & a = 1 \\ \cline{6-13}
& & & & & 3a & 3 & [\ell^a]^3 & R_{16,\ell} & S_3 & 3 & & > 3 \\ \cline{3-13}
 & & q^2 - 1 & A_2^{-\varepsilon}(q) & {1} & 2a & 2 & [\ell^a]^2 & 
R_{9,\ell} & 2 & 2 & & \geq 3 \\ \cline{5-13}
& & & & A_2^{-\varepsilon}(q), (2,1) & a & 1 & [\ell^a] & R_{2,\ell} & 1 & 1 & & \geq 3 \\ \cline{2-13}
	& 5, 4 & q^2 + \varepsilon q + 1 & A_2^{\varepsilon}(q) & {1} & 2a+2 & 3 &  (17,k) & 
R_{29} & 2 & 2 & & 3 \\
 & & & & & & & & R_{19} & \SL_2(3) & 1 \\ \cline{6-13}
& & & & & 2a & 3 & [\ell^a]^2 & R_{12,\ell} & S_3 & 3 & & > 3 \\ \cline{3-13}
& & q^2 - \varepsilon q + 1 & A_2^{-\varepsilon}(q) & {1} & a & 2 & [\ell^a] & R_{3,\ell} & 2 & 2 & & \geq 3 \\ \hline\hline
% \cline{5-13}
%& & & & A_2^{-\varepsilon}(q), (2,1) & 0 & 1 & 1 & 1 & 1 & 1 & & \geq 3 \\ \hline\hline
\end{array}
$$
\end{table}
\end{landscape}

\clearpage

\begin{landscape}
\begin{table}[h]
\caption{\label{14} The $\ell$-blocks of $F_4(q)$ of geometric type~$14$}
$$
\begin{array}{cccccccccccc}
\multicolumn{12}{c}{e = 1, 2} \\ \hline\hline
1 & 2 & 3 & 4 & 5 & 6 & 7 & 8 & 9 & 10 & 11 & 13 \\ 
i & k & |Z| & [C,C] & b & d(b) & l(b) & D(b) & R & \Out(b_R) & \cW & \ell 
\\ \hline\hline
 14 & 1, 4 & (q - \varepsilon)^2 & A_1(q)\tilde{A}_1(q) & {1} & 4a & 4 & [\ell^a]^4 & 
R_{18,\ell} & 2^2 & 4 & \geq 3 \\ \cline{3-12}
 & & (q + \varepsilon)^2 & A_1(q)\tilde{A}_1(q) & {1} & 2a & 4 & [\ell^a]^2 & 
R_{9,\ell} & 2^2 & 4 & \geq 3 \\ \cline{2-12}
 & 2 & q^2 - 1 & A_1(q)\tilde{A}_1(q) & {1} & 3a & 4 & [\ell^a]^3 &  
R_{16,\ell}, R_{17,\ell} & 2^2 & 4 & \geq 3 \\ \cline{2-12}
 & 3 & q^2 - 1 & A_1(q)\tilde{A}_1(q) & {1} & 3a & 4 & [\ell^a]^3 &  
R_{17,\ell}, R_{16,\ell} & 2^2 & 4 & \geq 3 \\ \hline\hline
\end{array}
$$
\end{table}
\end{landscape}

\clearpage

\begin{landscape}
\begin{table}[h]
\caption{\label{15} The $\ell$-blocks of $F_4(q)$ of geometric type~$15$}
$$
\begin{array}{cccccccccccc} 
\multicolumn{12}{c}{e = 1, 2} \\ \hline\hline
1 & 2 & 3 & 4 & 5 & 6 & 7 & 8 & 9 & 10 & 11 & 13 \\ 
i & k & |Z| & [C,C] & b & d(b) & l(b) & D(b) & R & \Out(b_R) & \cW & \ell
\\ \hline\hline
 15 & 1, 3 & (q - \varepsilon)^2 & C_2(q) &  {1} & 4a & 5 & [\ell^a]^4 & 
R_{18,\ell} & D_8 & 5 & \geq 3 \\ \cline{5-12}
& & & & B_2(q), \zeta_e & 2a & 1 & [\ell^a]^2 & R_{10,\ell} & 1 & 1 & \geq 3 \\ \cline{3-12}
 & & (q + \varepsilon)^2 & C_2(q) & {1} & 2a & 5 &  [\ell^a]^2 & 
R_{10,\ell} & D_8 & 5 & \geq 3 \\ \cline{2-12}
%& & & & B_2(q), \zeta_e & 0 & 1 & 1 & 1 & 1 & 1 & \geq 3 \\ \cline{2-12}
& 2/4 & q^2 - 1 & C_2(q) & {1} & 3a & 5 & [\ell^a]^3 & R_{17,\ell}/R_{16,\ell} & D_8 & 5 & \geq 3 \\ \cline{5-12}
& & & & B_2(q), \zeta_2 & a & 1 & [\ell^a] & R_{3,\ell}/R_{2,\ell} & 1 & 1 & \geq 3 \\ \cline{2-12}
& 5 & q^2 + 1 & C_2(q) & {1} & 2a & 5 & [\ell^a]^2 & R_{10,\ell} & D_8 & 5 & \geq 3 \\ \hline\hline
%\cline{5-12}
%& & & & B_2(q), \zeta_e & 0 & 1 & 1 & 1 & 1 & 1 & \geq 3 \\ \hline\hline
\end{array}
$$
\end{table}
\end{landscape}

\clearpage

\begin{landscape}
\begin{table}[h]
\caption{\label{16} The $\ell$-blocks of $F_4(q)$ of geometric type~$16$ 
(exist only when~$q$ is odd)}
$$
\begin{array}{cccccccccccc}
\multicolumn{12}{c}{e = 1, 2} \\ \hline\hline
1 & 2 & 3 & 4 & 5 & 6 & 7 & 8 & 9 & 10 & 11 & 13 \\ 
i & k & |Z| & [C,C] & b & d(b) & l(b) & D(b) & R & \Out(b_R) & \cW & \ell
\\ \hline\hline
 16 & 1, 9 & 2(q - \varepsilon)^2 & A_1(q)A_1(q) & {1} & 4a & 4 & [\ell^a]^4 & 
R_{18,\ell} & 2^2 & 4 & \geq 3 \\ \cline{3-12}
 & & 2(q + \varepsilon)^2 & A_1(q)A_1(q) & {1} & 2a & 4 & [\ell^a]^2 & 
R_{10,\ell} & 2^2 & 4 & \geq 3 \\ \cline{2-12}
 & 2/ 6 & 2(q^2 - 1) & A_1(q)A_1(q) & {1} & 3a & 4 & [\ell^a]^3 & 
R_{16,\ell}/R_{17,\ell} & 2^2 & 4 & \geq 3 \\ \cline{2-12}
 & 5 & 2(q^2 + 1) & A_1(q)A_1(q) & {1} & 2a & 4 & [\ell^a]^2 & 
R_{10,\ell} & 2^2 & 4 & \geq 3 \\ \cline{2-12}
 & 3, 10 & 2(q - \varepsilon)^2 & A_1(q^2) & {1} & 3a & 2 & [\ell^a]^3 & 
R_{16,\ell} & 2 & 2 & \geq 3 \\ \cline{3-12}
& & 2(q + \varepsilon)^2 & A_1(q^2) & {1} & a & 2 & [\ell^a] & R_{2,\ell} & 2 & 2 & \geq 3 \\ \cline{2-12}
& 4/ 7 & 2(q^2 - 1) & A_1(q^2) & {1} & 2a & 2 & [\ell^a]^2 & R_{10,\ell}/R_{9,\ell} & 2 & 2 & \geq 3 \\ \cline{2-12}
& 8 & 2(q^2 + 1) & A_1(q^2) & {1} & a & 2 & [\ell^a] & R_{2,\ell} & 2 & 2 & \geq 3 \\ \hline\hline
\end{array}
$$
\end{table}
\end{landscape}

\clearpage

\begin{landscape}
\begin{table}[h]
\caption{\label{17} The $\ell$-blocks of $F_4(q)$ of geometric type~$17$}
$$
\begin{array}{ccccccccccccc} 
\multicolumn{13}{c}{e = 1, 2} \\ \hline\hline
1 & 2 & 3 & 4 & 5 & 6 & 7 & 8 & 9 & 10 & 11 & 12 & 13 \\ 
i & k & |Z| & [C,C] & b & d(b) & l(b) & D(b) & R & \Out(b_R) & \cW & \text{Rem} & \ell 
\\ \hline\hline
17 & 1, 6 & (q - \varepsilon)^2 & \tilde{A}_2^{\varepsilon}(q) & {1} & 4a+1 & 3 & (13,k) & 
R_{34} & 2 & 2 & & 3 \\
 & & & & & & & & R_{26} & \SL_2(3) & 1 \\ \cline{6-13}
& & & & & 4a & 3 & [\ell^a]^4 & R_{18,\ell} & S_3 & 3 & & > 3 \\ \cline{3-13}
& & (q + \varepsilon)^2 & \tilde{A}_2^{-\varepsilon}(q) & {1} & a & 2 & [\ell^a] & R_{2,\ell} & 2 & 2 & & \geq 3 \\ \cline{2-13}
%& & & & \tilde{A}_2^{-\varepsilon}(q), (2,1) & 0 & 1 & 1 & 1 & 1 & 1 & & \geq 3 \\ \cline{2-13} 
& 3, 2 & q^2 - 1 & \tilde{A}_2^{\varepsilon}(q) & {1} & 3a+1 & 3 & (13,k') & 
R_{32} & 2 & 2 & & 3 \\
 & & & & & & & & R_{24} & \SL_2(3) & 1 & a \geq 2 \\ 
 & & & & & & & & R_{20} & \SL_2(3) & 1 & a = 1 \\ \cline{6-13}
 & & & & & 3a & 3 & [\ell^a]^3 & R_{17,\ell} & S_3 & 3 & & > 3 \\ \cline{3-13}
 & & q^2 - 1 & \tilde{A}_2^{-\varepsilon}(q) & {1} & 2a & 2 & [\ell^a]^2 & 
R_{9,\ell} & 2 & 2 & & \geq 3 \\ \cline{5-13}
& & & & \tilde{A}_2^{-\varepsilon}(q), (2,1) & a & 1 & [\ell^a] & R_{3,\ell} & 1 & 1 & & \geq 3 \\ \cline{2-13}
& 5, 4 & q^2 + \varepsilon q + 1 & \tilde{A}_2^{\varepsilon}(q) & {1} & 2a+2 & 3 &  (13,k) & 
R_{30} & 2 & 2 & & 3 \\
 & & & & & & & & R_{20} & \SL_2(3) & 1 \\ \cline{6-13}
& & & & & 2a & 3 & [\ell^a]^2 & R_{11,\ell} & S_3 & 3 & & > 3 \\ \cline{3-13}
& & q^2 - \varepsilon q + 1 & \tilde{A}_2^{-\varepsilon}(q) & {1} & a & 2 & [\ell^a] & R_{2,\ell} & 2 & 2 & & \geq 3 \\ \hline\hline
%\cline{5-13}
%& & & & \tilde{A}_2^{-\varepsilon}(q), (2,1) & 0 & 1 & 1 & 1 & 1 & 1 & & \geq 3 \\ \hline\hline
\end{array}
$$
\end{table}
\end{landscape}

\clearpage

\begin{landscape}
\begin{table}[h]
\caption{\label{18} The $\ell$-blocks of $F_4(q)$ of geometric type~$18$}
$$
\begin{array}{cccccccccccc}
\multicolumn{12}{c}{e = 1, 2} \\ \hline\hline
1 & 2 & 3 & 4 & 5 & 6 & 7 & 8 & 9 & 10 & 11 & 13 \\ 
i & k & |Z| & [C,C] & b & d(b) & l(b) & D(b) & R & \Out(b_R) & \cW & \ell
\\ \hline\hline
18 & 1, 10 & (q - \varepsilon)^3 & A_1(q) & {1} & 4a & 2 & [\ell^a]^4 & R_{18,\ell} & 2 & 2 & \geq 3 \\ \cline{3-12}
& & (q + \varepsilon)^3 & A_1(q) & {1} & a & 2 & [\ell^a] & R_{3,\ell} & 2 & 2 & \geq 3 \\ \cline{2-12}
 & 2, 8/ 9, 5 & (q^2 - 1)(q - \varepsilon) & A_1(q) & {1} & 3a & 2 & [\ell^a]^3 & 
R_{16,\ell}/R_{17,\ell} & 2 & 2 & \geq 3 \\ \cline{3-12}
 & & (q^2 - 1)(q + \varepsilon) & A_1(q) & {1} & 2a & 2 &  [\ell^a]^2 & 
R_{9,\ell}/R_{10,\ell} & 2 & 2 & \geq 3 \\ \cline{2-12}
 & 3, 7 & q^3 - \varepsilon & A_1(q) & {1} & 2a + 1 & 2 & [3^{a+1}] \times [3^a] & 
R_{12} & 2 & 2 & 3 \\ \cline{6-12}
& & & & & 2a & 2 & [\ell^a]^2 & R_{12,\ell} & 2 & 2 & > 3 \\ \cline{3-12}
& & q^3 + \varepsilon & A_1(q) & {1} & a & 2 & [\ell^a] & R_{3,\ell} & 2 & 2 & \geq 3 \\ \cline{2-12}
 & 6, 4 & (q^2 + 1)(q - \varepsilon) & A_1(q) & {1} & 2a & 2 & [\ell^a]^2 & 
R_{10,\ell} & 2 & 2 & \geq 3 \\ \cline{3-12}
& & (q^2 + 1)(q + \varepsilon) & A_1(q) & {1} & a & 2 & [\ell^a] & R_{3,\ell} & 2 & 2 & \geq 3 \\ \hline\hline
\end{array}
$$
\end{table}
\end{landscape}

\clearpage

\begin{landscape}
\begin{table}[h]
\caption{\label{19} The $\ell$-blocks of $F_4(q)$ of geometric type~$19$}
$$
\begin{array}{cccccccccccc} 
\multicolumn{12}{c}{e = 1, 2} \\ \hline\hline
1 & 2 & 3 & 4 & 5 & 6 & 7 & 8 & 9 & 10 & 11 & 13 \\ 
i & k & |Z| & [C,C] & b & d(b) & l(b) & D(b) & R & \Out(b_R) & \cW & \ell
\\ \hline\hline
19 & 1, 10 & (q - \varepsilon)^3 & \tilde{A}_1(q) & {1} & 4a & 2 & [\ell^a]^4 & R_{18,\ell} & 2 & 2 & \geq 3 \\ \cline{3-12}
& & (q + \varepsilon)^3 & \tilde{A}_1(q) & {1} & a & 2 & [\ell^a] & R_{2,\ell} & 2 & 2 & \geq 3 \\ \cline{2-12}
 & 2, 5/ 3, 4 & (q^2 - 1)(q - \varepsilon) & \tilde{A}_1(q) & {1} & 3a & 2 & [\ell^a]^3 & 
R_{16,\ell}/R_{17,\ell} & 2 & 2 & \geq 3 \\ \cline{3-12}
 & & (q^2 - 1)(q + \varepsilon) & \tilde{A}_1(q) & {1} & 2a & 2 &  [\ell^a]^2 & 
R_{10,\ell}/R_{9,\ell} & 2 & 2 & \geq 3 \\ \cline{2-12}
 & 8, 9 & q^3 - \varepsilon & \tilde{A}_1(q) & {1} & 2a + 1 & 2 & [3^{a+1}] \times [3^a] & 
R_{11} & 2 & 2 & 3 \\ \cline{6-12}
& & & & & 2a & 2 & [\ell^a]^2 & R_{11,\ell} & 2 & 2 & > 3 \\ \cline{3-12}
& & q^3 + \varepsilon & \tilde{A}_1(q) & {1} & a & 2 & [\ell^a] & R_{2,\ell} & 2 & 2 & \geq 3 \\ \cline{2-12}
 & 6, 7 & (q^2 + 1)(q - \varepsilon) & \tilde{A}_1(q) & {1} & 2a & 2 & [\ell^a]^2 & 
R_{10,\ell} & 2 & 2 & \geq 3 \\ \cline{3-12}
& & (q^2 + 1)(q + \varepsilon) & \tilde{A}_1(q) & {1} & a & 2 & [\ell^a] & R_{2,\ell} & 2 & 2 & \geq 3 \\ \hline\hline
\end{array}
$$
\end{table}
\end{landscape}

\clearpage

\begin{landscape}
\addtocounter{table}{1}
\begin{table}[h]
	\caption{\label{21} The faithful $\ell$-blocks of $\hat{G} = 2.F_4(2)$ of non-cyclic defect}
$$
\begin{array}{cccccc}
\multicolumn{6}{c}{\ell = 3} \\ \hline\hline
 6 & 7 & 8 & 9 & 10' & 11 \\
d(b) & l(b) & D(b) & R & N_{\hat{G}}(R)/R & \cW 
\\ \hline\hline
6 & 17 & (1,1) & R_{38} & 2^4 & 8 \\
 & & & R_{37} & 2 \times \GL_2(3) & 2  \\
 & & & R_{15} & 2 \times \SL_3(3) & 1 \\
 & & & R_{18} &  2 \times W(F_4) & 4 \\
 & & & R_{21} & 2\dot{\ }[(Q_8\times Q_8)\colon\!S_3] & 2 \\ \hline
 2 & 5 & 3^2 & R_{10} & 2 \times D_8 \times \Sp_4(2) & 5 \\ \hline
 2 & 5 & 3^2 & R_{10} & 2 \times D_8 \times \Sp_4(2)  & 5 \\ \hline\hline
\end{array}
$$

\vspace*{0.6cm}

$$
\begin{array}{cccccc}
\multicolumn{6}{c}{\ell = 5} \\ \hline\hline
 6 & 7 & 8 & 9 & 10' & 11 \\
	d(b) & l(b) & D(b) & R & N_{\hat{G}}(R)/R & \cW
\\ \hline\hline
2 & 16 & [5]^2 & [5]^2 & 2 \times \SL_2(3)\colon\![4] & 16 \\ \hline\hline
\end{array}
$$

\vspace*{0.6cm}

$$
\begin{array}{cccccc}
\multicolumn{6}{c}{\ell = 7} \\ \hline\hline
6 & 7 & 8 & 9 & 10' & 11 \\
	d(b) & l(b) & D(b) & R & N_{\hat{G}}(R)/R & \cW
\\ \hline\hline
2 & 21 & [7]^2 & R_{8} & 2 \times \SL_2(3) \times 3 & 21 \\ \hline\hline
\end{array}
$$
\end{table}
\end{landscape}

\clearpage

\subsection{Maximal tori}
\label{AppendixClassTypes}
This subsection contains information on the maximal tori of~$F_4(q)$. The 
$G$-conjugacy classes of the $F$-stable maximal tori of~$\mathbf{G}$ are in
bijection with the conjugacy classes of~$W$. For each such conjugacy class, 
Table~\ref{MaximalTori} gives the order of a representative $w$, the order of 
the centralizer $|C_W(w)|$, the $2$- and $3$-power maps on the conjugacy 
classes, the action of the automorphism~$\dagger$ of~$W$ (which swaps $s_1$ 
with $s_4$ and $s_2$ with $s_3$), as well as the names of the conjugacy classes 
following \cite{CWeyl} and the structure of the corresponding maximal tori. The 
first column of Table~\ref{MaximalTori} numbers the conjugacy classes. The 
reflections corresponding to the long respectively short roots lie in conjugacy 
class~$12$, respectively~$17$.

\begin{table}
\caption{\label{MaximalTori} Maximal tori of~$G$}
$$
\begin{array}{rrrrrrrl} \\ \hline\hline
\multicolumn{1}{c}{\text{No.}} & 
\multicolumn{1}{c}{|w|} & 
\multicolumn{1}{c}{|C_W(w)|} & 
\multicolumn{1}{c}{\text{2p}} &
\multicolumn{1}{c}{\text{3p}} &
\multicolumn{1}{c}{\dagger} &
\multicolumn{1}{c}{\text{Name}} & 
\multicolumn{1}{c}{\text{Structure}} \\ \hline\hline
 1 & 1 & 1152 & 1 & 1 & 1 &A_0       &  [\Phi_1(q)]^4 \\
 2 & 2 & 1152 & 1 & 2 & 2 &4A_1      &  [\Phi_2(q)]^4 \\
 3 & 2 &   64 & 1 & 3 & 3 &2A_1      &  [d] \times [\Phi_1(q)\Phi_2(q)/d] \times [\Phi_1(q)\Phi_2(q)] \\
 4 & 3 &   36 & 4 & 1 & 7 &A_2       &  [\Phi_1(q)] \times [\Phi_1(q)\Phi_3(q)] \\
 5 & 6 &   36 & 4 & 2 & 8 &D_4       &  [\Phi_2(q)] \times [\Phi_2(q)\Phi_6(q)] \\
 6 & 4 &   96 & 2 & 6 & 6 &D_4(a_1)  &  [\Phi_4(q)]^2 \\
 7 & 3 &   36 & 7 & 1 & 4 &\tilde A_2      &  [\Phi_1(q)] \times [\Phi_1(q)\Phi_3(q)] \\
 8 & 6 &   36 & 7 & 2 & 5 &C_3+A_1   &  [\Phi_2(q)] \times [\Phi_2(q)\Phi_6(q)]  \\
 9 & 3 &   72 & 9 & 1 & 9 &A_2+\tilde A_2  &  [\Phi_3(q)]^2  \\
10 & 6 &   72 & 9 & 2 &10 &F_4(a_1)  &  [\Phi_6(q)]^2 \\
11 &12 &   12 &10 & 6 &11 &F_4       &  [\Phi_{12}(q)] \\
12 & 2 &   96 & 1 &12 &17 &A_1       &  [\Phi_1(q)]^2 \times [\Phi_1(q)\Phi_2(q)] \\
13 & 2 &   96 & 1 &13 &18 &3A_1      &  [\Phi_2(q)]^2 \times [\Phi_1(q)\Phi_2(q)] \\
14 & 6 &   12 & 7 &12 &19 &\tilde A_2+A_1  &  [\Phi_1(q)\Phi_2(q)\Phi_3(q)] \\
15 & 6 &   12 & 7 &13 &20 &C_3       &  [\Phi_1(q)\Phi_2(q)\Phi_6(q)] \\
16 & 4 &   16 & 3 &16 &21 &A_3       &  [d] \times [\Phi_1(q)\Phi_2(q)\Phi_4(q)/d] \\
17 & 2 &   96 & 1 &17 &12 &\tilde A_1      &  [\Phi_1(q)]^2 \times [\Phi_1(q)\Phi_2(q)] \\
18 & 2 &   96 & 1 &18 &13 &2A_1+\tilde A_1 &  [\Phi_2(q)]^2 \times [\Phi_1(q)\Phi_2(q)] \\
19 & 6 &   12 & 4 &17 &14 &A_2+\tilde A_1  &  [\Phi_1(q)\Phi_2(q)\Phi_3(q)] \\
20 & 6 &   12 & 4 &18 &15 &B_3       &  [\Phi_1(q)\Phi_2(q)\Phi_6(q)] \\
21 & 4 &   16 & 3 &21 &16 &B_2+A_1   &  [d] \times [\Phi_1(q)\Phi_2(q)\Phi_4(q)]/d] \\
22 & 2 &   16 & 1 &22 &22 &A_1+\tilde A_1  &  [\Phi_1(q)\Phi_2(q)]^2 \\
23 & 4 &   32 & 3 &23 &23 &B_2       &  [\Phi_1(q)] \times [\Phi_1(q)\Phi_4(q)] \\
24 & 4 &   32 & 3 &24 &24 &A_3+\tilde A_1  &  [\Phi_2(q)] \times [\Phi_2(q)\Phi_4(q)] \\
25 & 8 &    8 & 6 &25 &25 &B_4       &  [\Phi_8(q)]
\\ \hline\hline
\end{array}
$$
\end{table}

\clearpage

\setlength{\extrarowheight}{0.5ex}

\subsection{The construction of the centralizers}
\label{ConstructionOfCentralizers}
Recall that $F = F_1^f$ with~$F_1$ as in Subsection~\ref{SetupF4}. In this 
subsection we assume that~$m$ is a positive integer dividing~$f$. 
Let~$\mathbf{T}_0$ denote a $1$-$F_1$-split maximal torus of~$\mathbf{G}$, and 
let $W := N_{\mathbf{G}}( \mathbf{T}_0 )/\mathbf{T}_0$.

Recall from Subsection~\ref{ClassTypes} that the $G$-class types of semisimple
elements of~$G$ are numbered by pairs of integers $(i,k)$, with $1 \leq i 
\leq 20$ and~$k$ depending on~$i$. The integer~$i$ labels certain subsets 
$\Gamma_i \subseteq \Sigma$ with $\Gamma_1 = \{\alpha_1, \ldots, \alpha_4 \}$ 
and $\Gamma_{20} = \emptyset$, and such that $\Gamma_i$ is a base of 
$\bar{\Gamma}_i$ for $1 \leq i \leq 20$. Moreover,~$k$ labels the conjugacy 
classes of $\Stab_W( \Gamma_i )$. 

Let us now explain the contents of Table~\ref{CT}, referring to
Subsections~\ref{SetupF4} and~\ref{ClassTypes} for the notation. For each pair 
$(i,k)$ as above, with 
$2 \leq i \leq 19$, this table lists several pairs of elements $v, w \in W$;
the choice of these pairs depends on two parameters, namely on certain 
congruence classes of~$f/m$ and on integers $e \in \{ 1, 2, 3, 4, 6 \}$. The 
first column of Table~\ref{CT} contains the sets~$\Gamma_i$ as list of integers 
$j_1, \ldots , j_c$ if $\Gamma_i = \{ \alpha_{j_1} , \ldots , \alpha_{j_c} \}$, 
and the second column contains~$(i,k)$. Column~$3$ and~$4$ give the values 
of~$e$ and of the congruence class of~$f/m$, respectively, where a hyphen 
indicates that there is no condition on $f/m$. Columns $5$--$7$ of 
Table~\ref{CT} contain the elements $v$, $v^{f/m}$ and~$w$, respectively, where 
we use the following conventions. Let $\Gamma$ be one of the $\Gamma_i$. 
As~$\Gamma$ is a base of $\bar{\Gamma}$, we have $W_{\Gamma} = 
W_{\bar{\Gamma}}$, and thus~$W_{\Gamma}$ is a Weyl group. The longest element 
of $W_{\Gamma}$ is denoted by $w_{\Gamma}$ and $w_0$ is the longest element 
of~$W$. Also, $s_j \in W$ denotes the reflection corresponding to the root 
$\alpha_j$ and $r_j$ denotes a representative of the $W$-conjugacy class with 
number~$j$, according to Table~\ref{MaximalTori}. For $w \in W$, we denote 
by~$w'$ any primitive power of~$w$, i.e.\ any power of~$w$ of the same order 
as~$w$.

The content of the last two columns of  Table~\ref{CT} and further properties
of the data listed are collected in the following remark, whose assertions
are easily checked with CHEVIE.
\enlargethispage{5ex}

\addtocounter{thm}{3}
\begin{rem}
\label{ProofCT}
{\rm 
Let $\Gamma \subseteq \Sigma$, $i, k \in \mathbb{Z}$, 
$e \in \{ 1, 2, 3, 4, 6 \}$ and $v, w \in W$ such that some row of 
{\rm Table~\ref{CT}} has the values
$(\Gamma, (i,k) , -, -, v, -, w, - , -)$. 
Put $C_{\Gamma}( v^{f/m}w ) := C_W( v^{f/m}w ) \cap W_\Gamma =
\C_{W_{\Gamma}}( v^{f/m}w )$. Then the following holds.

{\rm (i)} We have $v \in W_{\Gamma}.\Stab_W( \Gamma )$, 
and the centralizers of semisimple elements of~$G$ of $G$-class type
$(i,k)$ are of $F$-type $(\bar{\Gamma}, [v])$. 
%for the definition of $F$-type see Subsection~\ref{MaximalRank}.

{\rm (ii)} We have $w \in W_{\Gamma} \leq W$ and the number, according to 
Table~\ref{MaximalTori}, of the conjugacy class of $v^{f/m}w$ in~$W$ is as given 
in the column of {\rm Table~\ref{CT}} headed with~$\text{cl}$.

{\rm (iii)} The elements~$v$ and~$w$ commute.

{\rm (iv)} The element~$v$ centralizes $C_{\Gamma}( v^{f/m}w )$, unless~$f/m$ is 
even and $(i,k) \in \{ (12,2), (12,4), (16,k), k \in \{ 3, 4, 7, 8, 10 \} \}$.

{\rm (v)} The structure of $C_{\Gamma}( v^{f/m}w )$ is as given in the last
column of {\rm Table~\ref{CT}}.
}\hfill{$\Box$}
\end{rem}

\clearpage 

\begin{table}[h]
\caption{\label{CT} Construction of Centralizers}
$$
\begin{array}{ccrcccccc} 
\Gamma & (i,k) & e & f/m & v & v^{f/m} & w & \text{cl} & C_{\Gamma}( v^{f/m}w ) 
\\ \hline\hline
 1, 2, 3, 48 & (2,1) & 1 & - & 1 & 1 & 1 & 1 & W(B_4) \\
 & & 2 & - & 1 & 1 & w_0 & 2 & W(B_4) \\ 
 & & 4 & - & 1 & 1 & r_6 & 6 & [4]^2.2 \\ \hline
 1, 2, 4, 48 & (3,1) & 1 & - & 1 & 1 & 1 & 1 & W(A_3) \times 2 \\
 & & 2 & - & 1 & 1 & s_1s_4s_{22} & 18 & D_8 \times 2 \\
 & (3,2) & 1 &  \not\equiv 0(2) & s_{17} & s_{17} & 1 & 17 & D_8 \times 2 \\
 & & 2 & \not\equiv 0(2) & w_0 & w_0 & 1 & 2 & W(A_3) \times 2 \\ \hline
 1, 3, 4, 48 & (4,1) & 1 & - & 1 & 1 & 1 & 1 & S_3 \times S_3 \\
 & & 2 & - & 1 & 1 & s_7s_{23} & 22 & 2^2 \\
 & & 3 & - & 1 & 1 & r_9 & 9 & 3^2 \\
 & (4,2) & 1 &  \not\equiv 0(2) & s_{12}s_{11} & s_{12}s_{11} & 1 & 22 & 2^2 \\
 & & 2 & \not\equiv 0(2) & w_0 & w_0 & 1 & 2 & S_3 \times S_3 \\ 
 & & 6 & \not\equiv 0(2) & w_0 & w_0 & r_{9} & 10 & 3^2 \\ \hline
 2, 3, 4, 48 & (5,1) & 1 & - & 1 & 1 & 1 & 1 & W(C_3) \times 2 \\
 & & 2 & - & 1 & 1 & w_0 & 2 & W(C_3) \times 2 \\ \hline
 1, 2, 3 & (6,1) & 1 & - & 1 & 1 & 1 & 1 & W(B_3) \\
 & & 2 & - & 1 & 1 & w_{\Gamma} & 18 & \\
 & (6,2) & 1 & \equiv 0(2) & w_0 & 1 & 1 & 1 & \\
 & & 2 &  \equiv 0(2) & s_{21} & 1 & w_{\Gamma} & 18 & \\
 & & 1 &  \not\equiv 0(2) & s_{21} & s_{21} & 1 & 17 & \\
 & & 2 & \not\equiv 0(2) & w_0 & w_0 & 1 & 2 & \\ \hline
 1, 2, 4 & (7,1) & 1 & - & 1 & 1 & 1 & 1 & 2 \times S_3 \\
 & & 2 & - & 1 & 1 & w_{\Gamma} & 22 & 2^2 \\
 & (7,2) & 1 & \equiv 0(2) & w_0 & 1 & 1 & 1 & 2 \times S_3 \\
 & & 2 &  \equiv 0(2) & s_{15}s_{23} & 1 & w_{\Gamma} & 22 & 2^2 \\
 & & 1 &  \not\equiv 0(2) & s_{15}s_{23} & s_{15}s_{23} & 1 & 22 & 2^2 \\
 & & 2 & \not\equiv 0(2) & w_0 & w_0 & 1 & 2 & 2 \times S_3 \\ \hline\hline
\end{array}
$$
\end{table}

\clearpage

\addtocounter{table}{-1}
\begin{table}[h]
\caption{Construction of Centralizers (continued)}
$$
\begin{array}{ccrcccccc}
	\Gamma & (i,k) & e & f/m & v & v^{f/m} & w & \text{cl}      & C_{\Gamma}( v^{f/m}w ) 
\\ \hline\hline
 1, 2, 48 & (8,1) & 1 & - & 1 & 1 & 1 & 1 & W(A_3) \\
 & & 2 & - & 1 & 1 & s_1s_{22} & 3 & D_8 \\
 & (8,4) & 1 & \equiv 0(2) & w_0 & 1 & 1 & 1 & W(A_3) \\
 & & 2 &  \equiv 0(2) & s_4s_{17} & 1 & s_1s_{22} & 3 & D_8 \\
 & & 1 &  \not\equiv 0(2) & s_4s_{17} & s_4s_{17} & 1 & 3 & D_8 \\
 & & 2 & \not\equiv 0(2) & w_0 & w_0 & 1 & 2 & W(A_3) \\ 
 & (8,2) & 1 & \equiv 0(2) & s_4 & 1 & 1 & 1 & W(A_3) \\
 & & 2 &  \equiv 0(2) & s_4 & 1 & s_1s_{22} & 3 & D_8 \\
 & & 1 &  \not\equiv 0(2) & s_4 & s_4 & 1 & 17 & W(A_3) \\
 & & 2 & \not\equiv 0(2) & w_0s_{17} & w_0 s_{17}& 1 & 18 & D_8 \\
 & (8,3) & 1 & \equiv 0(2) & w_0s_{4} & 1 & 1 & 1 & W(A_3) \\
 & & 2 &  \equiv 0(2) & s_{17} & 1 & s_1s_{22} & 3 & D_8 \\
 & & 1 &  \not\equiv 0(2) & s_{17} & s_{17} & 1 & 17 & D_8 \\
 & & 2 & \not\equiv 0(2) & w_0s_{4} & w_0 s_{4}& 1 & 18 & W(A_3) \\ \hline
 1, 3, 4 & (9,1) & 1 & - & 1 & 1 & 1 & 1 & S_3 \times 2 \\
 & & 2 & - & 1 & 1 & w_{\Gamma} & 22 & 2^2 \\
 & (9,2) & 1 & \equiv 0(2) & w_0 & 1 & 1 & 1 & S_3 \times 2 \\
 & & 2 &  \equiv 0(2) & s_{20}s_{19}& 1 & w_{\Gamma} & 22 & 2^2 \\
 & & 1 &  \not\equiv 0(2) & s_{20}s_{19} & s_{20}s_{19} & 1 & 22 & 2^2 \\
 & & 2 & \not\equiv 0(2) & w_0 & w_0 & 1 & 2 & S_3 \times 2 \\ \hline
 2, 3, 4 & (10,1) & 1 & - & 1 & 1 & 1 & 1 & W(C_3) \\
 & & 2 & - & 1 & 1 & w_{\Gamma} & 13 & \\
 & (10,2) & 1 & \equiv 0(2) & w_0 & 1 & 1 & 1 & \\
 & & 2 &  \equiv 0(2) & s_{24} & 1 & w_{\Gamma} & 13 & \\
 & & 1 &  \not\equiv 0(2) & s_{24} & s_{24} & 1 & 12 & \\
 & & 2 & \not\equiv 0(2) & w_0 & w_0 & 1 & 2 & \\ \hline\hline
\end{array}
$$
\end{table}

\clearpage

\addtocounter{table}{-1}
\begin{table}[h]
\caption{Construction of Centralizers (continued)}
$$
\begin{array}{ccrcccccc}
	\Gamma & (i,k) & e & f/m & v & v^{f/m} & w & \text{cl}      & C_{\Gamma}( v^{f/m}w ) 
\\ \hline\hline
 2, 3, 48 & (11,1) & 1 & - & 1 & 1 & 1 & 1 & D_8 \times 2 \\
 & & 2 & - & 1 & 1 & w_{\Gamma} & 13 & \\
 & (11,2) & 1 & \equiv 0(2) & w_0 & 1 & 1 & 1 & \\
 & & 2 &  \equiv 0(2) & s_{16} & 1 & w_{\Gamma} & 13 & \\
 & & 1 &  \not\equiv 0(2) & s_{16} & s_{16} & 1 & 12 & \\
 & & 2 & \not\equiv 0(2) & w_0 & w_0 & 1 & 2 & \\ \hline
 2, 4, 48 & (12,1) & 1 & - & 1 & 1 & 1 & 1 & 2^3 \\
 & & 2 & - & 1 & 1 & w_{\Gamma} & 18 & \\
 & (12,3) & 1 & \equiv 0(2) & w_0 & 1 & 1 & 1 & \\
 & & 2 &  \equiv 0(2) & s_{13} & 1 & w_{\Gamma} & 18 & \\
 & & 1 &  \not\equiv 0(2) & s_{13} & s_{13} & 1 & 17 & \\
 & & 2 & \not\equiv 0(2) & w_0 & w_0 & 1 & 2 & \\ 
 & (12,2) & 1 & \equiv 0(2) & s_{17} & 1 & 1 & 1 & 2^2 \\
 & & 2 &  \equiv 0(2) & s_{17} & 1 & w_{\Gamma} & 18 & \\
 & & 1 &  \not\equiv 0(2) & s_{17} & s_{17} & 1 & 17 & \\
 & & 2 & \not\equiv 0(2) & w_0s_{13}s_{17} & w_0s_{13}s_{17}& 1 & 3 & \\
 & (12,4) & 1 & \equiv 0(2) & s_{13}s_{17} & 1 & 1 & 1 & \\
 & & 2 &  \equiv 0(2) & s_{13}s_{17} & 1 & w_{\Gamma} & 18 & \\
 & & 1 &  \not\equiv 0(2) & s_{13}s_{17} & s_{13}s_{17} & 1 & 3 & \\
 & & 2 & \not\equiv 0(2) & w_0s_{17} & w_0 s_{17}& 1 & 18 & \\ \hline\hline
\end{array}
$$
\end{table}

\clearpage

\addtocounter{table}{-1}
\begin{table}[h]
\caption{Construction of Centralizers (continued)}
$$
\begin{array}{ccrcccccc}
	\Gamma & (i,k) & e & f/m & v & v^{f/m} & w & \text{cl}      & C_{\Gamma}( v^{f/m}w ) 
\\ \hline\hline
 1, 2 & (13,1) & 1 & - & 1 & 1 & 1 & 1 & S_3 \\
 & & 2 & - & 1 & 1 & w_{\Gamma} & 12 & 2 \\
 & (13,6) & 1 & \equiv 0(2) & w_0 & 1 & 1 & 1 & S_3 \\
 & & 2 &  \equiv 0(2) & w_0w_{\Gamma} & 1 & w_{\Gamma} & 12 & 2 \\
 & & 1 &  \not\equiv 0(2) & w_0w_{\Gamma} & w_0w_{\Gamma} & 1 & 13 & 2 \\
 & & 2 & \not\equiv 0(2) & w_0 & w_0 & 1 & 2 & S_3 \\ 
 & (13,2) & 1 & \equiv 0(2) & s_4 & 1 & 1 & 1 & S_3 \\
 & & 2 &  \equiv 0(2) & s_4 & 1 & w_{\Gamma} & 12 & 2 \\
 & & 1 &  \not\equiv 0(2) & s_4 & s_4 & 1 & 17 & S_3 \\
 & & 2 & \not\equiv 0(2) & s_4w_{\Gamma} & s_4w_{\Gamma} & 1 & 22 & 2 \\
 & (13,3) & 1 & \equiv 0(2) & w_0s_4 & 1 & 1 & 1 & S_3 \\
 & & 2 &  \equiv 0(2) & w_0s_4 & 1 & w_{\Gamma} & 12 & 2 \\
 & & 1 &  \not\equiv 0(2) & w_0s_4w_{\Gamma} & w_0s_4w_{\Gamma} & 1 & 22 & 2 \\
 & & 2 & \not\equiv 0(2) & w_0s_4 & w_0s_4 & 1 & 18 & S_3 \\ 
 & (13,5) & 1 & \equiv 0(3) & s_4s_{19} & 1 & 1 & 1 & S_3 \\
 & & 2 &  \equiv 0(3) & s_4s_{19} & 1 & w_{\Gamma} & 12 & 2 \\
 & & 1 &  \not\equiv 0(3) & s_4s_{19} & (s_4s_{19})' & 1 & 7 & S_3 \\
 & & 2 & \not\equiv 0(3) & s_4s_{19}w_{\Gamma} & (s_4s_{19}w_{\Gamma})' & 1 & 14 & 2 \\
 & (13,4) & 1 & \equiv 0(6) & w_0s_4s_{19} & 1 & 1 & 1 & S_3 \\
 & & 2 &  \equiv 0(6) & w_0s_4s_{19} & 1 & w_{\Gamma} & 12 & 2 \\
 & & 1 &  \equiv 0(2)\not\equiv 0(3) & w_0s_4s_{19} & (s_4s_{19})' & 1 & 7 & S_3 \\
 & & 2 & \equiv 0(2)\not\equiv 0(3) & w_0s_4s_{19}w_{\Gamma} & (s_4s_{19}w_{\Gamma})' & 1 & 14 & 2 \\
 & & 1 &  \not\equiv 0(2)\equiv 0(3) & w_0s_4s_{19}w_{\Gamma} & w_0w_{\Gamma} & 1& 13 & 2 \\
 & & 2 & \not\equiv 0(2)\equiv 0(3) & w_0s_4s_{19} & w_0 & 1 & 2 & S_3 \\
 & & 1 & \gcd(6,f/m) = 1 & w_0s_4s_{19}w_{\Gamma} & (w_0s_4s_{19}w_{\Gamma})' & 1 & 15 & 2 \\
 & & 2 & \gcd(6,f/m) = 1 & w_0s_4s_{19} & (w_0s_4s_{19})' & 1 & 8 & S_3 
 \\ \hline\hline
\end{array}
$$
\end{table}

\clearpage

\addtocounter{table}{-1}
\begin{table}[h]
\caption{Construction of Centralizers (continued)}
$$
\begin{array}{ccrcccccc}
	\Gamma & (i,k) & e & f/m & v & v^{f/m} & w & \text{cl}      & C_{\Gamma}( v^{f/m}w ) 
\\ \hline\hline
 1, 4 & (14,1) & 1 & - & 1 & 1 & 1 & 1 & 2^2 \\
 & & 2 & - & 1 & 1 & w_{\Gamma} & 22 & \\
 & (14,4) & 1 & \equiv 0(2) & s_{17}s_{22} & 1 & 1 & 1 & \\
 & & 2 &  \equiv 0(2) & s_{17}s_{22} & 1 & w_{\Gamma} & 22 & \\
 & & 1 &  \not\equiv 0(2) & s_{17}s_{22} & s_{17}s_{22} & 1 & 22 & \\
 & & 2 & \not\equiv 0(2) & w_0 & w_0 & 1 & 2 & \\
 & (14,2) & 1 & \equiv 0(2) & s_{17} & 1 & 1 & 1 & \\
 & & 2 &  \equiv 0(2) & s_{17} & 1 & w_{\Gamma} & 22 & \\
 & & 1 &  \not\equiv 0(2) & s_{17} & s_{17} & 1 & 17 & \\
 & & 2 & \not\equiv 0(2) & w_0s_{22} & w_0s_{22} & 1 & 13 & \\
 & (14,3) & 1 & \equiv 0(2) & s_{22} & 1 & 1 & 1 & \\
 & & 2 &  \equiv 0(2) & s_{22} & 1 & w_{\Gamma} & 22 & \\
 & & 1 &  \not\equiv 0(2) & s_{22} & s_{22} & 1 & 12 & \\
 & & 2 & \not\equiv 0(2) & w_0s_{17} & w_0 s_{17}& 1 & 18 & \\ \hline
 2, 3 & (15,1) & 1 & - & 1 & 1 & 1 & 1 & D_8 \\
 & & 2 & - & 1 & 1 & w_{\Gamma} & 3 & \\
 & (15,3) & 1 & \equiv 0(2) & w_0w_{\Gamma} & 1 & 1 & 1 & \\
 & & 2 &  \equiv 0(2) & w_0w_{\Gamma} & 1 & w_{\Gamma} & 3 & \\
 & & 1 &  \not\equiv 0(2) & w_0w_{\Gamma} & w_0w_{\Gamma} & 1 & 3 & \\
 & & 2 & \not\equiv 0(2) & w_0 & w_0 & 1 & 2 & \\
 & (15,2) & 1 & \equiv 0(2) & s_{16} & 1 & 1 & 1 & \\
 & & 2 &  \equiv 0(2) & s_{16} & 1 & w_{\Gamma} & 3 & \\
 & & 1 &  \not\equiv 0(2) & s_{16} & s_{16} & 1 & 12 & \\
 & & 2 & \not\equiv 0(2) & w_0s_{16} & w_0s_{16} & 1 & 13 & \\
 & (15,4) & 1 & \equiv 0(2) & s_{8} & 1 & 1 & 1 & \\
 & & 2 &  \equiv 0(2) & s_{8} & 1 & w_{\Gamma} & 3 & \\
 & & 1 &  \not\equiv 0(2) & s_{8} & s_{8} & 1 & 17 & \\
 & & 2 & \not\equiv 0(2) & w_0s_{8} & w_0 s_{8}& 1 & 18 & \\ 
 & (15,5) & 1 & \equiv 0(4) & s_8s_{16} & 1 & 1 & 1 & \\
 & & 2 &  \equiv 0(4) & s_8s_{16} & 1 & w_{\Gamma} & 3 & \\
 & & 1 & \equiv 2(4) & s_8s_{16} & (s_8s_{16})^2 & 1 & 3 & \\
 & & 2 &  \equiv 2(4) & s_8s_{16} & (s_8s_{16})^2  & w_{\Gamma} & 2 & \\
 & & 1 &  \not\equiv 0(2) & s_8s_{16} & (s_8s_{16})' & 1 & 23 & \\
 & & 2 & \not\equiv 0(2) & w_0s_8s_{16} & (w_0s_8s_{16})' & 1 & 24 & \\ 
 \hline\hline
\end{array}
$$
\end{table}

\clearpage

\addtocounter{table}{-1}
\begin{table}[h]
\caption{Construction of Centralizers (continued)}
$$
\begin{array}{ccrcccccc}
	\Gamma & (i,k) & e & f/m & v & v^{f/m} & w & \text{cl}      & C_{\Gamma}( v^{f/m}w ) 
\\ \hline\hline
 2, 48 & (16,1) & 1 & - & 1 & 1 & 1 & 1 & 2^2 \\
 & & 2 & - & 1 & 1 & w_{\Gamma} & 3 & \\
 & (16,9) & 1 & \equiv 0(2) & (s_4s_9)^2 & 1 & 1 & 1 & \\
 & & 2 &  \equiv 0(2) & (s_4s_9)^2 & 1 & w_{\Gamma} & 3 & \\
 & & 1 &  \not\equiv 0(2) & (s_4s_9)^2 & (s_4s_9)^2 & 1 & 3 & \\
 & & 2 & \not\equiv 0(2) & w_0 & w_0 & 1 & 2 & \\
 & (16,2) & 1 & \equiv 0(2) & s_{4} & 1 & 1 & 1 & \\
 & & 2 &  \equiv 0(2) & s_{4} & 1 & w_{\Gamma} & 3 & \\
 & & 1 &  \not\equiv 0(2) & s_{4} & s_{4} & 1 & 17 & \\
 & & 2 & \not\equiv 0(2) & w_0s_{4} & w_0s_{4} & 1 & 18 & \\
 & (16,6) & 1 & \equiv 0(2) & s_{9} & 1 & 1 & 1 & \\
 & & 2 &  \equiv 0(2) & s_{9} & 1 & w_{\Gamma} & 3 & \\
 & & 1 &  \not\equiv 0(2) & s_{9} & s_{9} & 1 & 12 & \\
 & & 2 & \not\equiv 0(2) & w_0s_{9} & w_0 s_{9}& 1 & 13 & \\
 & (16,5) & 1 & \equiv 0(4) & s_4s_{9} & 1 & 1 & 1 & \\
 & & 2 &  \equiv 0(4) & s_4s_{9} & 1 & w_{\Gamma} & 3 & \\
 & & 1 & \equiv 2(4) & s_4s_{9} & (s_4s_{9})^2 & 1 & 3 & \\
 & & 2 &  \equiv 2(4) & s_4s_{9} & (s_4s_{9})^2  & w_{\Gamma} & 2 & \\
 & & 1 &  \not\equiv 0(2) & s_4s_{9} & (s_4s_{9})' & 1 & 23 & \\
 & & 2 & \not\equiv 0(2) & w_0s_4s_{9} & (w_0s_4s_{9})' & 1 & 24 & \\
 \hline\hline
\end{array}
$$
\end{table}

\clearpage

\addtocounter{table}{-1}
\begin{table}[h]
\caption{Construction of Centralizers (continued)}
$$
\begin{array}{ccrcccccc}
	\Gamma & (i,k) & e & f/m & v & v^{f/m} & w & \text{cl}      & C_{\Gamma}( v^{f/m}w ) 
\\ \hline\hline
 & (16,4) & 1 & \equiv 0(2) & s_{4}s_{17} & 1 & 1 & 1 & 2^2 \\
	& & 2 &  \equiv 0(2) & s_{4}s_{17} & 1 & w_{\Gamma} & 3 & 2^2 \\
 & & 1 &  \not\equiv 0(2) & s_{4}s_{17} & s_{4}s_{17} & 1 & 3 & 2 \\
 & & 2 & \not\equiv 0(2) & s_{4}s_{17} & s_{4}s_{17} & 1 & 3 & 2 \\
 & (16,7) & 1 & \equiv 0(2) & s_{9}s_{17} & 1 & 1 & 1 & 2^2 \\
 & & 2 &  \equiv 0(2) & s_{9}s_{17} & 1 & w_{\Gamma} & 3 & 2^2 \\
 & & 1 &  \not\equiv 0(2) & s_{9}s_{17} & s_{9}s_{17} & 1 & 22 & 2 \\
 & & 2 & \not\equiv 0(2) & s_{9}s_{17} & s_{9}s_{17} & 1 & 22 & 2 \\
 & (16,3) & 1 & \equiv 0(2) & s_{17} & 1 & 1 & 1 & 2^2 \\
 & & 2 &  \equiv 0(2) & s_{17} & 1 & w_{\Gamma} & 3 & 2^2 \\
 & & 1 &  \not\equiv 0(2) & s_{17} & s_{17} & 1 & 17 & 2 \\
 & & 2 & \not\equiv 0(2) & s_{17} & s_{17} & 1 & 17 & 2 \\
 & (16,10) & 1 & \equiv 0(2) & w_0s_{17} & 1 & 1 & 1 & 2^2 \\
 & & 2 &  \equiv 0(2) & w_0s_{17} & 1 & w_{\Gamma} & 3 & 2^2 \\
 & & 1 &  \not\equiv 0(2) & w_0s_{17} & w_0s_{17} & 1 & 18 & 2 \\
 & & 2 & \not\equiv 0(2) & w_0s_{17} & w_0s_{17} & 1 & 18 & 2 \\
 & (16,8) & 1 & \equiv 0(4) & s_4s_{9}s_{17} & 1 & 1 & 1 & 2^2 \\
 & & 2 &  \equiv 0(4) & s_4s_{9}s_{17} & 1 & w_{\Gamma} & 3 & 2^2 \\
 & & 1 & \equiv 2(4) & s_4s_{9}s_{17} & (s_4s_{9})^2 & 1 & 3 & 2^2 \\
 & & 2 &  \equiv 2(4) & s_4s_{9}s_{17} & (s_4s_{9})^2  & w_{\Gamma} & 2 & 2^2 \\
 & & 1 &  \not\equiv 0(2) & s_4s_{9}s_{17} & (s_4s_{9}s_{17})' & 1 & 16 & 2 \\
 & & 2 & \not\equiv 0(2) & s_4s_{9}s_{17} & (s_4s_{9}s_{17})' & 1 & 16 & 2 \\
 \hline\hline
\end{array}
$$
\end{table}

\clearpage

\addtocounter{table}{-1}
\begin{table}[h]
\caption{Construction of Centralizers (continued)}
$$
\begin{array}{ccrcccccc}
	\Gamma & (i,k) & e & f/m & v & v^{f/m} & w & \text{cl}      & C_{\Gamma}( v^{f/m}w ) 
\\ \hline\hline
 3, 4 & (17,1) & 1 & - & 1 & 1 & 1 & 1 & S_3 \\
 & & 2 & - & 1 & 1 & w_{\Gamma} & 17 & 2 \\
 & (17,6) & 1 & \equiv 0(2) & w_0 & 1 & 1 & 1 & S_3 \\
 & & 2 &  \equiv 0(2) & w_0w_{\Gamma} & 1 & w_{\Gamma} & 17 & 2 \\
 & & 1 &  \not\equiv 0(2) & w_0w_{\Gamma} & w_0w_{\Gamma} & 1 & 18 & 2 \\
 & & 2 & \not\equiv 0(2) & w_0 & w_0 & 1 & 2 & S_3 \\ 
 & (17,3) & 1 & \equiv 0(2) & s_1 & 1 & 1 & 1 & S_3 \\
 & & 2 &  \equiv 0(2) & s_1 & 1 & w_{\Gamma} & 17 & 2 \\
 & & 1 &  \not\equiv 0(2) & s_1 & s_1 & 1 & 12 & S_3 \\
 & & 2 & \not\equiv 0(2) & s_1w_{\Gamma} & s_1w_{\Gamma} & 1 & 22 & 2 \\
 & (17,2) & 1 & \equiv 0(2) & w_0s_1 & 1 & 1 & 1 & S_3 \\
 & & 2 &  \equiv 0(2) & w_0s_1 & 1 & w_{\Gamma} & 17 & 2 \\
 & & 1 &  \not\equiv 0(2) & w_0s_1w_{\Gamma} & w_0s_1w_{\Gamma} & 1 & 22 & 2 \\
 & & 2 & \not\equiv 0(2) & w_0s_1 & w_0s_1 & 1 & 13 & S_3 \\ 
 & (17,5) & 1 & \equiv 0(3) & s_1s_{23} & 1 & 1 & 1 & S_3 \\
 & & 2 &  \equiv 0(3) & s_1s_{23} & 1 & w_{\Gamma} & 17 & 2 \\
 & & 1 &  \not\equiv 0(3) & s_1s_{23} & (s_1s_{23})' & 1 & 4 & S_3 \\
 & & 2 & \not\equiv 0(3) & s_1s_{23}w_{\Gamma} & (s_1s_{23}w_{\Gamma})' & 1 & 19 & 2 \\
 & (17,4) & 1 & \equiv 0(6) & w_0s_1s_{23} & 1 & 1 & 1 & S_3 \\
 & & 2 &  \equiv 0(6) & w_0s_1s_{23} & 1 & w_{\Gamma} & 17 & 2 \\
 & & 1 &  \equiv 0(2)\not\equiv 0(3) & w_0s_1s_{23} & (s_1s_{23})' & 1 & 4 & S_3 \\
 & & 2 & \equiv 0(2)\not\equiv 0(3) & w_0s_1s_{23}w_{\Gamma} & (s_1s_{23}w_{\Gamma})' & 1 & 19 & 2 \\
 & & 1 &  \not\equiv 0(2)\equiv 0(3) & w_0s_1s_{23}w_{\Gamma} & w_0w_{\Gamma} & 1& 18 & 2 \\
 & & 2 & \not\equiv 0(2)\equiv 0(3) & w_0s_1s_{23} & w_0 & 1 & 2 & S_3 \\
 & & 1 & \gcd(6,f/m) = 1 & w_0s_1s_{23}w_{\Gamma} & (w_0s_1s_{23}w_{\Gamma})' & 1 & 20 & 2 \\
 & & 2 & \gcd(6,f/m) = 1 & w_0s_1s_{23} & (w_0s_1s_{23})' & 1 & 5 & S_3 
 \\ \hline\hline
\end{array}
$$
\end{table}

\clearpage

\addtocounter{table}{-1}
\begin{table}[h]
\caption{Construction of Centralizers (continued)}
$$
\begin{array}{ccrcccccc}
	\Gamma & (i,k) & e & f/m & v & v^{f/m} & w & \text{cl}      & C_{\Gamma}( v^{f/m}w ) 
\\ \hline\hline
 1 & (18,1) & 1 & - & 1 & 1 & 1 & 1 & 2 \\
 & & 2 & - & 1 & 1 & w_{\Gamma} & 12 & \\
 & (18,10) & 1 & \equiv 0(2) & w_0w_{\Gamma} & 1 & 1 & 1 & \\
 & & 2 &  \equiv 0(2) & w_0w_{\Gamma} & 1 & w_{\Gamma} & 12 & \\
 & & 1 &  \not\equiv 0(2) & w_0w_{\Gamma} & w_0w_{\Gamma} & 1 & 13 & \\
 & & 2 & \not\equiv 0(2) & w_0 & w_0 & 1 & 2 & \\ 
 & (18,2) & 1 & \equiv 0(2) & s_4 & 1 & 1 & 1 & \\
 & & 2 &  \equiv 0(2) & s_4 & 1 & w_{\Gamma} & 12 & \\
 & & 1 &  \not\equiv 0(2) & s_4 & s_4 & 1 & 17 & \\
 & & 2 & \not\equiv 0(2) & s_4w_{\Gamma} & s_4w_{\Gamma} & 1 & 22 & \\
 & (18,8) & 1 & \equiv 0(2) & w_0s_4 & 1 & 1 & 1 & \\
 & & 2 &  \equiv 0(2) & w_0s_4 & 1 & w_{\Gamma} & 12 & \\
 & & 1 &  \not\equiv 0(2) & w_0s_4w_{\Gamma} & w_0s_4w_{\Gamma} & 1 & 22 & \\
 & & 2 & \not\equiv 0(2) & w_0s_4 & w_0s_4 & 1 & 18 & \\ 
 & (18,9) & 1 & \equiv 0(2) & s_{14} & 1 & 1 & 1 & \\
 & & 2 &  \equiv 0(2) & s_{14} & 1 & w_{\Gamma} & 12 & \\
 & & 1 &  \not\equiv 0(2) & s_{14} & s_{14} & 1 & 12 & \\
 & & 2 & \not\equiv 0(2) & s_{14}w_{\Gamma} & s_{14}w_{\Gamma} & 1 & 3 & \\
 & (18,5) & 1 & \equiv 0(2) & w_0s_{14} & 1 & 1 & 1 & \\
 & & 2 &  \equiv 0(2) & w_0s_{14} & 1 & w_{\Gamma} & 12 & \\
 & & 1 &  \not\equiv 0(2) & w_0s_{14}w_{\Gamma} & w_0s_{14}w_{\Gamma} & 1 & 3 & \\
 & & 2 & \not\equiv 0(2) & w_0s_{14} & w_0s_{14} & 1 & 13 & \\
 \hline\hline
\end{array}
$$
\end{table}

\clearpage

\addtocounter{table}{-1}
\begin{table}[h]
\caption{Construction of Centralizers (continued)}
$$
\begin{array}{ccrcccccc}
	\Gamma & (i,k) & e & f/m & v & v^{f/m} & w & \text{cl}      & C_{\Gamma}( v^{f/m}w ) 
\\ \hline\hline
 & (18,6) & 1 & \equiv 0(4) & s_{14}s_{4} & 1 & 1 & 1 & 2 \\
 & & 2 &  \equiv 0(4) & s_{14}s_{4} & 1 & w_{\Gamma} & 12 & \\
 & & 1 & \equiv 2(4) & s_{14}s_{4} & (s_{14}s_{4})^2 & 1 & 3 & \\
 & & 2 &  \equiv 2(4) & s_{14}s_{4} & (s_{14}s_{4})^2  & w_{\Gamma} & 13 & \\
 & & 1 &  \not\equiv 0(2) & s_{14}s_{4} & (s_{14}s_{4})' & 1 & 23 & \\
 & & 2 & \not\equiv 0(2) & s_{14}s_{4}w_{\Gamma} & (s_{14}s_{4}w_{\Gamma})' & 1 & 21 & \\
 & (18,4) & 1 & \equiv 0({4}) & w_0s_{14}s_{4} & 1 & 1 & 1 & \\
 & & 2 &  \equiv 0(4) & w_0s_{14}s_{4} & 1 & w_{\Gamma} & 12 & \\
 & & 1 & \equiv 2(4) & w_0s_{14}s_{4} & (s_{14}s_{4})^2 & 1 & 3 & \\
 & & 2 &  \equiv 2(4) & w_0s_{14}s_{4} & (s_{14}s_{4})^2  & w_{\Gamma} & 13 & \\
 & & 1 &  \not\equiv 0(2) & w_0s_{14}s_{4}w_{\Gamma} & (w_0s_{14}s_{4}w_{\Gamma})' & 1 & 21 & \\
 & & 2 & \not\equiv 0(2) & w_0s_{14}s_{4} & (w_0s_{14}s_{4})' & 1 & 24 & \\
 & (18,3) & 1 & \equiv 0(3) & s_4s_{3} & 1 & 1 & 1 & \\
 & & 2 &  \equiv 0(3) & s_4s_{3} & 1 & w_{\Gamma} & 12 &  \\
 & & 1 &  \not\equiv 0(3) & s_4s_{3} & (s_4s_{3})' & 1 & 7 &  \\
 & & 2 & \not\equiv 0(3) & s_4s_{3}w_{\Gamma} & (s_4s_{3}w_{\Gamma})' & 1 & 14 & \\
 & (18,7) & 1 & \equiv 0(6) & w_0s_4s_{3} & 1 & 1 & 1 & \\
 & & 2 &  \equiv 0(6) & w_0s_4s_{3} & 1 & w_{\Gamma} & 12 & \\
 & & 1 &  \equiv 0(2)\not\equiv 0(3) & w_0s_4s_{3} & (s_4s_{3})' & 1 & 7 & \\
 & & 2 & \equiv 0(2)\not\equiv 0(3) & w_0s_4s_{3}w_{\Gamma} & (s_4s_{3}w_{\Gamma})' & 1 & 14 & \\
 & & 1 &  \not\equiv 0(2)\equiv 0(3) & w_0s_4s_{3}w_{\Gamma} & w_0w_{\Gamma} & 1& 13 & \\
 & & 2 & \not\equiv 0(2)\equiv 0(3) & w_0s_4s_{3} & w_0 & 1 & 2 & \\
 & & 1 & \gcd(6,f/m) = 1 & w_0s_4s_{3}w_{\Gamma} & (w_0s_4s_{3}w_{\Gamma})' & 1 & 15 & \\
 & & 2 & \gcd(6,f/m) = 1 & w_0s_4s_{3} & (w_0s_4s_{3})' & 1 & 8 & \\
 \hline\hline
\end{array}
$$
\end{table}

\clearpage

\addtocounter{table}{-1}
\begin{table}[h]
\caption{Construction of Centralizers (continued)}
$$
\begin{array}{ccrcccccc}
	\Gamma & (i,k) & e & f/m & v & v^{f/m} & w & \text{cl}      & C_{\Gamma}( v^{f/m}w ) 
\\ \hline\hline
 4 & (19,1) & 1 & - & 1 & 1 & 1 & 1 & 2 \\
 & & 2 & - & 1 & 1 & w_{\Gamma} & 17 & \\
 & (19,10) & 1 & \equiv 0(2) & w_0w_{\Gamma} & 1 & 1 & 1 & \\
 & & 2 &  \equiv 0(2) & w_0w_{\Gamma} & 1 & w_{\Gamma} & 17 & \\
 & & 1 &  \not\equiv 0(2) & w_0w_{\Gamma} & w_0w_{\Gamma} & 1 & 18 & \\
 & & 2 & \not\equiv 0(2) & w_0 & w_0 & 1 & 2 & \\ 
 & (19,2) & 1 & \equiv 0(2) & s_{13} & 1 & 1 & 1 & \\
 & & 2 &  \equiv 0(2) & s_{13} & 1 & w_{\Gamma} & 17 & \\
 & & 1 &  \not\equiv 0(2) & s_{13} & s_{13} & 1 & 17 & \\
 & & 2 & \not\equiv 0(2) & s_{13}w_{\Gamma} & s_{13}w_{\Gamma} & 1 & 3 & \\
 & (19,5) & 1 & \equiv 0(2) & w_0s_{13} & 1 & 1 & 1 & \\
 & & 2 &  \equiv 0(2) & w_0s_{13} & 1 & w_{\Gamma} & 17 & \\
 & & 1 &  \not\equiv 0(2) & w_0s_{13}w_{\Gamma} & w_0s_{13}w_{\Gamma} & 1 & 3 & \\
 & & 2 & \not\equiv 0(2) & w_0s_{13} & w_0s_{13} & 1 & 18 & \\ 
 & (19,3) & 1 & \equiv 0(2) & s_{1} & 1 & 1 & 1 & \\
 & & 2 &  \equiv 0(2) & s_{1} & 1 & w_{\Gamma} & 17 & \\
 & & 1 &  \not\equiv 0(2) & s_{1} & s_{1} & 1 & 12 & \\
 & & 2 & \not\equiv 0(2) & s_{1}w_{\Gamma} & s_{1}w_{\Gamma} & 1 & 22 & \\
 & (19,4) & 1 & \equiv 0(2) & w_0s_{1} & 1 & 1 & 1 & \\
 & & 2 &  \equiv 0(2) & w_0s_{1} & 1 & w_{\Gamma} & 17 & \\
 & & 1 &  \not\equiv 0(2) & w_0s_{1}w_{\Gamma} & w_0s_{1}w_{\Gamma} & 1 & 22 & \\
 & & 2 & \not\equiv 0(2) & w_0s_{1} & w_0s_{1} & 1 & 13 & \\
 \hline\hline
\end{array}
$$
\end{table}

\clearpage

\addtocounter{table}{-1}
\begin{table}[h]
\caption{Construction of Centralizers (continued)}
$$
\begin{array}{ccrcccccc}
	\Gamma & (i,k) & e & f/m & v & v^{f/m} & w & \text{cl}      & C_{\Gamma}( v^{f/m}w ) 
\\ \hline\hline
 & (19,6) & 1 & \equiv 0(4) & s_{13}s_{1} & 1 & 1 & 1 & 2 \\
 & & 2 &  \equiv 0(4) & s_{13}s_{1} & 1 & w_{\Gamma} & 17 & \\
 & & 1 & \equiv 2(4) & s_{13}s_{1} & (s_{13}s_{1})^2 & 1 & 3 & \\
 & & 2 &  \equiv 2(4) & s_{13}s_{1} & (s_{13}s_{1})^2  & w_{\Gamma} & 18 & \\
 & & 1 &  \not\equiv 0(2) & s_{13}s_{1} & (s_{13}s_{1})' & 1 & 23 & \\
 & & 2 & \not\equiv 0(2) & s_{13}s_{1}w_{\Gamma} & (s_{13}s_{1}w_{\Gamma})' & 1 & 16 & \\
 & (19,7) & 1 & \equiv 0({1}) & w_0s_{13}s_{1} & 1 & 1 & 1 & \\
 & & 2 &  \equiv 0(4) & w_0s_{13}s_{1} & 1 & w_{\Gamma} & 17 & \\
 & & 1 & \equiv 2(4) & w_0s_{13}s_{1} & (s_{13}s_{1})^2 & 1 & 3 & \\
 & & 2 &  \equiv 2(4) & w_0s_{13}s_{1} & (s_{13}s_{1})^2  & w_{\Gamma} & 18 & \\
 & & 1 &  \not\equiv 0(2) & w_0s_{13}s_{1}w_{\Gamma} & (w_0s_{13}s_{1}w_{\Gamma})' & 1 & 16 & \\
 & & 2 & \not\equiv 0(2) & w_0s_{13}s_{1} & (w_0s_{13}s_{1})' & 1 & 24 & \\
 & (19,8) & 1 & \equiv 0(3) & s_1s_{2} & 1 & 1 & 1 & \\
 & & 2 &  \equiv 0(3) & s_1s_{2} & 1 & w_{\Gamma} & 17 &  \\
 & & 1 &  \not\equiv 0(3) & s_1s_{2} & (s_1s_{2})' & 1 & 4 &  \\
 & & 2 & \not\equiv 0(3) & s_1s_{2}w_{\Gamma} & (s_1s_{2}w_{\Gamma})' & 1 & 19 & \\
 & (19,9) & 1 & \equiv 0(6) & w_0s_1s_{2} & 1 & 1 & 1 & \\
 & & 2 &  \equiv 0(6) & w_0s_1s_{2} & 1 & w_{\Gamma} & 17 & \\
 & & 1 &  \equiv 0(2)\not\equiv 0(3) & w_0s_1s_{2} & (s_1s_{2})' & 1 & 4 & \\
 & & 2 & \equiv 0(2)\not\equiv 0(3) & w_0s_1s_{2}w_{\Gamma} & (s_1s_{2}w_{\Gamma})' & 1 & 19 & \\
 & & 1 &  \not\equiv 0(2)\equiv 0(3) & w_0s_1s_{2}w_{\Gamma} & w_0w_{\Gamma} & 1& 18 & \\
 & & 2 & \not\equiv 0(2)\equiv 0(3) & w_0s_1s_{2} & w_0 & 1 & 2 & \\
 & & 1 & \gcd(6,f/m) = 1 & w_0s_1s_{2}w_{\Gamma} & (w_0s_1s_{2}w_{\Gamma})' & 1 & 20 & \\
 & & 2 & \gcd(6,f/m) = 1 & w_0s_1s_{2} & (w_0s_1s_{2})' & 1 & 5 & \\
 \hline\hline
\end{array}
$$
\end{table}

\clearpage

\addtocounter{subsection}{1}
\subsection{Fusion of maximal tori and central elements of centralizers}

We say that an $F$-stable maximal torus~$\mathbf{T}$ of~$\mathbf{G}$ 
\textit{fuses} into an $F$-stable closed connected reductive 
subgroup~$\mathbf{M} \leq \mathbf{G}$ of maximal rank, if some $G$-conjugate 
of~$\mathbf{T}$ lies in~$\mathbf{M}$.
Table~\ref{cc}, contains, for each class type $(i,k)$, the $F$-stable maximal
tori of~$\mathbf{G}$, identified by their number, fusing into $\mathbf{M}_{i,k}$.
It also gives the class types of those elements $z \in Z( M_{i,k} )$ which
belong to a different class type. This table is an excerpt of~\cite{LL}.

\begin{landscape}
\begin{longtable}{rp{12.00cm}p{05.00cm}}
\caption{\label{cc} Tori and other class types in centralizers}
\mbox{}\\ \hline\hline
$(i,k)$ & \multicolumn{1}{c}{\text{Tori in $\mathbf{M}_{i,k}$}} & \multicolumn{1}{c}{$(i',k') \in Z(M_{i,k})$}  \\ \hline
	(1,1) & \{ 1, \ldots , 25 \}
     \\ \myhline
(2,1) & 1, 2, 3, 3, 4, 5, 6, 12, 13, 16, 16, 17, 18, 19, 20, 21, 22,
    23, 24, 25 &
     \\ \myhline
(3,1) & 1, 3, 4, 12, 16, 17, 18, 19, 22, 24 &
     \mbox{(2, 1)} \\ \myhline
(3,2) & 2, 3, 5, 13, 16, 17, 18, 20, 22, 23 &
     \mbox{(2, 1)} \\ \myhline
(4,1) & 1, 4, 7, 9, 12, 14, 17, 19, 22 &
     \\ \myhline
(4,2) & 2, 5, 8, 10, 13, 15, 18, 20, 22 &
     \\ \myhline
(5,1) & 1, 2, 3, 3, 7, 8, 12, 12, 13, 13, 14, 15, 17, 18, 21, 21, 22,
    22, 23, 24 &
     \\ \myhline
(6,1) & 1, 3, 4, 12, 16, 17, 18, 20, 22, 23 &
     \mbox{(2, 1)} \\ \myhline
(6,2) & 2, 3, 5, 13, 16, 17, 18, 19, 22, 24 &
     \mbox{(2, 1)} \\ \myhline
(7,1) & 1, 4, 12, 17, 19, 22 &
     \mbox{(2,1)}, \mbox{(3,1)}, \mbox{(4, 1)} \\ \myhline
(7,2) & 2, 5, 13, 18, 20, 22 &
     \mbox{(2,1)}, \mbox{(3,2)}, \mbox{(4, 2)} \\ \myhline
(8,1) & 1, 3, 4, 12, 16 &
     \mbox{(2,1)}, \mbox{(3,1)}, \mbox{(6, 1)} \\ \myhline
(8,2) & 17, 18, 19, 22, 24 &
     \mbox{(2,1)}, \mbox{(3,1)}, \mbox{(6, 2)} \\ \myhline
(8,3) & 17, 18, 20, 22, 23 &
     \mbox{(2,1)}, \mbox{(3,2)}, \mbox{(6, 1)} \\ \myhline
(8,4) & 2, 3, 5, 13, 16 &
     \mbox{(2,1)}, \mbox{(3,2)}, \mbox{(6, 2)} \\ \myhline
(9,1) & 1, 7, 12, 14, 17, 22 &
     \mbox{(4,1)}, \mbox{(5, 1)} \\ \myhline
(9,2) & 2, 8, 13, 15, 18, 22 &
     \mbox{(4,2)}, \mbox{(5, 1)} \\ \myhline
	(10,1) & 1, 3, 7, 12, 13, 15, 17, 21, 22, 23 &
     \mbox{(5, 1)} \\ \myhline
(10,2) & 2, 3, 8, 12, 13, 14, 18, 21, 22, 24 &
     \mbox{(5, 1)} \\ \myhline
\end{longtable}
\end{landscape}

\addtocounter{table}{-1}
\begin{landscape}
\begin{longtable}{rp{06.00cm}p{12.00cm}}
\caption{continued}
\mbox{}\\ \hline\hline
$(i,k)$ & \multicolumn{1}{c}{\text{Tori in $\mathbf{M}_{i,k}$}} & \multicolumn{1}{c}{$(i',k') \in Z(M_{i,k})$}  \\ \hline
(11,1) & 1, 3, 3, 12, 12, 13, 17, 21, 22, 23 &
     \mbox{(2,1)}, \mbox{(5,1)}, \mbox{(10, 1)} \\ \myhline
(11,2) & 2, 3, 3, 12, 13, 13, 18, 21, 22, 24 &
     \mbox{(2,1)}, \mbox{(5,1)}, \mbox{(10, 2)} \\ \myhline
(12,1) & 1, 3, 12, 12, 17, 18, 22, 22 &
     \mbox{(2,1)}, \mbox{(3,1)}, \mbox{(5,1)},
    \mbox{(6, 1)} \\ \myhline
(12,2) & 3, 16, 17, 23 &
     \mbox{(2,1)}, \mbox{(3,2)}, \mbox{(6, 1)} \\ \myhline
(12,3) & 2, 3, 13, 13, 17, 18, 22, 22 &
     \mbox{(2,1)}, \mbox{(3,2)}, \mbox{(5,1)},
    \mbox{(6, 2)} \\ \myhline
(12,4) & 3, 16, 18, 24 &
     \mbox{(2,1)}, \mbox{(3,1)}, \mbox{(6, 2)} \\ \myhline
(13,1) & 1, 4, 12 &
     \mbox{(2,1)}, \mbox{(3,1)}, \mbox{(4,1)},
    \mbox{(6,1)}, \mbox{(7,1)}, \mbox{(8, 1)} \\ \myhline
(13,2) & 17, 19, 22 &
     \mbox{(2,1)}, \mbox{(3,1)}, \mbox{(4,1)},
    \mbox{(6,2)}, \mbox{(7,1)}, \mbox{(8, 2)} \\ \myhline
(13,3) & 18, 20, 22 &
     \mbox{(2,1)}, \mbox{(3,2)}, \mbox{(4,2)},
    \mbox{(6,1)}, \mbox{(7,2)}, \mbox{(8, 3)} \\ \myhline
(13,4) & 8, 10, 15 &
     \mbox{(4, 2)} \\ \myhline
(13,5) & 7, 9, 14 &
     \mbox{(4, 1)} \\ \myhline
(13,6) & 2, 5, 13 &
     \mbox{(2,1)}, \mbox{(3,2)}, \mbox{(4,2)},
    \mbox{(6,2)}, \mbox{(7,2)}, \mbox{(8, 4)} \\ \myhline
(14,1) & 1, 12, 17, 22 &
     \mbox{(2,1)}, \mbox{(3,1)}, \mbox{(4,1)},
    \mbox{(5,1)}, \mbox{(6,1)}, \mbox{(7,1)}, \mbox{(9,1)},
    \mbox{(10,1)}, \mbox{(11,1)}, \mbox{(12, 1)} \\ \myhline
(14,2) & 3, 13, 17, 22 &
     \mbox{(2,1)}, \mbox{(3,2)}, \mbox{(5,1)},
    \mbox{(6,2)}, \mbox{(10,1)}, \mbox{(11,1)}, \mbox{(12, 3)} \\ \myhline
(14,3) & 3, 12, 18, 22 &
     \mbox{(2,1)}, \mbox{(3,1)}, \mbox{(5,1)},
    \mbox{(6,1)}, \mbox{(10,2)}, \mbox{(11,2)}, \mbox{(12, 1)} \\ \myhline
(14,4) & 2, 13, 18, 22 &
     \mbox{(2,1)}, \mbox{(3,2)}, \mbox{(4,2)},
    \mbox{(5,1)}, \mbox{(6,2)}, \mbox{(7,2)}, \mbox{(9,2)},
    \mbox{(10,2)}, \mbox{(11,2)}, \mbox{(12, 3)} \\ \myhline
	(15,1) & 1, 3, 12, 17, 23 &
     \mbox{(2,1)}, \mbox{(5,1)}, \mbox{(6,1)},
    \mbox{(10,1)}, \mbox{(11, 1)} \\ \myhline
(15,2) & 3, 12, 13, 21, 22 &
     \mbox{(2,1)}, \mbox{(5,1)}, \mbox{(10,1)},
    \mbox{(10,2)}, \mbox{(11,1)}, \mbox{(11, 2)} \\ \myhline
\end{longtable}
\end{landscape}

\addtocounter{table}{-1}
\begin{landscape}
\begin{longtable}{rp{03.00cm}p{14.00cm}}
\caption{continued}
\mbox{}\\ \hline\hline
$(i,k)$ & \multicolumn{1}{c}{\text{Tori in $\mathbf{M}_{i,k}$}} & \multicolumn{1}{c}{$(i',k') \in Z(M_{i,k})$}  \\ \hline
(15,3) & 2, 3, 13, 18, 24 &
     \mbox{(2,1)}, \mbox{(5,1)}, \mbox{(6,2)},
    \mbox{(10,2)}, \mbox{(11, 2)} \\ \myhline
(15,4) & 3, 16, 17, 18, 22 &
     \mbox{(2,1)}, \mbox{(6,1)}, \mbox{(6, 2)} \\ \myhline
(15,5) & 6, 16, 21, 23, 24 &
     \mbox{(2, 1)} \\ \myhline
(16,1) & 1, 3, 12, 12 &
     \mbox{(2,1)}, \mbox{(3,1)}, \mbox{(5,1)},
    \mbox{(6,1)}, \mbox{(8,1)}, \mbox{(10,1)}, \mbox{(11,1)},
    \mbox{(12,1)}, \mbox{(15, 1)} \\ \myhline
(16,2) & 17, 18, 22, 22 &
     \mbox{(2,1)}, \mbox{(3,1)}, \mbox{(3,2)},
    \mbox{(5,1)}, \mbox{(6,1)}, \mbox{(6,2)}, \mbox{(8,2)},
    \mbox{(8,3)}, \mbox{(12,1)}, \mbox{(12,3)}, \mbox{(15, 4)} \\ \myhline
(16,3) & 17, 23 &
     \mbox{(2,1)}, \mbox{(3,2)}, \mbox{(5,1)},
    \mbox{(6,1)}, \mbox{(8,3)}, \mbox{(10,1)}, \mbox{(11,1)},
    \mbox{(12,2)}, \mbox{(15, 1)} \\ \myhline
(16,4) & 3, 16 &
     \mbox{(2,1)}, \mbox{(3,1)}, \mbox{(3,2)},
    \mbox{(6,1)}, \mbox{(6,2)}, \mbox{(8,1)}, \mbox{(8,4)},
    \mbox{(12,2)}, \mbox{(12,4)}, \mbox{(15, 4)} \\ \myhline
(16,5) & 21, 21, 23, 24 &
     \mbox{(2,1)}, \mbox{(5,1)}, \mbox{(15, 5)} \\ \myhline
(16,6) & 3, 3, 12, 13 &
     \mbox{(2,1)}, \mbox{(5,1)}, \mbox{(10,1)},
    \mbox{(10,2)}, \mbox{(11,1)}, \mbox{(11,2)}, \mbox{(15, 2)} \\ \myhline
(16,7) & 21, 22 &
     \mbox{(2,1)}, \mbox{(5,1)}, \mbox{(10,1)},
    \mbox{(10,2)}, \mbox{(11,1)}, \mbox{(11,2)}, \mbox{(15, 2)} \\ \myhline
(16,8) & 6, 16 &
     \mbox{(2,1)}, \mbox{(15, 5)} \\ \myhline
(16,9) & 2, 3, 13, 13 &
     \mbox{(2,1)}, \mbox{(3,2)}, \mbox{(5,1)},
    \mbox{(6,2)}, \mbox{(8,4)}, \mbox{(10,2)}, \mbox{(11,2)},
    \mbox{(12,3)}, \mbox{(15, 3)} \\ \myhline
(16,10) & 18, 24 &
     \mbox{(2,1)}, \mbox{(3,1)}, \mbox{(5,1)},
    \mbox{(6,2)}, \mbox{(8,2)}, \mbox{(10,2)}, \mbox{(11,2)},
    \mbox{(12,4)}, \mbox{(15, 3)} \\ \myhline
(17,1) & 1, 7, 17 &
     \mbox{(4,1)}, \mbox{(5,1)}, \mbox{(9,1)},
    \mbox{(10, 1)} \\ \myhline
(17,2) & 13, 15, 22 &
     \mbox{(4,2)}, \mbox{(5,1)}, \mbox{(9,2)},
    \mbox{(10, 1)} \\ \myhline
(17,3) & 12, 14, 22 &
     \mbox{(4,1)}, \mbox{(5,1)}, \mbox{(9,1)},
    \mbox{(10, 2)} \\ \myhline
(17,4) & 5, 10, 20 &
     \mbox{(4, 2)} \\ \myhline
	(17,5) & 4, 9, 19 &
     \mbox{(4, 1)} \\ \myhline
(17,6) & 2, 8, 18 &
     \mbox{(4,2)}, \mbox{(5,1)}, \mbox{(9,2)},
    \mbox{(10, 2)} \\ \myhline
\end{longtable}
\end{landscape}

\addtocounter{table}{-1}
\begin{landscape}
\begin{longtable}{rp{01.50cm}p{15.50cm}}
\caption{continued}
\mbox{}\\ \hline\hline
$(i,k)$ & \multicolumn{1}{c}{\text{Tori in $\mathbf{M}_{i,k}$}} & \multicolumn{1}{c}{$(i',k') \in Z(M_{i,k})$}  \\ \hline
(18,1) & 1, 12 &
     \mbox{(2,1)}, \mbox{(3,1)}, \mbox{(4,1)},
    \mbox{(5,1)}, \mbox{(6,1)}, \mbox{(7,1)}, \mbox{(8,1)},
    \mbox{(9,1)}, \mbox{(10,1)}, \mbox{(11,1)}, \mbox{(12,1)},
    \mbox{(13,1)}, \mbox{(14,1)}, \mbox{(15,1)}, \mbox{(16, 1)} \\ \myhline
(18,2) & 17, 22 &
     \mbox{(2,1)}, \mbox{(3,1)}, \mbox{(3,2)},
    \mbox{(4,1)}, \mbox{(5,1)}, \mbox{(6,1)}, \mbox{(6,2)},
    \mbox{(7,1)}, \mbox{(8,2)}, \mbox{(8,3)}, \mbox{(9,1)},
    \mbox{(10,1)}, \mbox{(11,1)}, \mbox{(12,1)}, \mbox{(12,3)},
    \mbox{(13,2)}, \mbox{(14,1)}, \mbox{(14,2)}, \mbox{(15,4)},
    \mbox{(16, 2)} \\ \myhline
(18,3) & 7, 14 &
     \mbox{(4,1)}, \mbox{(5,1)}, \mbox{(9,1)},
    \mbox{(13, 5)} \\ \myhline
(18,4) & 21, 24 &
     \mbox{(2,1)}, \mbox{(5,1)}, \mbox{(10,2)},
    \mbox{(11,2)}, \mbox{(15,5)}, \mbox{(16, 5)} \\ \myhline
(18,5) & 3, 13 &
     \mbox{(2,1)}, \mbox{(3,2)}, \mbox{(5,1)},
    \mbox{(6,2)}, \mbox{(8,4)}, \mbox{(10,1)}, \mbox{(10,2)},
    \mbox{(11,1)}, \mbox{(11,2)}, \mbox{(12,3)}, \mbox{(14,2)},
    \mbox{(15,2)}, \mbox{(15,3)}, \mbox{(16,6)}, \mbox{(16, 9)} \\ \myhline
(18,6) & 21, 23 &
     \mbox{(2,1)}, \mbox{(5,1)}, \mbox{(10,1)},
    \mbox{(11,1)}, \mbox{(15,5)}, \mbox{(16, 5)} \\ \myhline
(18,7) & 8, 15 &
     \mbox{(4,2)}, \mbox{(5,1)}, \mbox{(9,2)},
    \mbox{(13, 4)} \\ \myhline
(18,8) & 18, 22 &
     \mbox{(2,1)}, \mbox{(3,1)}, \mbox{(3,2)},
    \mbox{(4,2)}, \mbox{(5,1)}, \mbox{(6,1)}, \mbox{(6,2)},
    \mbox{(7,2)}, \mbox{(8,2)}, \mbox{(8,3)}, \mbox{(9,2)},
    \mbox{(10,2)}, \mbox{(11,2)}, \mbox{(12,1)}, \mbox{(12,3)},
    \mbox{(13,3)}, \mbox{(14,3)}, \mbox{(14,4)}, \mbox{(15,4)},
    \mbox{(16, 2)} \\ \myhline
(18,9) & 3, 12 &
     \mbox{(2,1)}, \mbox{(3,1)}, \mbox{(5,1)},
    \mbox{(6,1)}, \mbox{(8,1)}, \mbox{(10,1)}, \mbox{(10,2)},
    \mbox{(11,1)}, \mbox{(11,2)}, \mbox{(12,1)}, \mbox{(14,3)},
    \mbox{(15,1)}, \mbox{(15,2)}, \mbox{(16,1)}, \mbox{(16, 6)} \\ \myhline
(18,10) & 2, 13 &
     \mbox{(2,1)}, \mbox{(3,2)}, \mbox{(4,2)},
    \mbox{(5,1)}, \mbox{(6,2)}, \mbox{(7,2)}, \mbox{(8,4)},
    \mbox{(9,2)}, \mbox{(10,2)}, \mbox{(11,2)}, \mbox{(12,3)},
    \mbox{(13,6)}, \mbox{(14,4)}, \mbox{(15,3)}, \mbox{(16, 9)} \\ \myhline
	\end{longtable}
\end{landscape}

\addtocounter{table}{-1}
\begin{landscape}
\begin{longtable}{rp{01.50cm}p{15.50cm}}
\caption{continued}
\mbox{}\\ \hline\hline
$(i,k)$ & \multicolumn{1}{c}{\text{Tori in $\mathbf{M}_{i,k}$}} & \multicolumn{1}{c}{$(i',k') \in Z(M_{i,k})$}  \\ \hline
(19,1) & 1, 17 &
     \mbox{(2,1)}, \mbox{(3,1)}, \mbox{(4,1)},
    \mbox{(5,1)}, \mbox{(6,1)}, \mbox{(7,1)}, \mbox{(9,1)},
    \mbox{(10,1)}, \mbox{(11,1)}, \mbox{(12,1)}, \mbox{(14,1)},
    \mbox{(15,1)}, \mbox{(17, 1)} \\ \myhline
(19,2) & 3, 17 &
     \mbox{(2,1)}, \mbox{(3,2)}, \mbox{(5,1)},
    \mbox{(6,1)}, \mbox{(6,2)}, \mbox{(10,1)}, \mbox{(11,1)},
    \mbox{(12,2)}, \mbox{(12,3)}, \mbox{(14,2)}, \mbox{(15,1)},
    \mbox{(15, 4)} \\ \myhline
(19,3) & 12, 22 &
     \mbox{(2,1)}, \mbox{(3,1)}, \mbox{(4,1)},
    \mbox{(5,1)}, \mbox{(6,1)}, \mbox{(7,1)}, \mbox{(9,1)},
    \mbox{(10,1)}, \mbox{(10,2)}, \mbox{(11,1)}, \mbox{(11,2)},
    \mbox{(12,1)}, \mbox{(14,1)}, \mbox{(14,3)}, \mbox{(15,2)},
    \mbox{(17, 3)} \\ \myhline
(19,4) & 13, 22 &
     \mbox{(2,1)}, \mbox{(3,2)}, \mbox{(4,2)},
    \mbox{(5,1)}, \mbox{(6,2)}, \mbox{(7,2)}, \mbox{(9,2)},
    \mbox{(10,1)}, \mbox{(10,2)}, \mbox{(11,1)}, \mbox{(11,2)},
    \mbox{(12,3)}, \mbox{(14,2)}, \mbox{(14,4)}, \mbox{(15,2)},
    \mbox{(17, 2)} \\ \myhline
(19,5) & 3, 18 &
     \mbox{(2,1)}, \mbox{(3,1)}, \mbox{(5,1)},
    \mbox{(6,1)}, \mbox{(6,2)}, \mbox{(10,2)}, \mbox{(11,2)},
    \mbox{(12,1)}, \mbox{(12,4)}, \mbox{(14,3)}, \mbox{(15,3)},
    \mbox{(15, 4)} \\ \myhline
(19,6) & 16, 23 &
     \mbox{(2,1)}, \mbox{(3,2)}, \mbox{(6,1)},
    \mbox{(12,2)}, \mbox{(15, 5)} \\ \myhline
(19,7) & 16, 24 &
     \mbox{(2,1)}, \mbox{(3,1)}, \mbox{(6,2)},
    \mbox{(12,4)}, \mbox{(15, 5)} \\ \myhline
(19,8) & 4, 19 &
     \mbox{(2,1)}, \mbox{(3,1)}, \mbox{(4,1)},
    \mbox{(7,1)}, \mbox{(17, 5)} \\ \myhline
(19,9) & 5, 20 &
     \mbox{(2,1)}, \mbox{(3,2)}, \mbox{(4,2)},
    \mbox{(7,2)}, \mbox{(17, 4)} \\ \myhline
(19,10) & 2, 18 &
     \mbox{(2,1)}, \mbox{(3,2)}, \mbox{(4,2)},
    \mbox{(5,1)}, \mbox{(6,2)}, \mbox{(7,2)}, \mbox{(9,2)},
    \mbox{(10,2)}, \mbox{(11,2)}, \mbox{(12,3)}, \mbox{(14,4)},
    \mbox{(15,3)}, \mbox{(17, 6)} \\ \myhline
\end{longtable}
\end{landscape}

\clearpage

\clearpage

\subsection{Sylow $3$-subgroups of centralizers}
\label{Sylow3SubgroupsInCentralizers}

Table~\ref{SylowsInCentralizersI} collects information on the $G$-class types
and their centralizers. The numbering of the class types follows the numbering
in the Tables in~\cite{LL}, with slight adjustments explained in
Subsection~\ref{ClassTypes}. For each class type with label $(i,k)$, we let
$\mathbf{M} := \mathbf{M}_{i,k}$ denote a corresponding $F$-stable regular
subgroup of maximal rank of~$\mathbf{G}$, such that $\mathbf{M}_{i,k}$ is
$G$-conjugate to the centralizers in~$\mathbf{G}$ of the elements of $G$-class
type $(i,k)$; see also Subsection~\ref{ConstructionOfCentralizers}. 
Table~\ref{SylowsInCentralizersI} gives the rough structure
of~$M$ following the conventions introduced in Subsection~\ref{Groups}
and~\ref{CentralProducts}. The structure of these groups can be determined with 
the methods utilized in Subsection~\ref{SomeSplitLeviSubgroups}. The structure 
of the groups for $i \in \{ 4, 7, 9, 13, 17 \}$ assumes $3 \nmid q$. (In case 
$3 \mid q$, these groups are direct products; see Propositions~\ref{C3C} 
and~\ref{NormalizerM13}.) We also give, in case $3 \nmid q$, a representative
of the conjugacy class of radical $3$-subgroups of~$G$ containing a Sylow 
$3$-subgroup of~$M$; cf.\ Proposition~\ref{SylowStructure} and 
Table~\ref{tabradicalageq2}. Finally, the last 
column of this table gives the conditions for the existence of a semisimple 
element $s \in G$ with $C_{\mathbf{G}}( s ) = \mathbf{M}$. The other notational
conventions used in the table have been in effect throughout this paper. The 
integer $\varepsilon \in \{ 1, -1 \}$ is defined by $3 \mid q - \varepsilon$,
provided $3 \nmid q$. Also $d = \gcd( 2, q - 1 ) \in \{ 1, 2 \}$. The symbol
$\circ_1$ denotes the direct product of two groups. Any number in square 
brackets denotes a cyclic group of that order. The notational conventions used 
to distinguish the cases $\varepsilon = 1$ and $\varepsilon = -1$, are the same 
as in Tables \ref{1}--\ref{19}; see the introduction to 
Subsection~\ref{ellBlocksAnTheirInvariants}. The information given in 
Table~\ref{SylowsInCentralizersI}, apart from that on the Sylow $3$-subgroups, 
is contained in~\cite{LL}.

\begin{table}
\caption{\label{SylowsInCentralizersI} Sylow $3$-subgroups of centralizers}
$$
\begin{array}{lllll} \\ \hline\hline
\multicolumn{1}{c}{i} & \multicolumn{1}{c}{k} & \multicolumn{1}{c}{M} 
& \multicolumn{1}{c}{\Syl_3(M)} & \multicolumn{1}{c}{\text{Condition}} \\ \hline\hline
1 & 1 & F_4( q ) & R_{38} \\ \hline
2 & 1 & \Spin_9( q ) & R_{34} & q \text{\ odd} \\ \hline
3 & 1 & (\SL_4^\varepsilon(q)\circ_2\SL_2(q)).2 & R_{34} & 4 \mid q - 1 \\
  &      & (\SL_4^{-\varepsilon}(q)\circ_2\SL_2(q)).2 & R_{17} & 4 \mid q - 1 \\ \hline
3 & 2 & (\SL_4^\varepsilon(q)\circ_2\SL_2(q)).2 & R_{34} & 4 \mid q + 1 \\
  &      & (\SL_4^{-\varepsilon}(q)\circ_2\SL_2(q)).2 & R_{17} & 4 \mid q + 1 \\ \hline
4 & 1, 2 & (\SL_3^\varepsilon(q) \circ_3 \SL_3^\varepsilon(q)).3 & R_{38} &  \\ \hline
5 & 1    & (\Sp_6(q)\circ_2\SL_2(q)).2 & R_{33} & q \text{\ odd} \\ \hline
6 & 1, 2 & ([q-\varepsilon]\circ_d \Spin_7(q)).d & R_{34} & (q,k) \neq (2,1) \\
&      & ([q+\varepsilon]\circ_d \Spin_7(q)).d & R_{32} & (q,k) \neq (2,1) \\ \hline
7 & 1, 2 & (\SL_3^\varepsilon(q) \circ_{3} \GL^\varepsilon_2(q)).3 & R_{34} & q \neq 2, (q,k) \neq (4,1) \\
&      & (\SL_3^{-\varepsilon}(q) \circ_{3} \GL^{-\varepsilon}_2(q)).3 & R_9 & q \neq 2, (q,k) \neq (4,1) \\ \hline
8 & 1, 4 & ([q-\varepsilon]\circ_2 \SL_4^\varepsilon(q)).2 & R_{34} & q \text{\ odd} \\
&      & ([q+\varepsilon]\circ_2 \SL_4^{-\varepsilon}(q)).2 & R_{10} & q \text{\ odd} \\ \cline{2-5}
& 2, 3 & ([q+\varepsilon]\circ_2 \SL_4^{\varepsilon}(q)).2  & R_{32} & q \text{\ odd} \\
&      & ([q-\varepsilon]\circ_2 \SL_4^{-\varepsilon}(q)).2  & R_{17} & q \text{\ odd} \\ \hline
9 & 1, 2 & (\GL^\varepsilon_2(q) \circ_{3} \SL_3^\varepsilon(q)).3 & R_{33} & q \neq 2, (q,k) \neq (4,1) \\
&      & (\GL^{-\varepsilon}_2(q) \circ_{3} (\SL_3^{-\varepsilon}(q)).3 & R_9 & q \neq 2, (q,k) \neq (4,1)  \\ \hline
10& 1, 2 & ([q-\varepsilon]\circ_d \Sp_6(q)).d & R_{33} & (q,k) \neq (2,1) \\
&      & ([q+\varepsilon]\circ_d \Sp_6(q)).d & R_{31} & (q,k) \neq (2,1) \\ \hline
11& 1, 2 & ([q-\varepsilon]\circ_2 (\Sp_4(q)\times \SL_2(q))).2 & R_{18} & q \text{\ odd} \\
&      & ([q+\varepsilon]\circ_2 (\Sp_4(q)\times \SL_2(q))).2 & R_{16} & q \text{\ odd} \\ \hline
	12& 1, 3 & ([q-\varepsilon]\circ_{2} (\SL_2(q)^2 \circ_2 \SL_2( q )).2).2 & R_{18} & q \text{\ odd} \\
&      & ([q+\varepsilon]\circ_{2} (\SL_2(q)^2 \circ_2 \SL_2( q )).2).2 & R_{17} & q \text{\ odd} \\ \cline{2-5}
& 2, 4 & ([q-\varepsilon]\circ_{2} (\SL_2(q^2)\circ_2 \SL_2(q)).2).2 & R_{17} & q \text{\ odd} \\ 
&      & ([q+\varepsilon]\circ_{2} (\SL_2(q^2)\circ_2 \SL_2(q)).2).2 & R_{10} & q \text{\ odd} \\ \hline
13& 1, 6 & ([q-\varepsilon]^2\circ_{3} \SL_3^\varepsilon(q)).3 & R_{34} & q \neq 2, 4 \\ 
&      & ([q+\varepsilon]^2\circ_{3} \SL_3^{-\varepsilon}(q)).3 & R_{2} & q \neq 2, 4 \\ \cline{2-5}
& 2, 3 & ([q^2-1]\circ_{3} \SL_3^\varepsilon(q)).3 & R_{32} & q \neq 2 \\
&      & ([q^2-1]\circ_{3} \SL_3^{-\varepsilon}(q)).3 & R_9  & q \neq 2\\ \cline{2-5}
& 5, 4 & ([q^2+\varepsilon q+1]\circ_{3} \SL_3^{\varepsilon}(q)).3 & R_{30} & (q,k) \neq (2,4) \\ 
&      & ([q^2-\varepsilon q+1]\circ_{3} \SL_3^{-\varepsilon}(q)).3 & R_{2} & (q,k) \neq (2,4)  
\\ \hline\hline
\end{array}
$$
\end{table}

\addtocounter{table}{-1}

\begin{table}
\caption{\label{SylowsInCentralizersII} Sylow $3$-subgroups of centralizers (continued)}
$$
\begin{array}{lllll} \\ \hline\hline
\multicolumn{1}{c}{i} & \multicolumn{1}{c}{k} & \multicolumn{1}{c}{M} 
& \multicolumn{1}{c}{\Syl_3(M)} & \multicolumn{1}{c}{\text{Condition}} \\ \hline\hline
14& 1, 4 & ([q-\varepsilon]^2\circ_{d^2}  \SL_2(q)^2).d^2 & R_{18} & q \neq 2, 4 \\
&      & ([q+\varepsilon]^2\circ_{d^2} \SL_2(q)^2).d^2 & R_{9} & q \neq 2, 4\\ \cline{2-5}
& 2    & ([q^2-1]\circ_{d} \SL_2(q)^2).d & R_{17}, R_{16} & q \neq 2 \\ \cline{2-5}
& 3    & ([q^2-1]\circ_{d} \SL_2(q)^2).d & R_{16}, R_{17} & q \neq 2 \\ \hline
15& 1, 3 & ([q-\varepsilon]^2\circ_{d} \Sp_4(q)).d & R_{18} & q \neq 2, (q,k) \neq (4,1)  \\
&      & ([q+\varepsilon]^2\circ_{d} \Sp_4(q)).d & R_{10} & q \neq 2, (q,k) \neq (4,1)  \\ \cline{2-5}
& 2/ 4 & ([q^2-1]\circ_{d} \Sp_4(q)).d & R_{16}/R_{17} & q \neq 2 \\ \cline{2-5}
& 5    & ([q^2+1]\circ_{d} \Sp_4(q)).d & R_{10} \\  \hline
16& 1, 9 & ([q-\varepsilon]^2\circ_{2} \SL_2(q)^2).2 & R_{18} & q \text{\ odd} \\
&      & ([q+\varepsilon]^2\circ_{2} \SL_2(q)^2).2 & R_{10} & q \text{\ odd} \\ \cline{2-5}
& 2/ 6 & ([q^2-1]\circ_{2} \SL_2(q)^2).2 & R_{17}/R_{16} & q \text{\ odd} \\ \cline{2-5}
& 5    & ([q^2+1]\circ_{2} \SL_2(q)^2).2 & R_{10} & q \text{\ odd} \\ \cline{2-5}
& 3, 10& ([q-\varepsilon]^2\circ_2 \SL_2(q^2)).2 & R_{17} & q \text{\ odd} \\ 
&      & ([q+\varepsilon]^2\circ_2 \SL_2(q^2)).2 & R_{3} & q \text{\ odd} \\ \cline{2-5}
& 4/ 7 & ([q^2-1]\circ_2 \SL_2(q^2)).2 & R_{10}/R_9 & q \text{\ odd} \\ \cline{2-5}
& 8    & ([q^2+1] \circ_2 \SL_2(q^2)).2 & R_3 & q \text{\ odd} \\ \hline
17& 1, 6 & ([q-\varepsilon]^2\circ_{3} \SL_3^\varepsilon(q)).3 & R_{33} & q \neq 2, 4 \\ 
&      & ([q+\varepsilon]^2\circ_{3} \SL_3^{-\varepsilon}(q)).3 & R_{3} & q \neq 2, 4 \\ \cline{2-5}
& 3, 2 & ([q^2-1]\circ_{3} \SL_3^\varepsilon(q)).3 & R_{31} & q \neq 2 \\
&      & ([q^2-1]\circ_{3} \SL_3^{-\varepsilon}(q)).3 & R_9 & q \neq 2 \\ \cline{2-5}
& 5, 4 & ([q^2+\varepsilon q+1]\circ_{3} \SL_3^\varepsilon(q)).3 & R_{29} & (q,k) \neq (2,4) \\
&      & ([q^2-\varepsilon q+1]\circ_{3} \SL_3^{-\varepsilon}(q)).3 & R_{3} & (q,k) \neq (2,4) \\ \hline
18& 1, 10& ([q-\varepsilon]^3\circ_{d} \SL_2(q)).d & R_{18} & q \neq 2, 4 \\ 
&      & ([q+\varepsilon]^3\circ_{d} \SL_2(q)).d & R_{2} & q \neq 2, 4 \\ \cline{2-5}
& 2, 8/ 9, 5& (([q^2-1]\times [q-\varepsilon])\circ_{d} \SL_2(q)).d & R_{17}/R_{16} & q \neq 2, (q,k) \neq (4,2), (4,9) \\ 
&           & (([q^2-1]\times [q+\varepsilon])\circ_{d} \SL_2(q)).d & R_9/R_{10} & q \neq 2, (q,k) \neq (4,2), (4,9) \\ \cline{2-5}
& 3, 7 & ([q^3-\varepsilon]\circ_{d} \SL_2(q)).d & R_{11} & (q,k) \neq (2,3) \\ 
&      & ([q^3+\varepsilon]\circ_{d} \SL_2(q)).d & R_{2} & (q,k) \neq (2,3) \\ \cline{2-5}
& 6, 4 & (([q-\varepsilon]\times [q^2+1])\circ_{d} \SL_2(q)).d & R_{10} & (q,k) \neq (2,6) \\ 
&      & (([q+\varepsilon]\times [q^2+1])\circ_{d} \SL_2(q)).d & R_{2} & (q,k) \neq (2,6) 
\\ \hline\hline
\end{array}
$$
\end{table}

\addtocounter{table}{-1}

\begin{table}
\caption{\label{SylowsInCentralizersIII} Sylow $3$-subgroups of centralizers (continued)}
$$
\begin{array}{lllll} \\ \hline\hline
\multicolumn{1}{c}{i} & \multicolumn{1}{c}{k} & \multicolumn{1}{c}{M} 
& \multicolumn{1}{c}{\Syl_3(M)} & \multicolumn{1}{c}{\text{Condition}} \\ \hline\hline
19& 1, 10& ([q-\varepsilon]^3\circ_{d} \SL_2(q)).d & R_{18} & q \neq 2, 4 \\ 
&      & ([q+\varepsilon]^3\circ_{d} \SL_2(q)).d & R_{3} & q \neq 2, 4 \\ \cline{2-5}
& 2, 5/ 3, 4& (([q^2-1]\times [q-\varepsilon])\circ_{d} \SL_2(q)).d & R_{17}/R_{16} & q \neq 2 \\
&           & (([q^2-1]\times [q+\varepsilon])\circ_{d} \SL_2(q)).d & R_{10}/R_9 & q \neq 2 \\ \cline{2-5}
& 8, 9 & ([q^3-\varepsilon]\circ_{d} \SL_2(q)).d & R_{12} & (q,k) \neq (2,8) \\ 
&      & ([q^3+\varepsilon]\circ_{d} \SL_2(q)).d & R_{3} & (q,k) \neq (2,8) \\ \cline{2-5}
& 6, 7 & (([q-\varepsilon]\times [q^2+1])\circ_{d} \SL_2(q)).d & R_{10} & (q,k) \neq (2,6) \\ 
&      & (([q+\varepsilon]\times [q^2+1])\circ_{d} \SL_2(q)).d & R_{3} & (q,k) \neq (2,6) \\ \hline
20& 1, 2 & [q-\varepsilon]^4 & R_{18} & q \neq 2 \\ 
  &      & [q+\varepsilon]^4 & 1 & q \neq 2 \\ \cline{2-5}
  & 3    & [d] \times [(q^2 - 1)/d] \times [q^2-1] & R_{10} & q \neq 2 \\ \cline{2-5}
  & 4, 5/ 7, 8& [q-\varepsilon] \times [q^3-\varepsilon] & R_{12}/R_{11} & q \neq 2 \\ 
  &           & [q+\varepsilon] \times [q^3+\varepsilon] & 1 & q \neq 2 \\ \cline{2-5}
  & 6    & [q^2 + 1]^2 & 1 & q \neq 2 \\ \cline{2-5}
  & 9, 10& [q^2+\varepsilon q+1]^2 & R_8 & q \neq 2 \\
  &      & [q^2-\varepsilon q+1]^2 & 1 & q \neq 2 \\ \cline{2-5}
  & 11   & [q^4 - q^2 + 1] & 1 \\ \cline{2-5}
  & 12,13/17,18& [q-\varepsilon]^2 \times [q^2 - 1] & R_{16}/R_{17} & q \neq 2 \\ 
  &            & [q+\varepsilon]^2 \times [q^2 - 1] & R_{2}/R_{3} & q \neq 2 \\ \cline{2-5}
  & 14,15/19,20& [(q^3-\varepsilon)(q+\varepsilon)] & R_{6}/R_{7} & (q,k) \neq (2,15), (2,20) \\ 
  &            & [(q^3+\varepsilon)(q-\varepsilon)] & R_{2}/R_{3} & (q,k) \neq (2,15), (2,20) \\ \cline{2-5}
  & 16/ 21& [d] \times [(q^4 - 1)/d] & R_{3}/R_{2} & q \neq 2 \\ \cline{2-5}
  & 22   & [q^2-1]^2 & R_9 & q \neq 2 \\ \cline{2-5}
  & 23, 24& [(q^2+1)(q-\varepsilon)] \times [q-\varepsilon] & R_{10} & q \neq 2 \\
  &       & [(q^2+1)(q+\varepsilon)] \times [q+\varepsilon] & 1 & q \neq 2 \\ \cline{2-5}
  & 25    & [q^4 + 1] & 1 
\\ \hline\hline
\end{array}
$$
\end{table}

\clearpage

\clearpage

%\addtocounter{subsection}{1}
\subsection{The radical $3$-subgroups of~$F_4(q)$}
\label{secradF4Appendix}
The table in this subsection lists the non-trivial radical $3$-subgroups of 
$G = F_4( q )$
and some of their properties in case $3 \nmid q$. The displayed information has
been determined in~\cite{AH2} for odd~$q$, and in~\cite{AnDF4} for 
even~$q$. 

The parameters $\varepsilon$,~$e$, $a$ and $d$ have the usual meaning: 
$\varepsilon \in \{ -1, 1 \}$ is such that $3 \mid q - \varepsilon$, and
$3^a$ is the highest power of~$3$ dividing $ q - \varepsilon$; moreover,
$e = 1$ if $\varepsilon = 1$, and $e = 2$, otherwise, and
$d = \gcd( 2, q - 1 )$.

Let~$R$ be one of the listed radical $3$-subgroups. For easier reference, the 
first column of Table~\ref{tabradicalageq2} contains a name for~$R$, or rather 
for a representative of the $G$-conjugacy class containing~$R$. This numbering 
is different from the one given in~\cite{AnDF4}.
The second column gives the rough structure of the groups, using the conventions
naming cyclic groups, extensions and central products as introduced in 
Subsections~\ref{Groups} and~\ref{CentralProducts}. These conventions are also 
used in Columns $5$--$7$,
where $C_G(R)$, $N_G(R)$ and $\Out_G( R ) = N_G( R )/RC_G( R )$ are described. 
The groups $L^i \cong \SL_3^\varepsilon(q)$, $i = 1, 2$ have the same 
significance as in Subsection~\ref{secradF4}, and $T_i$ and~$D_i$,
$i = 1, 2$, denote a maximal $e$-split torus and a Sylow $3$-subgroup 
of~$L^i$, respectively, constructed in Subsection~\ref{secradLe}.

The \textit{characteristic} of~$R$ is defined to be the
elementary abelian $3$-subgroup $\Omega_1(Z(R))$, unless
$R \in_G \{ R_{29}\text{--}R_{34} \}$, where the characteristic is
defined to be $\Omega_1([R,R])$. In the latter cases, $[R,R]$ is abelian.
Recall that if~$Q$ is a finite $r$-group for some prime~$r$, then
$\Omega_1( Q )$ denotes the subgroup of~$Q$ generated by its elements of
order~$r$. This is elementary abelian if $Q$ is abelian.
The characteristics of the radical subgroups of~$G$ are worked out in
\cite[Lemma~$4.2$]{AH2} and in \cite[Proposition~$4.5$ and Theorem~$4.8$]{AnDF4},
respectively.
The fourth column of Table~\ref{tabradicalageq2} gives the names of the
characteristics in the classification of {\rm \cite[Table~$4$]{ADHE6}} (which
holds for all~$q$).
The third column of Table~\ref{tabradicalageq2} gives the types of the
characteristics. If~$E$ is a characteristic, the \textit{type} of~$E$ is
the triple $(i, j, k)$ of non-negative integers adding up to $|E| - 1$,
indicating that $E$ has exactly $i$, $j$ and $k$ non-trivial elements
lying in the conjugacy classes $3\A$, $3\B$ and $3\C$ of~$G$, respectively. 
The type of~$E$ is written as $3\A_{i}\B_{j}\C_{k}$ (with obvious modifications 
if one or two of the $i, j, k$ are equal to~$0$).  
Column~$8$ contains the conditions on~$q$ for which the corresponding radical
subgroup exists. Finally, Column~$9$ gives the name of the group~$R^\dagger$,
defined in Subsection~\ref{DualRadicalSubgroups}.

\begin{landscape}
{\small
\begin{table}[ht]
\hspace*{-1.5cm}
\scalebox{1.0}{
\begin{tabular}[ht]{lllllllll}
 & \multicolumn{1}{c}{$R$} & \multicolumn{1}{c}{$\mbox{\rm Type}$} & 
   \multicolumn{1}{c}{$\mbox{\rm Char}$}  & \multicolumn{1}{c}{$C_G(R)$}  
  &\multicolumn{1}{c}{$N_G(R)$} & \multicolumn{1}{c}{$\Out_G(R)$} & 
   \multicolumn{1}{c}{$\mbox{\rm Cond.}$} & \multicolumn{1}{c}{$R^\dagger$} \\ \hline\hline
$R_1$& $3$ & $3\C$ & $E_3$ & $(\SL_3^\varepsilon(q)\circ_3\SL_3^\varepsilon(q)).3$&
 $((\SL_3^\varepsilon(q)\circ_3 \SL_3^\varepsilon(q)).3).2$ & $2$ & $q \neq 2$ & $R_1$ \\
 
$R_2$& $3^a$ & $3\A$ & $E_1$ & $([q-\varepsilon] \circ_{d} \Sp_6(q)).d$
& $(([q-\varepsilon] \circ_{d} \Sp_6(q)).d).2$ & $2$ & & $R_3$ \\

$R_3$& $3^a$ & $3\B$ & $E_{2}$ & $([q-\varepsilon] \circ_{d} \Spin_7(q)).d$&
$(([q-\varepsilon] \circ_{d} \Spin_7(q)).d).2$ & $2$ & & $R_2$ \\

$R_4$& $3^a$ & $3\C$ & $E_{3}$ & $(\GL_2^\varepsilon(q)\circ_3\SL_3^\varepsilon(q)).3$&
 $((\GL_2^\varepsilon(q)\circ_3 \SL_3^\varepsilon(q)).3).2$ & $2$ & $a \geq 2$ & $R_5$ \\
 
$R_5$ & $3^a$ & $3\C$ & $E_{3}$ & $(\SL_3^\varepsilon(q)\circ_3 \GL_2^\varepsilon(q)).3$&
 $((\SL_3^\varepsilon(q)\circ_3 \GL_2^\varepsilon(q)).3).2$ & $2$ & $a \geq 2$ & $R_4$ \\
 
$R_6$& $3^{a+1}$ & $3\C$ & $E_{3}$ & $(\GL_2^\varepsilon(q) \circ_3 [q^2+\varepsilon q+1]).3$&
 $((\GL_2^\varepsilon(q) \circ_3 [q^2+\varepsilon q+1]).3).6$ & $6$ & $q \neq 2$ & $R_7$ \\
 
$R_7$ & $3^{a+1}$ & $3\C$ & $E_{3}$ & $([q^2+\varepsilon q+1] \circ_3 \GL_2^\varepsilon(q)).3$&
 $(([q^2+\varepsilon q+1] \circ_3 \GL_2^\varepsilon(q)).3).6$
  & $6$ & $q \neq 2$ & $R_6$ \\ %\hline

$R_8$ & $3^2$ & $(3\C^2)_2$ & $E_{9}$ & $[q^2+\varepsilon q+1]^2.3$ & $([q^2+\varepsilon q+1]^2.3).\SL_2(3)$
& $\SL_2( 3 )$ & $q \neq 2$ & $R_8$ \\

$R_9$ & $[3^a]^2$ & $3\A_2\B_2\C_4$ & $E_{5}$ & $(\GL_2^\varepsilon(q)\circ_{3}\GL_2^\varepsilon(q)).3$&
$((\GL_2^\varepsilon(q)\circ_{3}\GL_2^\varepsilon(q)).3).2^2$ & $2^2$ & $q \neq 2$ & $R_9$ \\

$R_{10}$& $[3^a]^2$ & $3\A_4\B_4$ & $E_{4}$ & $ ([q-\varepsilon]^2 \circ_{d} \Sp_4(q)).d$&
$(([q-\varepsilon]^2 \circ_{d} \Sp_4(q)).d).D_8$ & $D_8$ & & $R_{10}$ \\

$R_{11}$ & $[3^a] \times [3^{a+1}]$ & $3\A_6\C_2$ & $E_{6}$ &
$[q-\varepsilon] \times [q^3-\varepsilon]$ & $([q-\varepsilon] \times [q^3-\varepsilon]).(S_3\times 6)$
& $S_3 \times 6$ & $q \neq 2, 4, 8$ & $R_{12}$ \\

$R_{12}$ & $[3^{a+1}] \times [3^a]$ & $3\B_6 \C_2$ & $E_{7}$ & $[q^3-\varepsilon] \times [q-\varepsilon]$ 
& $([q^3-\varepsilon] \times [q-\varepsilon]).(6\times S_3)$ & $6 \times S_3$ & $q \neq 2, 4, 8$ & $R_{11}$ \\

$R_{13}$& $[3^a]^2$ & $3\A_6 \C_2$ & $E_{6}$ & $(T_1\circ_3 L^2).3$ & 
$((T_1 \circ_3 L^2).3).(S_3\times 2)$ & $S_3 \times 2$ & $q \neq 2, 4, 8$ & $R_{14}$ \\

$R_{14}$& $[3^a]^2$ & $3\B_6 \C_2$ & $E_{7}$ & $(L^1 \circ_3 T_2).3$ & 
$((L^1\circ_3 T_2).3).(2\times S_3)$ & $2 \times S_3$ & $q \neq 2, 4, 8$ & $R_{13}$ \\ %\hline

$R_{15}$&  $3^3$ & $(3\C^3)_1$ & $E_{13}$ & $3^3$ & $3^3.\SL_3(3)$ & $\SL_3(3)$ & & $R_{15}$ \\

$R_{16}$&  $[3^a]^3$ & $3\A_{6}\B_{12}\C_8$ & $E_{11}$ & $(\SL_2(q) \circ_d [q-\varepsilon]^3).d$&
$((\SL_2(q) \circ_d [q-\varepsilon]^3).d).2^3.S_3$ & $2^3.S_3$ & $q \neq 2$ & $R_{17}$ \\

$R_{17}$&  $[3^a]^3$ & $3\A_{12}\B_{6}\C_8$ & $E_{10}$ & $([q-\varepsilon]^3 \circ_d \SL_2(q)).d$&
$(([q-\varepsilon]^3 \circ_d \SL_2(q)).d).2^3.S_3$ & $2^3.S_3$ & $q \neq 2$ & $R_{16}$ \\  %\hline

$R_{18}$&  $[3^a]^4$ & $3\A_{24}\B_{24}\C_{32}$ & $E_{15}$ & $T = [q-\varepsilon]^4$ & $T.W(F_4)$ & $W(F_4)$ 
& & $R_{18}$ \\ \hline\hline
\end{tabular}}
\bigskip
\caption{\label{tabradicalageq2} Non-trivial radical $3$-subgroups of $G = F_4(q)$}
\end{table}
}
\end{landscape}

\addtocounter{table}{-1}

\begin{landscape}
{\small
\begin{table}[ht]
\hspace*{-1.5cm}
\scalebox{1.0}{
\begin{tabular}[ht]{lllllllll}
 & \multicolumn{1}{c}{$R$} & \multicolumn{1}{c}{$\mbox{\rm Type}$} & 
   \multicolumn{1}{c}{$\mbox{\rm Char}$} & \multicolumn{1}{c}{$C_G(R)$}
  & \multicolumn{1}{c}{$N_G(R)$} & \multicolumn{1}{c}{$\Out_G(R)$} &
  \multicolumn{1}{c}{$\mbox{\rm Cond.}$} & \multicolumn{1}{c}{$R^\dagger$} \\ \hline\hline
$R_{19}$&  $3_+^{1+2}$ & $3\C$ & $E_{3}$ & $L^1=\SL_3^\varepsilon(q)$&
$(L^1\circ_3 (3_+^{1+2}.\SL_2(3))).2$ & $\SL_2(3).2$ & $q \neq 2$ & $R_{19}$ \\

$R_{20}$&  $3_+^{1+2}$ & $3\C$ & $E_{3}$ & $L^2=\SL_3^\varepsilon(q)$&
$((3_+^{1+2}.\SL_2(3))\circ_3 L^2).2$ & $\SL_2(3).2$ & $q \neq 2$ & $R_{20}$ \\

$R_{21}$&  $3_+^{1+4}$ & $3\C$ & $E_3$ & $3$ & 
$((3_+^{1+2}.Q_8)\circ_3 (3_+^{1+2}.Q_8)).S_3$ & $(Q_8 \times Q_8).S_3$ & $a = 1$ & $R_{21}$ \\ 

$R_{21}$&  $3_+^{1+4}$ & $3\C$ & $E_{3}$ & $3$ &
$((3_+^{1+2}.\SL_2(3))\circ_3 (3_+^{1+2}.\SL_2(3))).2$ & $(\SL_2(3) \times \SL_2(3)).2$ & $a \geq 2$ & $R_{21}$ \\

$R_{22}$&  $3_+^{1+4}$ & $3\C$ & $E_{3}$ & $3$ &
$(3_+^{1+2}.\SL_2(3))\circ_3 (3_+^{1+2}.\SL_2(3))$ & $\SL_2(3) \times \SL_2(3)$ & $a \geq 2$ & $R_{22}$ \\ %\hline

$R_{23}$&  $[3^a]\circ_3 3_+^{1+2}$ & $3\C$ & $E_{3}$ & $\GL_2^\varepsilon(q) \leq L^1$&
$(\GL_2^\varepsilon(q)\circ_3 (3_+^{1+2}.\SL_2(3))).2$ & $\SL_2(3).2$ & $a \geq 2$ & $R_{24}$ \\

$R_{24}$&  $3_+^{1+2}\circ_3 [3^a]$ & $3\C$ & $E_{3}$ & $\GL_2^\varepsilon(q) \leq L^2$&
$((3_+^{1+2}.\SL_2(3))\circ_3 \GL_2^\varepsilon(q)).2$ & $\SL_2(3).2$ & $a \geq 2$ & $R_{23}$ \\

$R_{25}$&  $[3^a]^2\circ_3 3_+^{1+2}$ & $3\A_6\C_2$ & $E_{6}$ & $T_1 = [q-\varepsilon]^2$&
$((T_1.S_3)\circ_3(3_+^{1+2}.\SL_2(3)) ).2$ & $(S_3 \times \SL_2(3)).2$ & $q \neq 2, 4, 8$ & $R_{26}$ \\

$R_{26}$ & $3_+^{1+2}\circ_3 [3^a]^2$ & $3\B_6\C_2$ & $E_{7}$ & $T_2 = [q-\varepsilon]^2$&
$((3_+^{1+2}.\SL_2(3))\circ_3(T_2.S_3))).2$ & $(\SL_2(3) \times S_3).2$ & $q \neq 2, 4, 8$ & $R_{25}$ \\  %\hline

$R_{27}$&  $D_2=[3^a]^2.3$ & $3\C$ & $E_{3}$ & $L^1$ & $(L^1\circ_3 (D_2.2)).S_3$ & $2 \times S_3$ & $a \geq 2$ & $R_{27}$ \\

$R_{28}$&  $D_1=[3^a]^2.3$ & $3\C$ & $E_{3}$ & $L^2$ & $((D_1.2)\circ_3 L^2).S_3$ & $2 \times S_3$ & $a \geq 2$ & $R_{28}$ \\

$R_{29}$&  $D_2.3$ & $3\B_6\C_2$ & $E_{7}$ & $[q^2+\varepsilon q+1] \leq L^1 $&
$(([q^2+\varepsilon q+1] \circ_3 (D_2.2)).3).S_3$ & $2 \times S_3$ & $q \neq 2$ & $R_{30}$ \\

$R_{30}$&  $D_1.3$ & $3\A_6\C_2$ & $E_{6}$ & $[q^2+\varepsilon q+1] \leq L^2$&
$(((D_1.2) \circ_3 [q^2+\varepsilon q+1]).3).S_3$ & $2 \times S_3$ & $q \neq 2$ & $R_{29}$ \\

$R_{31}$&  $([3^a]\circ_3 D_2).3$ & $3\B_6\C_2$ & $E_{7}$ & $\GL_2^\varepsilon(q) \leq L^1$ &
$(\GL_2^\varepsilon(q)\circ_3 (D_2.2)).S_3$ & $2^2$ & $q \neq 2$ & $R_{32}$ \\

$R_{32}$&  $(D_1\circ_3 [3^a]).3$ & $3\A_6\C_2$ & $E_{6}$ & $\GL_2^\varepsilon(q) \leq L^2$ &
$((D_1.2)\circ_3 \GL_2^\varepsilon(q)).S_3$ & $2^2$ & $q \neq 2$ & $R_{31}$ \\ %\hline

$R_{33}$&  $([3^a]^2\circ_3 D_2).3$ & $3\B_6\C_2$ & $E_{7}$ & $T_1 = [q-\varepsilon]^2$&
$((T_1.2)\circ_3 (D_2.2)).S_3$ & $2^3$ & $q \neq 2, 4, 8$ & $R_{34}$ \\

$R_{34}$&  $(D_1\circ_3 [3^a]^2).3$ & $3\A_6\C_2$ & $E_{6}$ & $T_2 = [q-\varepsilon]^2$&
$((D_1.2)\circ_3 (T_2.2)).S_3$ & $2^3$ & $q \neq 2, 4, 8$ & $R_{33}$ \\

$R_{35}$&  $3_+^{1+2}\circ_3 D_2$ & $3\C$ & $E_{3}$ & $3$
& $((3_+^{1+2}.\SL_2(3)\circ_3(D_2.2)).2$ & $(\SL_2(3) \times 2).2$ & $a \geq 2$ & $R_{35}$ \\

$R_{36}$&  $D_1\circ_3 3_+^{1+2}$ & $3\C$ & $E_{3}$ & $3$ 
& $((D_1.2)\circ_3 (3_+^{1+2}.\SL_2(3))).2$ & $(2 \times \SL_2(3)).2$ & $a \geq 2$ & $R_{36}$ \\ %\hline

$R_{37}$&  $[3^a]^4.3$ & $(3\C^2)_1$ & $E_{8}$ & $3^2$ & $[3^a]^4.3.\SL_2(3).2$ & $\SL_2(3).2$ & & 
$R_{37}$ \\\hline

$R_{38}$&  $\Syl_3(G)$ & $3\C$ & $E_{3}$ & $3$ & $\Syl_3(G).2^3$ & $2^3$ & & $R_{38}$ \\\hline\hline

\end{tabular}}
\bigskip
\caption{Non-trivial radical $3$-subgroups of $G = F_4(q)$ (continued)}
\end{table}
}
\end{landscape}

\end{document}